\documentclass[runningheads,a4paper]{llncs}
\usepackage{graphicx}
\usepackage{amsmath,amssymb} % define this before the line numbering.
\usepackage{color}
\usepackage{amssymb}
\usepackage{graphicx}
\usepackage{url}
\usepackage{float}
\usepackage{comment}
\usepackage[utf8]{inputenc}

\graphicspath{{./}{chap3/}}

\DeclareMathOperator{\Span}{Span}

\DeclareMathOperator{\supp}{supp}

\DeclareMathOperator*{\argmin}{\arg\!\min}

\newcommand{\newbf}{}
\newcommand{\newbfa}{}
\newcommand{\newbfb}{}

\def\M{\mathcal{M}}
\def\F{\mathcal{F}}
\def\E{\mathcal{E}}

\def\R{\mathbb{R}}

\def\T{\mathbb{T}}
\def\T{\mathbb{T}}
\def\Z{\mathbb{Z}}

\def\e{\varepsilon}
\def\a{\beta}
\def\O{\Omega}
\def\la{\lambda}
\def\vphi{\varphi}

\begin{document}

%\pagestyle{headings}
%\mainmatter

\title{Adapted Basis for Non-Local\\ Reconstruction of Missing Spectrum}
\titlerunning{Adapted Basis for Non-Local Reconstruction of Missing Spectrum}
%\authorrunning{Khalid Jalalzai, Antonin Chambolle}
%\author{Khalid Jalalzai\thanks{\tt khalid.jalalzai@polytechnique.edu}
%, Antonin Chambolle\thanks{\tt antonin.chambolle@polytechnique.fr}
%}
\authorrunning{Antonin Chambolle, Khalid Jalalzai}
\author{Antonin Chambolle\thanks{\tt antonin.chambolle@cmap.polytechnique.fr}, Khalid Jalalzai\thanks{\tt khalid.jalalzai@polytechnique.edu}
}
\institute{CMAP, \'Ecole Polytechnique, CNRS}

\maketitle

\begin{abstract}
The object of this work is to design an adequate regularization for the problem of recovering missing Fourier coefficients, \newbfa{particularly in some
non standard situations were low frequency coefficients are lost.}
%\newbf{We suggest a general approach based on a new and simple Non-Local quadratic energy and (among other) we propose a technique to build an original patchwise similarity measure that is adapted to the missing spectrum.} 
\newbfa{In the framework of non-local regularization, we propose a technique to build an original patchwise similarity measure that is adapted to the missing spectrum. Then, a simple Non-Local quadratic energy is minimized.}
By construction, the similarity criterion is invariant under the corruption process so that the distance between two patches of the corrupted image is almost
exactly equal to the one computed on the clean image. We illustrate our method with experiments which show its efficiency, 
both in terms of speed and quality of the results, with respect to other common approaches. We show that the method is practical on synthetic examples which are built upon models of inverse scattering problems, synthetic aperture mirrors for spatial imaging or also medical imaging.
\end{abstract}

% \begin{comment}
\section{Introduction}

\begin{figure}[H]
\centering
\begin{minipage}[c]{0.49\linewidth}
\label{spectrum_scatt}
\centering
\includegraphics[width=0.6\linewidth]{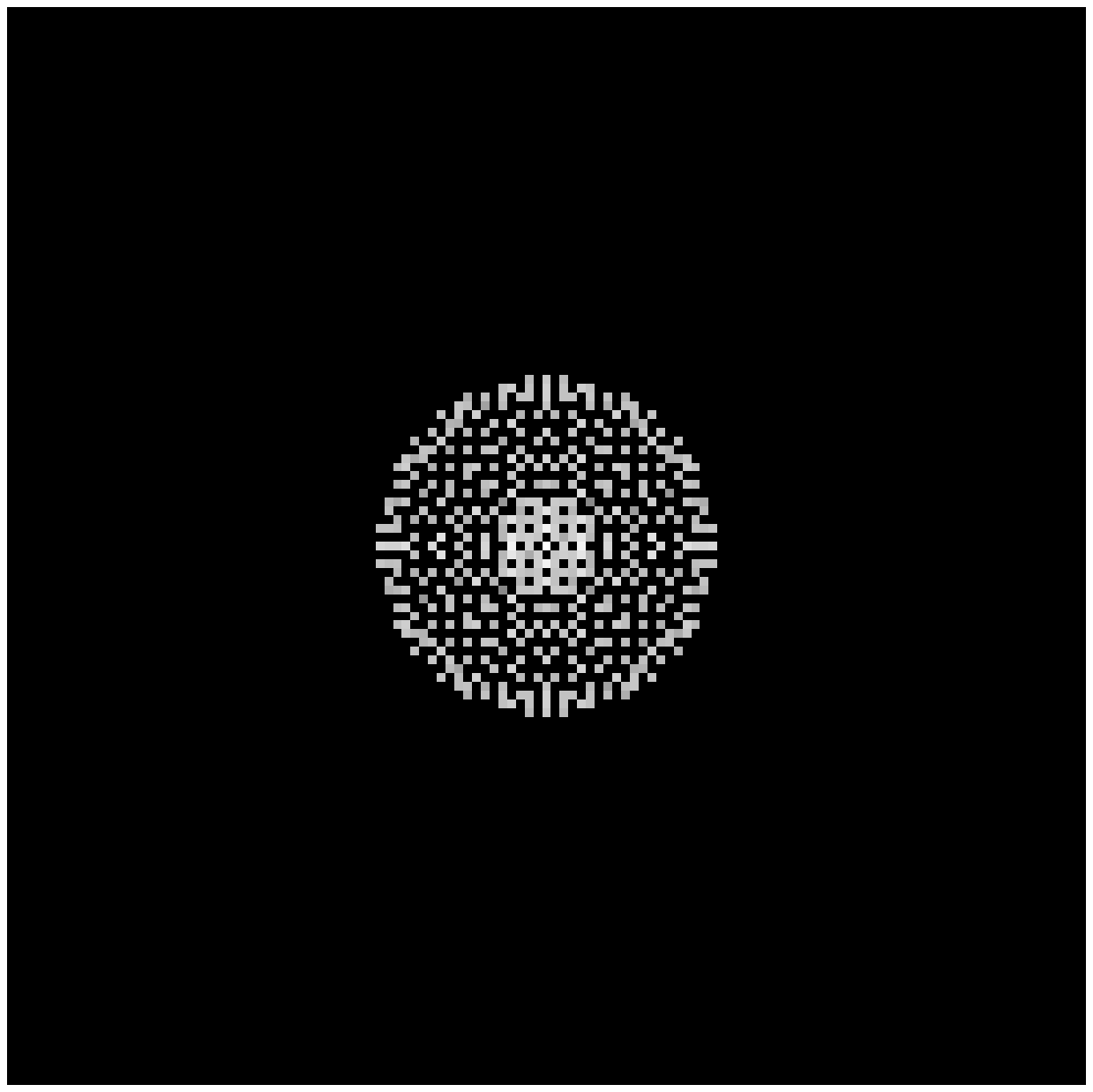}
\vspace{-0.1cm}
\caption%[]
{Spectrum obtained with 32 incident planewaves and 32 measurement directions.}
%{$\F(\chi_D)\chi_{B(0,2k)}$ with 32 incident waves and 32 measurement directions.}
\end{minipage}
\begin{minipage}[c]{0.49\linewidth}
\centering
% \vspace{-0.1cm}
\includegraphics[height=0.2\vsize]{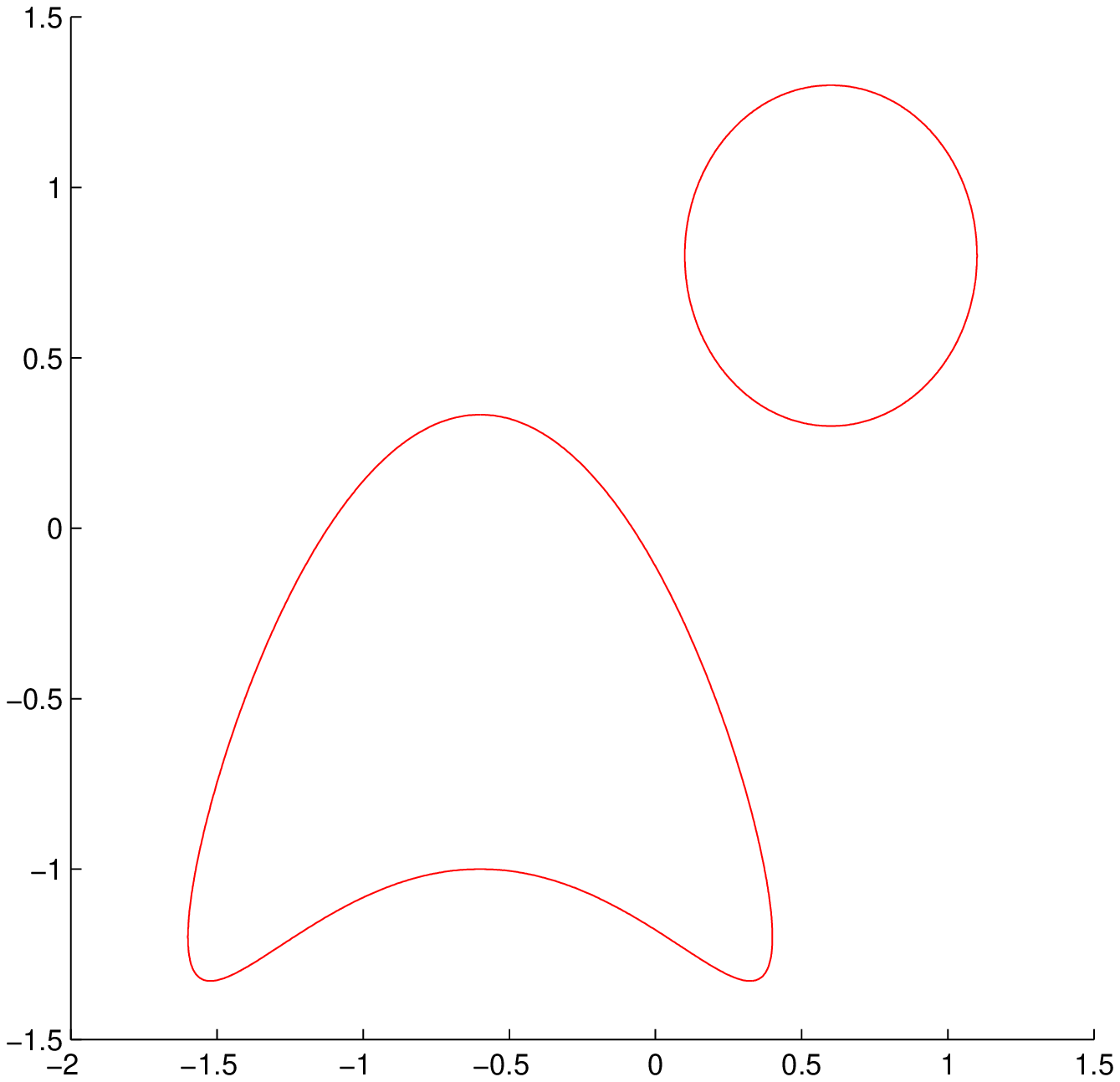}%trim=l b r t
\vspace{-0.1cm}
\caption%[]
{Scatterer $D$.}
% \vspace{0.5cm}
\label{scatterer}
\end{minipage}
\end{figure}

% \label{Motivation}
\label{sec_scatt}
In this paper, we consider the problem of retrieving missing Fourier coefficients of a raw data. This problem has important applications in various fields: 
% \noindent In particular, several important inverse problems can be casted into this framework:
%\vspace{-0.2cm}
\begin{itemize}
\item The zooming problem in image processing where one has to construct high frequency data from a low-resolution sample. This is of some importance nowadays for the transition of SD videos to HD ones. 
\item Aperture Synthesis for spatial imaging where the corruption process is given by a mask whose typical example is shown in Fig.~\ref{Atom_mask}. (For further details see \cite{Johnson}.)
\item The tomography problem for medical imaging or seismic imaging. In this case Fourier coefficients usually lie on straight lines that are either parallel or that cross at the origin.
\item %The scattering theory, which is concerned with the effects of static objects on traveling waves. 
The inverse scattering problem where one is interested in recovering the shape of an hidden object using electromagnetic or acoustic waves.
\end{itemize}
%\vspace{-0.2cm}
Let us consider the latter example (see \cite{Colton03} for details). 
%As an exemple, we decide to detail the latter application. The topic is the object of a recent survey \cite{Colton03} we refer to for the details. Scattering theory is concerned with the effects of static objects on traveling waves.
Let an \textit{incident acoustic plane wave} $u^i(x)=e^{ik x\cdot d}$ propagate in the direction of the unit vector $d$ in an isotropic medium. Here $k>0$ is the \textit{wave number}. In case there is an inhomogeneity $D$ (the hidden object), the wave will be ``scattered'' and give rise to another wave $u^s$. 
% The resulting wave $u=u^i+u^s$ satisfies in $\R^N$ the Helmholtz equation
% \begin{align}
% \Delta u +k^2(1-\chi_D)u=0 
% \end{align}
% and the \textit{scattered wave} $u^s$ satisfies the \textit{Sommerfeld radiation condition}
% \begin{align}
% \lim_{r\to +\infty} r\left(\frac{\partial u^s}{\partial r}-iku^s\right)=0,
% \end{align}
% where $r=|x|$ and the limit is uniform in the direction $\hat{x}=x/|x|$.
The latter has an asymptotic behavior
\begin{align}
u^s(x)=\frac{e^{ik|x|}}{|x|}u_\infty(\hat{x},d)+O\left(\frac{1}{|x|^2}\right),
\end{align}
where $u_\infty(\hat{x},d)$ is known as the \textit{far field pattern} and gives an idea of the behavior of the scattered wave at large distance. Here, $\hat{x}$ is the observation direction and $d$ is the incident wave direction. The \textit{direct problem} amounts to find $u_\infty(\hat{x},d)$ whereas the \textit{nonlinear inverse problem} takes the direct method as a starting point and asks what is the nature of the scatterer $D$ that gave rise to such a farfield. %This problem is well-posed since the knowledge of $u_\infty(\hat{x},d)$ determines uniquely the obstacle $D$. A significant development in this field is the introduction of the \textit{factorization method}. 
In what follows, we assume that $k$ is sufficiently small so that the \textit{Born approximation} implies that
\begin{align}
\label{Born_approx}
u_\infty(\hat{x},d)\approx\int_{\mathbb{R}^N}\chi_D(y)e^{-ik(\hat{x}-d)\cdot y} dy.
\end{align}
This means that $u_\infty(\hat{x},d)$ can be interpreted as Fourier coefficients in the ball of radius $2k$.
% \begin{align}
% u_\infty=\F(\chi_D)\chi_{B(0,2k)}.
% \end{align}
In practice, having only limited incident waves and measurements we end up with a sampling of the spectrum of $\chi_D(y)$ that is depicted in Fig.~\ref{spectrum_scatt}.
%\vspace{-0.3cm}
% \newpage
\vspace{-0.3cm}
%\textbf{In \cite{} the authors address the problem when these coefficients lie all along a manifold.}
In general non regular sampling of the spectrum is a difficult problem (see for instance \cite{Meyer09}). Though if $D$ is bounded, its Fourier transform is a $C^\infty$ function that can be interpolated on a uniform grid. This allows to use the \verb+fft+. %This is sometimes referred to as the \textit{regridding} procedure.

In \cite{Guichard,Moisan}, the authors consider the Total Variation (henceforth denoted $TV$) to recover missing Fourier coefficients. %The problem of recovering missing Fourier coefficients was also tackled in \cite{Moisan} where a weighted $TV$ term was considered.
Let us compare the performance of this approach with that of the factorization method \cite{Kirsch} thanks to a numerical test. We consider the object $D$ that is bounded by the red curves in Fig.~\ref{scatterer}. 
% \vspace{-0.4cm}
The inverse problem is then solved given 32 incident waves and 32 measurement directions%\footnote{I would like to thank Armin Lechleichter who provided me with the solution of the direct problem and the code for the factorization method.}
. The wave number is $k=3\pi$. The factorization method, which performs quite well for small values of $k$, yields quite poor results in this extreme case:%The respective results are depicted in Fig.~\ref{scat_result}.
\vspace{-0.3cm}
\begin{figure}[H]
\centering
\begin{minipage}[c]{\linewidth}
\centering
\includegraphics[height=0.32\linewidth]{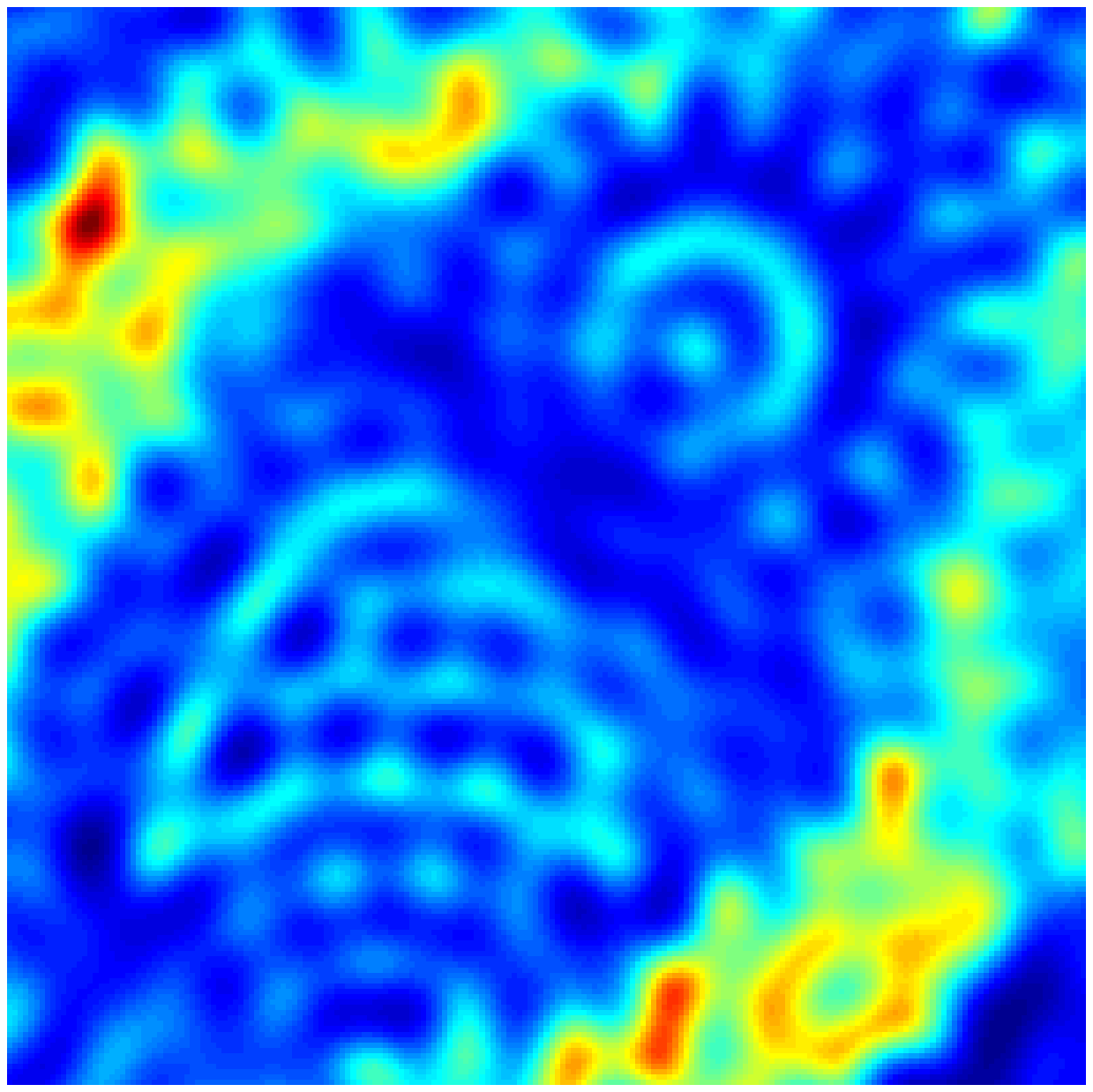}%trim=l b r t
\includegraphics[height=0.32\linewidth,clip=true,trim=50 50 50 55]{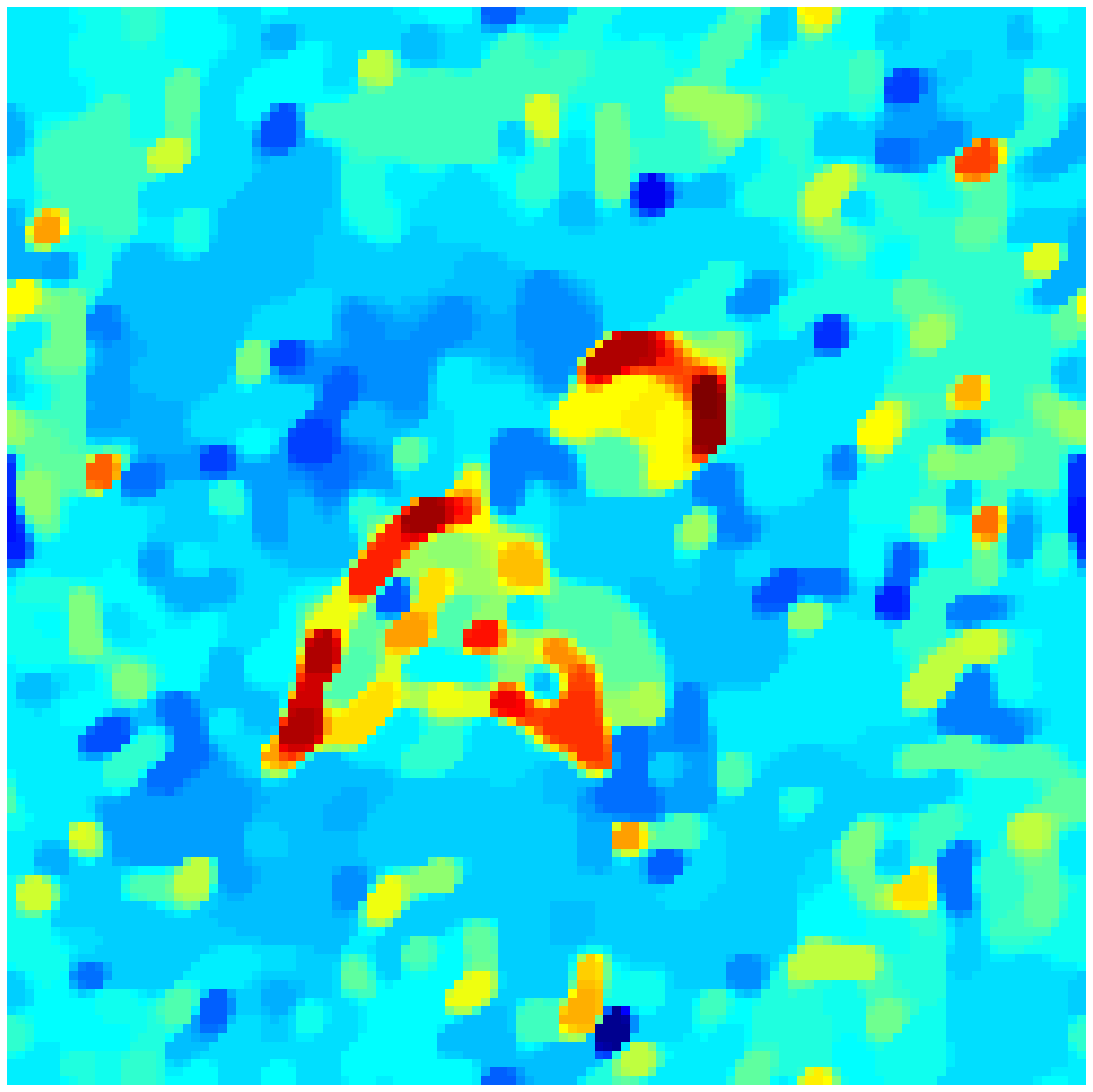}%trim=l b r t
\includegraphics[height=0.32\linewidth,clip=true,trim=50 50 50 55]{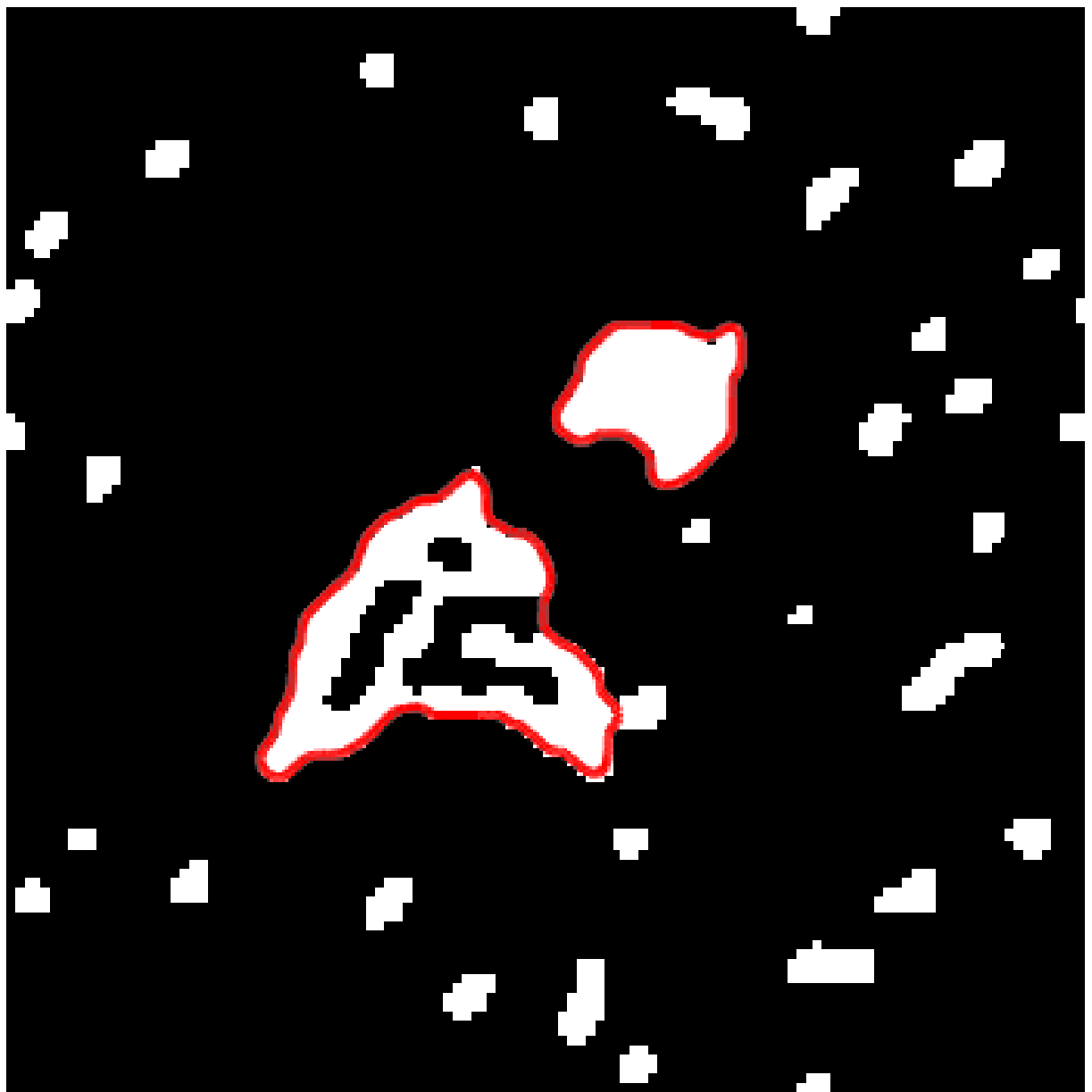}%trim=l b r t
\vspace{-0.1cm}
\caption%[]
% {From left to right, scatterer obtained by the \textit{factorization method} \cite{Kirsch} which is then thresholded, scatterer obtained by minimizing the Total Variation which is then thresholded. The final result is segmented \cite{Chan01} (red curve on the right).}
	{From left to right, scatterer obtained by the \textit{factorization method} \cite{Kirsch}, scatterer obtained by Total Variation minimization, the result is then thresholded and segmented \cite{Chan01} (red curve on the right).}
% \vspace{0.5cm}
\label{scat_result}
\end{minipage}
\end{figure}
 \vspace{-0.3cm}

The $TV$-based method yields really good results in this case. Indeed, from partial and noisy measurements, we were able to distinguish two objects that are merely separated by $3/4\la$ where $\la=2\pi/k$ is the \textit{wavelength}. This value is close to the theoretical diffraction limit that is $\la/2$. Note also it is not necessary to tune a Lagrange multiplier as it is usually done. Indeed, if one denotes $M$ the discrete set of points where $u_\infty$ is known, one can simply minimize the constrained problem
\begin{align}
\min_{{\F u_{|M}}=u_\infty}TV(u),
\end{align}
to get the restoration shown in Fig.~\ref{scat_result}.
Total Variation was already used for solving the nonlinear inverse scattering problem in \cite{Vandenberg}. The idea of considering the inverse problem within the Born approximation as an inpainting problem in the Fourier domain is discussed in \cite{Devaney}.\\
% \newpage

This example shows that using variational methods for the general spectrum interpolation problem may yield quite satisfactory results. We should also mention the not so different problem of restoring missing wavelet coefficients which is treated in \cite{Chan06,Chan09}. For such problems, can we do better than using $TV$? Some recent papers address the problem of solving general inverse problems in imaging by means of non-local methods. These methods obviously apply to the problem of retrieving missing Fourier coefficients.

\section{Prior Works}
Traditional methods in image processing are based on local properties of images (Wavelets, Total Variation). Recently, state of the art results were achieved for the denoising problem by Non-Local methods that exploit redundancies in images. The idea of using self-similarities in images was first exploited in \cite{Buades}. They proposed a filter that averages similar patterns of a noisy image $g=g_0+n$ defined on $\O\subset\R^2$ even though these self-similarities occur at large distance. The filter reads as follows
\begin{align}
\text{NLMeans}(g)(x)=\frac{1}{C(x)}\int_{\O} g(y)\exp\left(-\frac{{\|p_g(x)-p_g(y)\|}_2}{h}\right)dy,
\end{align}
where $p_g(x),\ p_g(y)$ are patches of $g$ and $h$ determines the selectivity of the similarity measure. 
% The idea of looking at such a term dates back to the filter 
% % \begin{align}
% % Y(g)(x)=\frac{1}{C(x)}\int_\O g(y)\exp\left(-\frac{\|g(x)-g(y)\|}{h}\right)dy
% % \end{align}
% proposed in \cite{Yaroslavsky} which merely compares gray levels.
Many modifications of this filter were considered for denoising purposes (adaptive $h$, adaptive window \cite{Kervrann}, %non-square patches \cite{Duval}, 
shape adaptive patches \cite{Dabov}).

% Though, an important drawback of the NLMeans filter is that the image is not processed if no similar patches are found. This phenomenon is referred to as the \textit{rare patch effect}.  To remedy this artifact, the authors of \cite{Louchet10} propose to mix local and non-local methods by preprocessing the patches with the total variation. The regularization parameter being set locally so that one gets enough similar patches.

The use of patches has been widely adopted in the image processing community and in the recent years, these non-local methods were extended to the study of general inverse problems (see \cite{Kindermann,Gilboa06,Gilboa,Peyre08,Faccolio}). Such a general inverse problem can be modeled as follows
%\begin{align}
$g=Ag_0+n.$
%\end{align}
Here $g_0$ is the original image defined on $\O\subset\R^2$, $A$ is a linear transformation and $n$ is a white Gaussian noise. The classical local methods were adapted by considering a \textit{pointwise Non-Local regularization} as introduced in \cite{Gilboa06} 
\begin{align}
J_w(u)=\int_{\O\times\O}{|u(x)-u(y)|}^\alpha w(x,y)\eta(x-y)dxdy
\end{align}
for some $\alpha\geq 1$ and where $w(x,y)=\exp\left(-\frac{{\|p_g(x)-p_g(y)\|}_2}{h}\right)$
is the weight function used in the NLMeans filter and is based on the Sum of Squared Differences (henceforth called SSD distance). As for $\eta$, it is a function of compact support centered at the origin. It indicates that patches that are too far from each other should not be taken into account. Indeed, it was observed that the Non-Local methods yield much better results if one seeks for similar patches in the support of $\eta$, referred to as the \textit{search window}
%% . In the end these methods are not fully Non-Local % and $h$ plays the role of a tuning paramter that 
(see also assumption $(ii)$ in Section \ref{num_ex_NL}). % and $C(x)$ is merely a normalizing factor.

Peyré et al. proposed in \cite{Peyre09} (see also \cite{Faccolio}) a \textit{patchwise Non-Local regularization}
\begin{align}
J_w(u)=\int_{\O\times\O}{{\|p_u(x)-p_u(y)\|}_2} w(x,y)\eta(x-y) dxdy\,.
\end{align}
%that is able to reconstruct entire patches. 
To get a restored image one then has to minimize the following Non-Local energy
\begin{align}
\E(u)=\frac{1}{2}{\|Au-g\|}_2^2+\la J_w(u).
\end{align}
% where for some $\alpha\geq 1$ the regularization term is given by
% \begin{align}
% J_w(u)=\int_{\O\times\O} {|u(x)-u(y)|}^\alpha w(x,y)\eta(x-y)dxdy
% \end{align}
It is further remarked that one can enhance the restoration by recomputing the weight $w$ regularly. % to improve the overall restoration.
% In these papers it was also suggested that the patch distance should be recomputed on the processed image  
In \cite{Faccolio,Peyre09}, the authors proposed a general framework where the weight $w(x,\cdot)$ is interpreted up to renormalization as a density of probability and is an unknown of the problem. Then the energy to be minimized
% \begin{align}\label{Energie_NL}
% \E(u,w)=\frac{1}{2}{\|Au-g\|}_2^2+\la J_w(u)-h^2\int_{\O\times\O} w\log w
% \end{align}
% which
 involves the potential energy of $w(x,\cdot)$, used to infer unknown probability distributions.
The resulting energy is not convex in $w$ and an alternate coordinate descent gives back the SSD-based weight but computed this time on the $u$ being processed. This actually corresponds to the aforementioned weight recomputation procedure.

\section{Contribution}

In this paper, we deal with the problem of restoring missing Fourier coefficients of an image. In particular we propose a technique to build a distance insensitive to the degradation of an image. The corruption process is supposed to be known. We show that using this distance between patches gives very good reconstruction results by solving a simple quadratic energy.

In the classical Non-Local approaches, a critical step is to compute the SSD distance between two patches of the corrupted image. This strategy dates back to the NLMeans~\cite{Buades}.% and to our knowledge no other similarity measure has been considered since then.
We explain how to replace this step by the computation of a simpler $\ell^2$-type distance that is really adapted to the problem of spectrum reconstruction and does not incorporate spurious information. The idea is to define atoms similar to Gabor filters to test whether two regions of the corrupted image are similar. These atoms should be as concentrated as possible and should not depend on the corruption process. This way two patches that are close in the clean image will be close in the corrupted image. We now first introduce our variational framework, which is relatively standard.

\section{A Non-Local Energy for the Problem}
Henceforth, we work in $\R^2$ which will is equipped with the 
norm $\max\{|x_1|,|x_2|\}$ so that ``balls'' are in fact squares.
%. Considering squares instead of balls can be handy in image processing.

\begin{figure}[H]
\centering
%\begin{minipage}[c]{\linewidth}
\begin{tabular}{cc}
Original image & Corrupted image\\
 \includegraphics[width=0.25\linewidth]{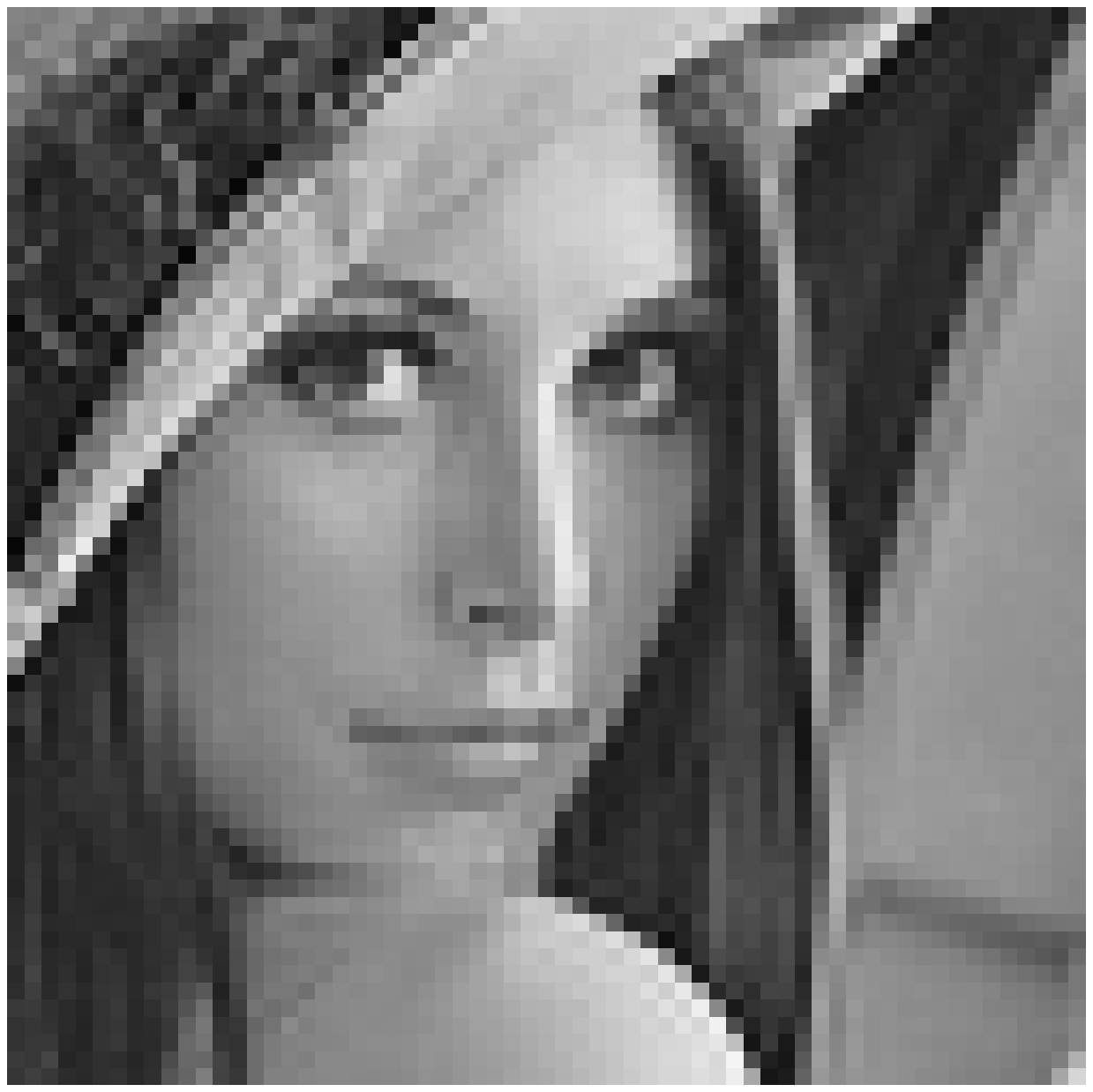}
&\includegraphics[width=0.25\linewidth]{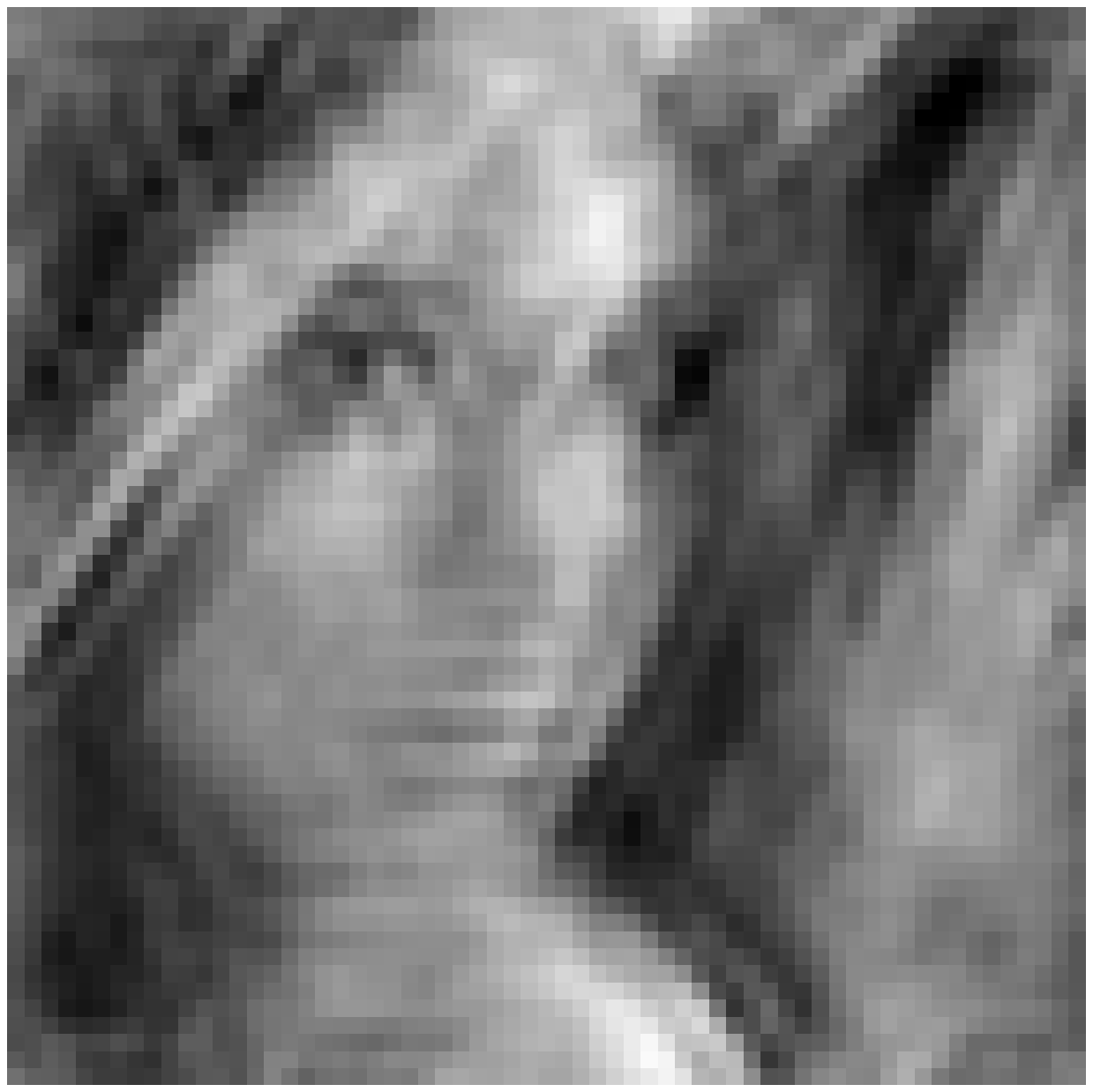}\\
 &PSNR=20.0dB
\end{tabular}
\begin{tabular}{ccc}
$\ell^2$ Oracle ($\alpha=1$)& $\ell^2$ Oracle ($\alpha=2$)&TV restored\\
%\includegraphics[width=0.25\linewidth]{perfect_lena_NL_7}
%&\includegraphics[width=0.25\linewidth]{perfect_lena_NL_7}
\includegraphics[width=0.25\linewidth]{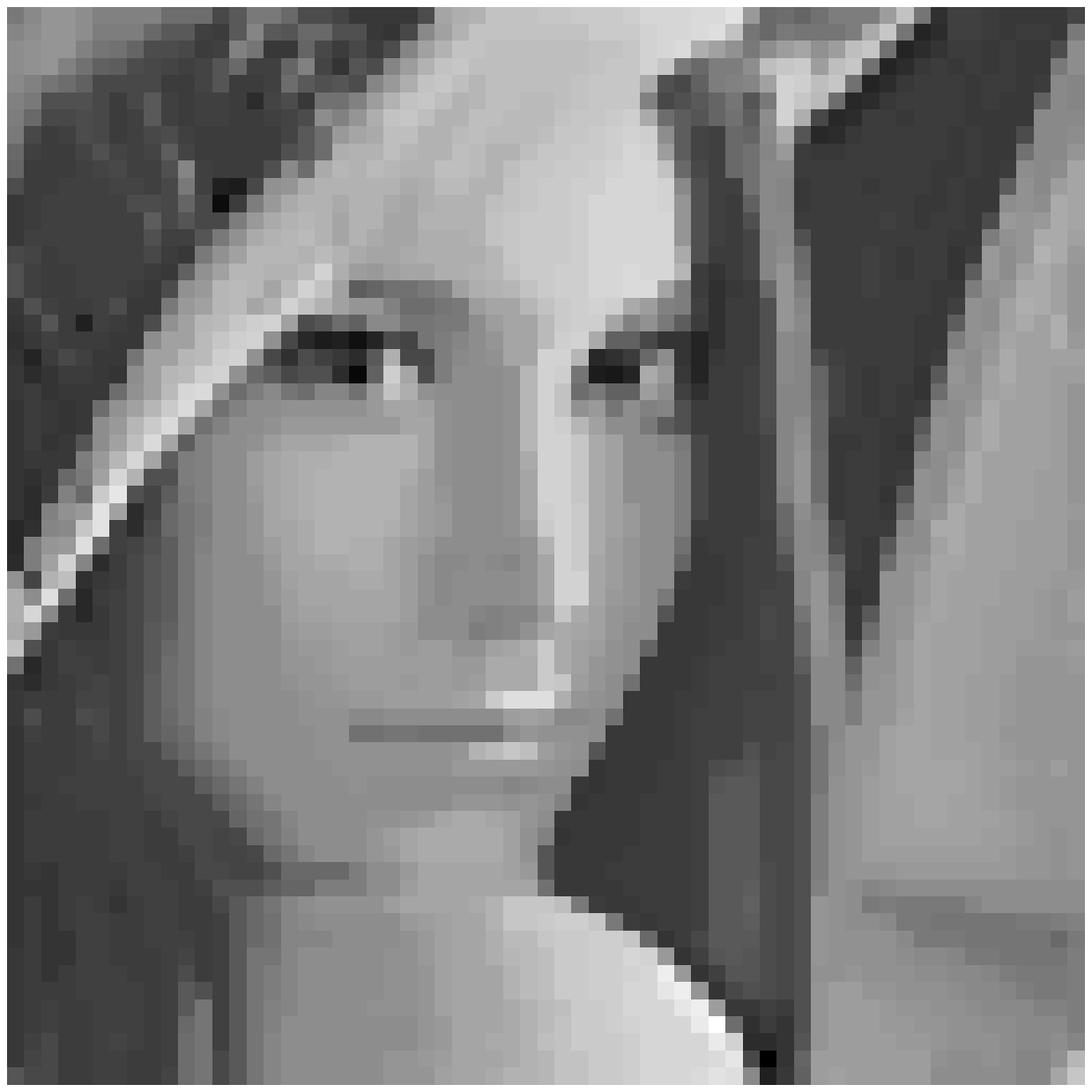}
&\includegraphics[width=0.25\linewidth]{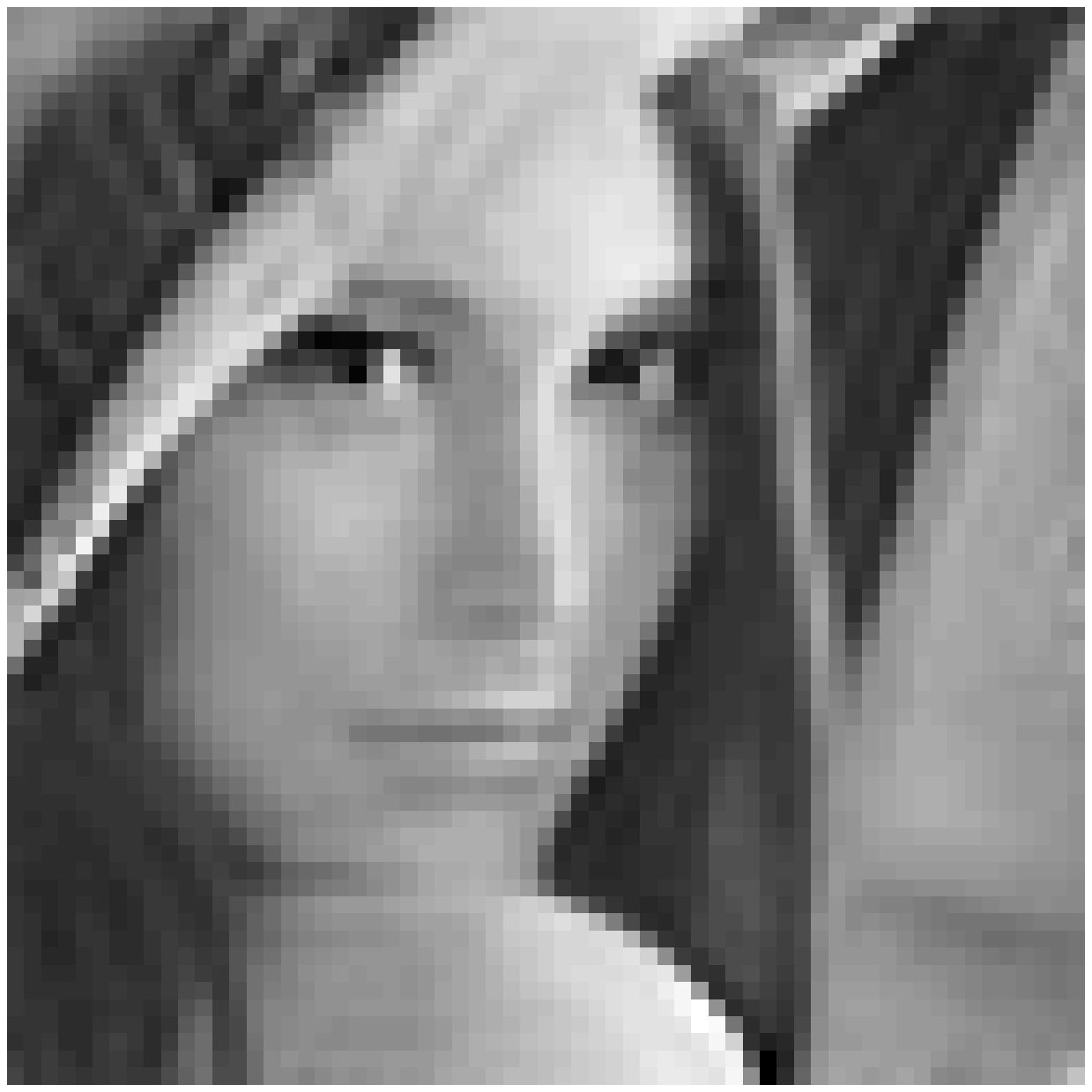}
&\includegraphics[width=0.25\linewidth]{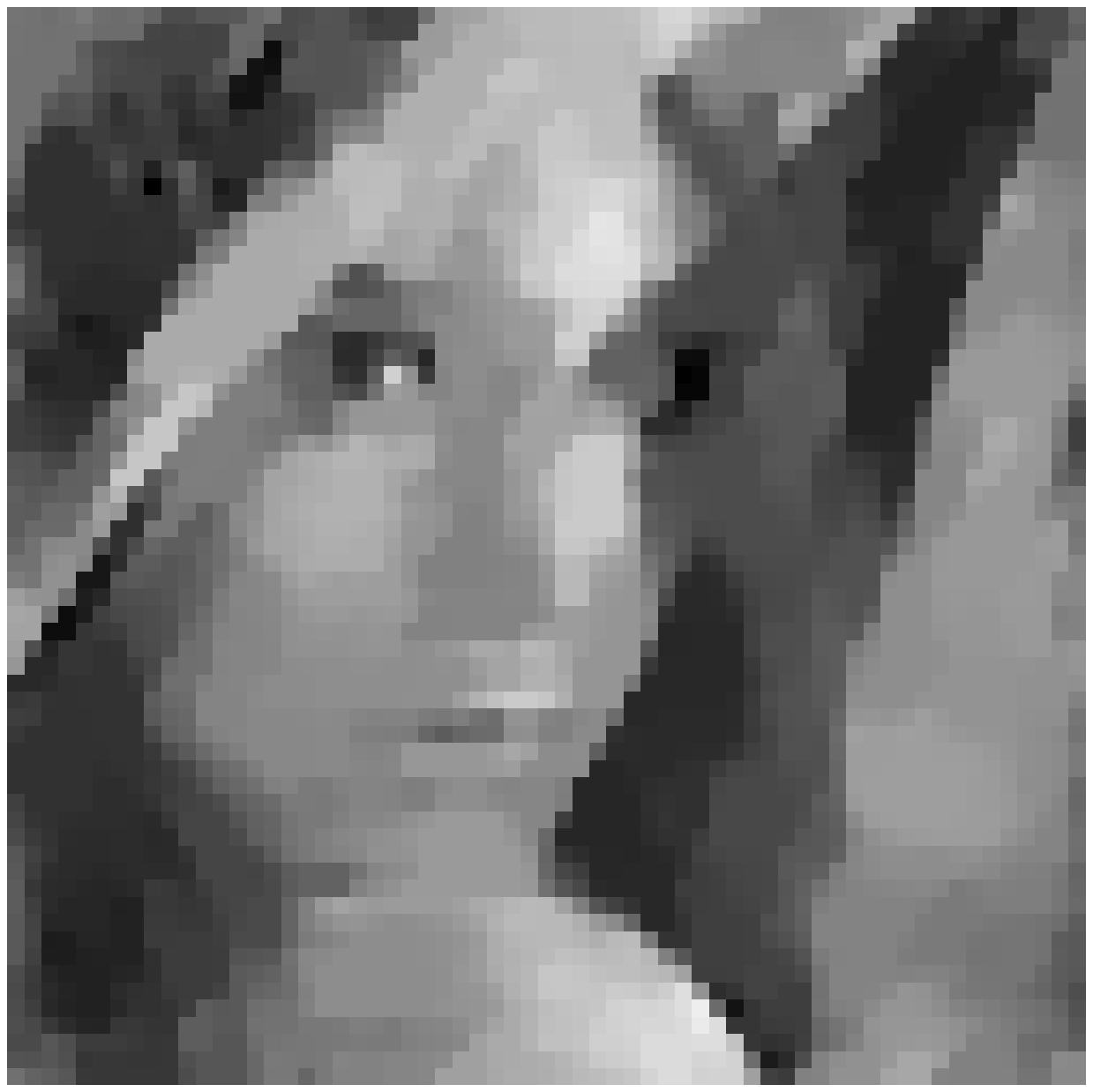}\\
 PSNR=25.6dB& PSNR=26.2dB& PSNR=20.9dB\\\\
\end{tabular}
% \end{minipage}
% \begin{minipage}{[c]{\linewidth}}

\begin{tabular}{cc}
Original image & Corrupted image\\
 \includegraphics[width=0.25\linewidth]{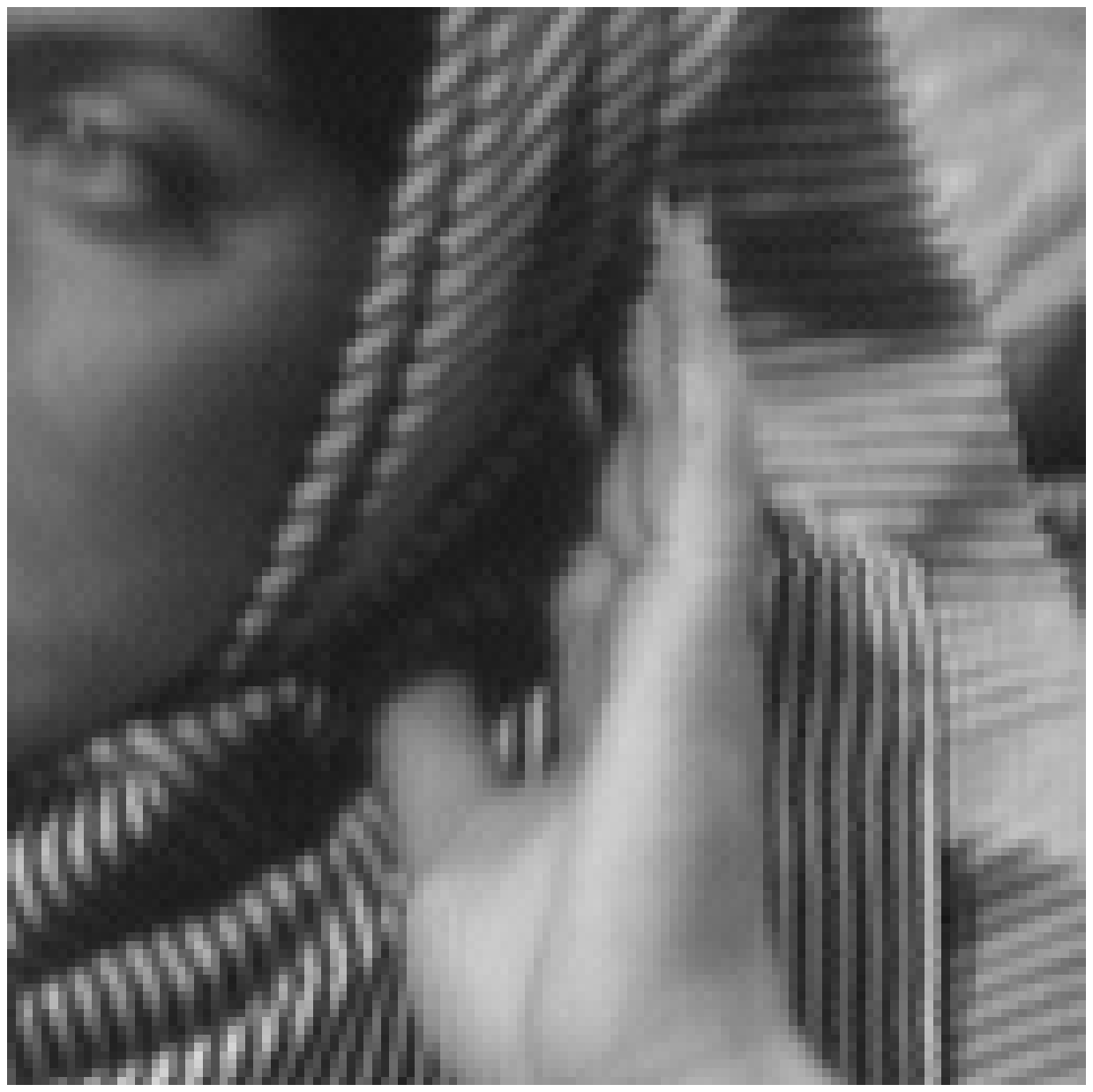}
&\includegraphics[width=0.25\linewidth]{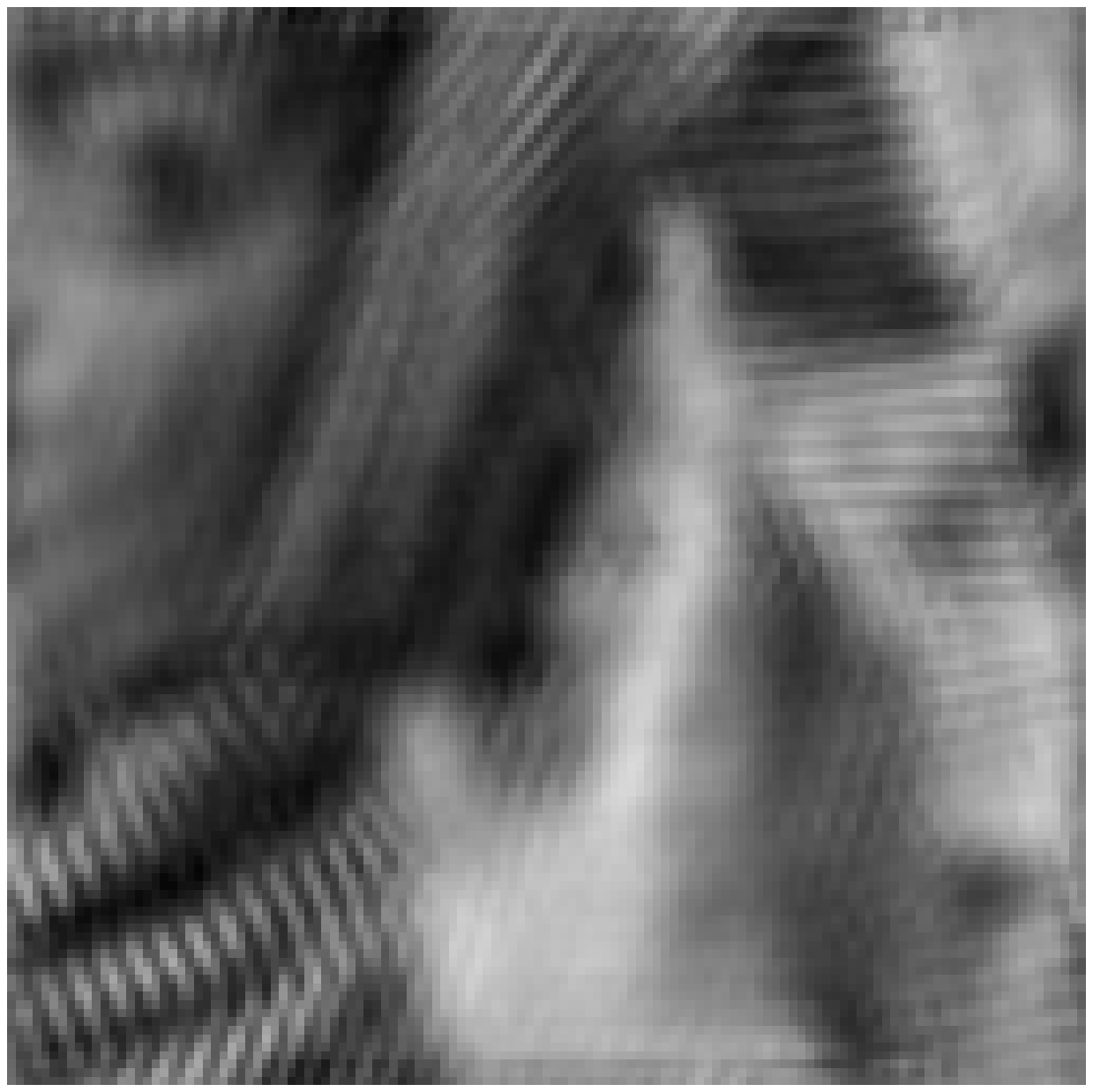}\\
 &PSNR=24.1dB
\end{tabular}
\begin{tabular}{ccc}
$\ell^2$ Oracle ($\alpha=1$)& $\ell^2$ Oracle ($\alpha=2$)&TV restored\\
%\includegraphics[width=0.25\linewidth]{perfect_barbara_NL_7}
%&\includegraphics[width=0.25\linewidth]{perfect_barbara_NL_7}
\includegraphics[width=0.25\linewidth]{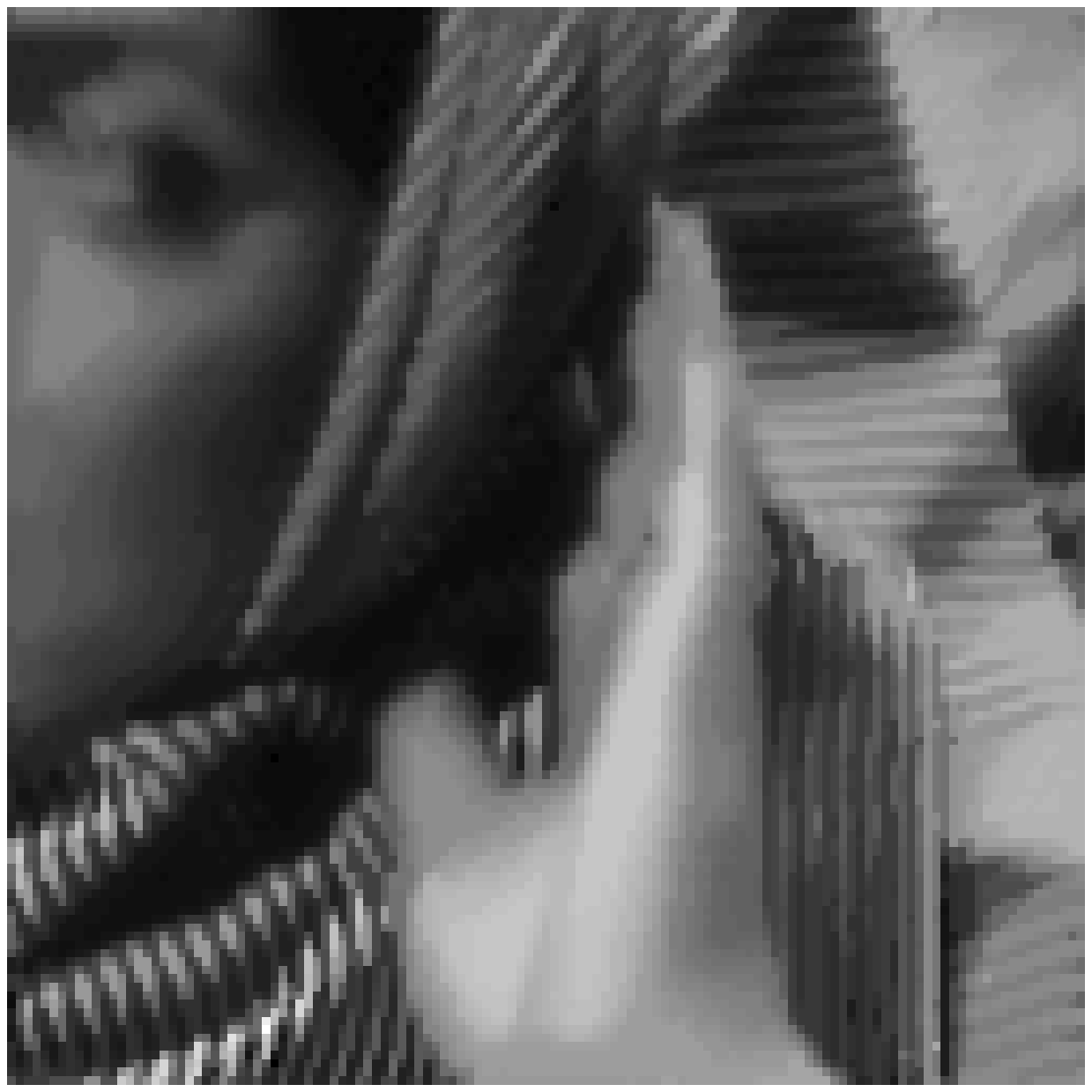}
&\includegraphics[width=0.25\linewidth]{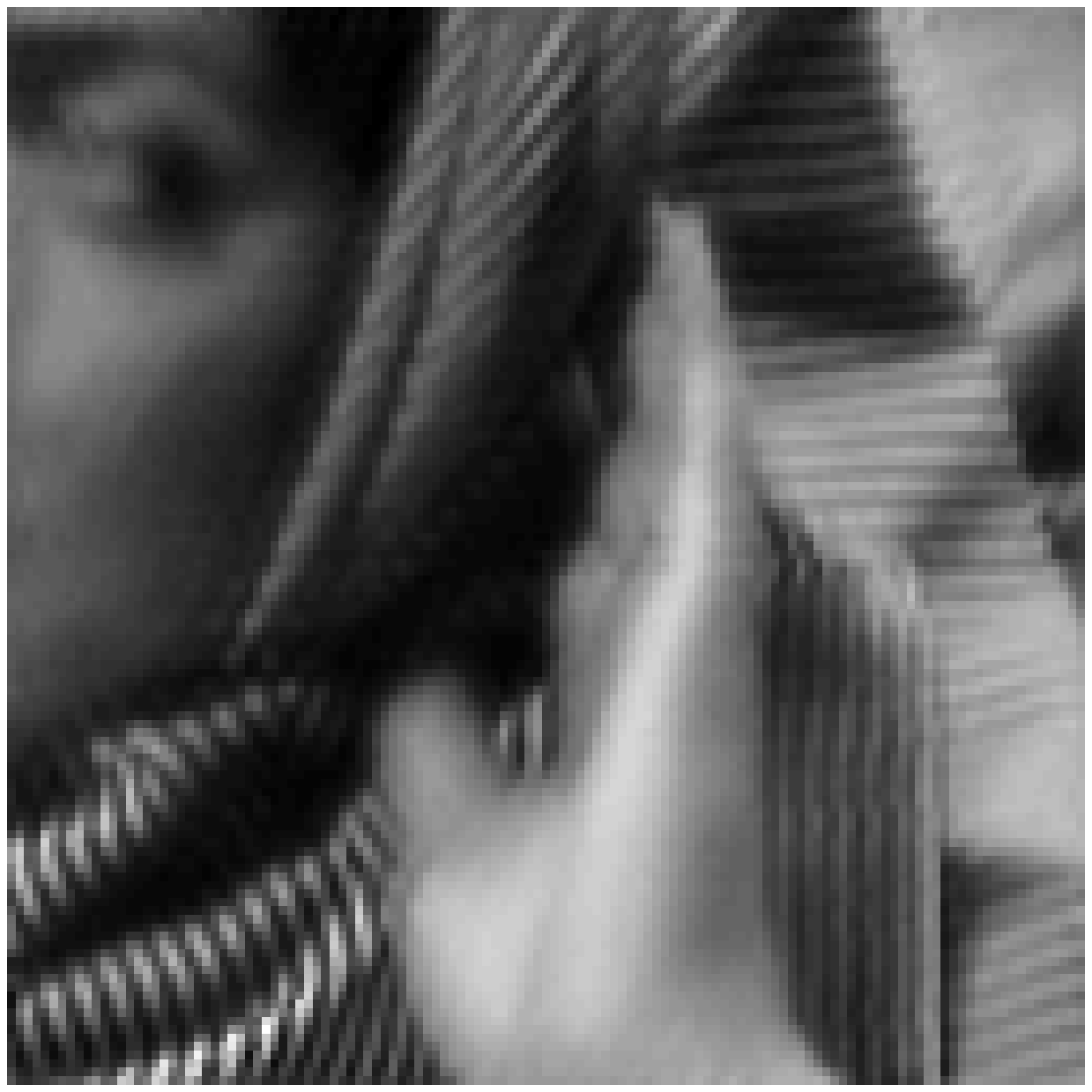}
&\includegraphics[width=0.25\linewidth]{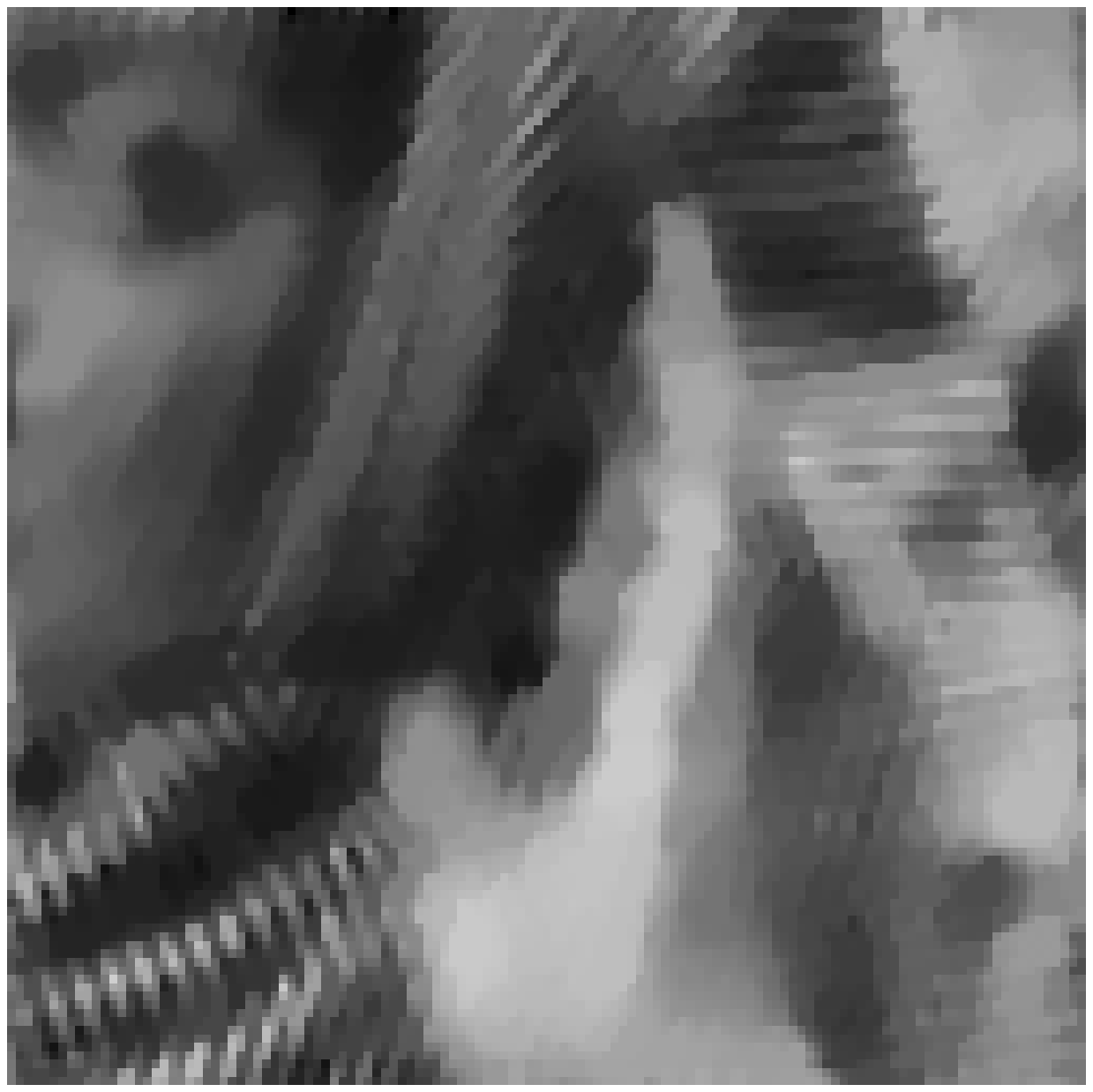}\\
 PSNR=28.9dB& PSNR=29.5dB& PSNR=24.5dB
\end{tabular}
%\begin{tabular}{cccc}
%Original image & Corrupted image & Oracle weight& TV restored\\
 %\includegraphics[width=0.25\linewidth]{perfect_barbara_NL_4}
%&\includegraphics[width=0.25\linewidth]{perfect_barbara_NL_6}
%&\includegraphics[width=0.25\linewidth]{perfect_barbara_NL_7}
%&\includegraphics[width=0.25\linewidth]{perfect_barbara_NL_8}\\
%&PSNR=24.1dB&  PSNR=29.8dB&  PSNR=24.5dB
%\end{tabular}
\vspace{-0.2cm}
\caption%[]
{Non-Local restoration thanks to the ideal distance. \newbf{Both cases $\alpha=1$ and $\alpha=2$ are considered.}}
% \end{minipage}
\label{simul_ideal}
\end{figure}

In this section, we first formalize the problem of reconstructing unknown coefficients of a Fourier series. Let us call ${M}\subset\Z^2$ the finite mask of points where the Fourier coefficients are known. \newbf{We shall assume that it is symmetric with respect to the origin.} Assuming that the image is periodic and defined on $\T=[0,1]^2$, one thus considers respectively the clean image and the corrupted images
\begin{align}
g_0(x)=\sum_{k\in\Z^2}c_k e^{-2i\pi k\cdot x},\  
g(x)=\sum_{k\in M}c_k e^{-2i\pi k\cdot x}.
\end{align}
We assume that not only the high frequencies but also middle range frequencies are lost in the corruption process.
%Our aim is to capture low frequencies and \textit{some} middle range frequencies but we discard the high frequencies. 
In other words, we keep the Fourier coefficients corresponding to the white
areas of the  mask $M$. The typical $M$ we will consider here
is the following
\vspace{-0.3cm}
\begin{figure}[H]
%    \begin{minipage}[c]{\linewidth}
\centering
     \includegraphics[width=0.35\hsize]{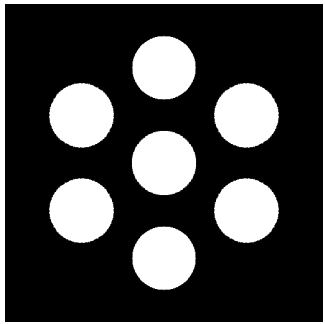}
\vspace{-0.1cm}
    \caption{Corruption mask $M$.}

% \end{minipage}
\label{Atom_mask}
\end{figure}
% \textbf{Give an example of degraded image}
\vspace{-0.3cm}
for which usual \textit{local} TV-based reconstruction techniques
totally fail at recovering the middle-range frequencies,
see for instance Fig.~\ref{simul_ideal}, right.

In the sequel, we introduce the two
subsets of $L^2$: $\M=\{\sum_{k\in M} c_k e^{-2i\pi k\cdot x}\}$
and its complement $\M^\perp=\{\sum_{k\not\in M} c_k e^{-2i\pi k\cdot x}\}$
where in both cases $(c_k)_k$ runs over all sequences in $\ell^2(\Z^2;\R)$,
which are even in $k$.
Let us denote $P_\M$ the orthogonal projection on $\M$. We recall that the Fourier transform on the torus $\F:L^2(\T)\to\ell^2(\Z^2)$ is an isometry and thus invertible so the corrupted data $g\in\M$ that one has to interpolate is obtained by
\begin{align}
g=\F^{-1}\circ P_{\mathcal{M}}\circ\F(g_0)=\F^{-1}(\chi_M\F(g_0)),
\end{align}
%The corruption process results into a perturbation that is similar in nature to the Gibbs phenomenon. %Note also that $g$ is real-valued as soon as $M$ is symmetric with respect to the origin.

We shall call \textit{patch} centered at $x_k\in\R^2$, denoted $p_k(g)$, the image $g(\cdot-x_k)\vphi$ where $\vphi$ is a test function with support $C(0,\frac{\rho}{2}):={\left[-\frac{\rho}{2},\frac{\rho}{2}\right]}^2$ and $\rho$ is the \textit{patch size}. 

For the moment let us assume that the original image is such that there are two distinct $x_k,x_\ell\in\R^2$ with
% \begin{align}
$g_0(x-x_k)\vphi(x)=g_0(x-x_\ell)\vphi(x),\ \forall x\in\R^2.$
% \end{align}
This is to say that two patches are similar in the original image. Then clearly,
% one expects that
 the \newbfa{part of the image $v$ corresponding to the missing frequencies}
is among the solutions of
\begin{align}\label{E_k_ell}
\min_{v\in\M^\perp}\int_\T\psi(x)^2|(g+v)(x-x_k)-(g+v)(x-x_\ell)|^\alpha dx
\end{align}
where $\psi$ is smooth and such that $\supp(\psi)\subset\supp(\vphi)=C(0,\frac{\rho}{2})$, and $\alpha\ge 1$ is fixed. %Note that we do not necessarily assume that $\vphi=\psi$.
This means that knowing that $p_k(g_0)$ and $p_\ell(g_0)$ are similar one can hope to get a restored spectrum by minimizing (\ref{E_k_ell}). The reconstruction is obviously not unique since modifying $v$ out of $(x_k+\supp(\psi))\cup(x_\ell+\supp(\psi))$ does not change the energy (and numerical simulations actually show that the reconstructed $v$ has the same support as $(x_k+\supp(\psi))\cup(x_\ell+\supp(\psi))$).
Thus, we have to take into account all the patches in the image to get a global reconstruction.

We are therefore led to consider the following problem
\begin{align}\label{E_NL}
\min_{v\in\M^\perp}
\int_{\T}\sum_{(k,\ell)\in I}\psi^2|(g+v)(x-x_k)-(g+v)(x-x_\ell)|^\alpha w(x_k,x_\ell)\,dx
\end{align}
where $w(x_k,x_\ell)=\exp\left(-\frac{\delta(x_k,x_\ell)}{h}\right)$ and $\delta(x_k,x_l)$, that is going to be defined in the sequel, should tell us from the corrupted image $g$ whether one had for the original image 
%\begin{align}
$p_{g_0}(x_k)\sim p_{g_0}(x_\ell).$
%\end{align}

\newbfa{\indent There are two interesting values for $\alpha$, namely $\alpha=1$ and $\alpha=2$.  Our first experiment (Fig.~\ref{simul_ideal}) shows 
that the reconstruction with~(\ref{E_NL}) for both these choices is almost
perfect if the weights $w(x_k,x_\ell)$ are derived from the 
``oracle'' distance
$\delta(x_k,x_\ell)={\|p_{g_0}(x_k)-p_{g_0}(x_\ell)\|}_2$ given by the
SSD distance of the patches of the (normally unknown) original image $g_0$.}

%To solve (\ref{E_NL}), the weight computation is a critical step as can be seen in the simulations of Fig.~\ref{simul_ideal} where we considered the \textit{oracle distance} 
% \begin{align}
%$\delta(x_k,x_\ell)={\|p_{g_0}(x_k)-p_{g_0}(x_\ell)\|}_2.$
% \end{align}

%\begin{figure}[H]
%\centering
%% \begin{minipage}[c]{\linewidth}
%\begin{tabular}{cccc}
%Original image & Corrupted image & Oracle weight& TV restored\\
 %\includegraphics[width=0.25\linewidth]{perfect_lena_NL_4}
%&\includegraphics[width=0.25\linewidth]{perfect_lena_NL_6}
%&\includegraphics[width=0.25\linewidth]{perfect_lena_NL_7}
%&\includegraphics[width=0.25\linewidth]{perfect_lena_NL_8}\\
%& PSNR=20.0dB& PSNR=26.5dB& PSNR=21.0dB\\\\
%% \end{tabular}
%% \end{minipage}
%% \begin{minipage}{[c]{\linewidth}}
%% \begin{tabular}{cccc}
%Original image & Corrupted image & Oracle weight& TV restored\\
 %\includegraphics[width=0.25\linewidth]{perfect_barbara_NL_4}
%&\includegraphics[width=0.25\linewidth]{perfect_barbara_NL_6}
%&\includegraphics[width=0.25\linewidth]{perfect_barbara_NL_7}
%&\includegraphics[width=0.25\linewidth]{perfect_barbara_NL_8}\\
%&PSNR=24.1dB&  PSNR=29.8dB&  PSNR=24.5dB
%\end{tabular}
%\vspace{-0.2cm}
%\caption%[]
%{Non-Local restoration thanks to the ideal distance.}
%% \end{minipage}
%\label{simul_ideal}
%\end{figure}
%\vspace{-0.4cm}

\vspace{-0.4cm}

\newbf{From these experiments, it seems clear that both values of $\alpha$ allowed us to produce a cleaner image. Due to the huge size of our problems, we will stick in general
to the case of the exponent $\alpha=2$ which \textit{de facto} excludes
the case of standard Non-Local Total Variation. We will see that this
choice leads to computationally tractable problems with excellent
results, for some masks (such as the one in Fig.~\ref{Atom_mask}),
provided the distance between similar patches is correctly estimated. For the considered mask (Fig.~\ref{Atom_mask}), standard Total Variation yield low quality reconstructions (right).}

%The "oracle" distance, which is an ideal that \newbf{cannot} be used in practice, leads to a very good restoration and shows that the Non-Local model,
%with exponent $\alpha\in\{1,2\}$, performs very well if one is able to define a good measure of similarity for patches.

\section{Construction of Atoms Adapted to the Corruption}
Our aim is to define a similarity measure for patches that is not modified through the corruption process $\F^{-1}\circ P_{\mathcal{M}}\circ\F$. This way, if two points are close in the original image they will remain close in the perturbed image.

To do so, our idea is to find a family of test functions ${(\phi_\a)}_\a$ dense in $\M$, that do not depend on $g$ and such that
\begin{align}
g*\phi_\a=g_0*\phi_\a,\ \forall \a.
\end{align}
The atoms ${(\phi_\a)}_\a$ should obviously depend on the mask $M$. Indeed, for any $g\in \M$,
\begin{align}
\langle g,\phi_\a\rangle&=\langle \F^{-1}(\chi_M \F(g)),\phi_\a\rangle=\langle g,\F^{-1}(\chi_M \F(\phi_\a))\rangle
\end{align}
which means that $\F\phi_\a=\chi_M\F\phi_\a$ and therefore $\supp(\F\phi_\a)\subset M$ for any $\a$. From the uncertainty principle (see \cite{Mallat}), one knows that $\phi_\a$ is not of compact support (since \newbf{M} is)
but  one could still expect it to be localized in space. This can be ensured by defining these functions as minimizers of the following parameterized problems  
\begin{align}
\phi_\a=\argmin\left\{\int_\O{|\phi(x)|}^2{|x|}^p_2dx,\phi\in\M,{\|\phi\|}_2= 1,\phi\perp\Span\{\phi_{\a'},\a'<\a\}\right\}
%\phi_\a=\argmin\left\{\int_\O{|\phi(x)|}^2{|x|}^p_2dx,\phi\in\M,\phi\perp\Span\{\phi_{\a'},\a'<\a\}\right\}
\end{align}
where we minimize the $p$-moments of $\phi$, for some $p>1$.

% \noindent In the end, we define for each $\phi_\a$
% \begin{align}
% \delta_{\a}(x_k,x_\ell)&=|g*\phi_\a(x_k)-g*\phi_\a(x_\ell)|\\
% &=|g_0*\phi_\a(x_k)-g_0*\phi_\a(x_\ell)|.
% \end{align}
% We can then consider a measure of similarity between patches that is of the form 
% \begin{align}
% \delta(x_k,x_l)={\left(\sum_{\a<\a_0}\delta_{\a}(x_k,x_\ell)^2\right)}^{\frac{1}{2}},
% \end{align}
We can then consider a measure of similarity that is of the form 
\begin{align}
\delta(x_k,x_l)={\left(\sum_{\a<\a_0}{|g*\phi_\a(x_k)-g*\phi_\a(x_\ell)|}^2\right)}^{\frac{1}{2}},
\end{align}
where $\a_0$ sets how localized the considered atoms $\phi_\a$ are. %Note that we could have used a general $\ell^q$ norm, $q>1$, but this choice does not seem to improve significantly the final results. 
Now, by definition, the similarity measure $\delta$ is invariant under the corruption process $\F^{-1}\circ P_{\mathcal{M}}\circ\F$. \newbf{In other words, two patches that were identical in the original image $g_0$ may not necessarily match according to the SSD but will match according to the atom-based distance, provided that the atoms are well
localized.}

For the mask $M$ of Fig.~\ref{Atom_mask}, the first atoms are actually really localized and can be used as test functions. 
% \vspace{-0.3cm}
% \begin{figure}[H]
% % \begin{minipage}[c]{\linewidth}
% % \centering
% % \includegraphics[width=\hsize]{Atomes_base_50_moment_20_crop}
% \centering
% \includegraphics[width=0.8\linewidth]{test3_crop}
% 
% \vspace{-0.2cm}
% \caption%[]
% {The 50 first atoms (zoomed in) with $p=4$ adapted to the $128\times128$ mask $M$.}
% % \end{minipage}
% \vspace{-0.6cm}
% \end{figure}
Let us have a closer look at the 7 first atoms, computed for $p=4$,
with their respective spectra:
% \vspace{-0.6cm}
\begin{figure}[H]
% \begin{minipage}[c]{\linewidth}
% \centering
% \includegraphics[width=\linewidth]{first_atoms}
% \vspace{-0.6cm}
% \caption%[]
% {Atoms $\phi_n$, $n=1,...,7$ with $p=2$ adapted to the $128\times128$ mask $M$,}
% \vspace{0.5cm}
% \end{minipage}

% \begin{minipage}[c]{\linewidth}
% \centering
\includegraphics[width=\linewidth]{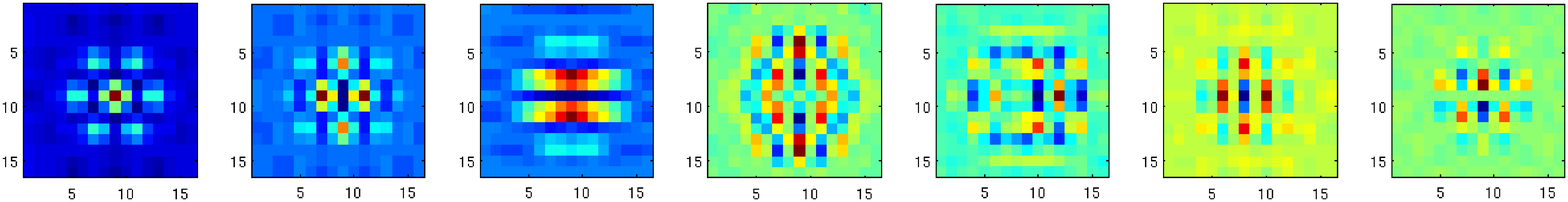}
% \includegraphics[width=\linewidth]{first_atoms_zoom_no_axis}
% \includegraphics[width=\linewidth]{test2_crop}

% \vspace{-0.6cm}
% \caption%[] 
% {Atoms $\phi_n$, $n=1,...,7$, zoomed in,}
% \vspace{0.5cm}
% \end{minipage}

% \begin{minipage}[c]{\linewidth}
% \centering
\includegraphics[width=\linewidth]{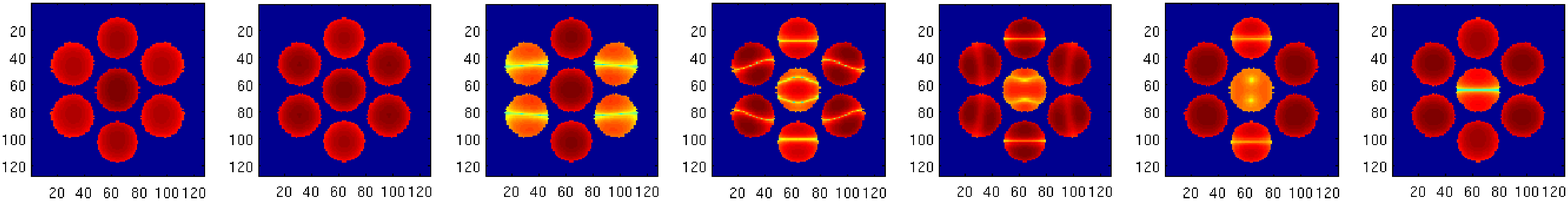}
% \includegraphics[width=\linewidth]{first_atoms_spectra_no_axis}
% \includegraphics[width=\linewidth]{test1_crop}
% \vspace{-0.6cm}
% \caption%[]
% {$\log(|\F(\phi_n)|)$, $n=1,...,7$.}
% \vspace{0.5cm}
% \end{minipage}
% \vspace{-0.4cm}
\caption{%First line: atoms ${(\phi_n)}_{n=1,...,7}$ with $p=2$ adapted to the $128\times128$ mask $M$.
 First line: Atoms ${(\phi_n)}_{n=1,...,7}$ with $p=4$ adapted to the $128\times128$ mask $M$ zoomed in.
 Second line: Their respective spectra $\log(|\F(\phi_n)|)$, $n=1,...,7$.}
\end{figure}

These atoms adapted to the mask $M$ may recall the reader of the Gabor atoms (see \cite{Mallat}). However, they have the advantage of having a prescribed spectrum and being as localized as possible. They are in \newbf{a} sense the optimal functions satisfying these two conditions
(\newbf{see Section~\ref{num_ex_NL} for a discussion on the numerical algorithms that can be used in the discrete setting to compute these atoms.})
\smallskip

\begin{figure}[H]
% \begin{minipage}[c]{\linewidth}
\centering
\includegraphics[width=0.18\linewidth]{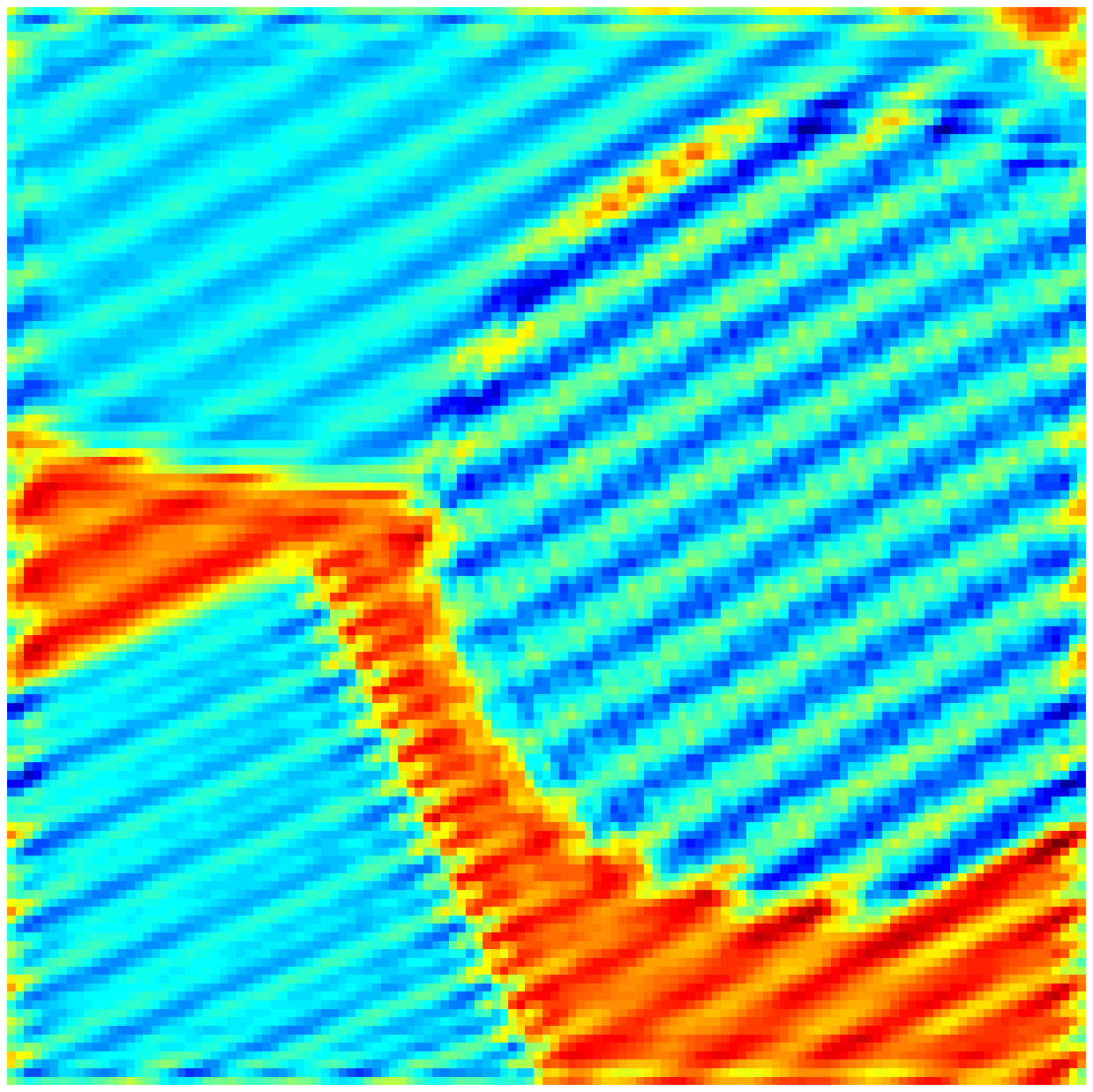}
\includegraphics[width=0.18\linewidth]{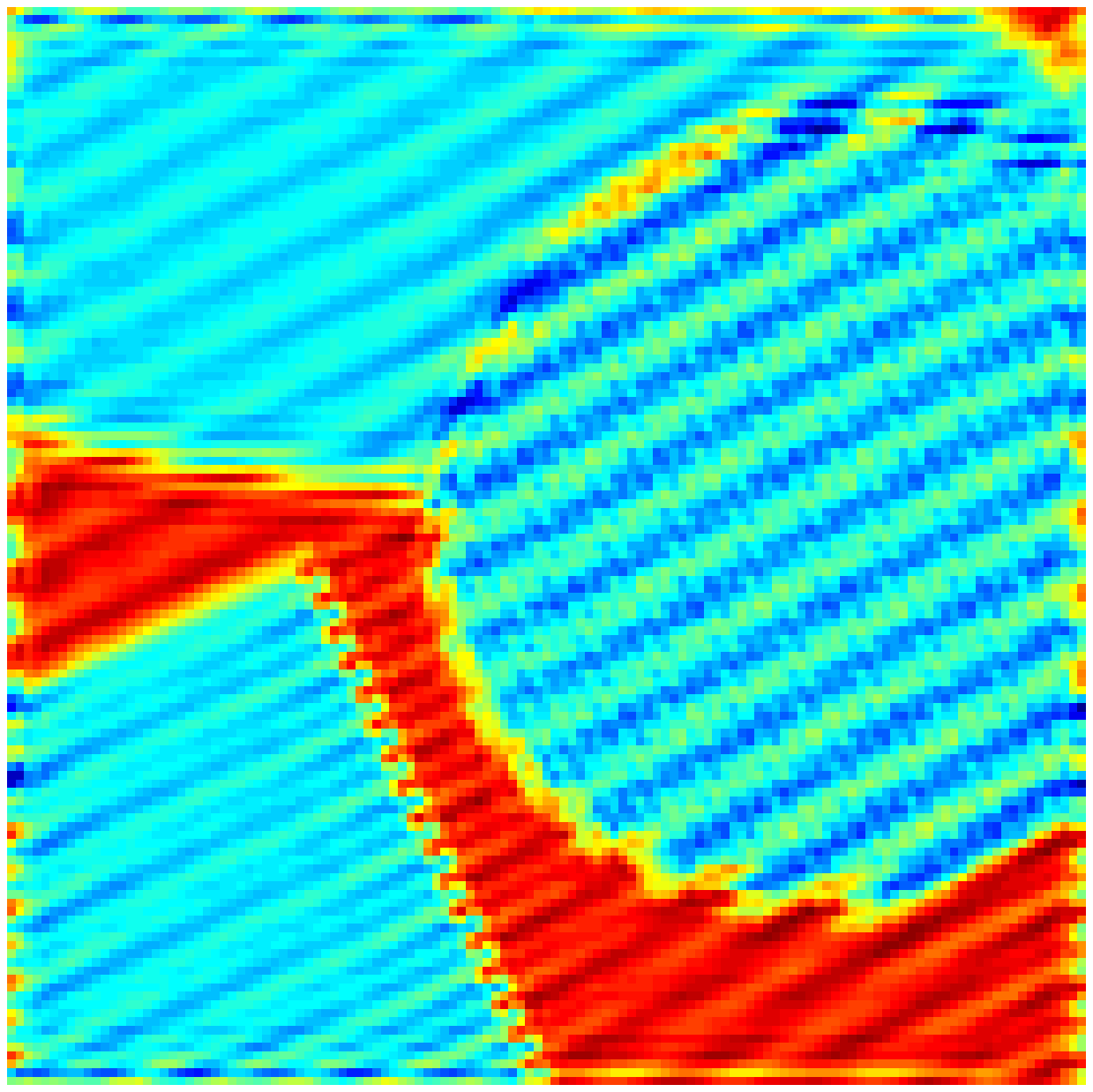}
\includegraphics[width=0.18\linewidth]{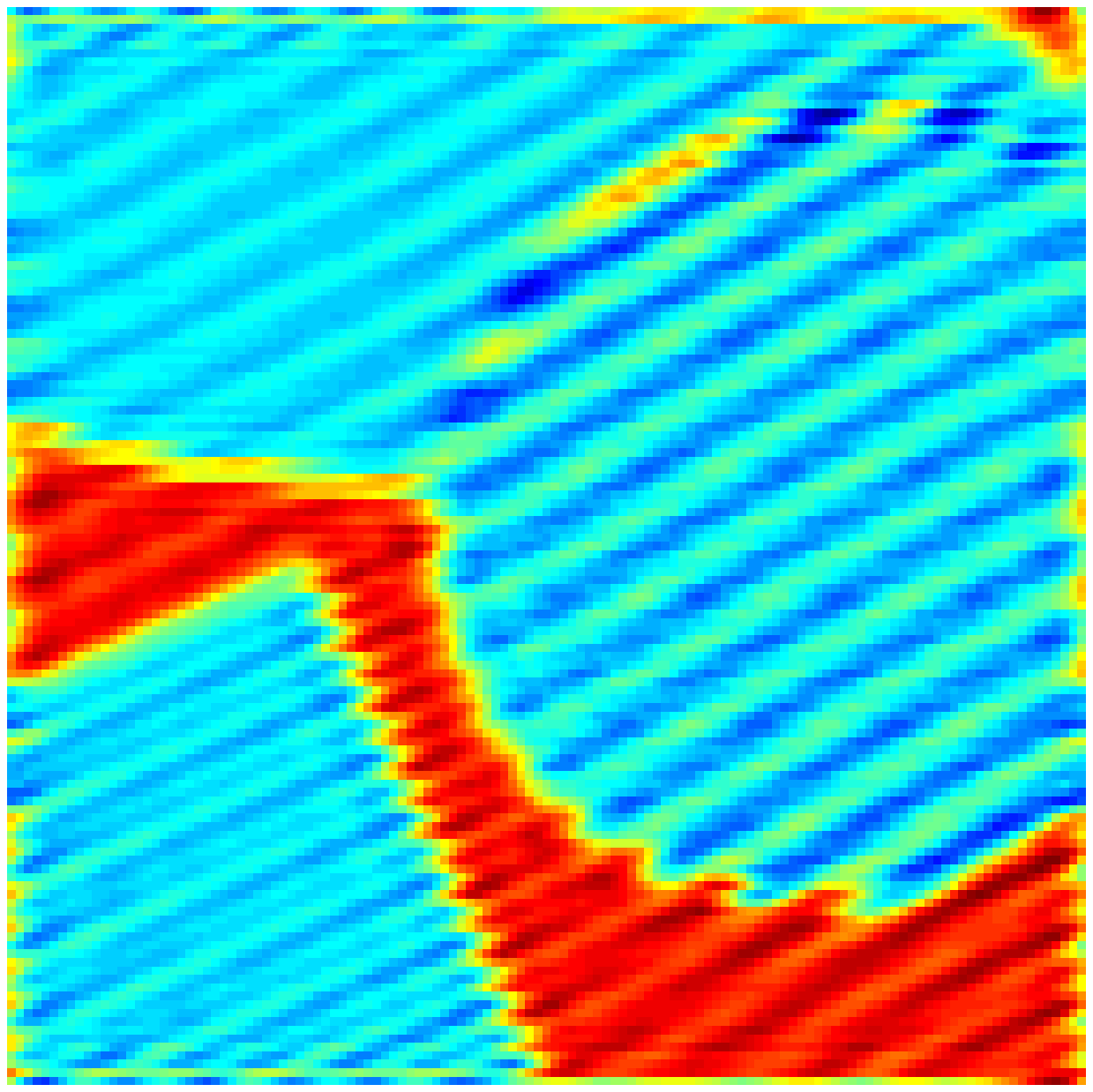}
\includegraphics[width=0.18\linewidth]{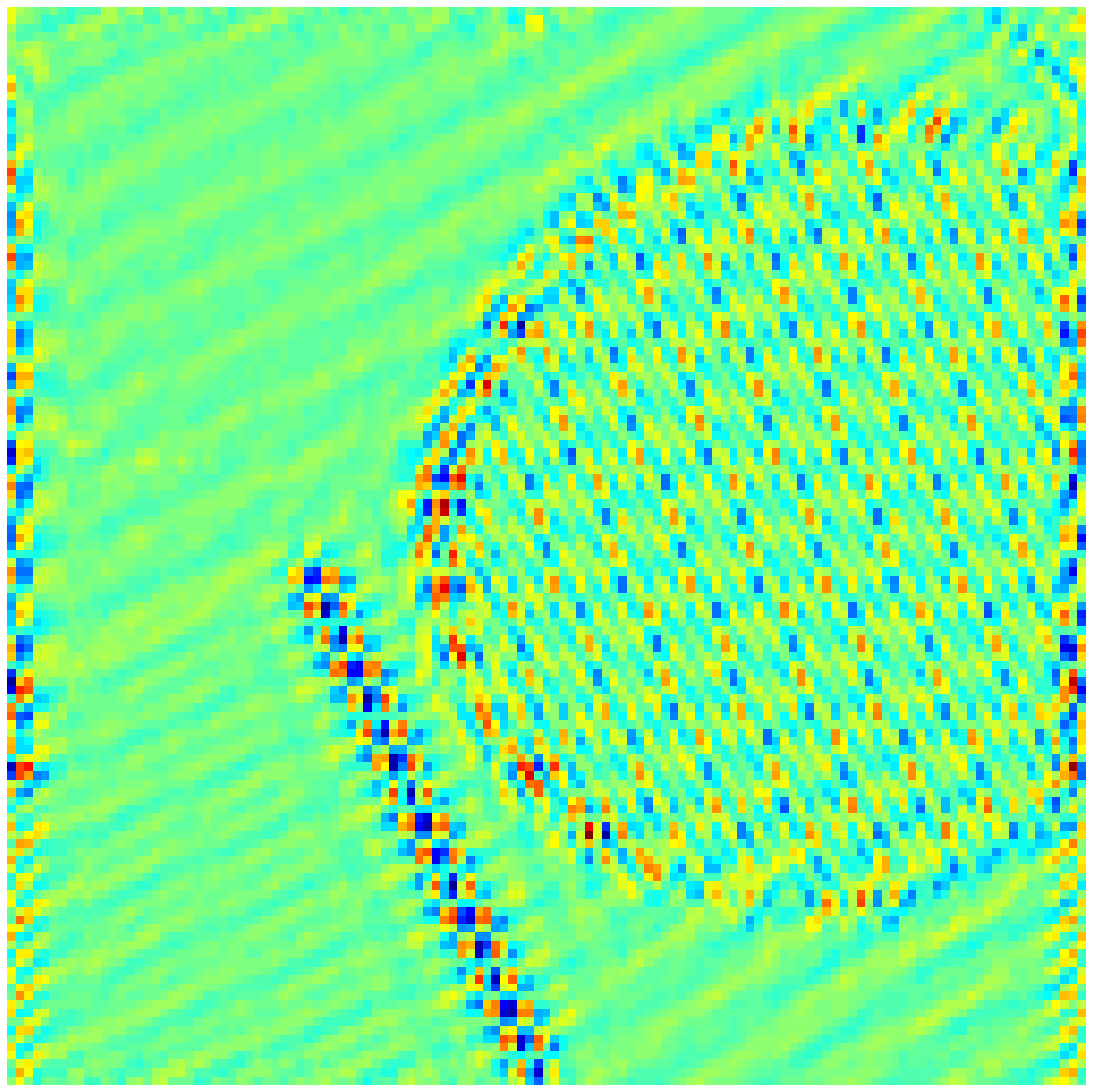}\\
\includegraphics[width=0.18\linewidth]{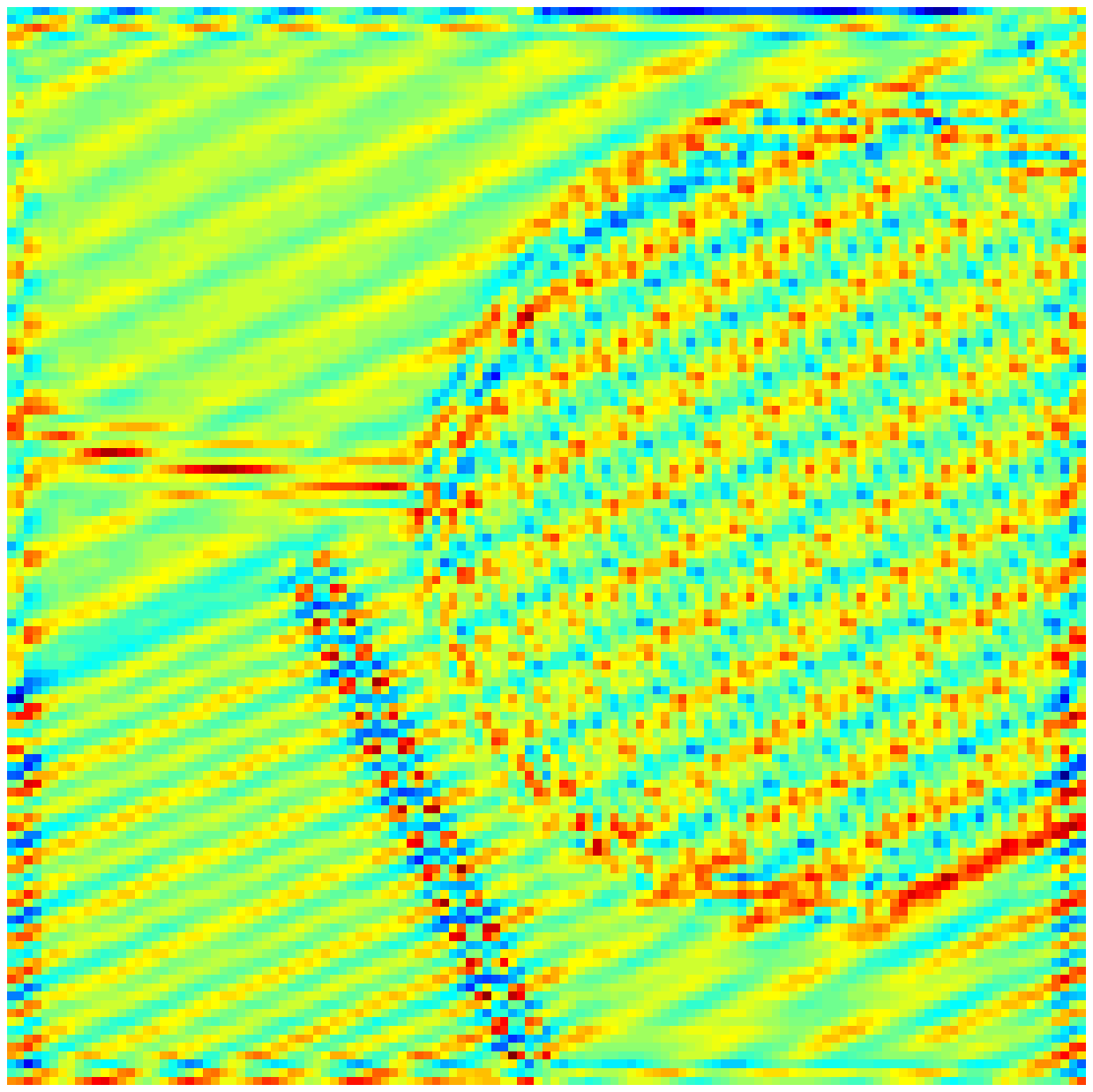}
\includegraphics[width=0.18\linewidth]{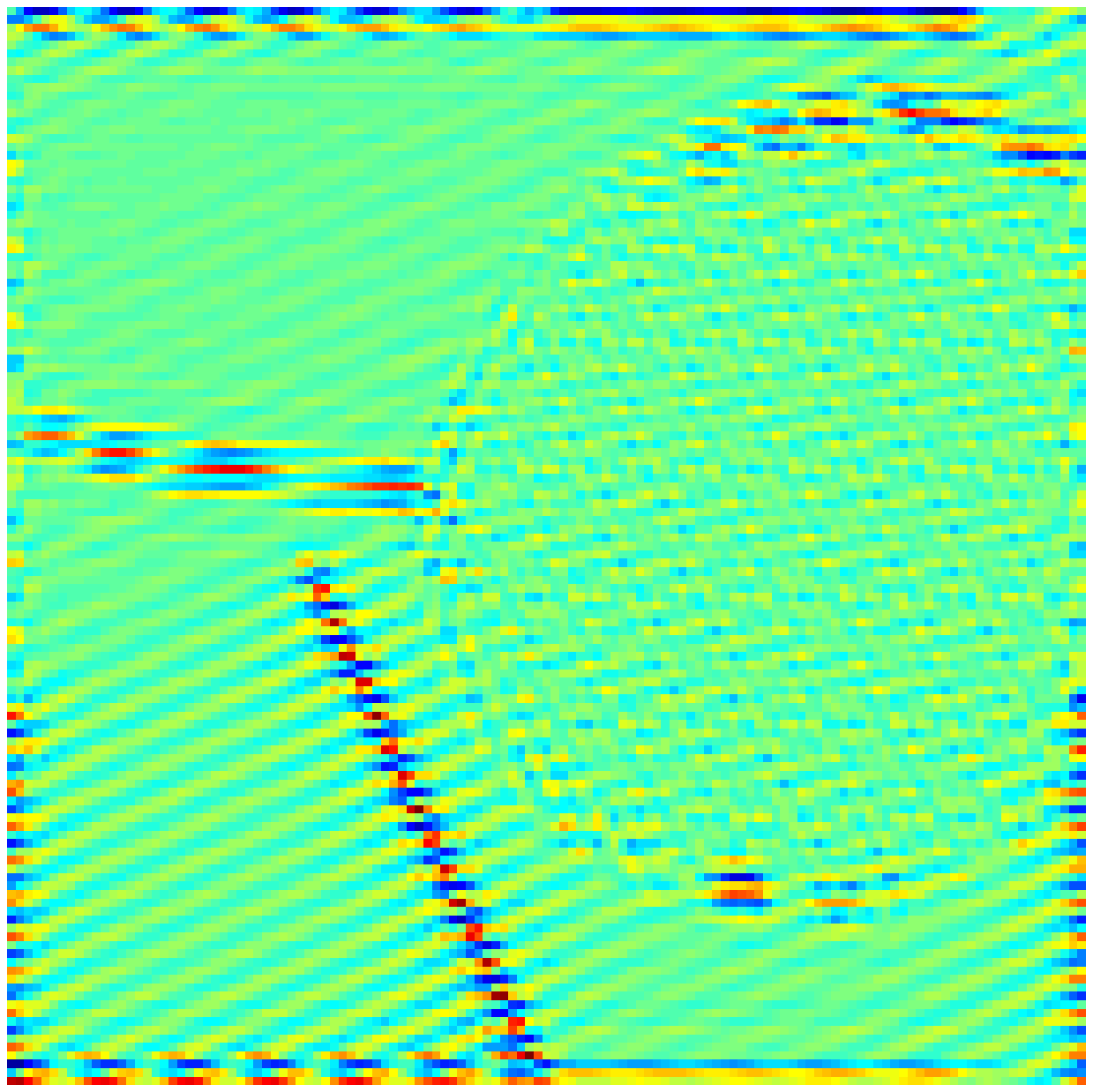}
\includegraphics[width=0.18\linewidth]{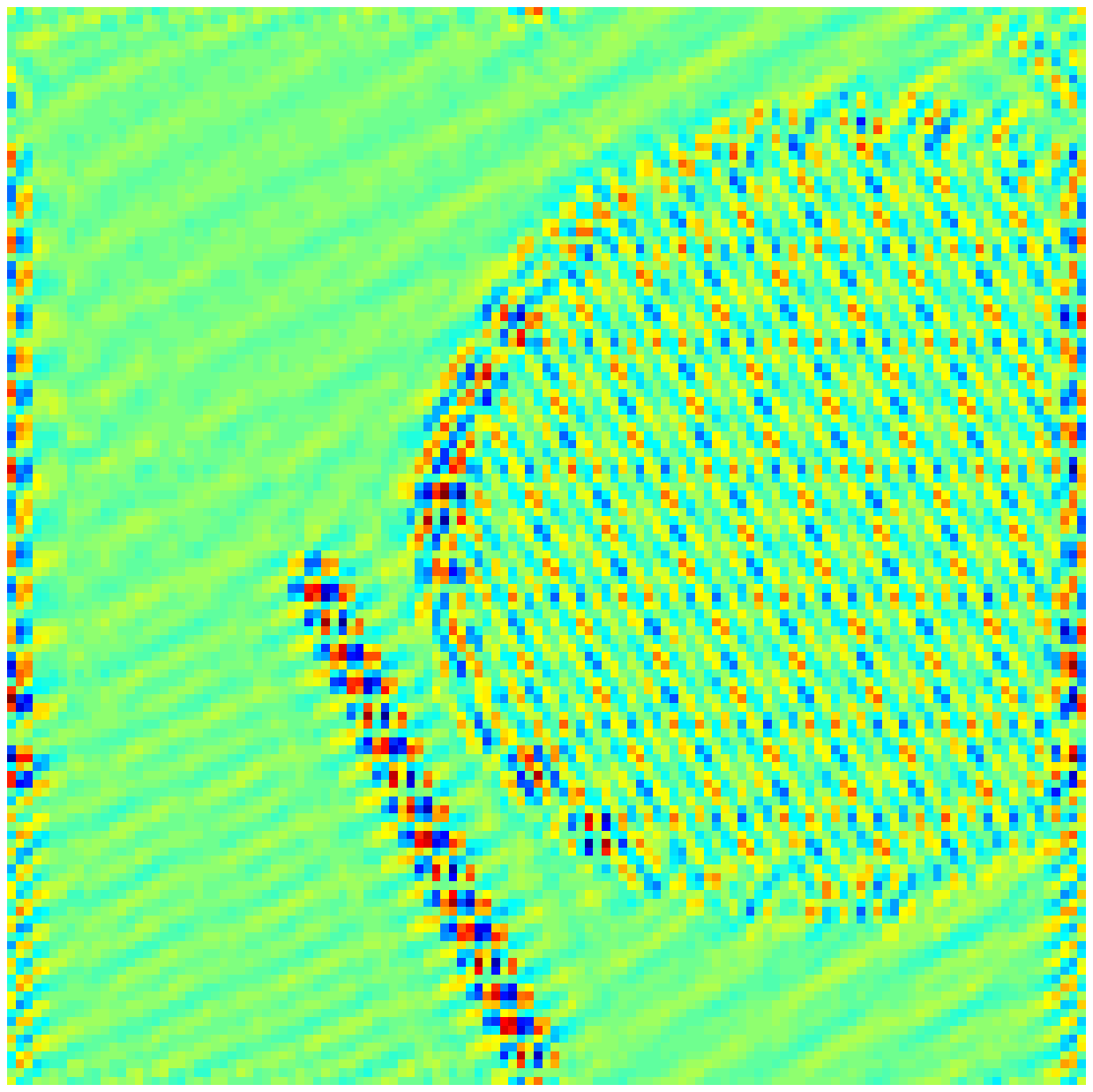}
\vspace{-0.2cm}
\caption%[]
{Filtered $g*\phi_n$, $n=1,...,7$.}
\label{filtered}
% \vspace{0.5cm}
% \end{minipage}
\end{figure}

Let us now consider the following synthetic image:
\vspace{-0.6cm}
\begin{figure}[H]
\begin{minipage}[c]{\linewidth}
\centering
\includegraphics[width=0.41\linewidth]{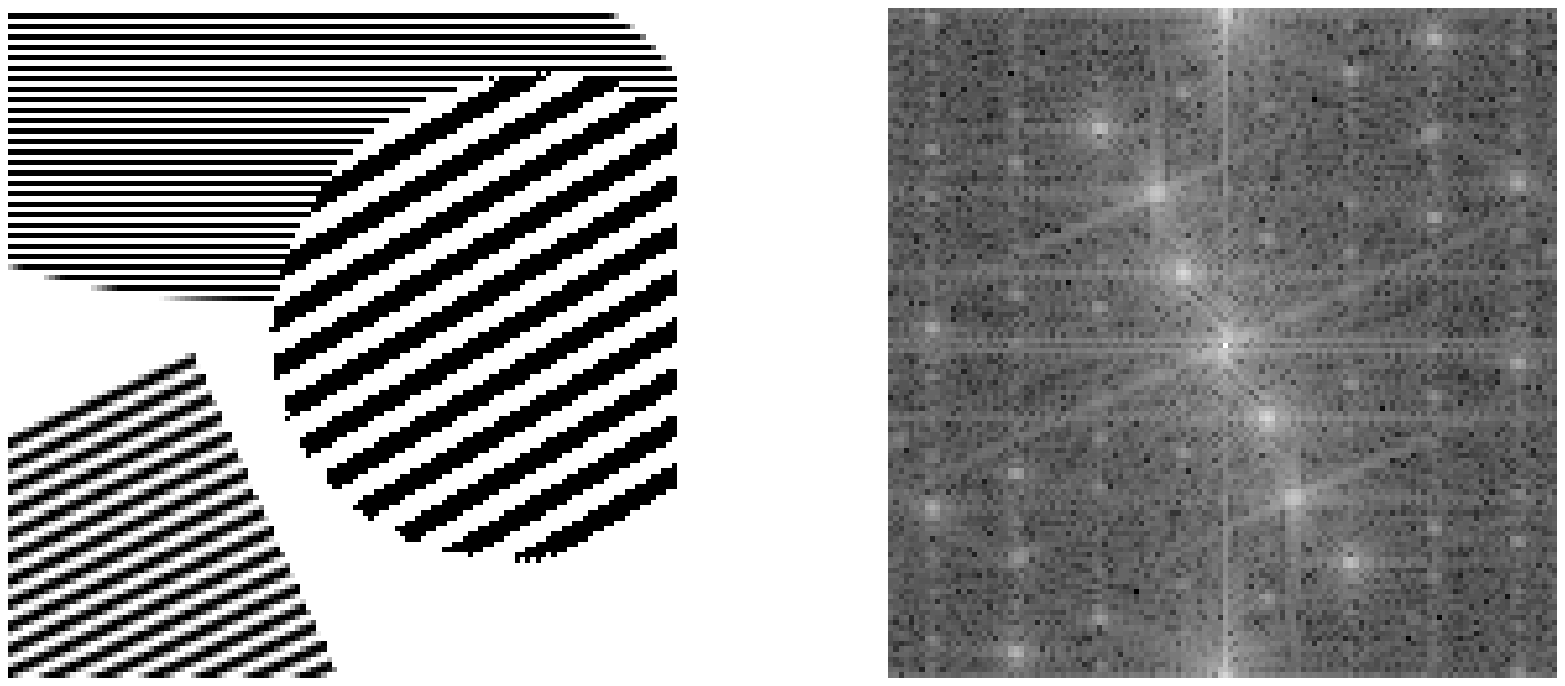}
\hspace{0.5cm}
\includegraphics[width=0.18\linewidth]{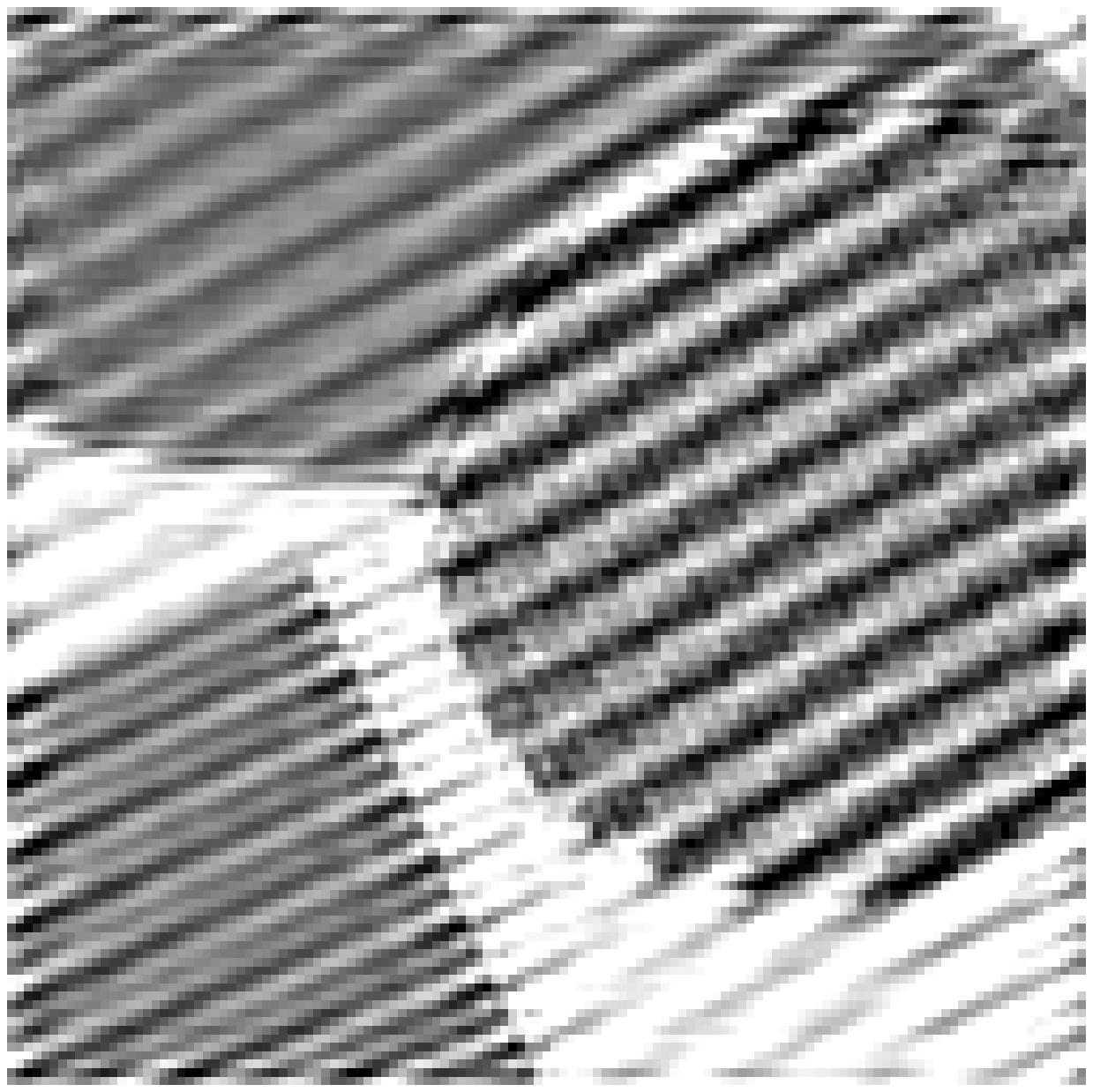}
\hspace{0.5cm}
\includegraphics[width=0.18\linewidth]{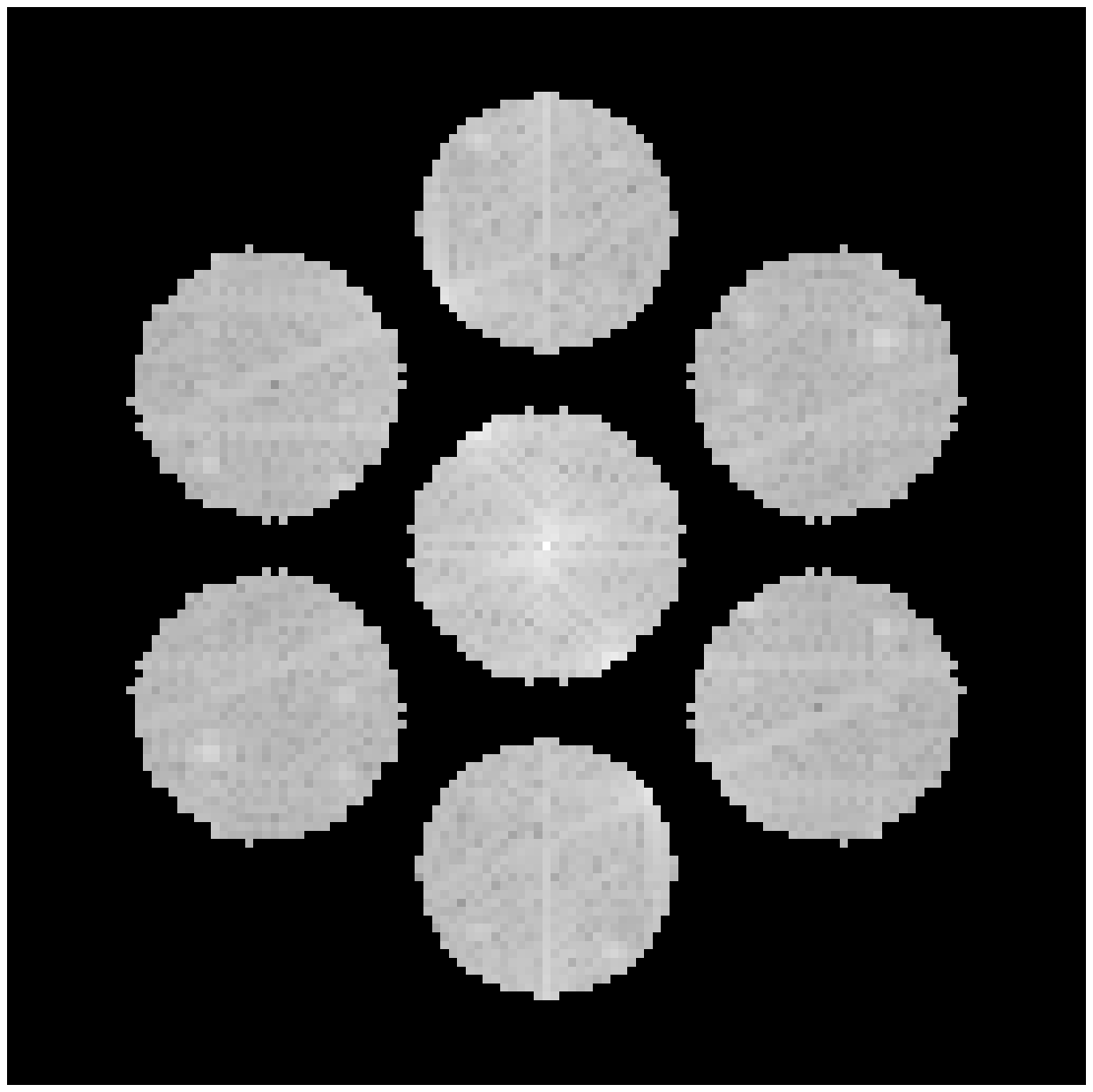}
\vspace{-0.2cm}
\caption%[]
{From left to right: the original $128\times128$ image $g_0$ and its spectrum $\log|\F(g_0)|$, the corrupted image $g$ and its spectrum $\log|\F(g)|$. %As can be seen, we got rid of the high frequencies.
}
% \vspace{0.5cm}
\end{minipage}
\end{figure}
\vspace{-0.6cm}
In Fig.~\ref{filtered}, we can see that these atoms behave as Gabor atoms by capturing different patterns in the corrupted image. The different regions of $g$ can thus be distinguished by analyzing the filtered $g*\phi_n$.
% Filtering this image by the atoms $\phi_n$, $n=1,...,7$ one gets:\\
\vspace{-0.6cm}
\vspace{-0.6cm}
% \begin{comment}
% \newpage
The previous atoms were computed for the mask considered in Fig.~\ref{Atom_mask}. However, our approach is quite general and can be adapted to any other corruption mask. As an example, here follow the atoms we get if we consider the mask that comes into play in the scattering problem:
% \begin{figure}
% % \begin{minipage}[c]{\linewidth}
% % \centering
% \includegraphics[width=\hsize]{scatt_50_atoms_axis_bis_crop}
% \vspace{-0.5cm}
% \caption%[]
% {The 50 first atoms (zoomed in) with $p=2$ adapted to the mask associated to the scattering problem.}
% % \end{minipage}
% \end{figure}
% % \vspace{-0.3cm}

%In the following we depict the very first atoms with their spectra:
\vspace{-.5cm}
\begin{figure}[H]
\includegraphics[width=\hsize]{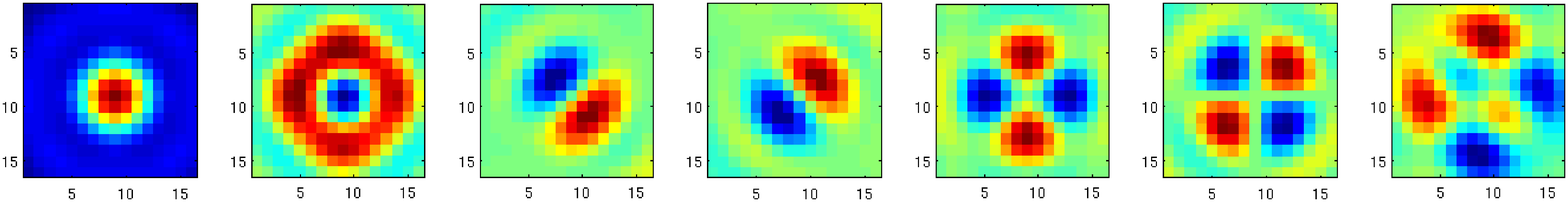}
\includegraphics[width=\hsize]{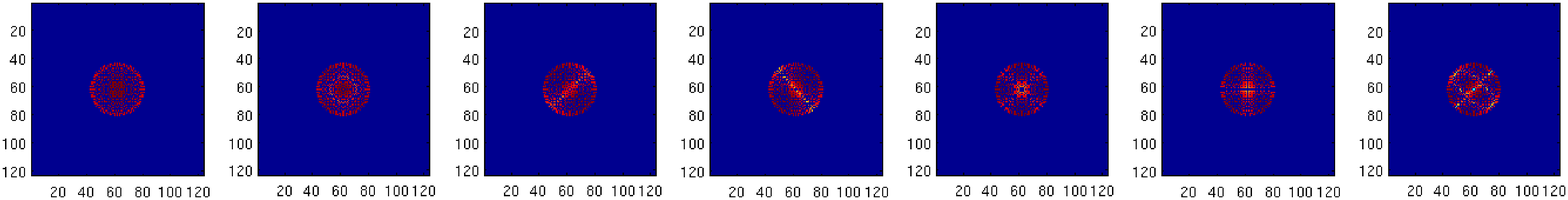}
\vspace{-0.5cm}
\caption{First line: Zoomed-in atoms ${(\phi_n)}_{n=1,...,7}$ with $p=4$ adapted to the scattering problem. Second line: Their respective spectra.}
\end{figure}
\vspace{-0.5cm}
\newbf{These atoms are also really localized and are about the same size as a classical image patch (about 10 pixels wide) which therefore allows us to correlate local features of an image within a Non-Local approach.}

\section{Distance Map Comparison}
We are going to compare the performance of the different similarity measures considered so far
% \begin{figure}
% \begin{minipage}[c]{\linewidth}
% \centering
% 
% % \vspace{-0.5cm}
% \caption%[]
% { Psychedelic \textit{Lena} image $g_0$ that and the corrupted $g$}
% % \vspace{0.5cm}
% \end{minipage}
% \end{figure}
by fixing one patch of the \textit{corrupted image} $g$ (indicated in \textcolor{green}{green} in the clean image) and identify the $13$ best matches in the \textit{corrupted image} (and indicated in \textcolor{red}{red} in the clean image) according to
%\vspace{-0.3cm}
\begin{itemize}
 \item[$1.$] The atom-based distance $\delta^1(x_k,x_\ell)={\left(\sum_{n=1}^7|g*\phi_n(x_k)-g*\phi_n(x_\ell)|^2\right)}^{\frac{1}{2}}$.
 \item[$2.$] The ideal distance $\delta^2(x_k,x_\ell)={\|p_{g_0}(x_k)-p_{g_0}(x_\ell)\|}_2$ based on the clean $g_0$.
\vspace{.15cm}
 \item[$3.$] The SSD $\delta^3(x_k,x_\ell)={\|p_g(x_k)-p_{g}(x_\ell)\|}_2$.
\end{itemize}
% \vspace{-0.5cm}
%Let us consider the following image $g_0$ and its corresponding $g$:
% Here is what we get:
% \textbf{Example 1:}
\vspace{-0.7cm}
\begin{figure}[H]
\centering
\begin{tabular}{cccc}
Original and corrupted & Atom-based % $\delta^1$ 
& Oracle % $\delta^2$ 
& SSD % $\delta^3$
\\
\includegraphics[width=0.24\linewidth]{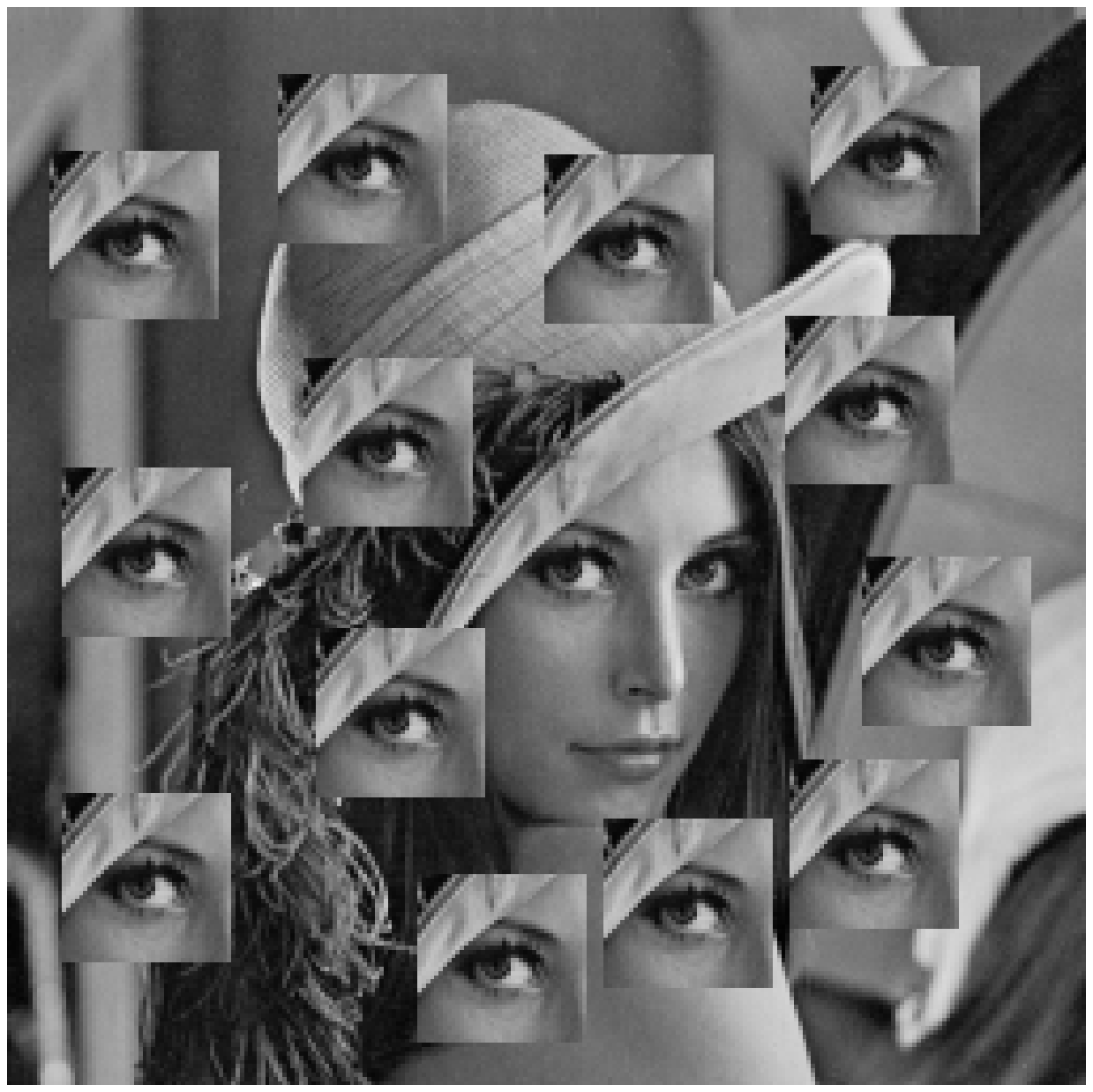}
&\includegraphics[width=0.24\linewidth]{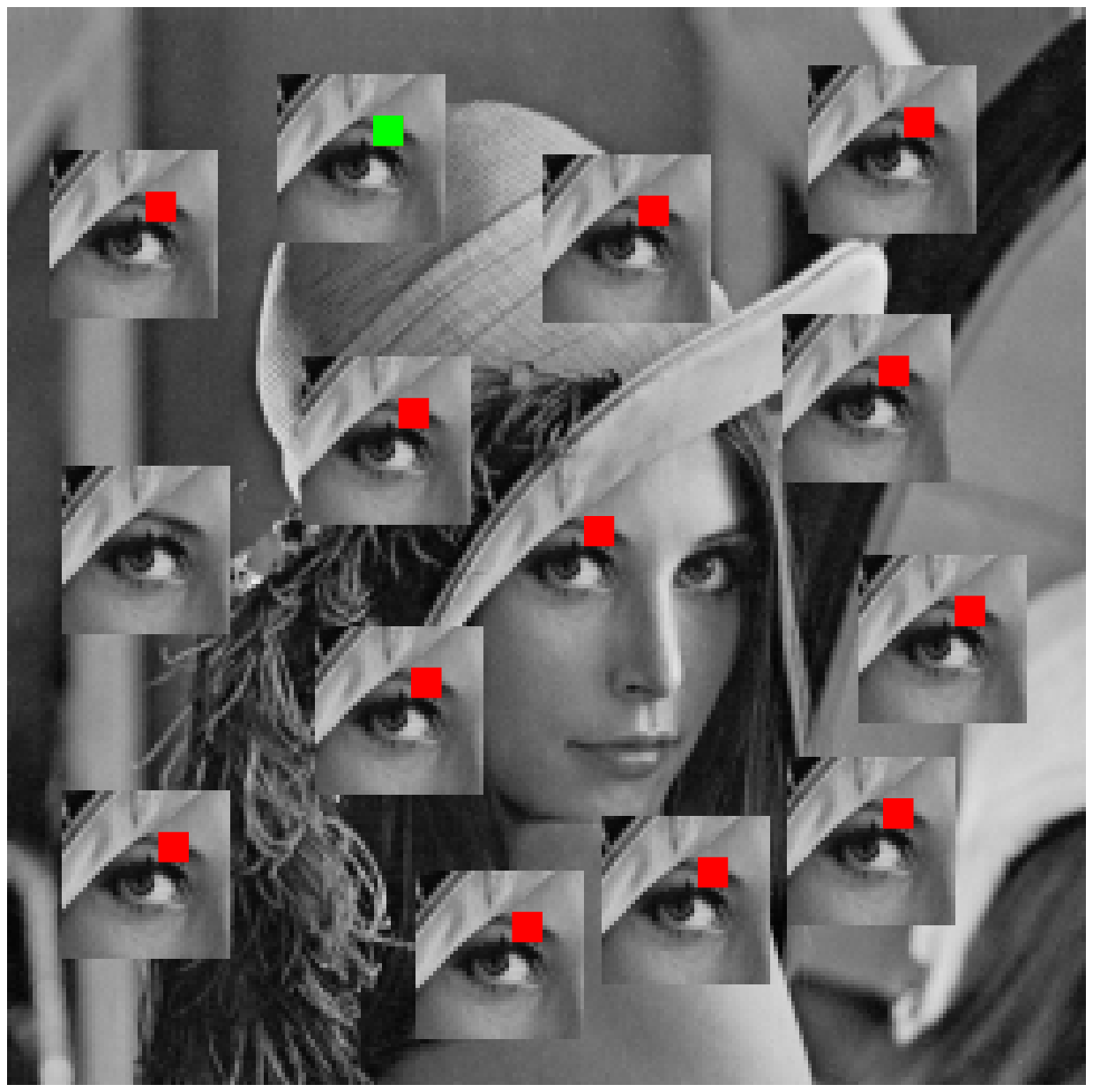}
&\includegraphics[width=0.24\linewidth]{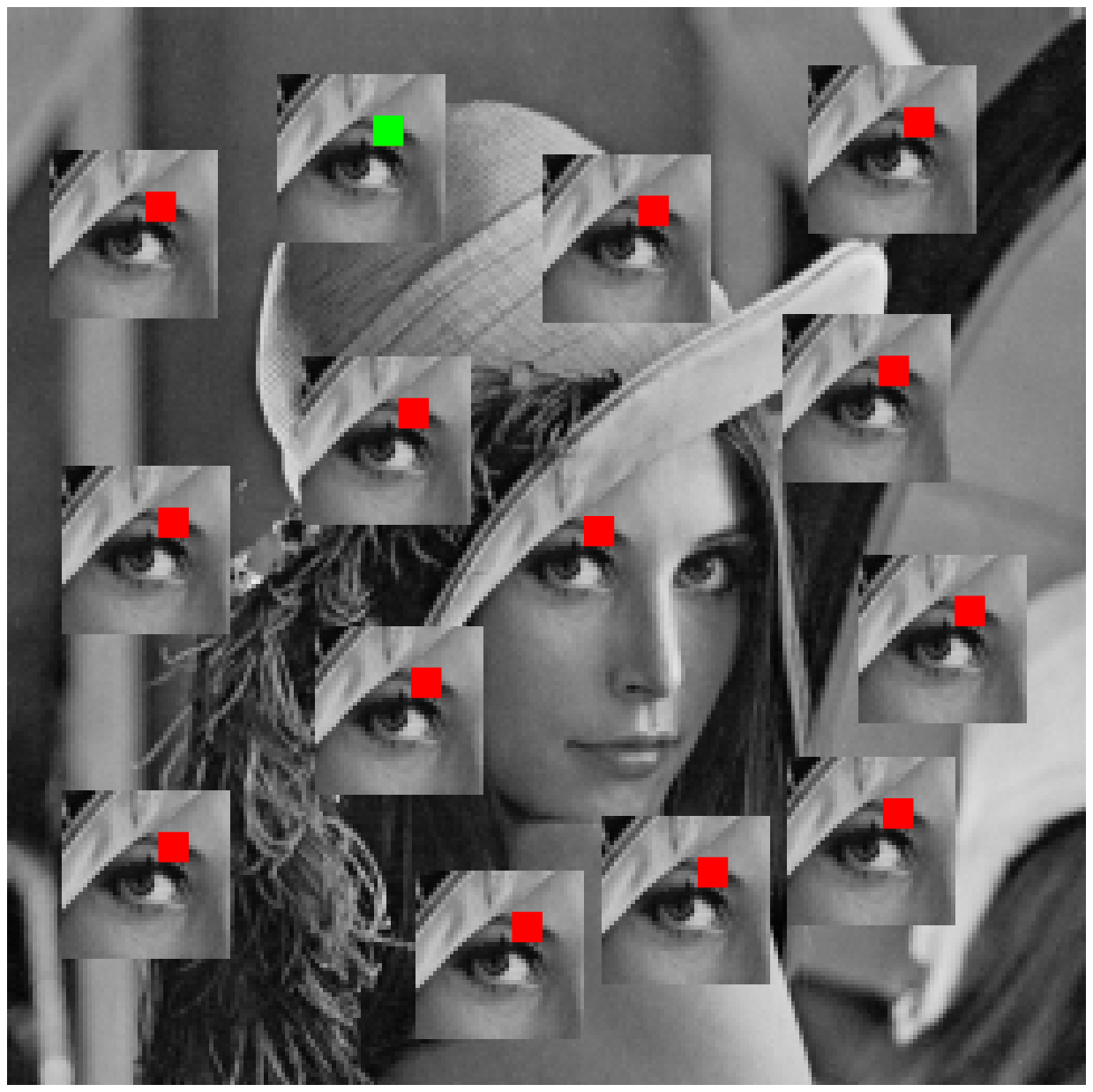}
&\includegraphics[width=0.24\linewidth]{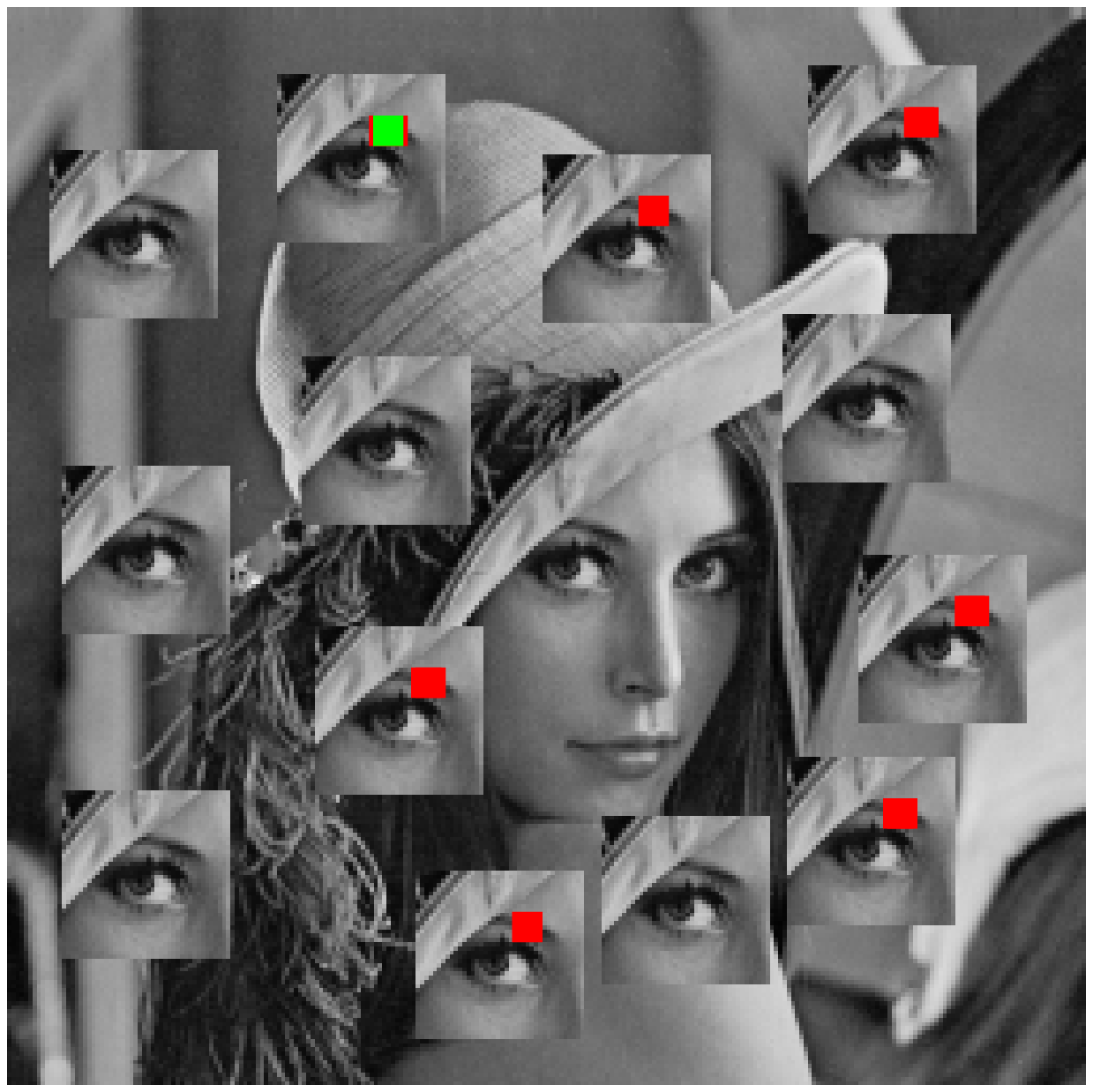}
\\
\includegraphics[width=0.24\linewidth]{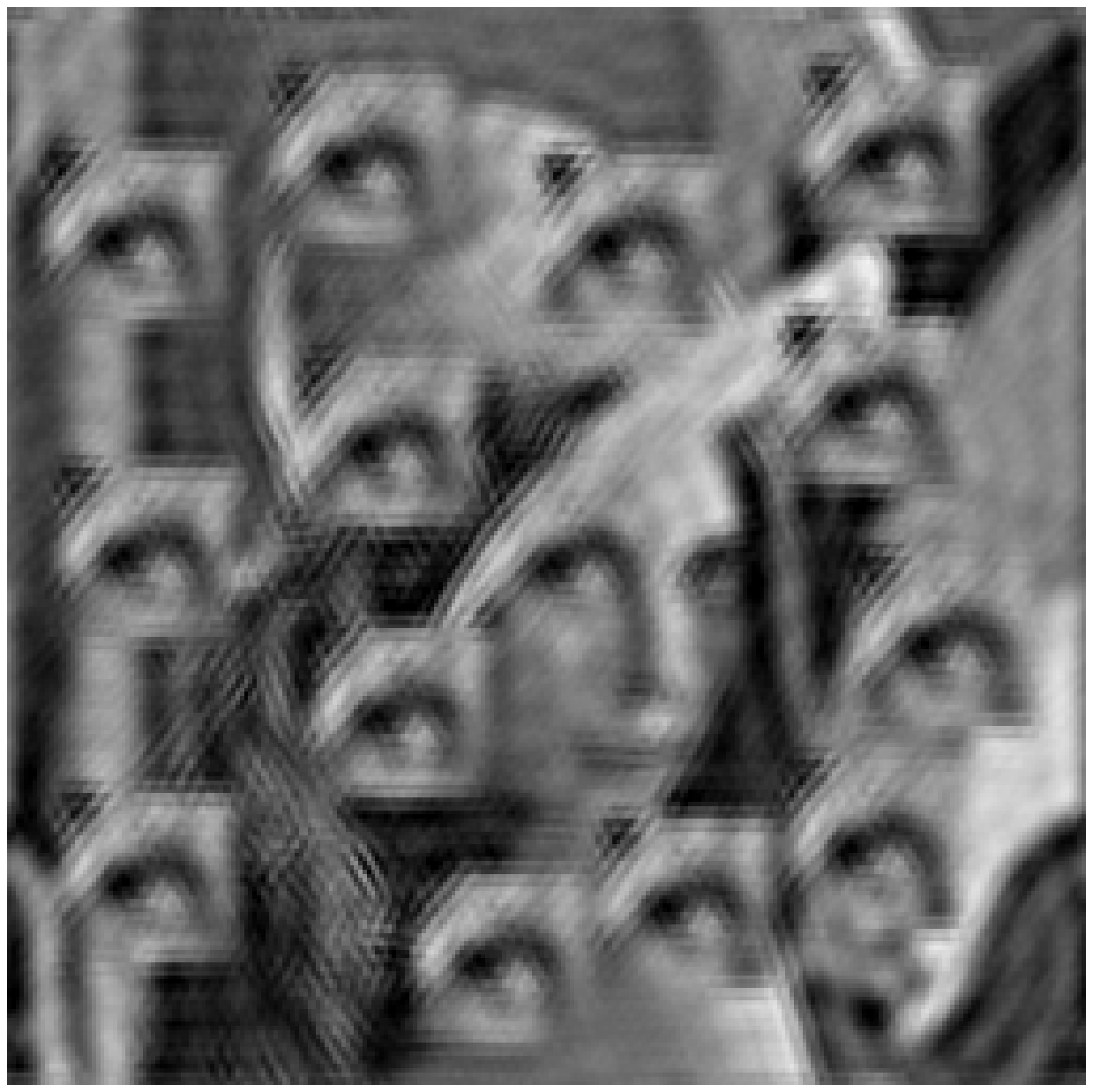}
&\includegraphics[width=0.24\linewidth]{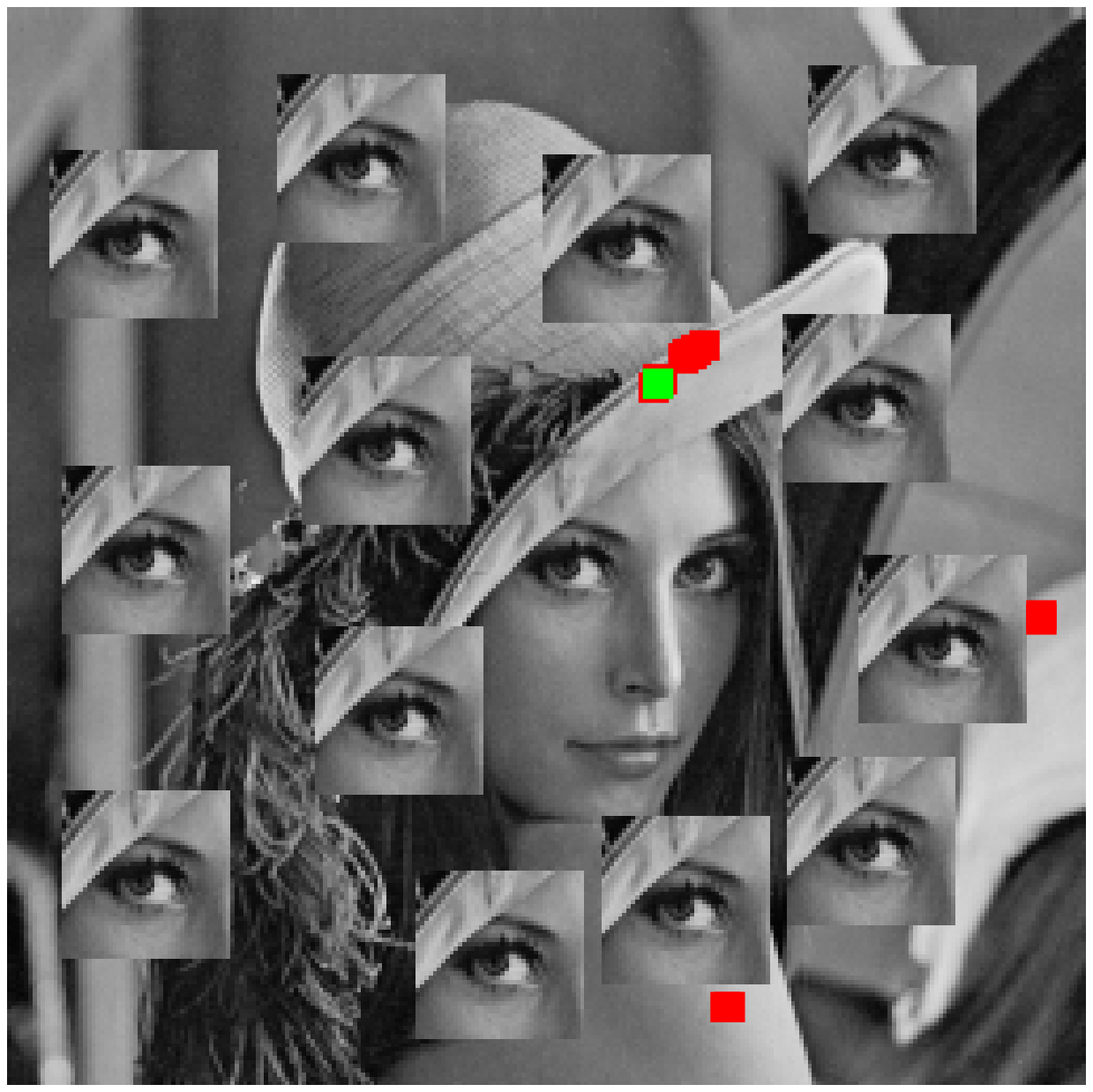}
&\includegraphics[width=0.24\linewidth]{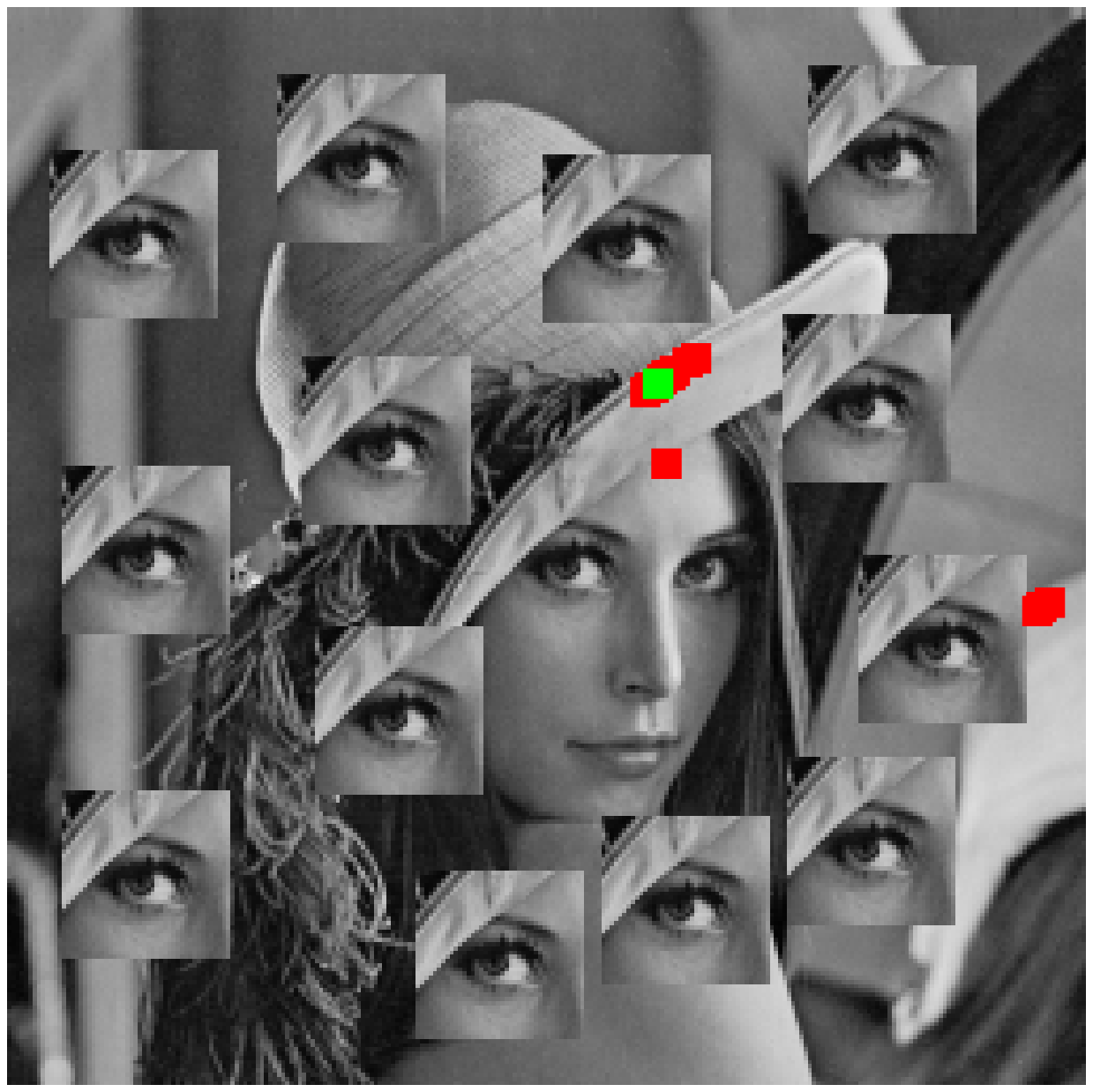}
&\includegraphics[width=0.24\linewidth]{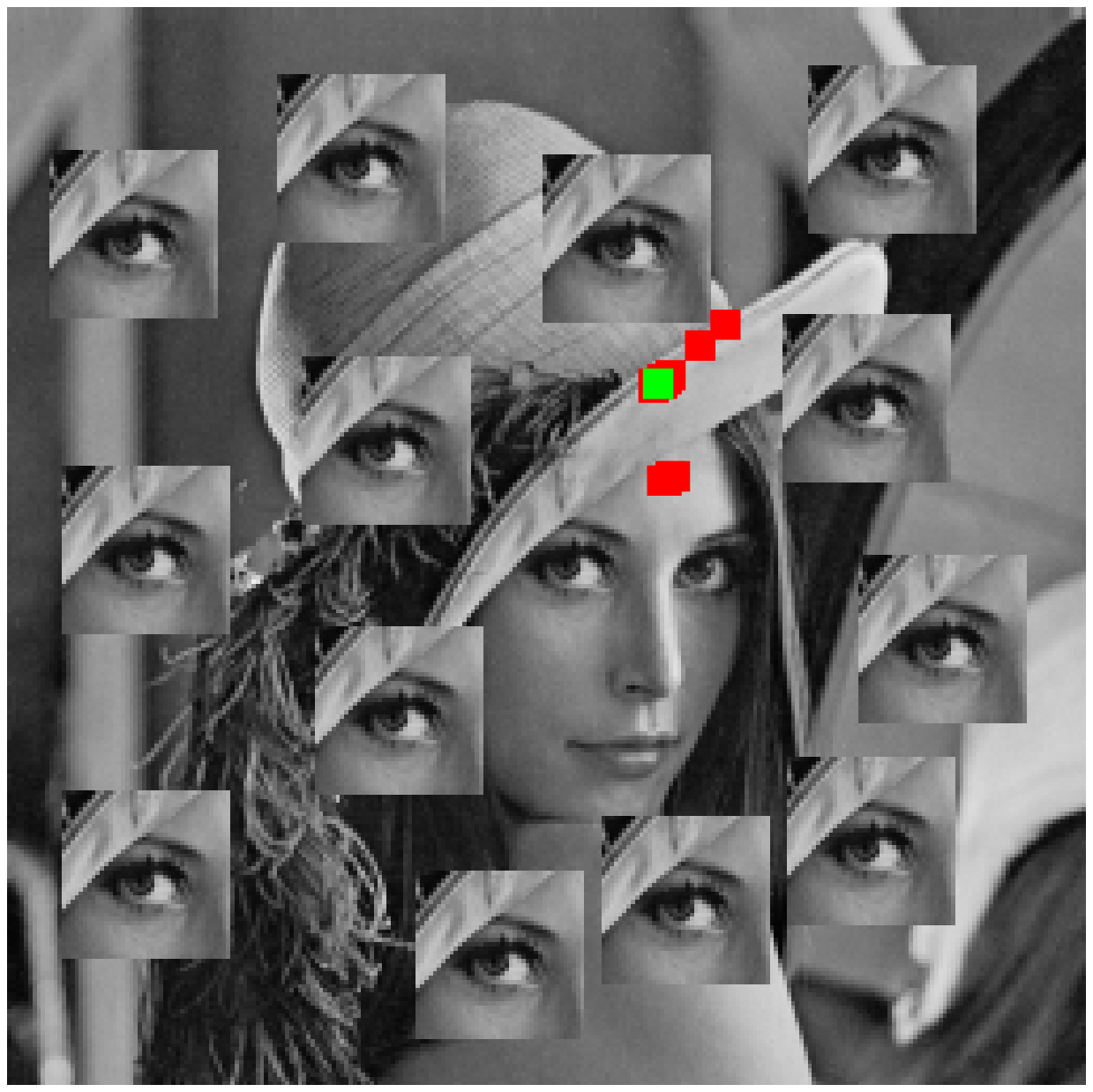}
\end{tabular}
\vspace{-0.3cm}
\caption{Best matches (in \textcolor{red}{red}) corresponding to a fixed patch (depicted in \textcolor{green}{green}) \newbf{in an image with lot of self-similarities.}}
\label{best_match}
\label{compar_test}
\end{figure}
\newbf{The two lines correspond to two different choices for the patch
or reference (first line: a small square on the eyebrow from the second
large eye patch,
second line: the edge of the hat towards the middle)
 The second and fourth columns are computed using the corrupted image $g$. The results are displayed on the clean image. In the third column we depicted the results obtained with the SSD computed with the clean image.}

These two examples suggest that the atom-based distance we designed performs
 better than the SSD in some cases: \newbf{indeed, the atom based results are closer to the ``oracle'' than the SSD which looses more matches}.

\section{Numerical Experiment}\label{num_ex_NL}
From now on, we are going to make important assumptions that let us drastically improve the complexity of the minimization problem (\ref{E_NL}) and save memory:
\begin{itemize}
\item[$(i)$] Two patches are unlikely to be the same if they are far from each other thus we can assume that for any fixed $x_k$, $x_\ell$ is a candidate if $|x_k-x_\ell|\leq\eta$ for some fixed $\eta$. This defines a neighborhood of candidates (also called \textit{window}).
\item[$(ii)$] Given $x_k,x_\ell$, minimizing (\ref{E_k_ell}) yields a $v_{k,\ell}$ whose support is actually $(x_k+\supp(\psi))\cup(x_\ell+\supp(\psi))$. Therefore, assuming that $|x_k-x_\ell|\geq \e$ for some $\e\leq\rho$, we can still get a global minimizer for (\ref{E_NL}). % by considering $w=\sum_{k,\ell} w_{k,\ell}$. 
Although, usually $\e=1$ pixel, assuming that $\e=7$ pixels for a patch size $\rho=7$ or $9$ pixels (to let the reconstructed patches to overlap) results in an acceleration of order $10$.
\item[$(iii)$] Once $w$ computed, we can for any fixed pixel $x_k$ keep only the $m_0$ best matches $x_\ell$ (see Fig.~\ref{best_match}). 
\end{itemize}
%Therefore in the sequel we are only going to consider that $I=\{(i,j),\e\leq|x_k-x_\ell|,|x_k-x_\ell|\leq\eta\}$.

%Let us remark that our first assumption can be improved by selecting regions in the image (not necessarily windows) where patches are likely to be the same as in \cite{Elad}.
% \newpage
\noindent\textbf{Weight computation:}
Notice that for a fixed mask $M$, the atoms ${(\phi_n)}_{n=1,...,n_0}$ can be computed in advance and stored. Then $g_n=\phi_n*g$ can be readily, computed once for all, using the \verb+fft+. Now, we recall that
\begin{align}
{\big(\delta^1(x_k,x_l)\big)}^2&=\sum_{n=1}^{n_0}{|g_n(x_k)-g_n(x_\ell)|}^2,\\ 
% \end{align}
% whereas
% \begin{align}
{\big(\delta^3(x_k,x_l)\big)}^2&=\sum_{i=1}^{\rho^2}{|{\big(p_g(x_k)\big)}_i-{\big(p_g(x_\ell)\big)}_i|}^2.
\end{align}

%\textbf{Discussion on complexity}:
 \newbf{From the previous definitions, it is readily seen that if $\rho=7$ or $9$ pixels is the patch size, then, as far as $(n_0\leq\rho^2)$\ \footnote{which will always be the case, otherwise it means that all atoms are essentially compactly supported which is possible only if the mask is full},
it is faster to compute $\delta^1$. In practice, this can result in an acceleration of order 10. Such a gain in complexity will be observed numerically in the sequel.}
\smallskip

\noindent\textbf{Algorithmic issues:} %The problem (\ref{E_NL}) involves a quadratic energy that can be minimized thanks to the conjugate gradient or also by Nesterov's algorithm . \\
{An interesting feature of our work is that most optimization problems which are solved are constrained quadratic minimizations. Thus, the best choice (which we used) is the conjugate gradient. We also used Beck-Teboulles's algorithm \cite{Beck} in some instances, since our same program was optionally implementing an additional TV-regularization, with no real difference. The assumptions $(i)-(iii)$ allow us to accelerate the computation.}
\smallskip

\begin{figure}[htb]
\centering
% \begin{minipage}[c]{\linewidth}
\begin{tabular}{cc}
Original $g_0$& Corrupted $g$\\
\includegraphics[width=0.19\linewidth]{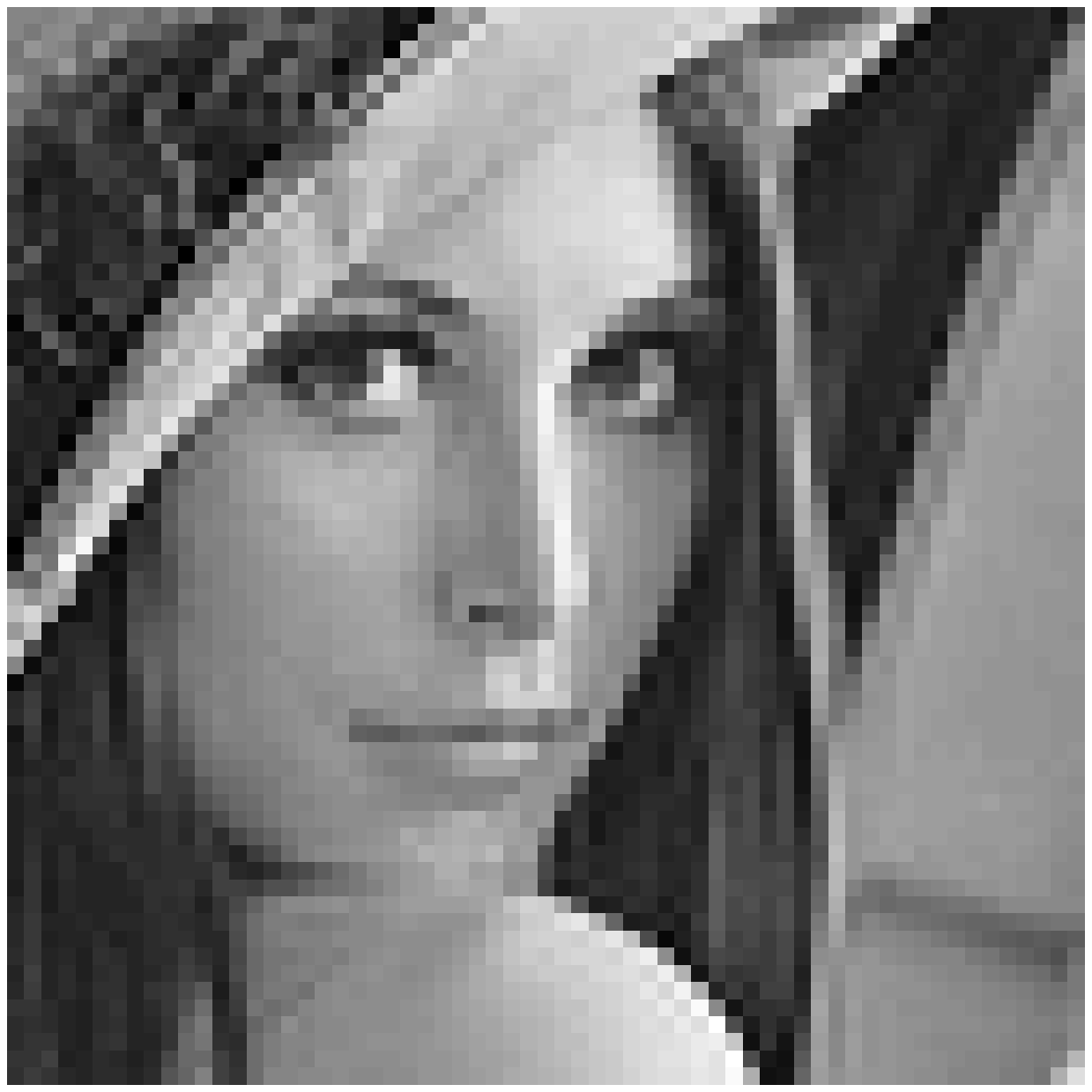}
&\includegraphics[width=0.19\linewidth]{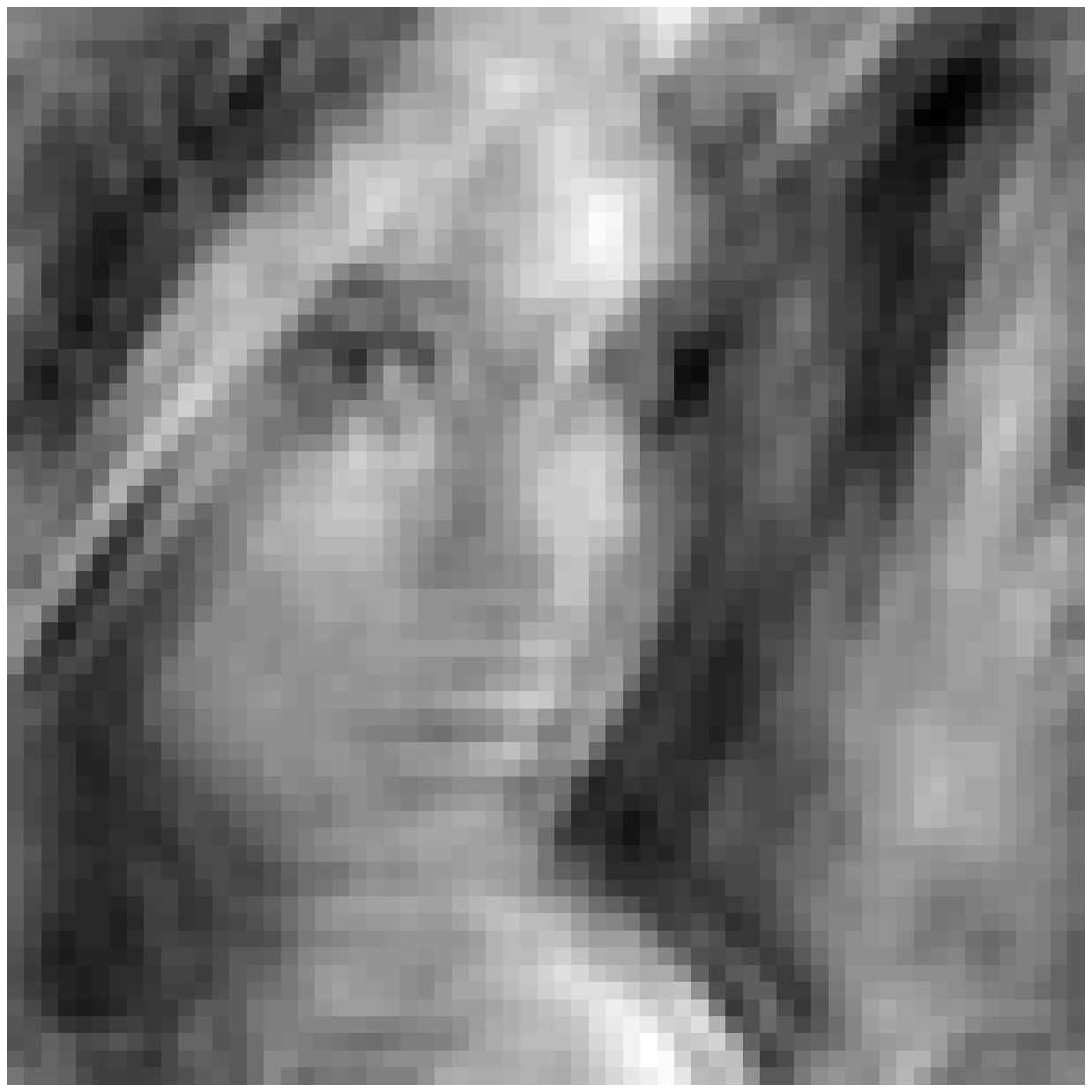}\\
&PSNR=20.0dB\\
\end{tabular}
\begin{tabular}{ccccc}
SSD - $\delta^3\ $& SSD - $\delta^3\ $ & NL-Atom - $\delta^1\ $ &NL-Atom - $\delta^1\ $ & \\
$(\alpha=1)$& $(\alpha=2)$ & $(\alpha=1)$ &$(\alpha=2)$ & TV restored\\
\includegraphics[width=0.19\linewidth]{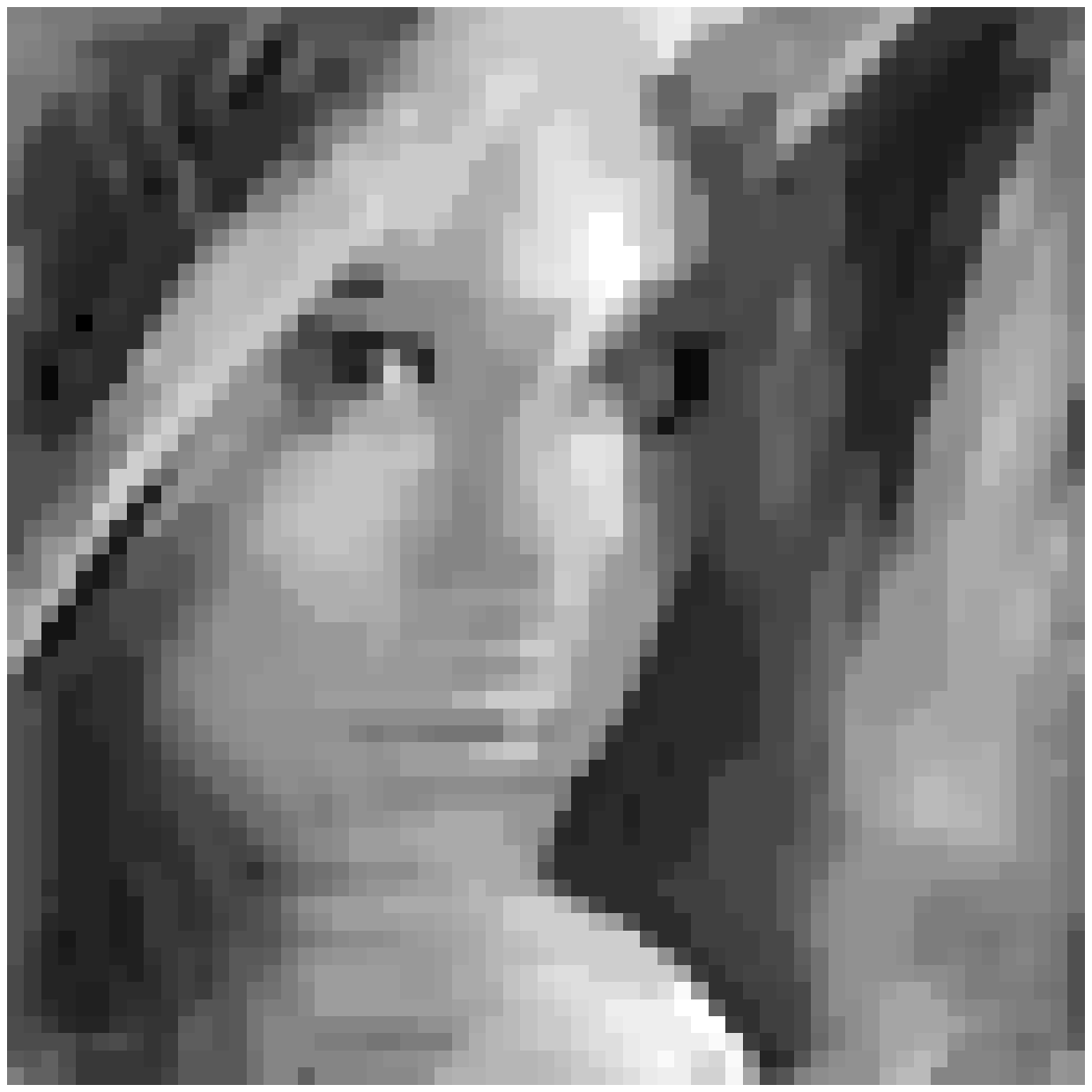}
&\includegraphics[width=0.19\linewidth]{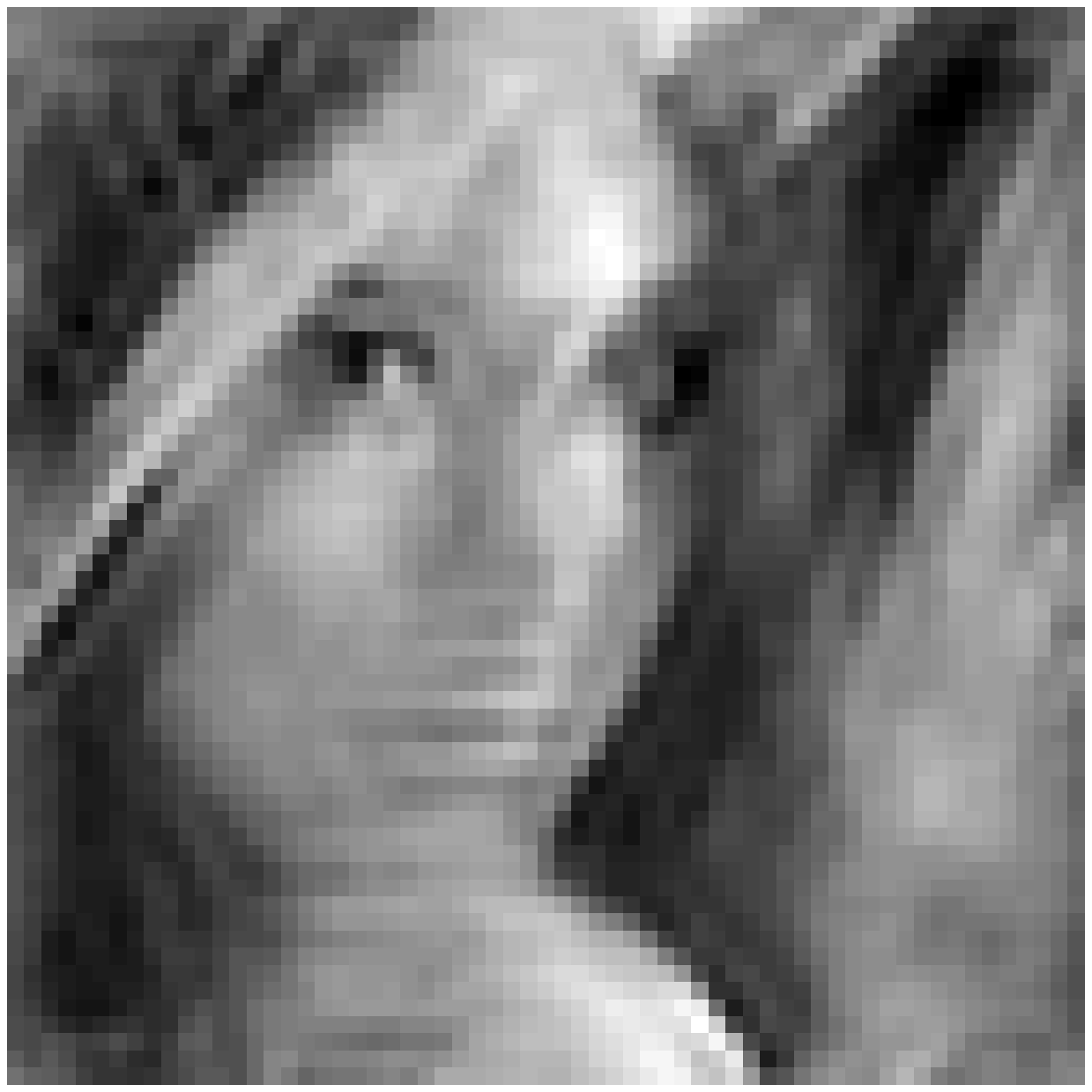}
&\includegraphics[width=0.19\linewidth]{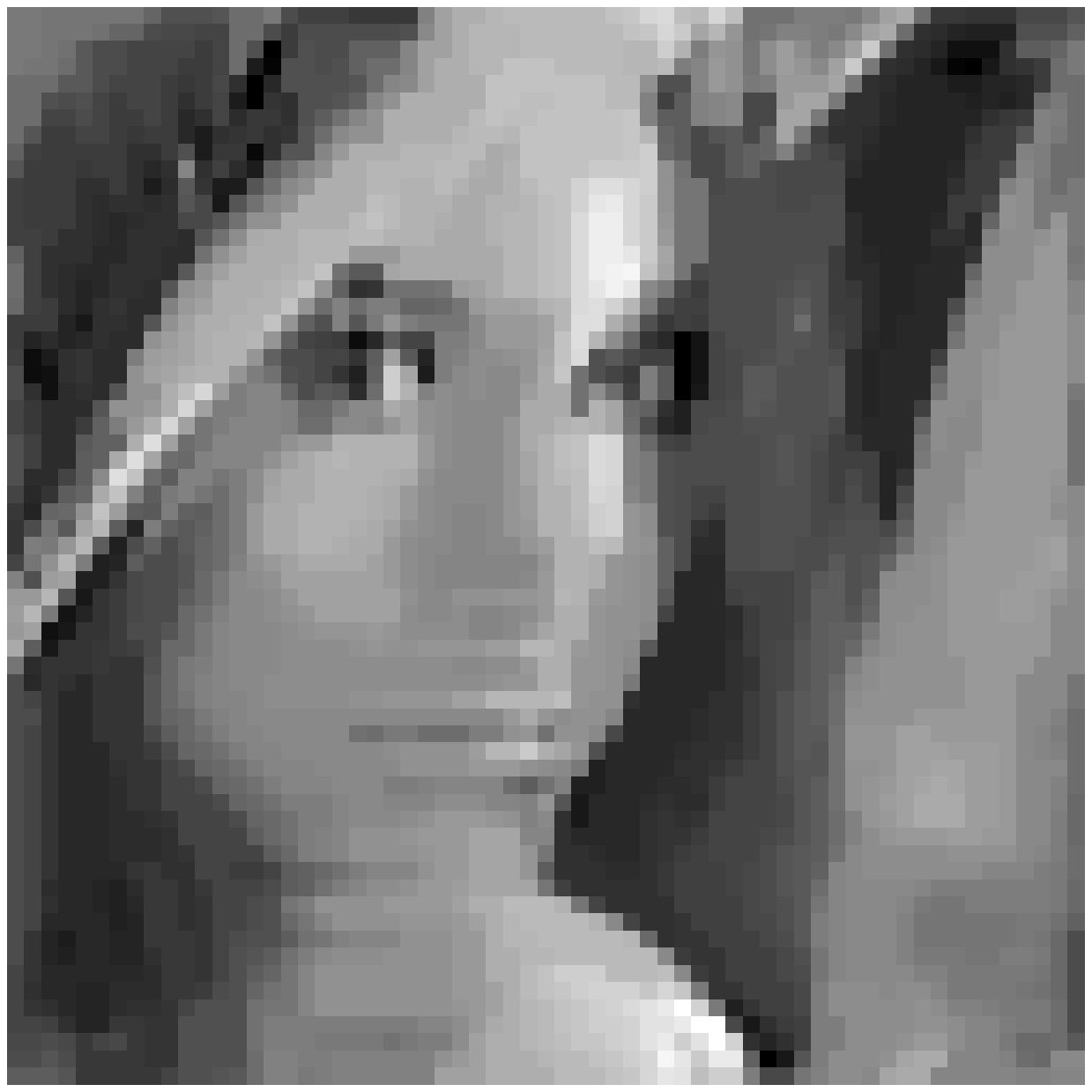}
&\includegraphics[width=0.19\linewidth]{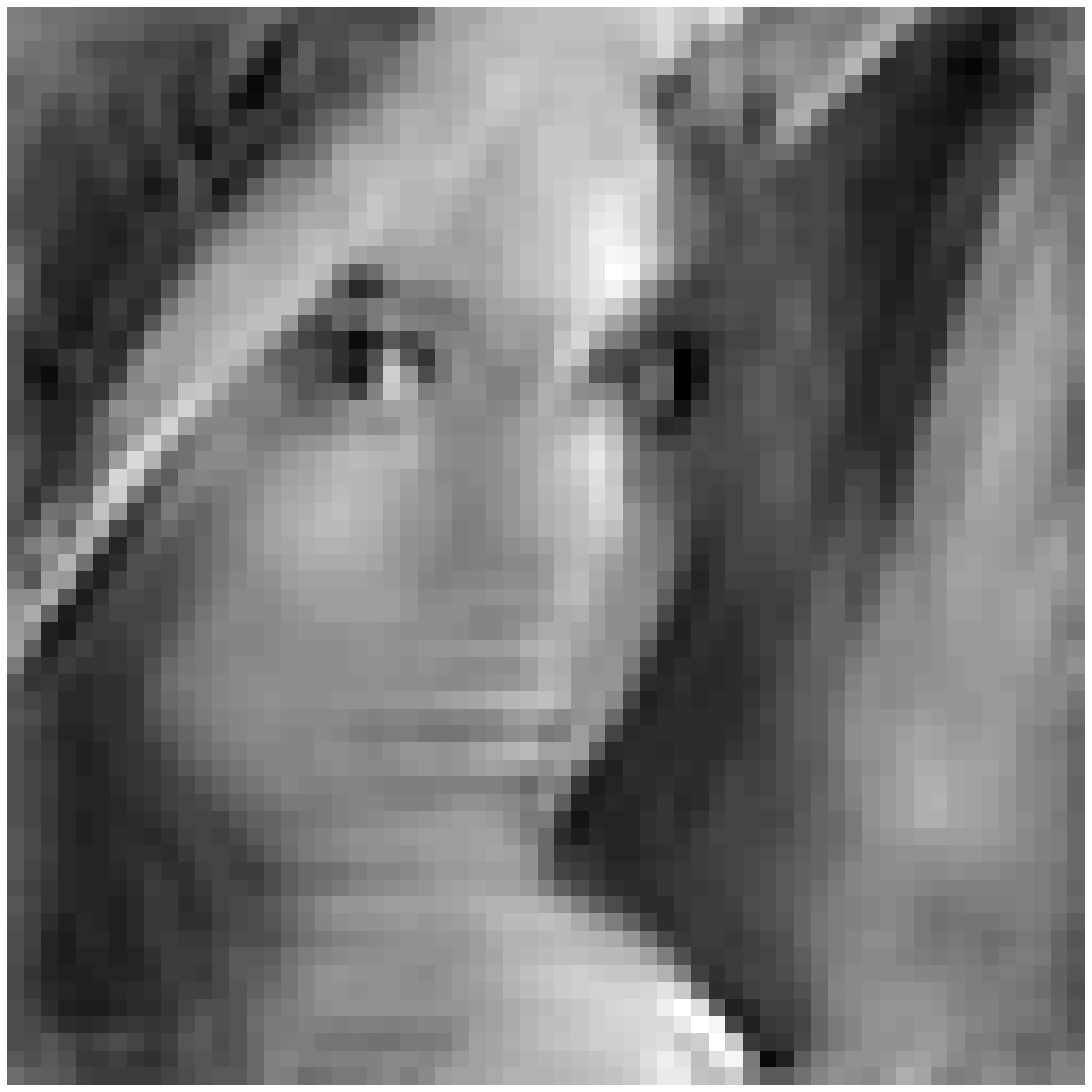}
&\includegraphics[width=0.19\linewidth]{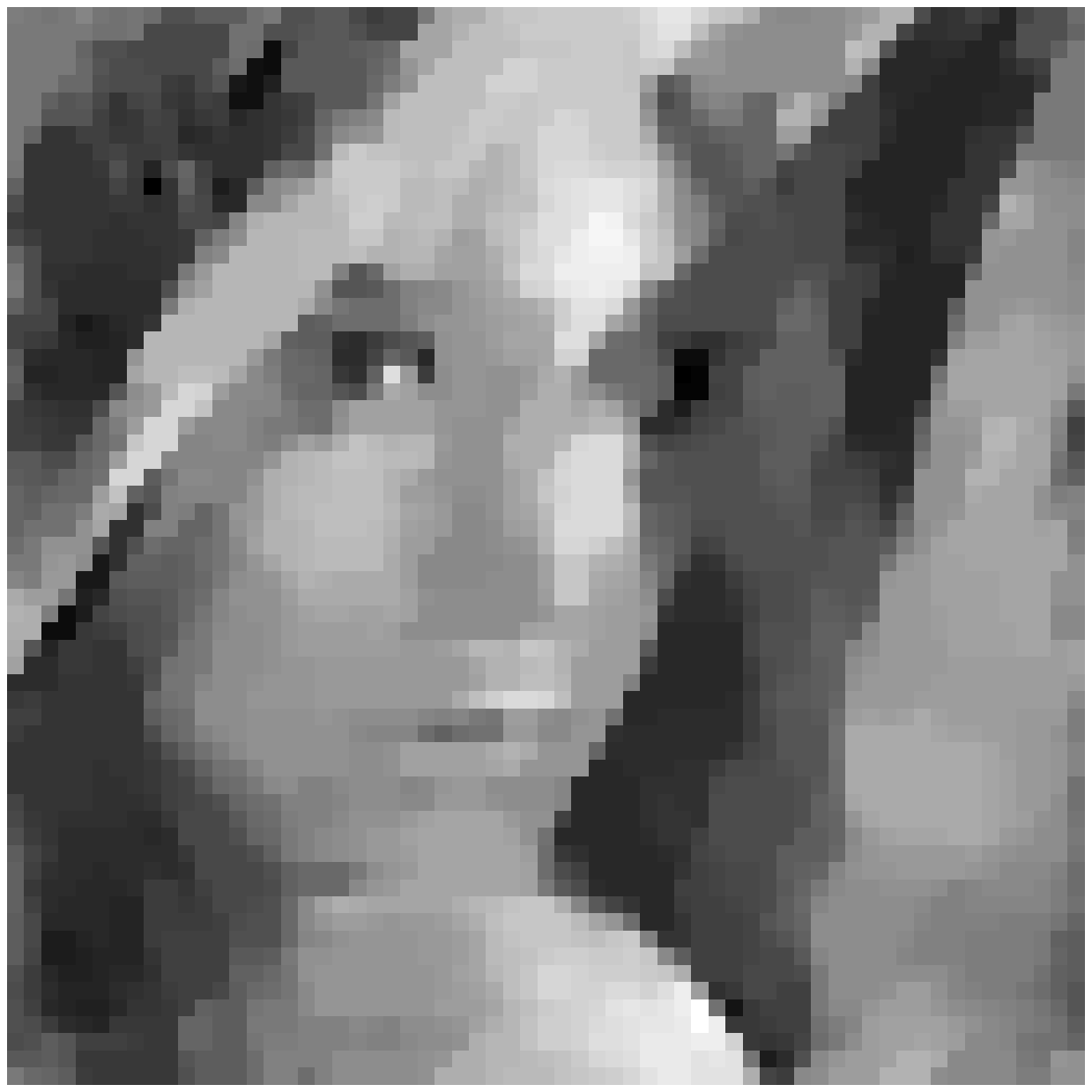}\\
%Corrupted image $u$ & $\delta^3$-NL restored & Atom-based restoration & TV restored\\
 PSNR=21.1dB &PSNR=20.7dB & PSNR=21.3dB & PSNR=21.1dB & PSNR=20.9dB\\
\includegraphics[width=0.19\linewidth]{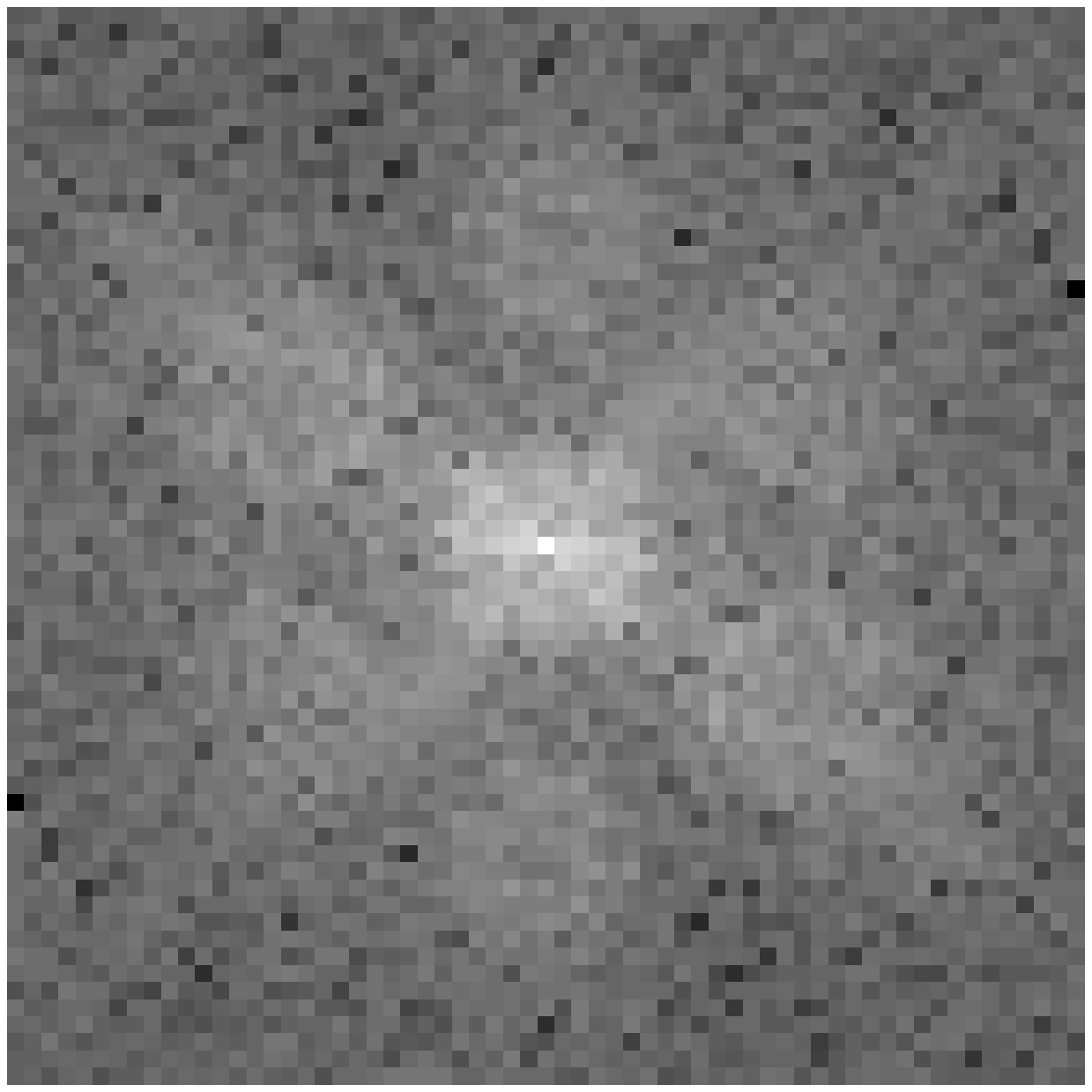}
&\includegraphics[width=0.19\linewidth]{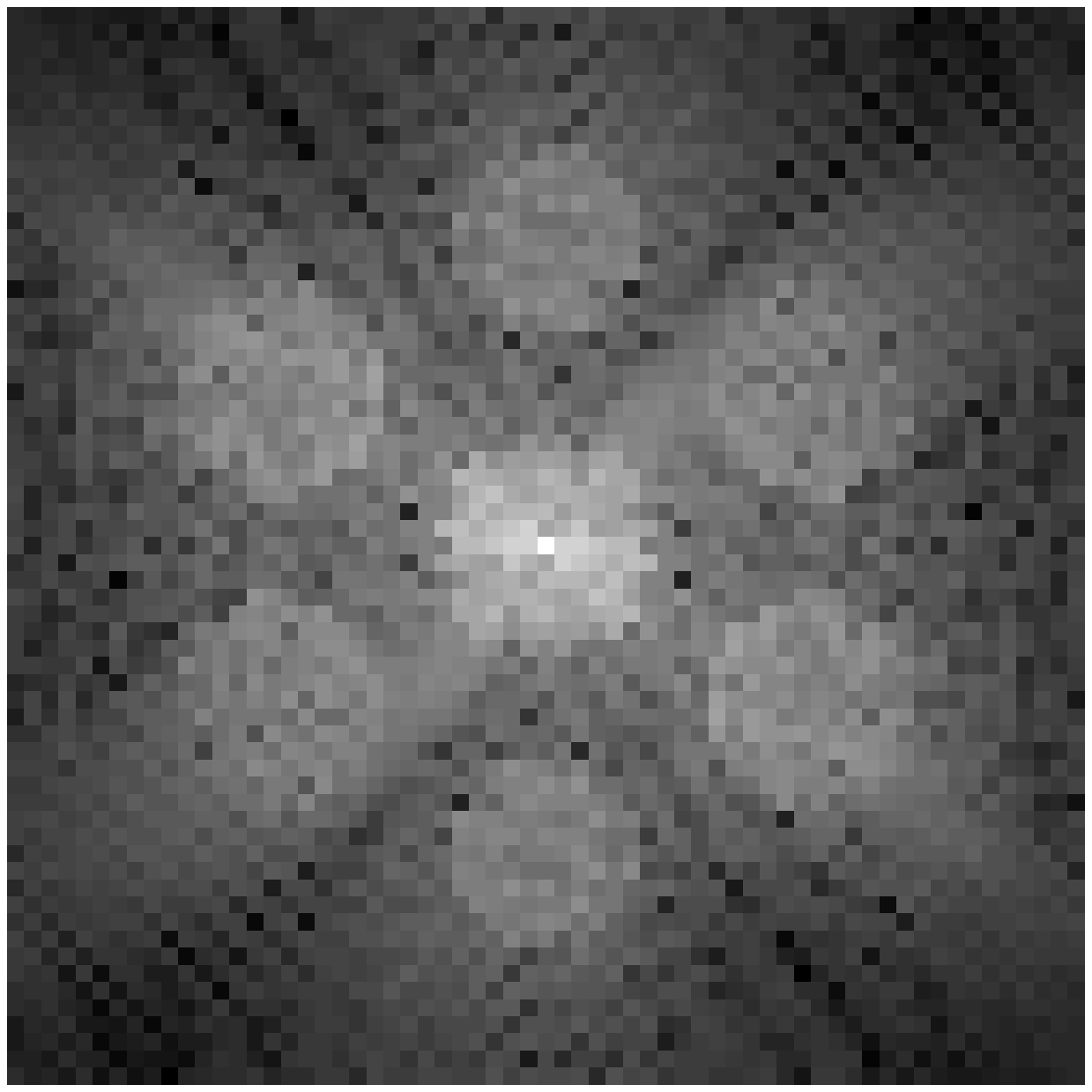}
&\includegraphics[width=0.19\linewidth]{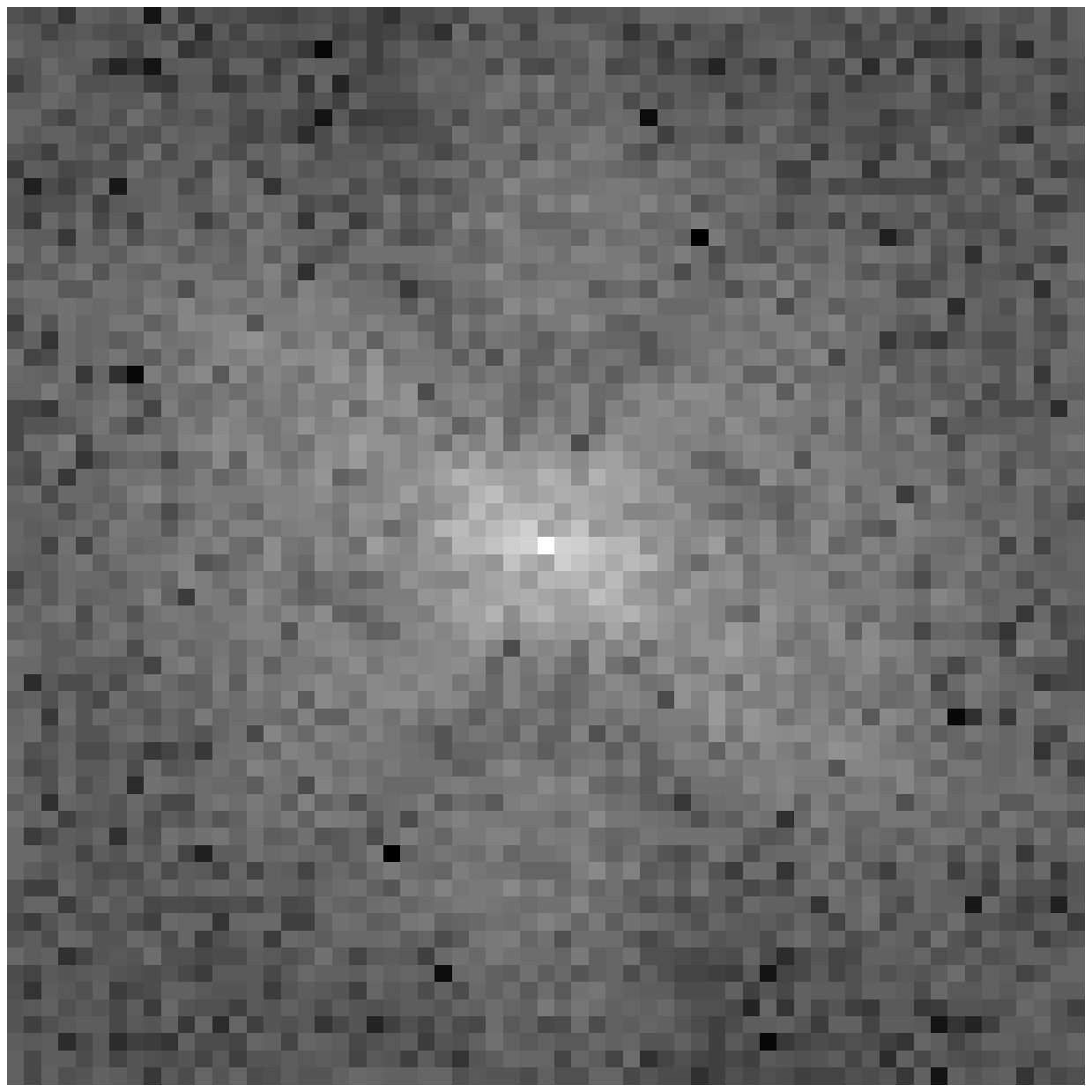}
&\includegraphics[width=0.19\linewidth]{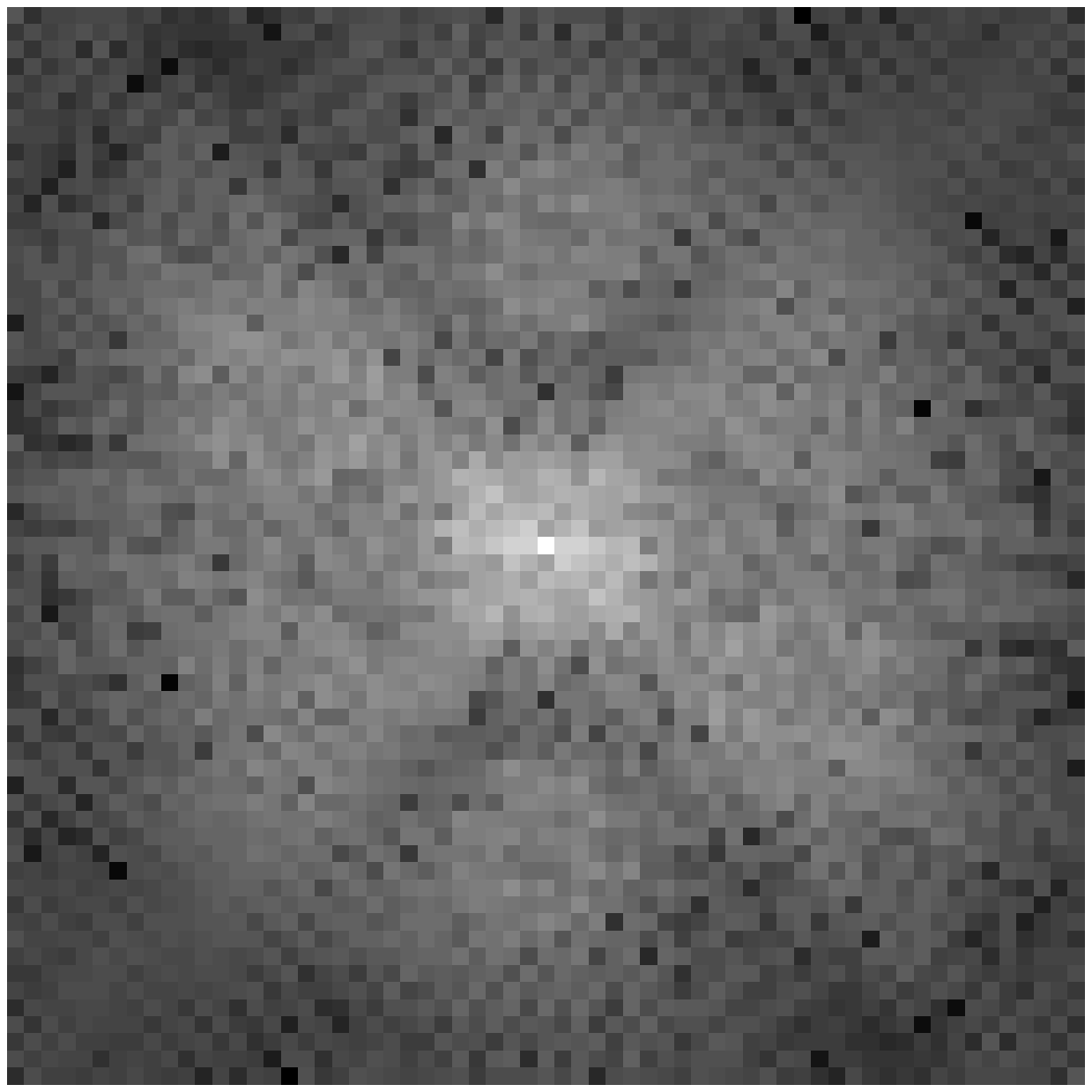}
&\includegraphics[width=0.19\linewidth]{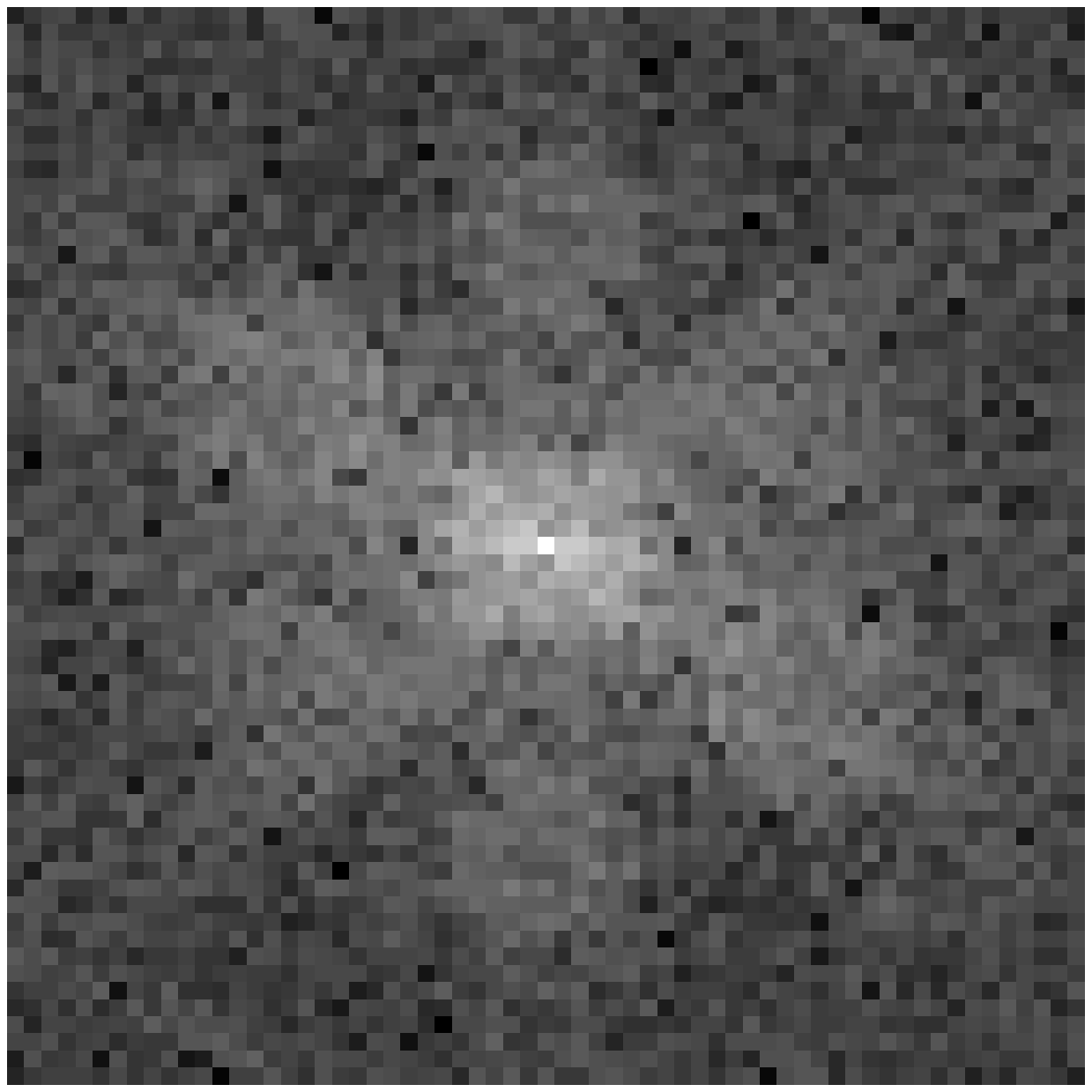}\\
%Spectrum of $u+v_{\delta^3}$ & Spectrum of $u+v_{\delta^1}$ & Spectrum of TV restored\\
\end{tabular}
\vspace{-0.4cm}
\caption%[]
{\newbfb{Corrupted and restored $64\times64$ \textit{Lena} and their respective spectra for exponent values $\alpha\in\{1,2\}$.}}\label{figLen}
% \end{minipage}
\end{figure}

\begin{figure}[H]
\centering
% \begin{minipage}[c]{\linewidth}
\begin{tabular}{ccccc}
Original $g_0$& Corrupted $g$ & SSD - $\delta^3$ & NL-Atom - $\delta^1$ & TV restored\\
 \includegraphics[width=0.19\linewidth]{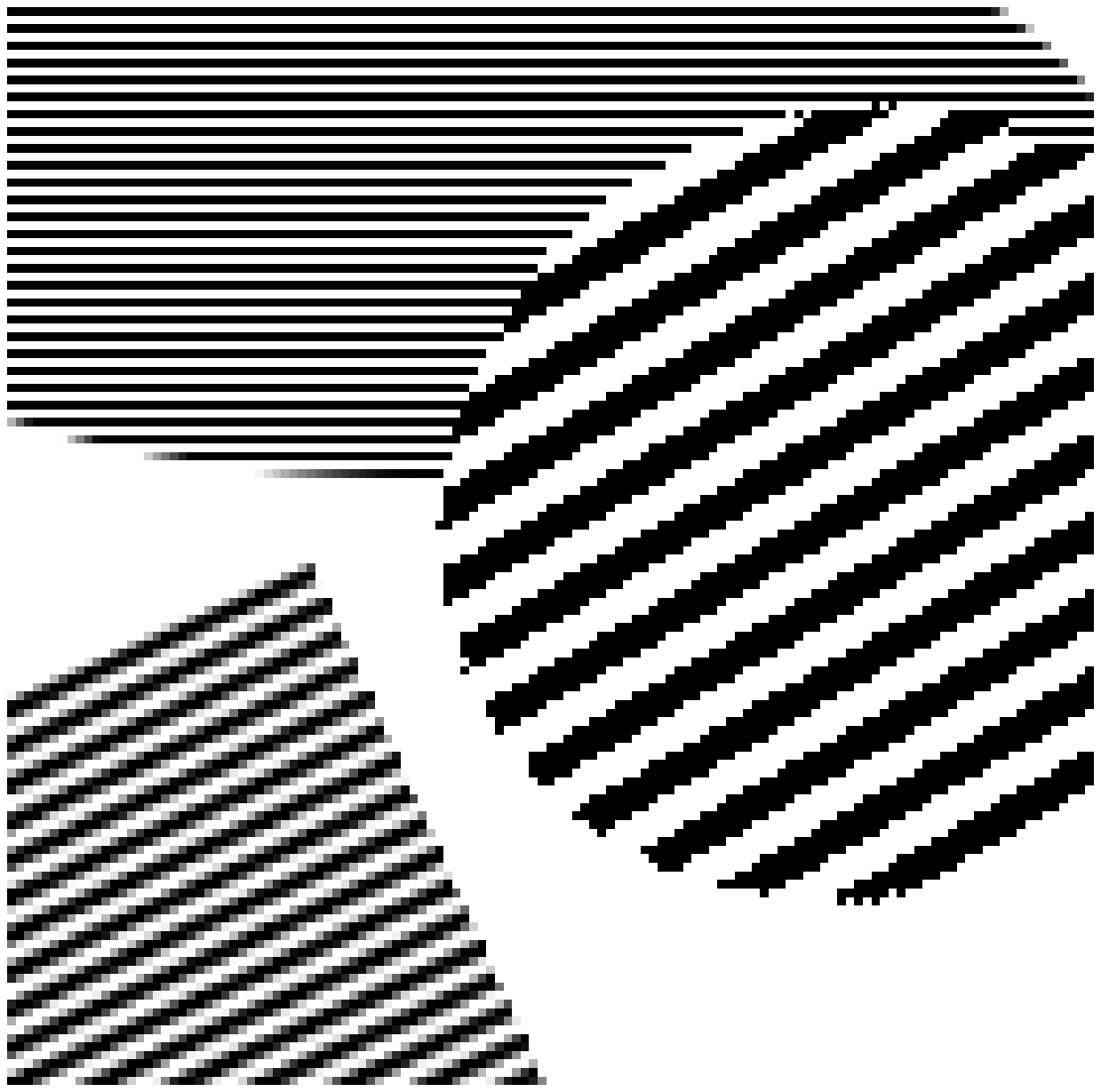}
&\includegraphics[width=0.19\linewidth]{rayures_poids_NL_means_NL_6}
&\includegraphics[width=0.19\linewidth]{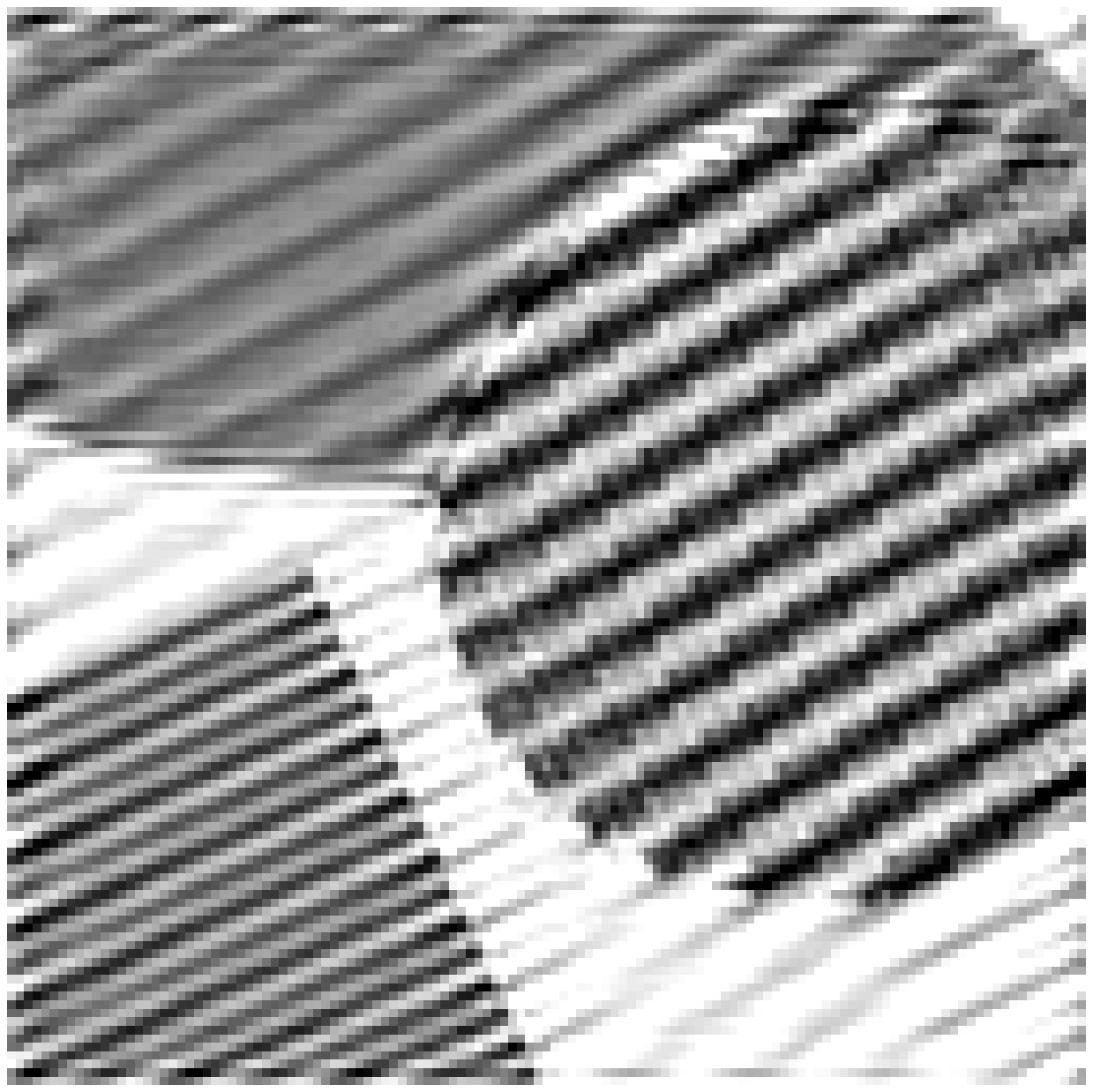}
&\includegraphics[width=0.19\linewidth]{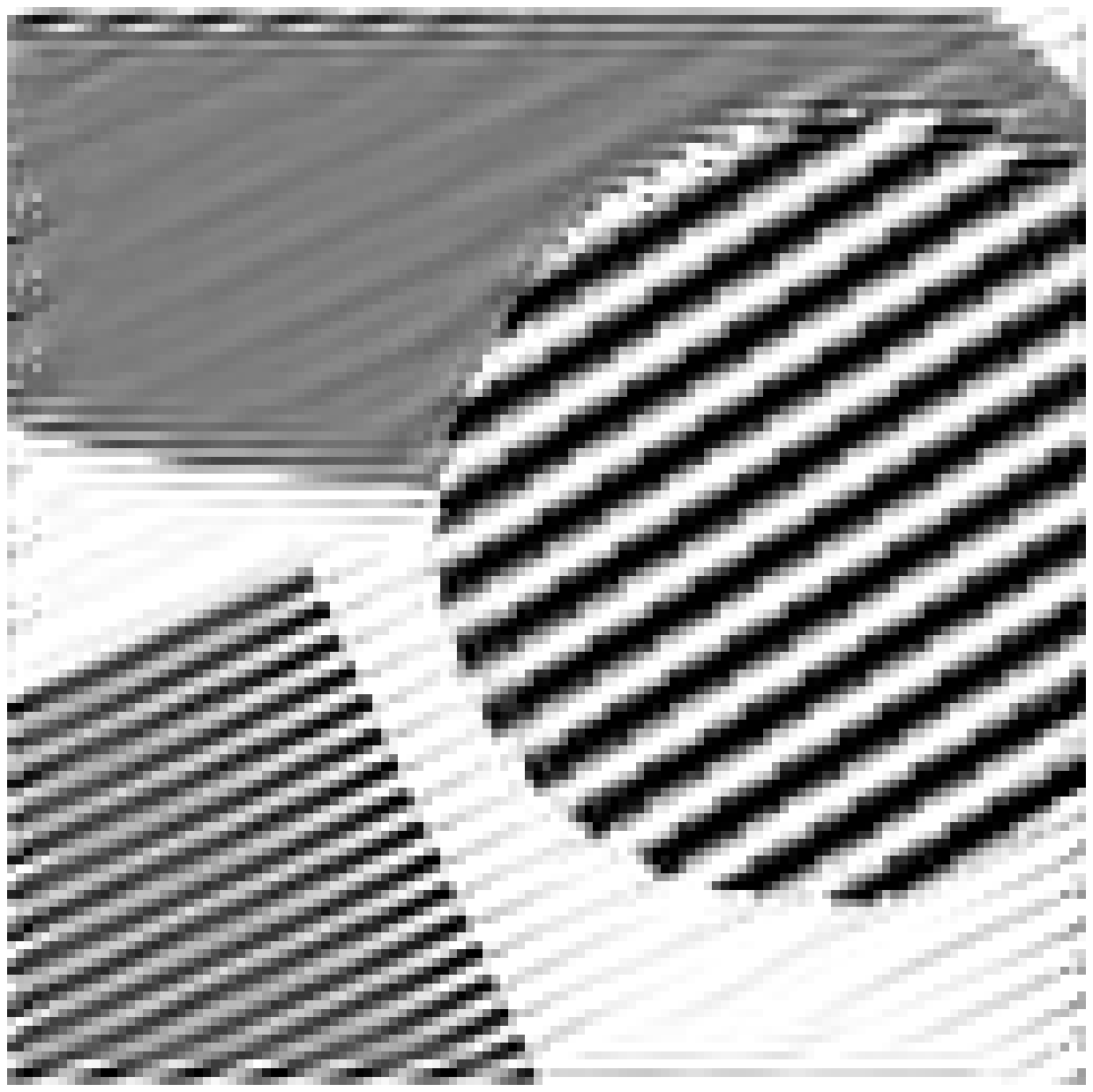}
&\includegraphics[width=0.19\linewidth]{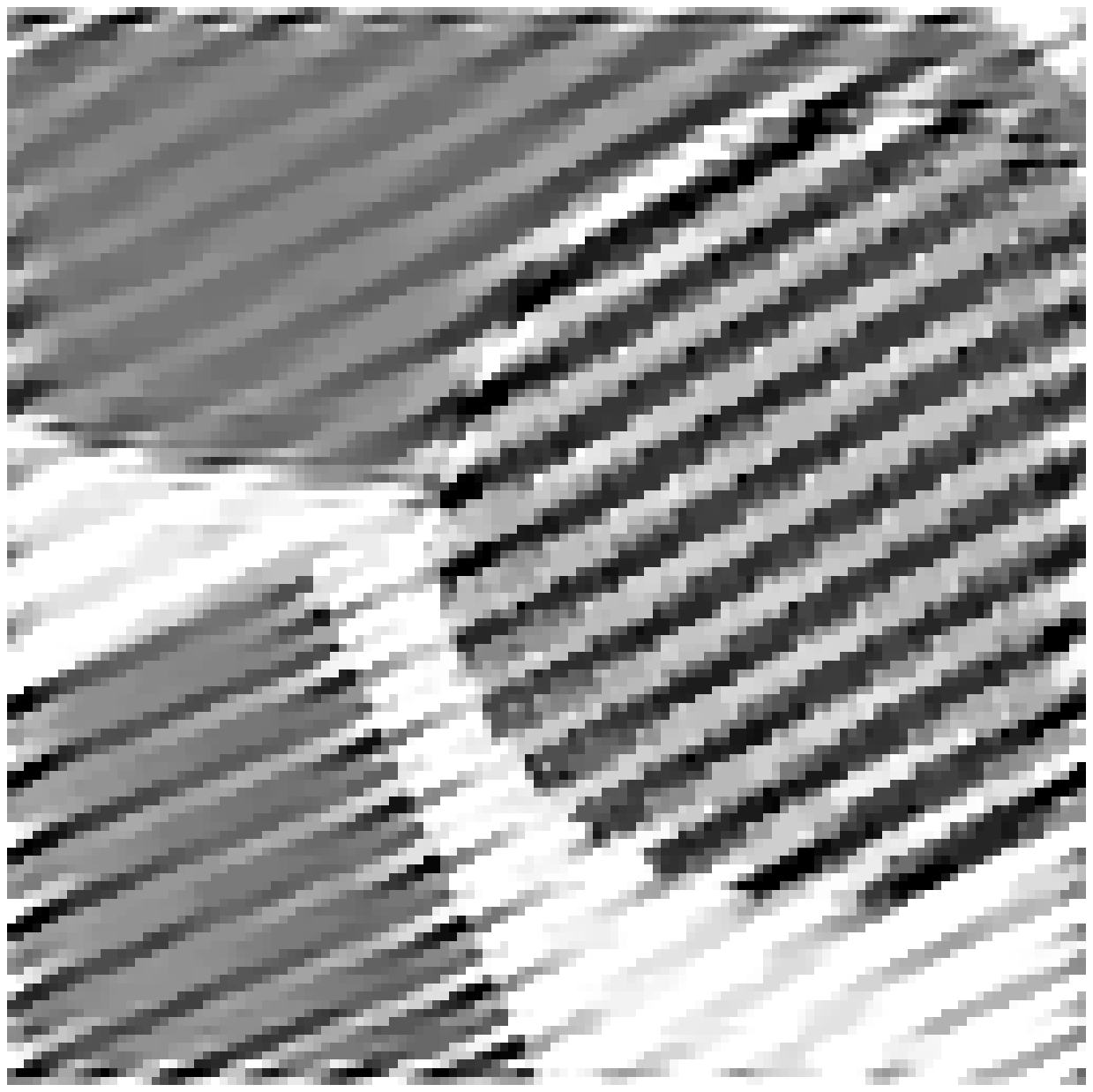}\\
%Corrupted image $u$ & $\delta^3$-NL restored & Atom-based restoration & TV restored\\
 &PSNR=8.5dB &PSNR=9.2dB & PSNR=10.6dB & PSNR=8.4dB\\
\includegraphics[width=0.19\linewidth]{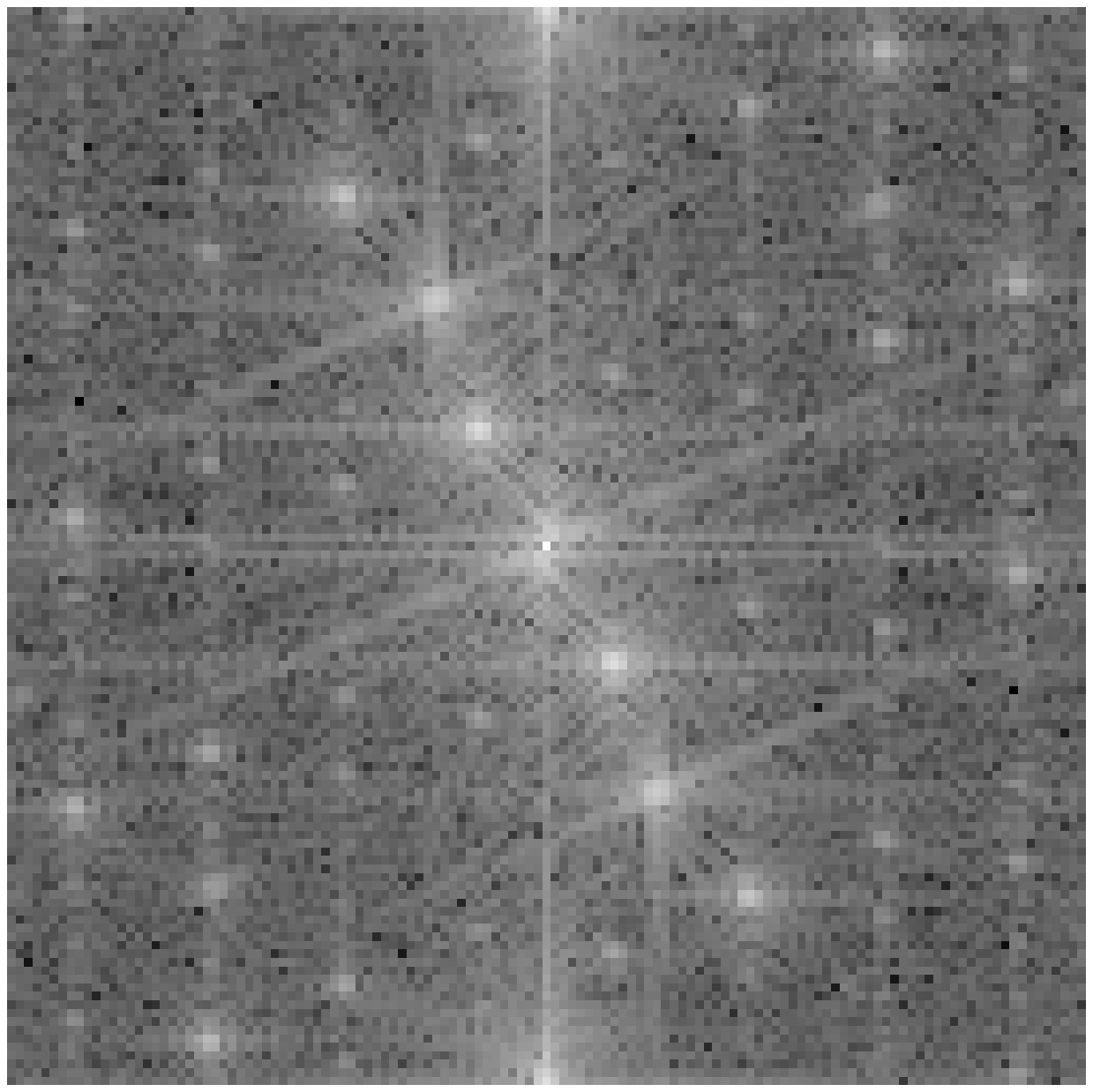}
&\includegraphics[width=0.19\linewidth]{rayures_poids_NL_means_NL_14}
&\includegraphics[width=0.19\linewidth]{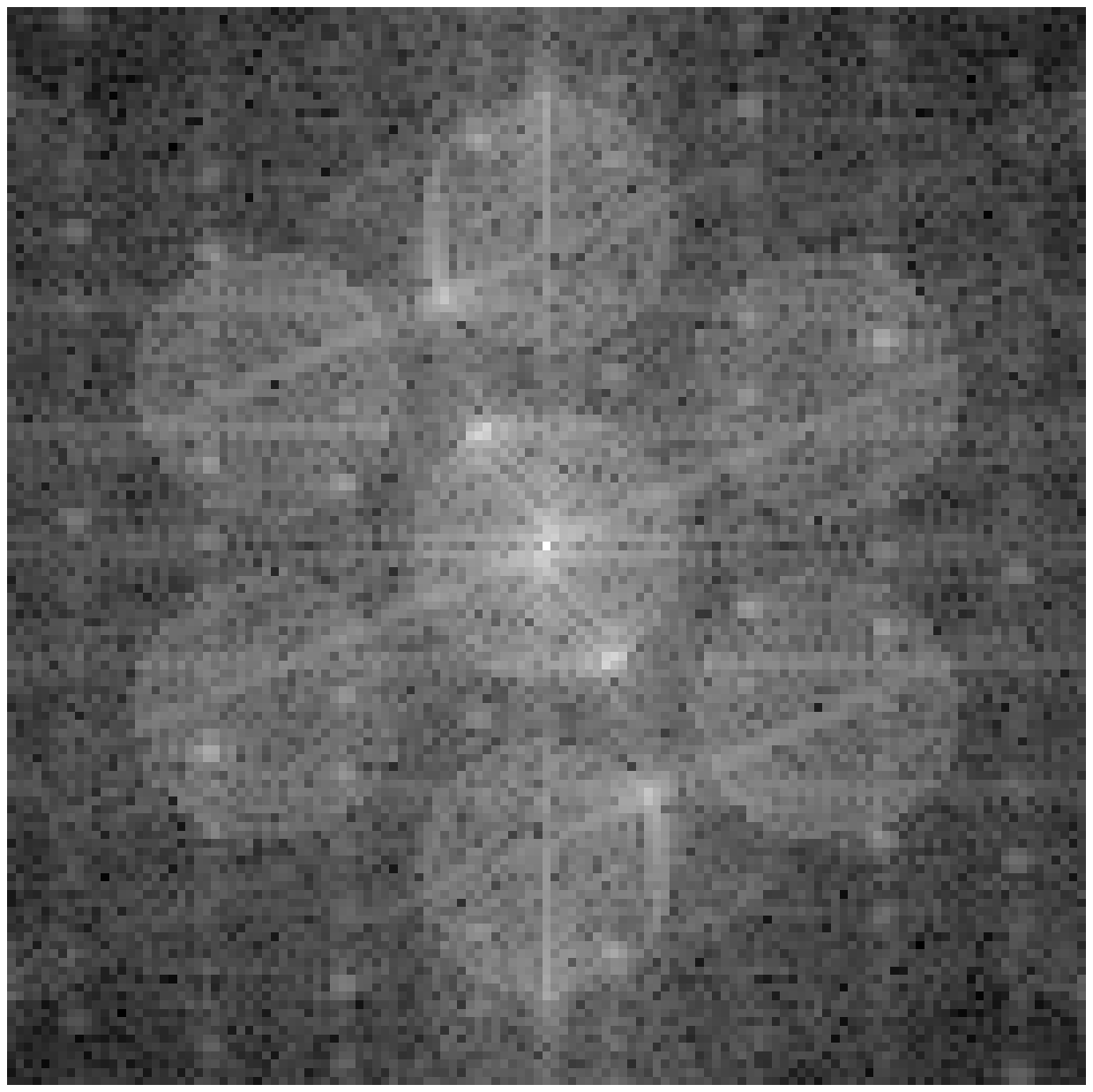}
&\includegraphics[width=0.19\linewidth]{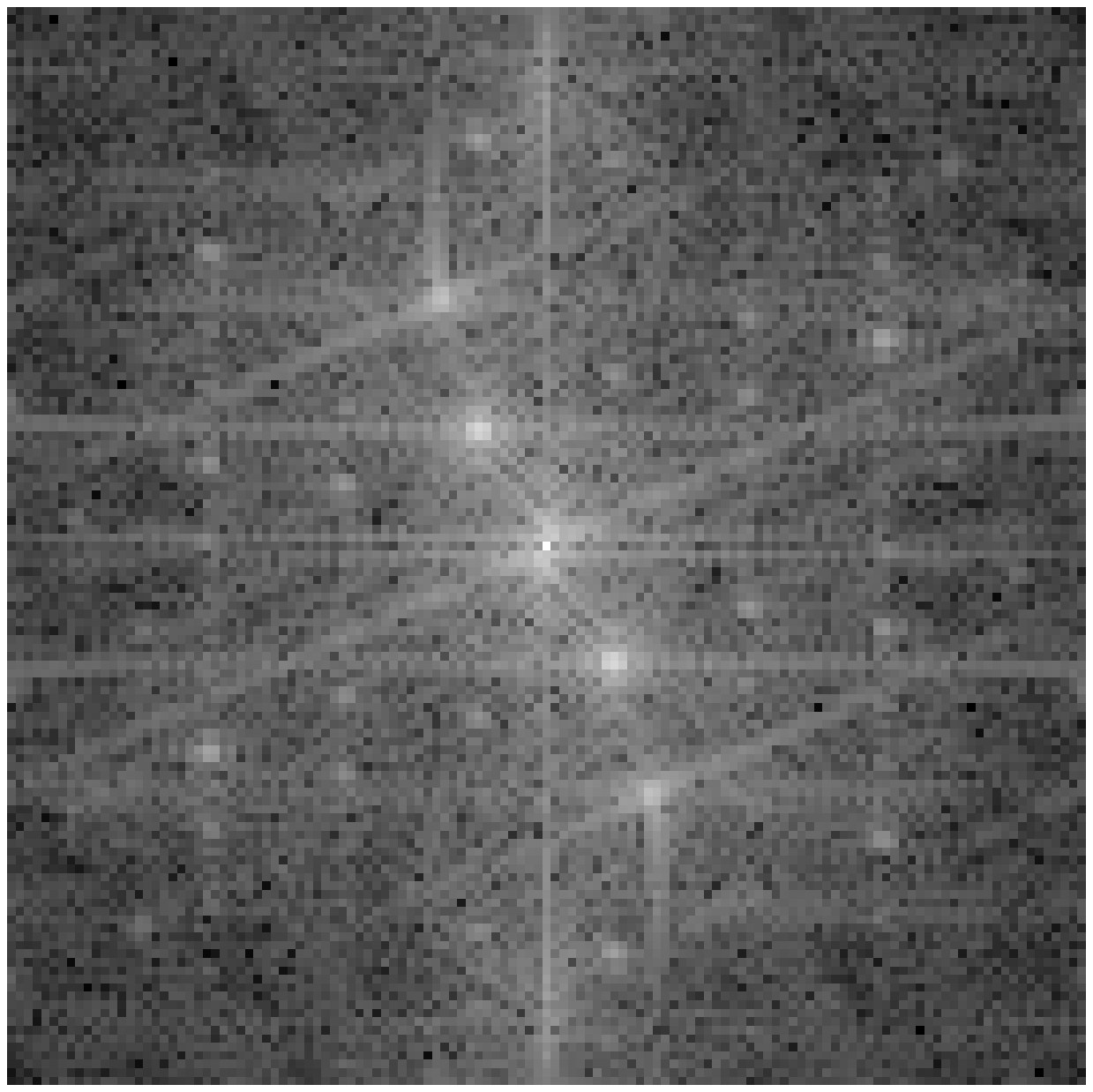}
&\includegraphics[width=0.19\linewidth]{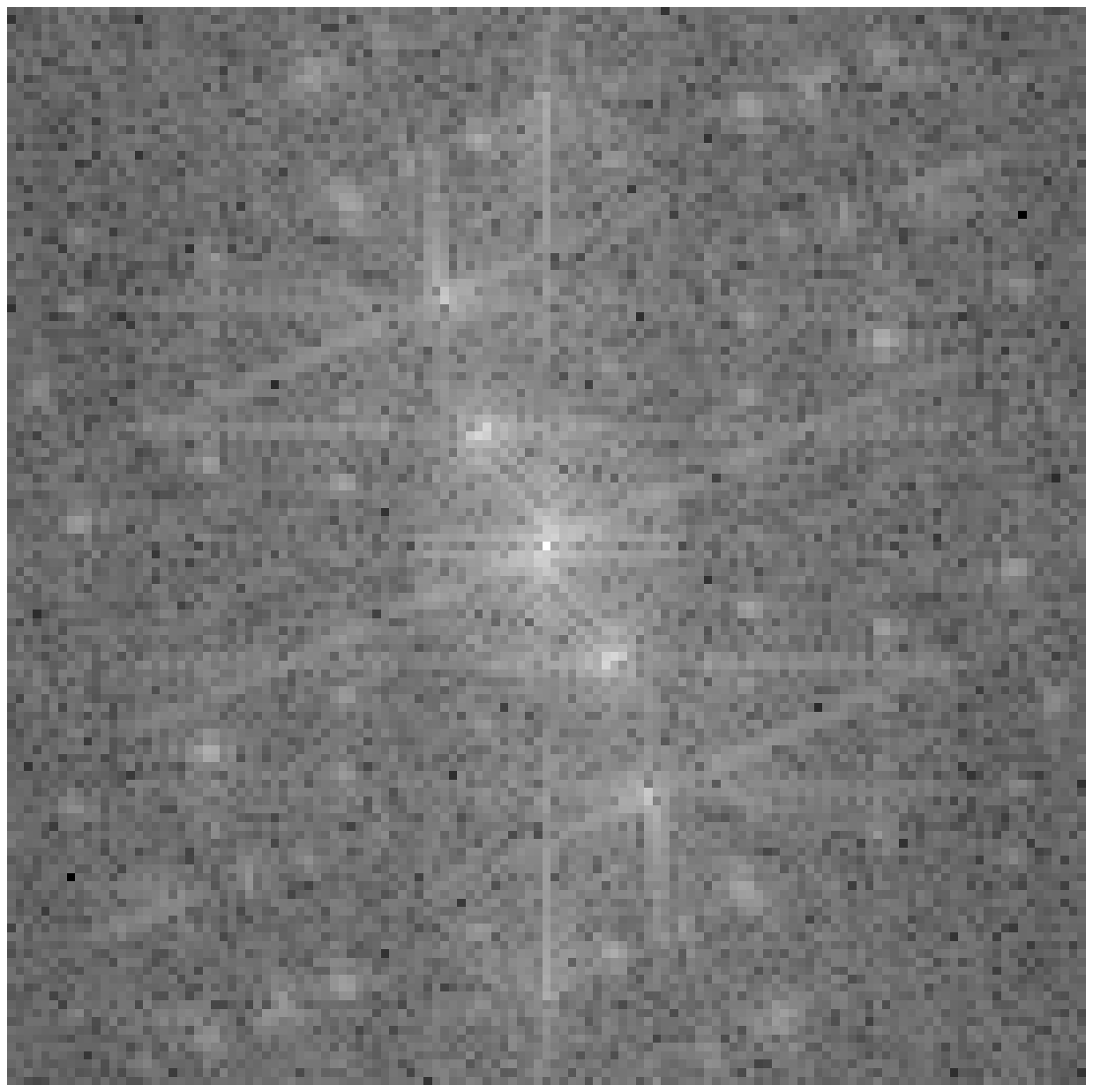}\\
%Spectrum of $u+v_{\delta^3}$ & Spectrum of $u+v_{\delta^1}$ & Spectrum of TV restored\\
\end{tabular}
\vspace{-0.4cm}
\caption%[]
{Corrupted and restored $128\times128$ toy example and their respective spectra.}\label{figToy}
% \end{minipage}
\end{figure}

\begin{figure}[H]
\centering
% \begin{minipage}[c]{\linewidth}
\begin{tabular}{ccccc}
Original $g_0$ & Corrupted $g$ & SSD - $\delta^3$ & NL-Atom - $\delta^1$ & TV restored\\
\includegraphics[width=0.19\linewidth]{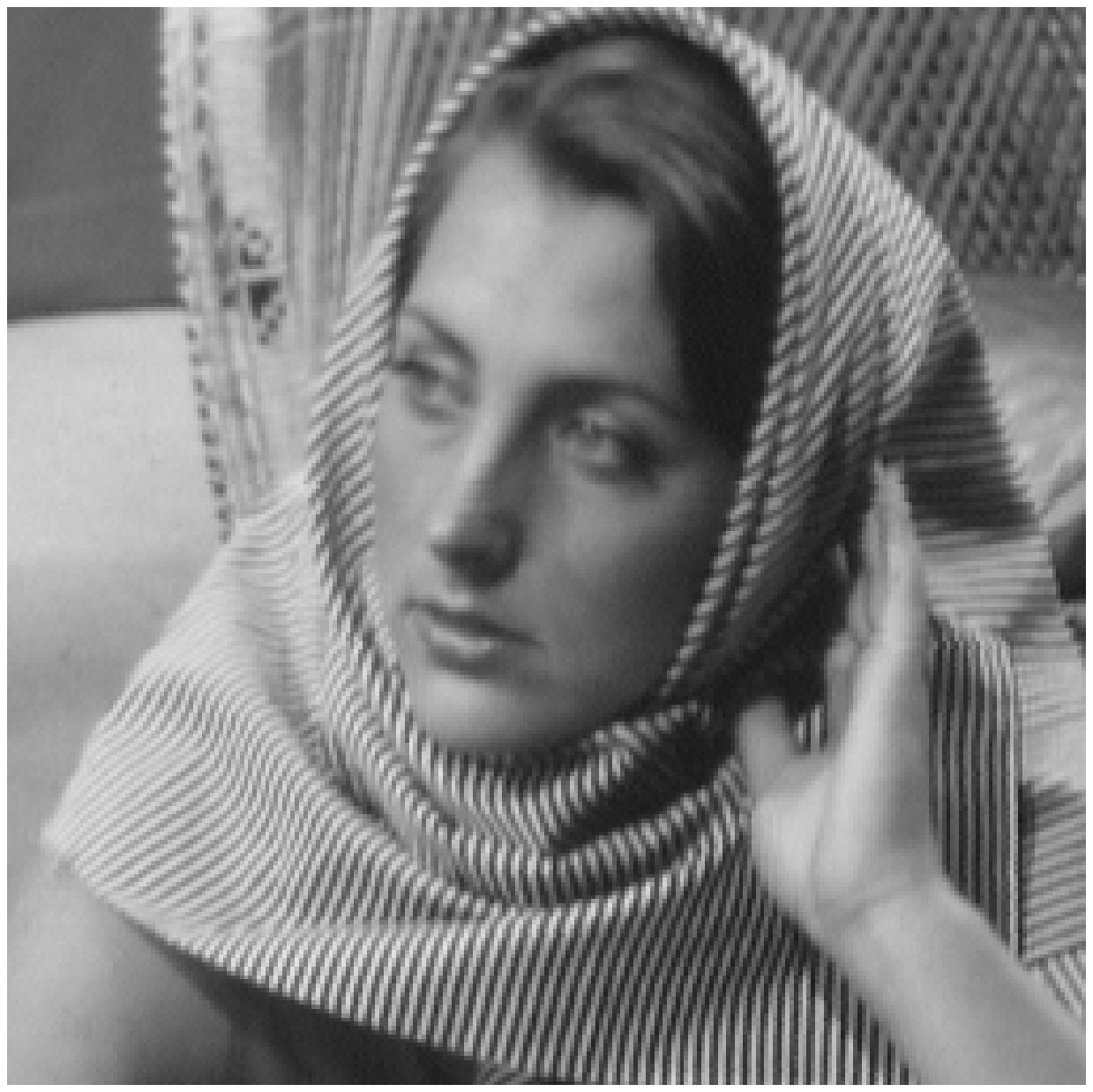}
&\includegraphics[width=0.19\linewidth]{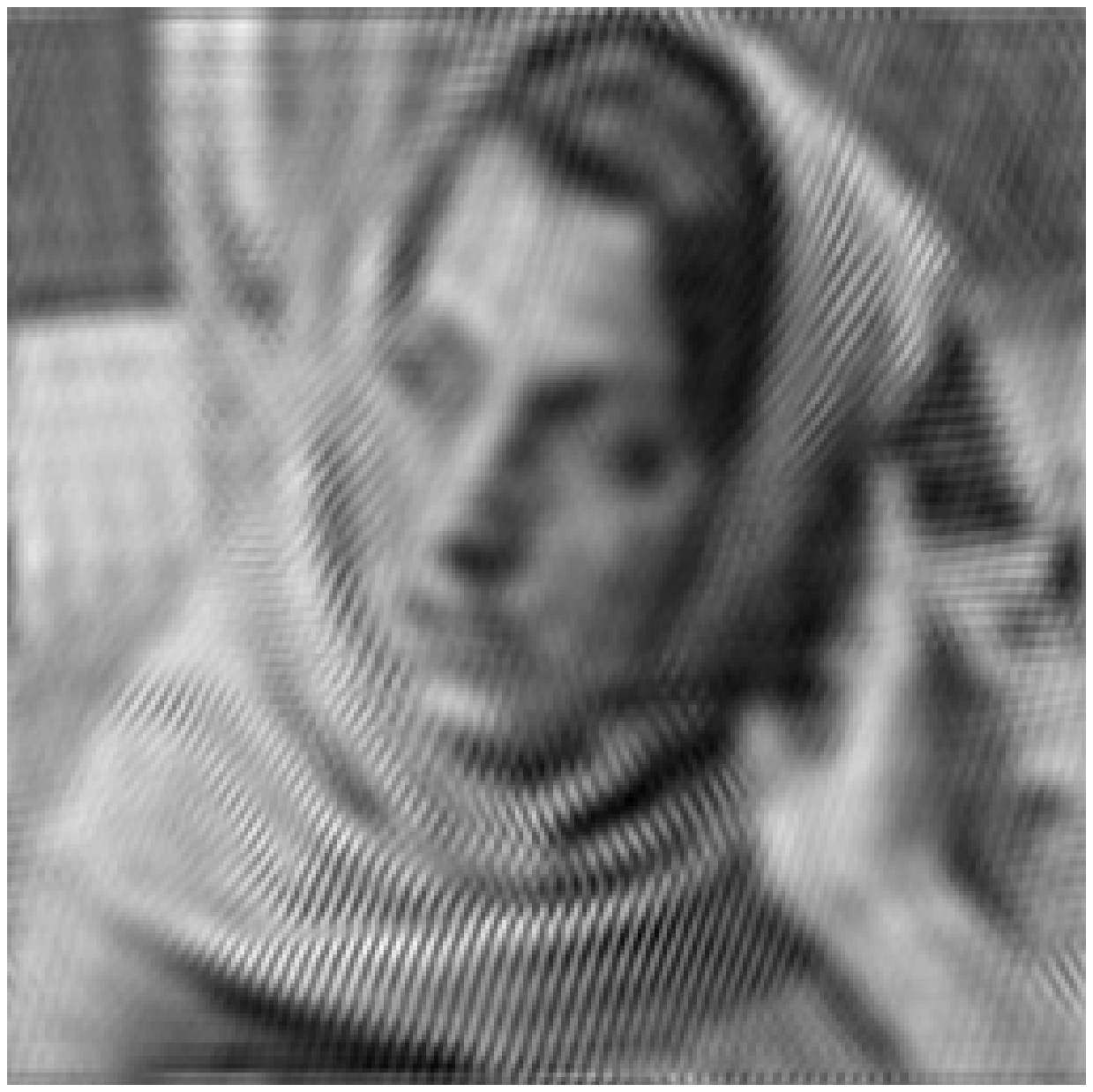}
&\includegraphics[width=0.19\linewidth]{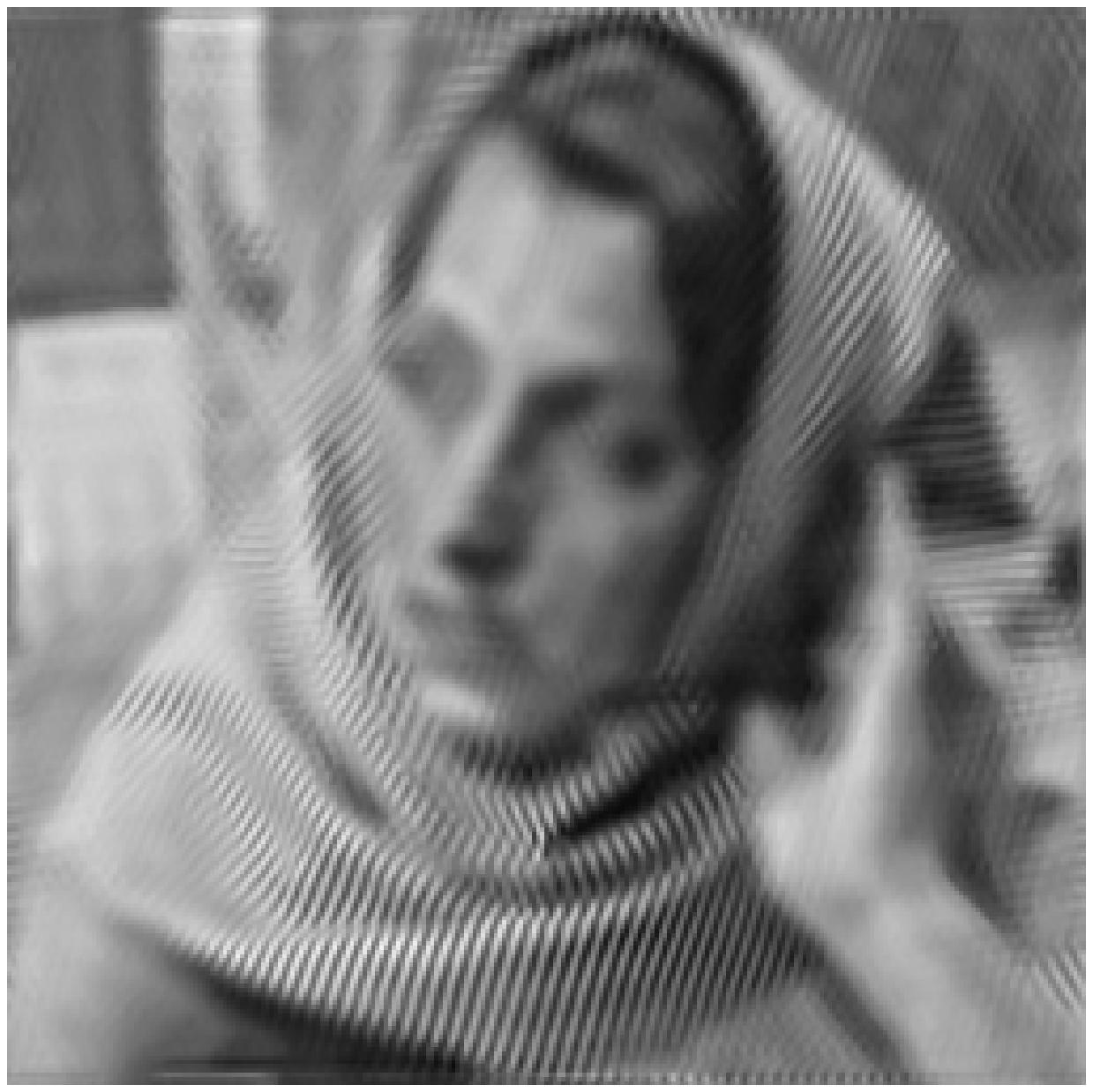}
&\includegraphics[width=0.19\linewidth]{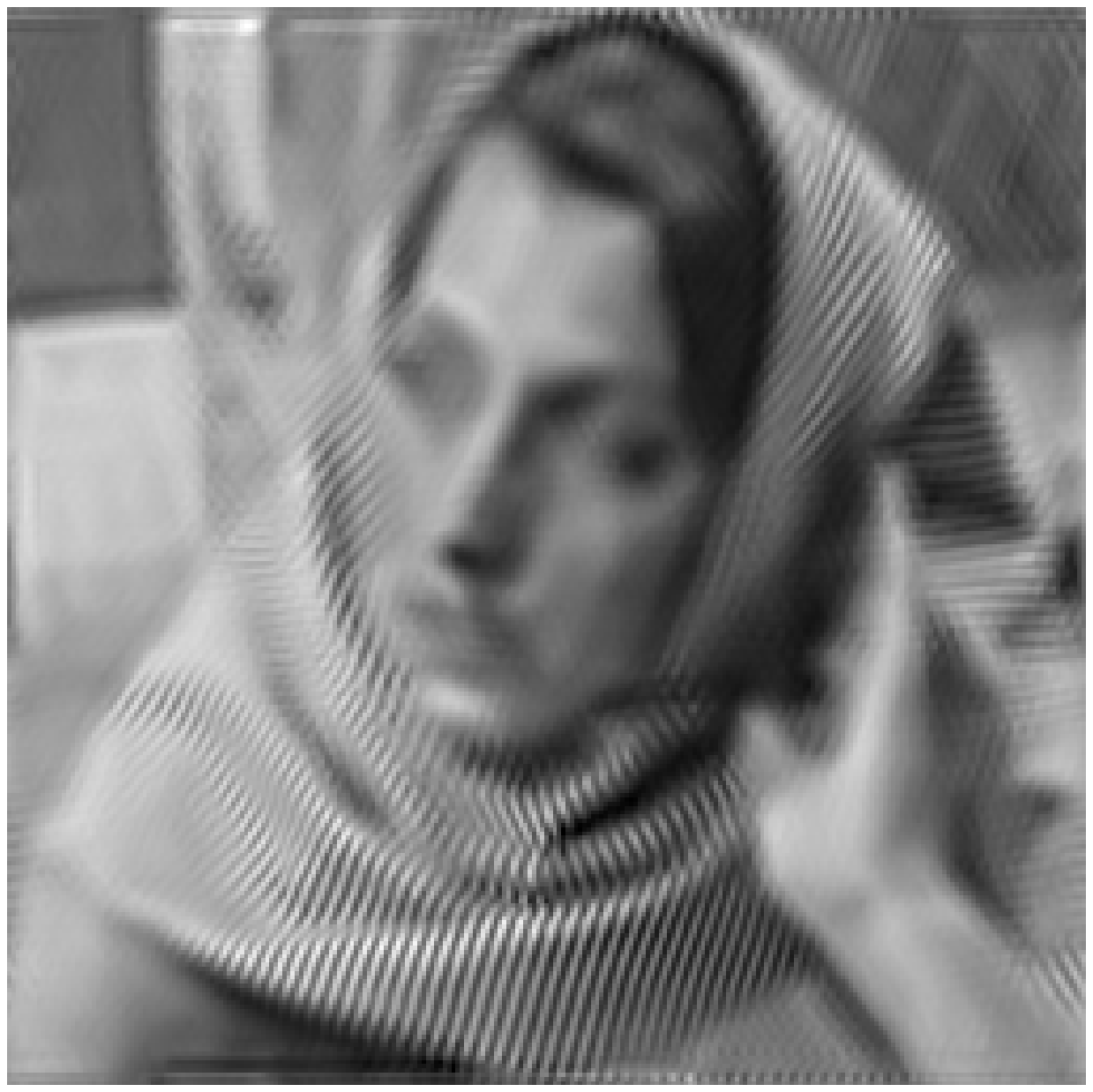}
&\includegraphics[width=0.19\linewidth]{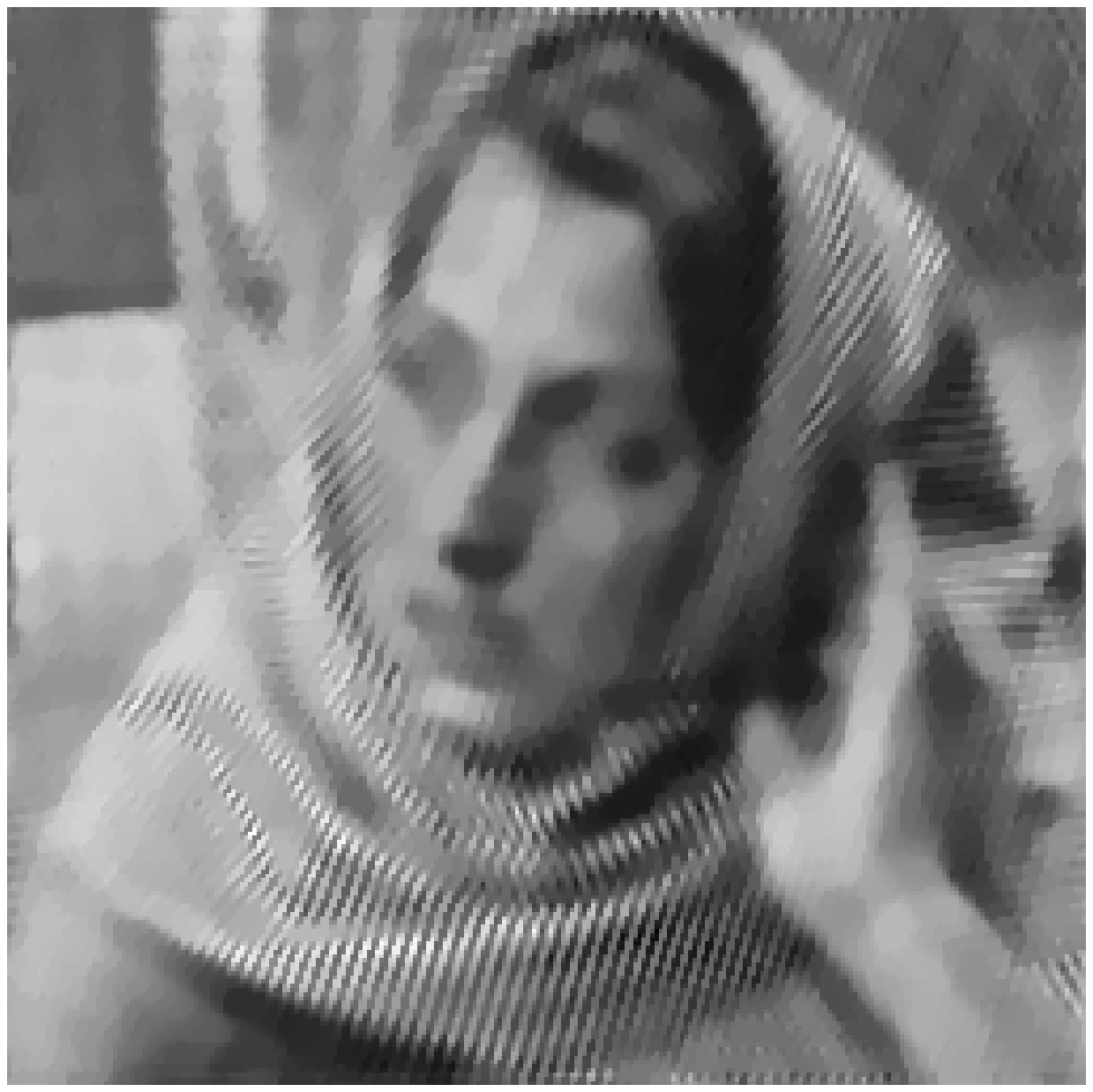}\\
%Corrupted image $u$ & $\delta^3$-NL restored & Atom-based restoration & TV restored\\
 & PSNR=24.0dB & PSNR=24.6dB & PSNR=24.7dB & PSNR=24.1dB\\
\includegraphics[width=0.19\linewidth]{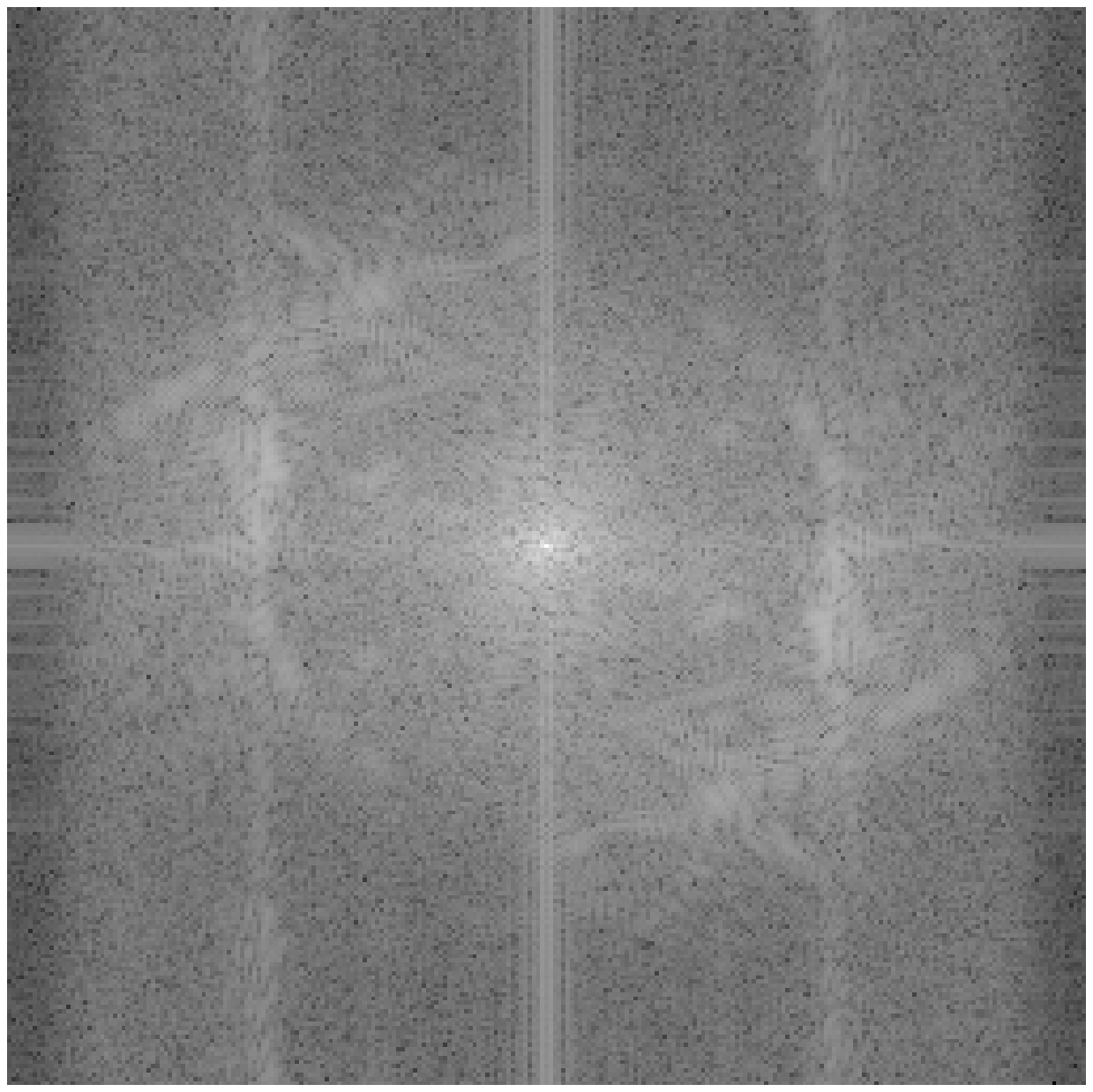}
&\includegraphics[width=0.19\linewidth]{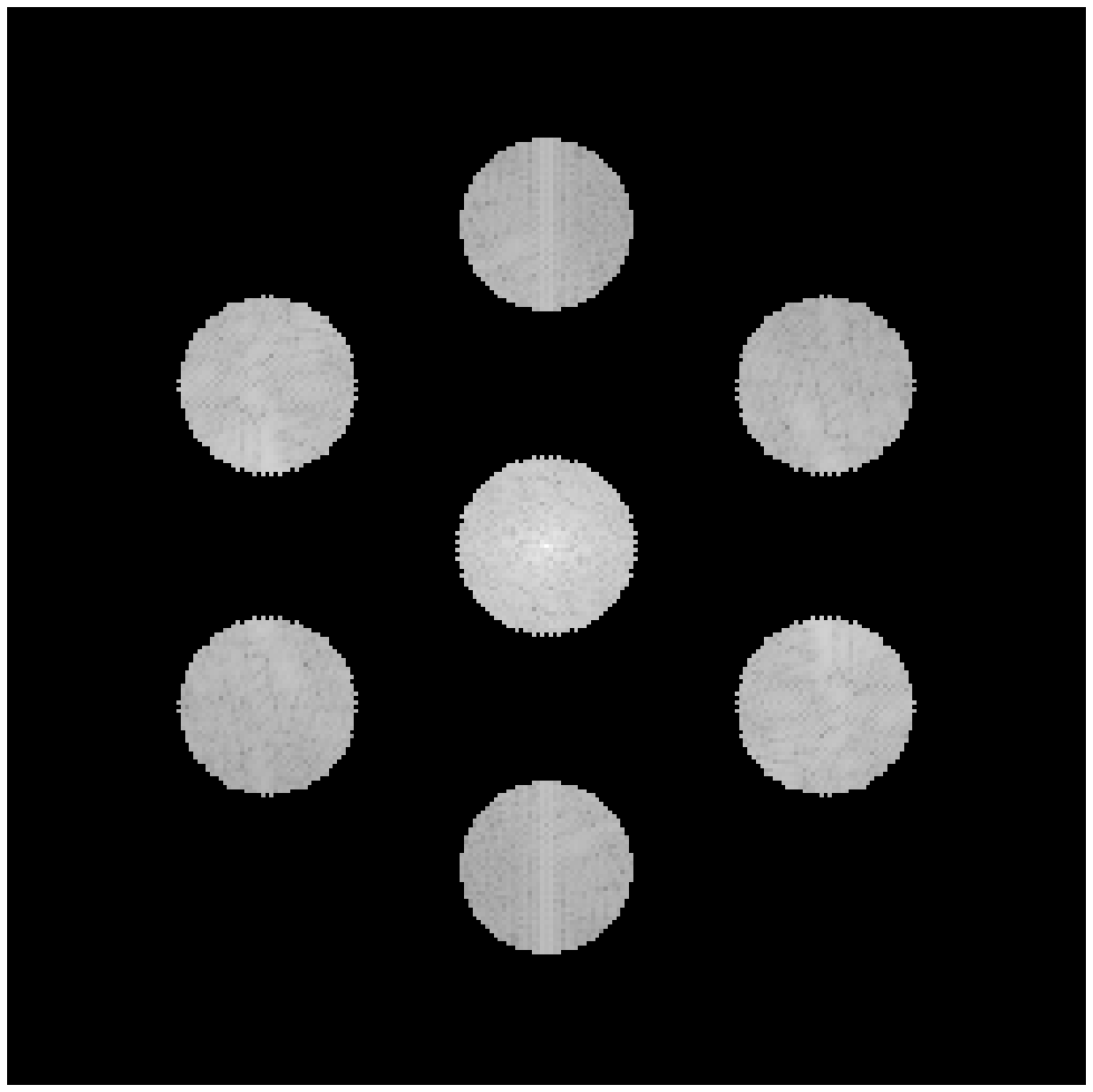}
&\includegraphics[width=0.19\linewidth]{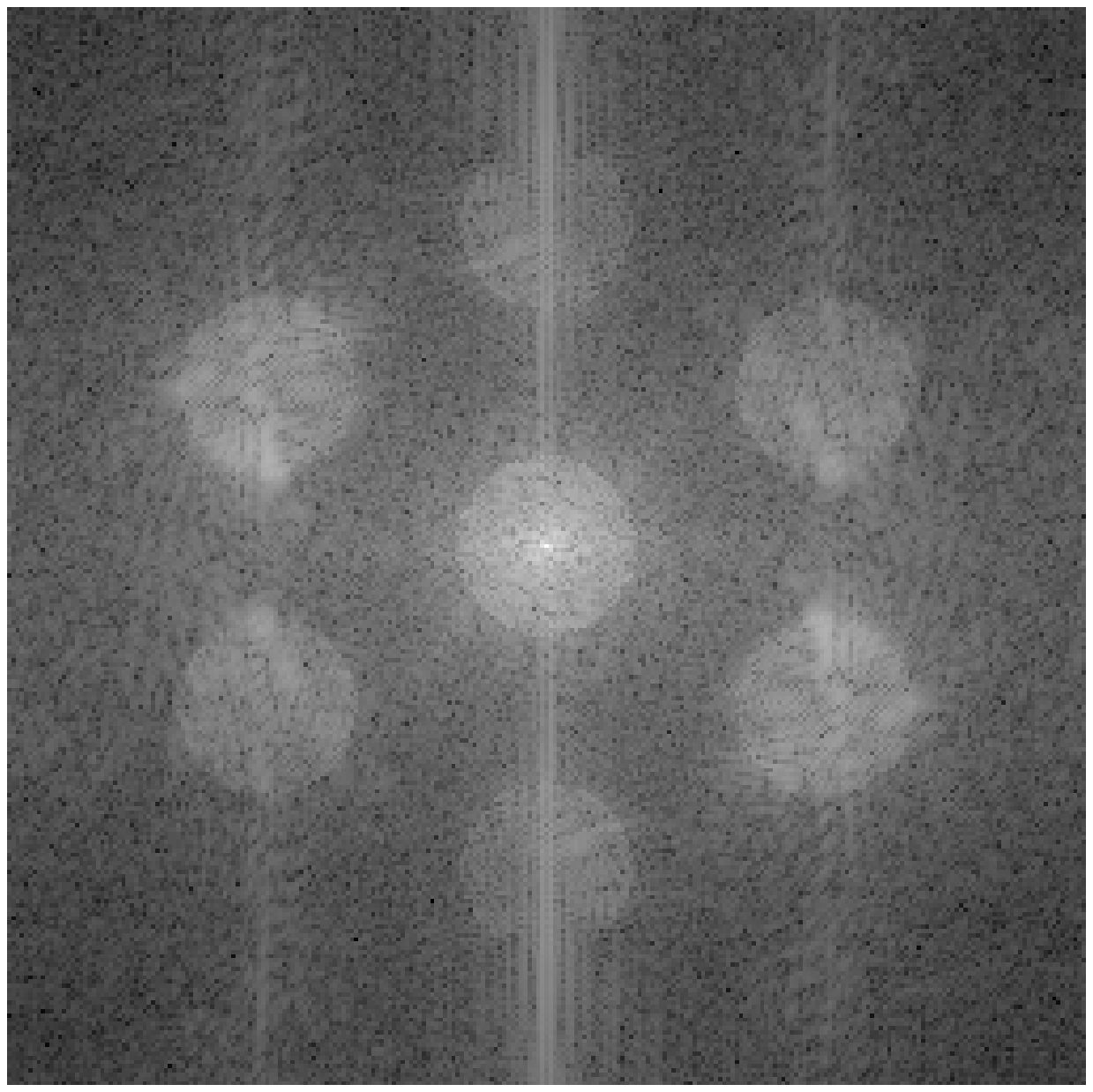}
&\includegraphics[width=0.19\linewidth]{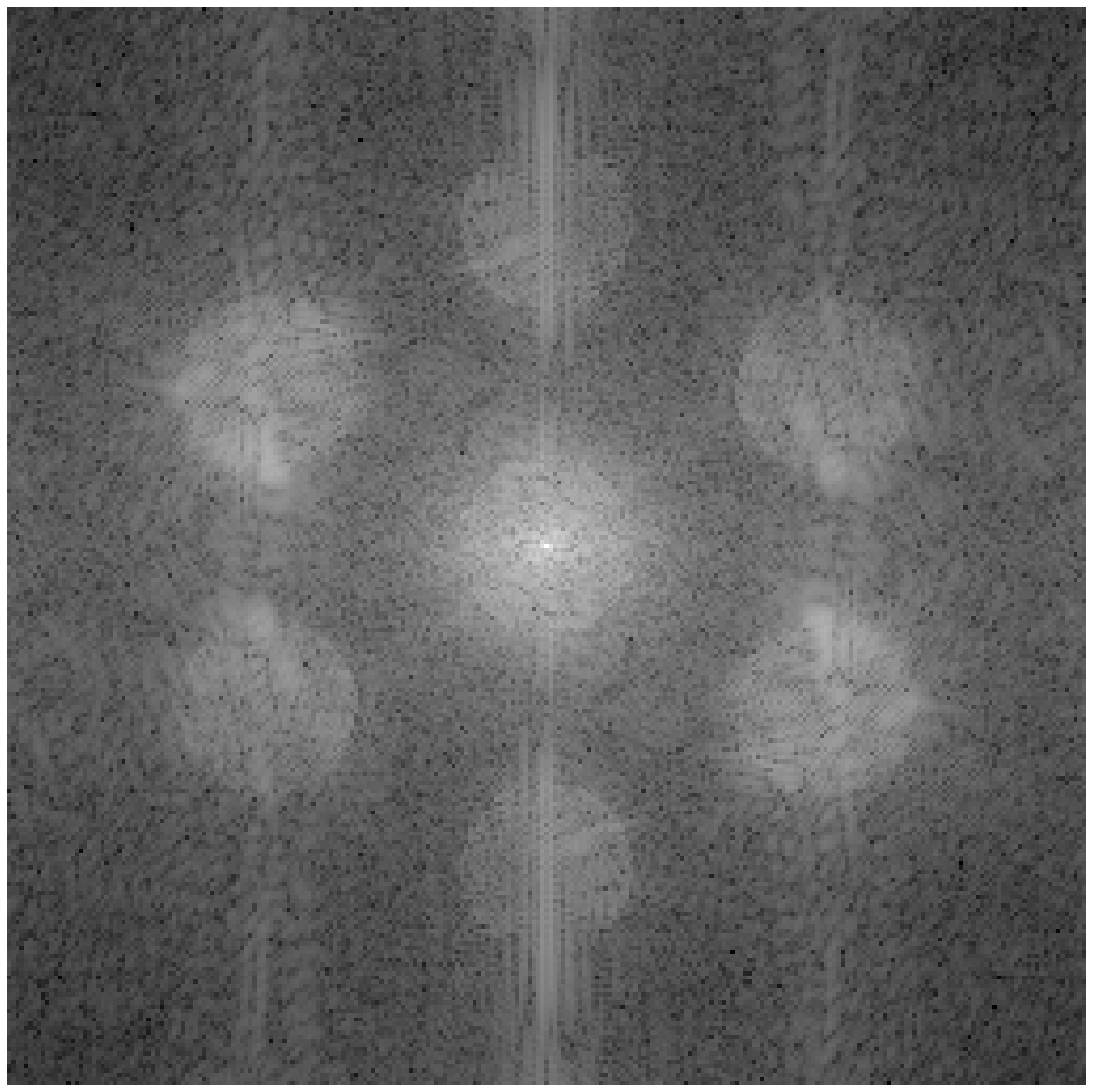}
&\includegraphics[width=0.19\linewidth]{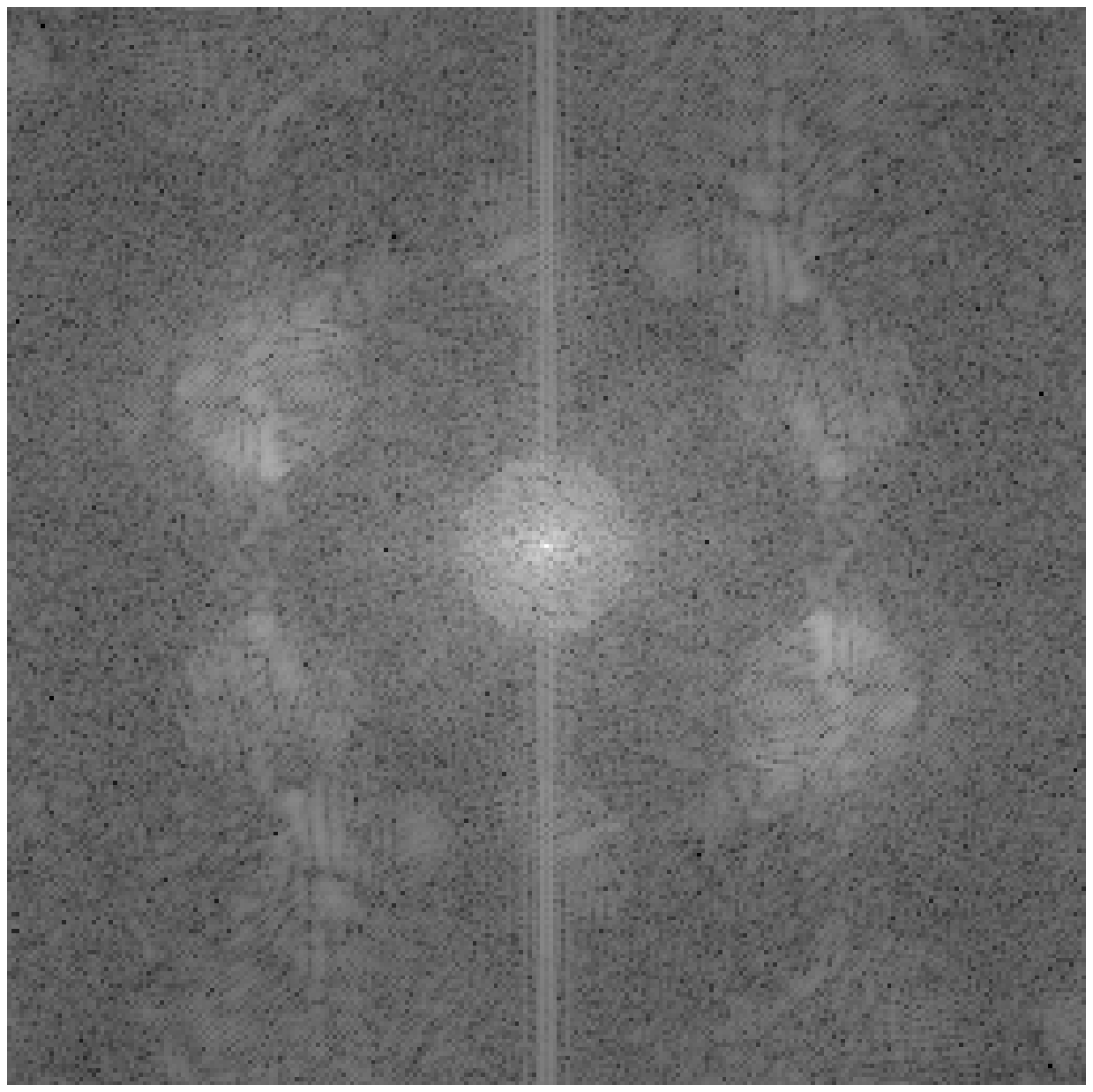}
%Spectrum of $u+v_{\delta^3}$ & Spectrum of $u+v_{\delta^1}$ & Spectrum of TV restored\\
\end{tabular}
\vspace{-0.4cm}
\caption%[]
{Corrupted and restored $256\times256$ crop of \textit{Barbara} and their respective spectra.}\label{figBarb}
% \end{minipage}
\end{figure}

\begin{figure}[H]
\centering
% \begin{minipage}[c]{\linewidth}
\begin{tabular}{ccccc}
Original $g_0$ & Corrupted $g$ & SSD - $\delta^3$ & NL-Atom - $\delta^1$ & TV restored\\
\includegraphics[width=0.19\linewidth]{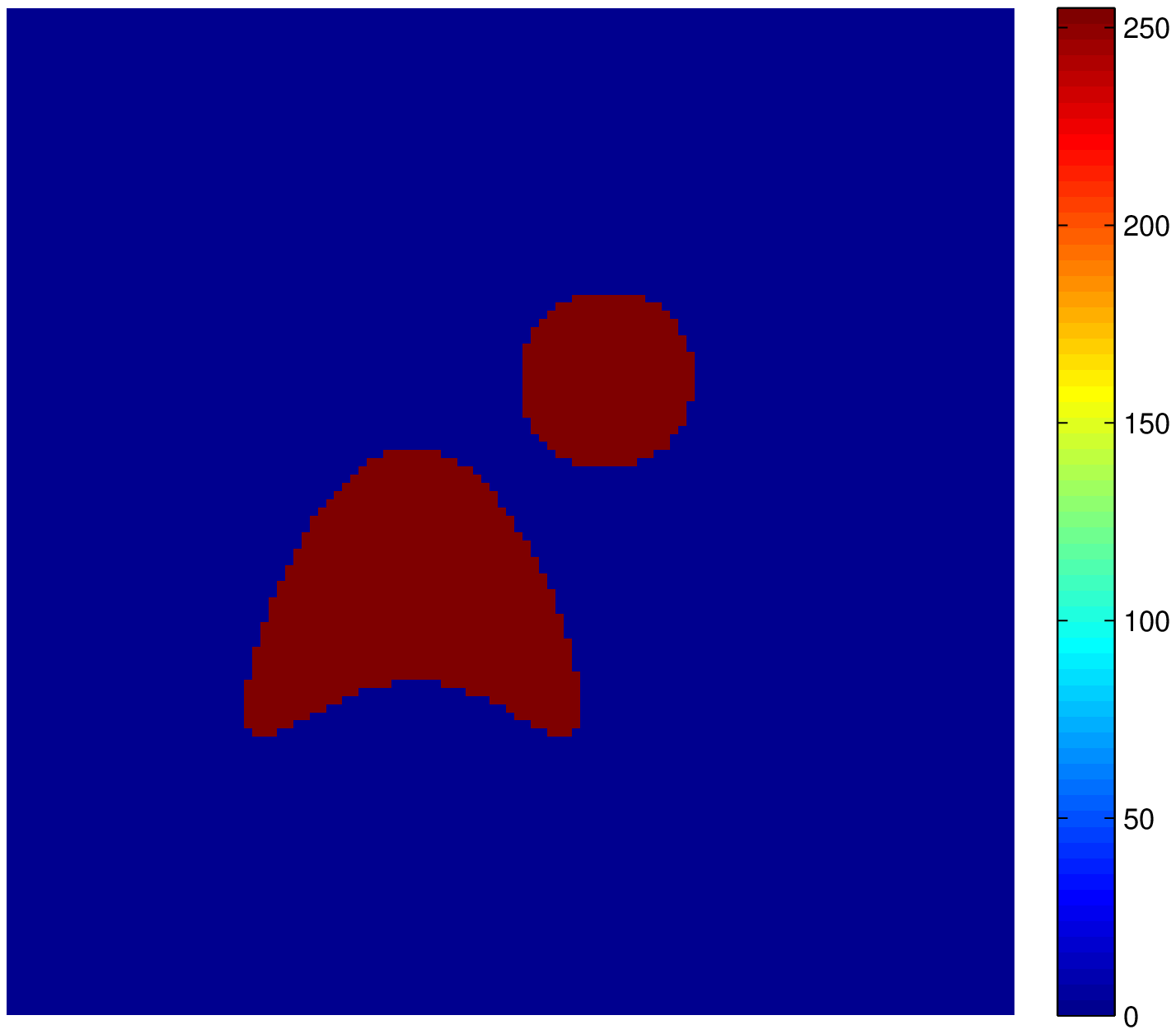}
&\includegraphics[width=0.19\linewidth]{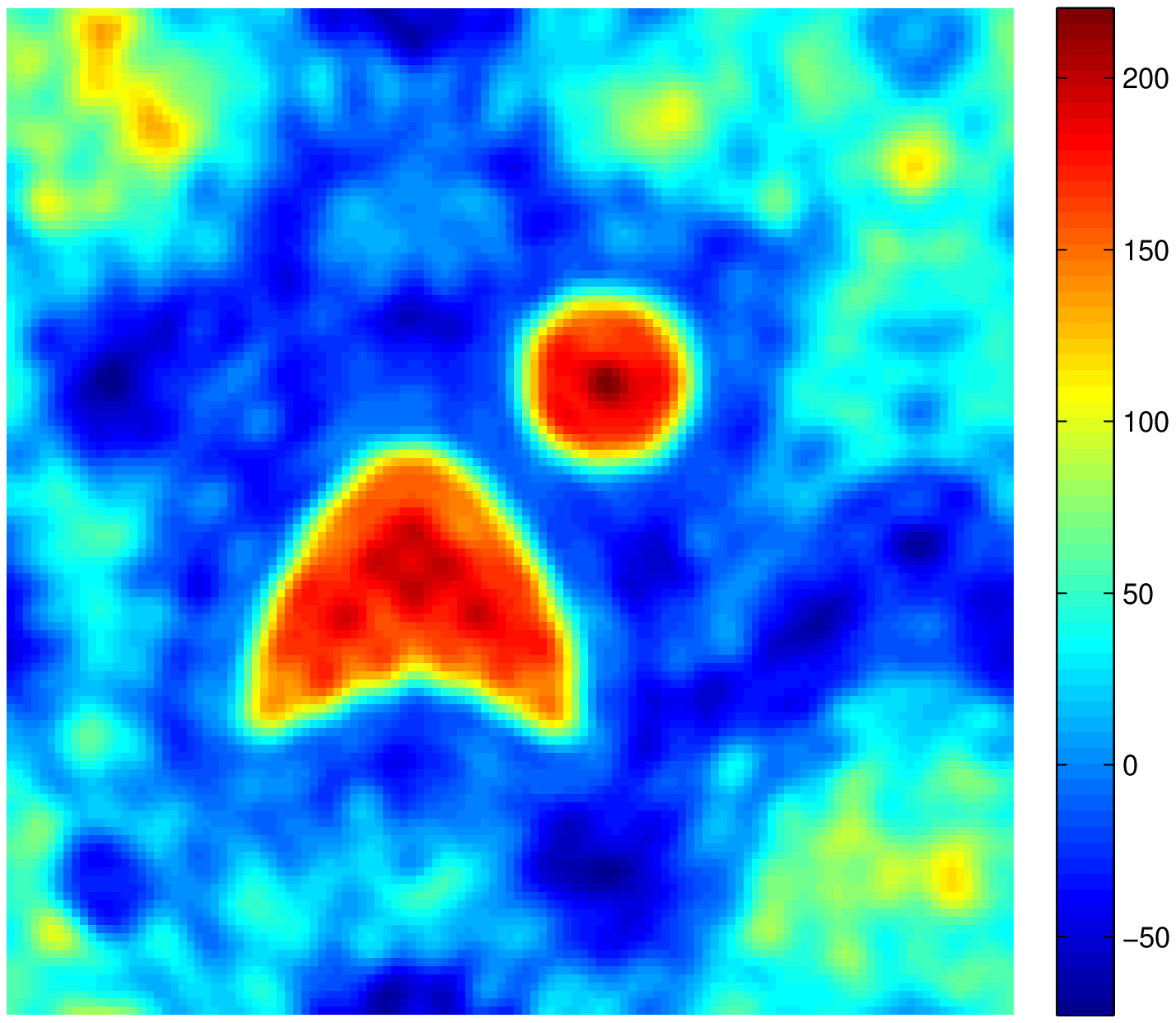}
&\includegraphics[width=0.19\linewidth]{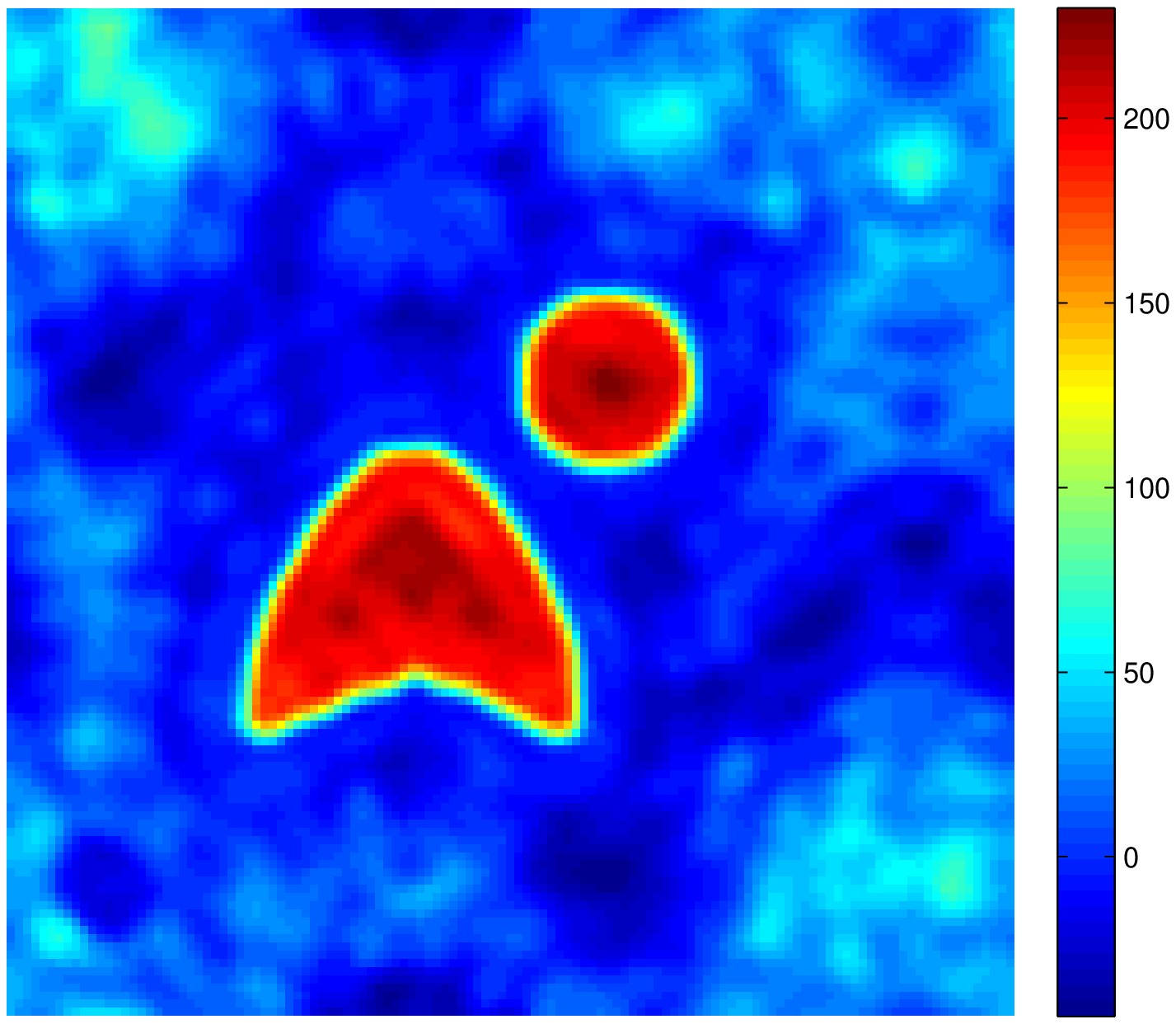}
&\includegraphics[width=0.19\linewidth]{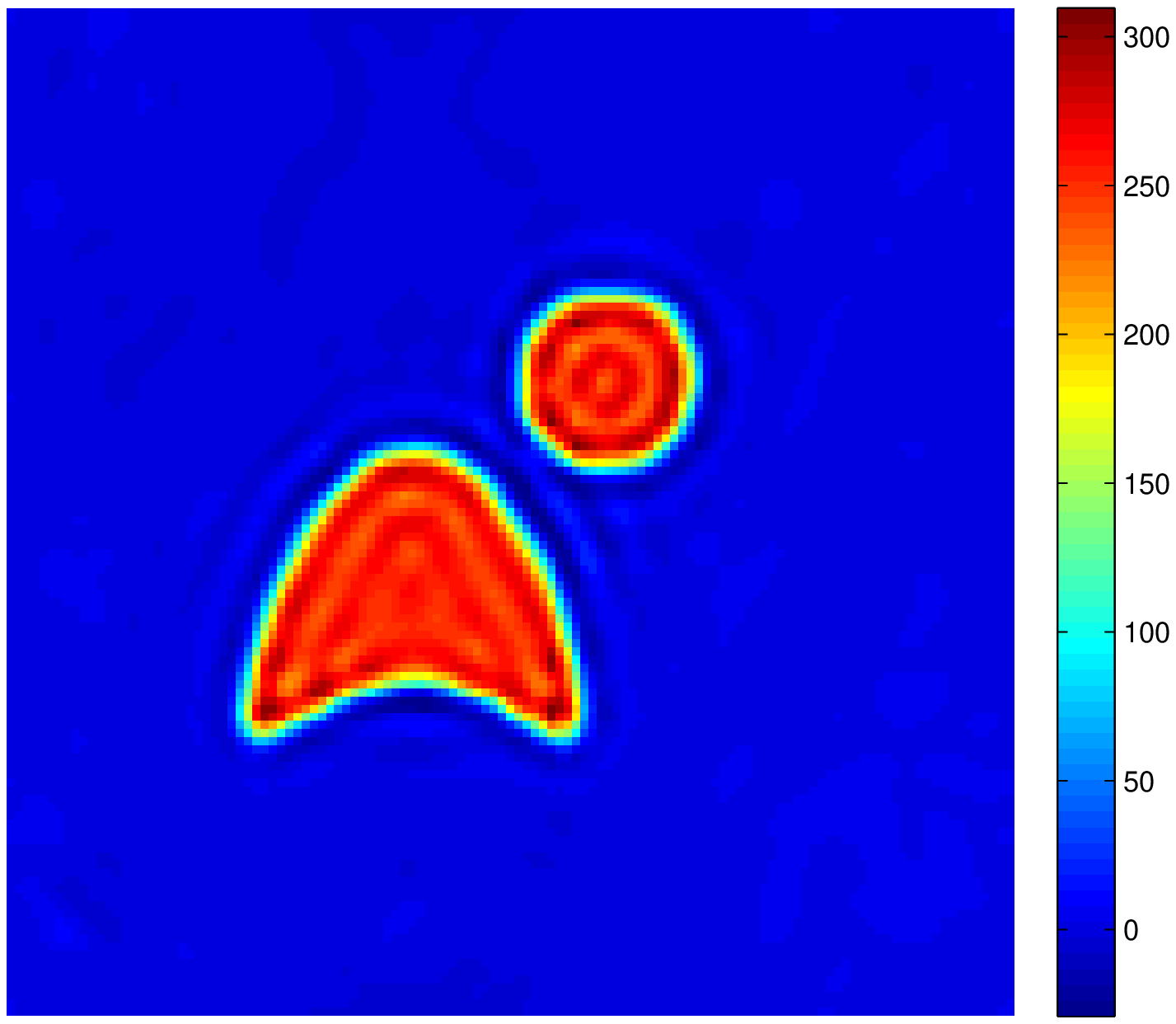}
&\includegraphics[width=0.19\linewidth]{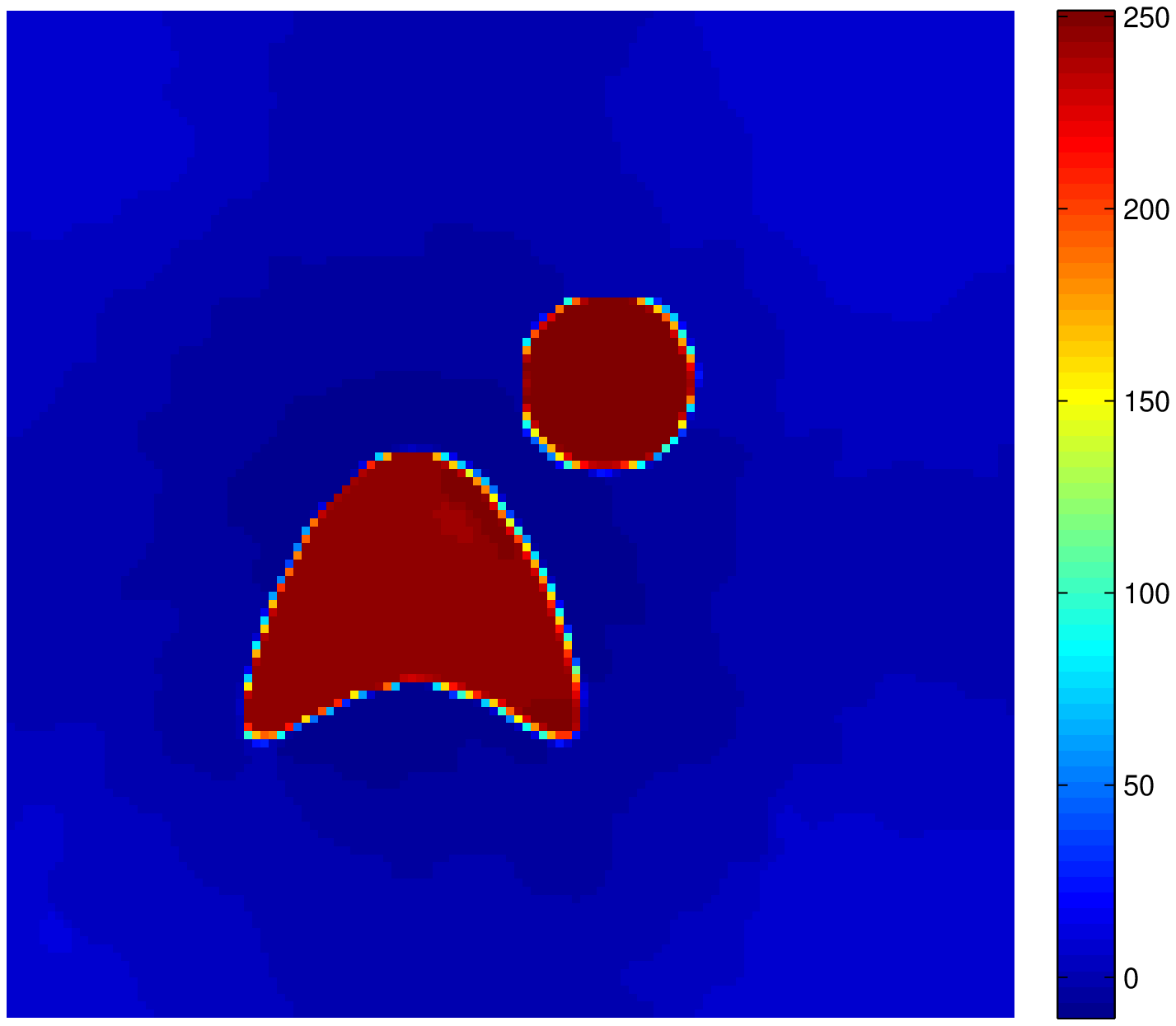}\\
%Corrupted image $u$ & $\delta^3$-NL restored & Atom-based restoration & TV restored\\
&PSNR=14.4dB &PSNR=18dB & PSNR=24.5dB & PSNR=29.1dB\\
\includegraphics[width=0.19\linewidth]{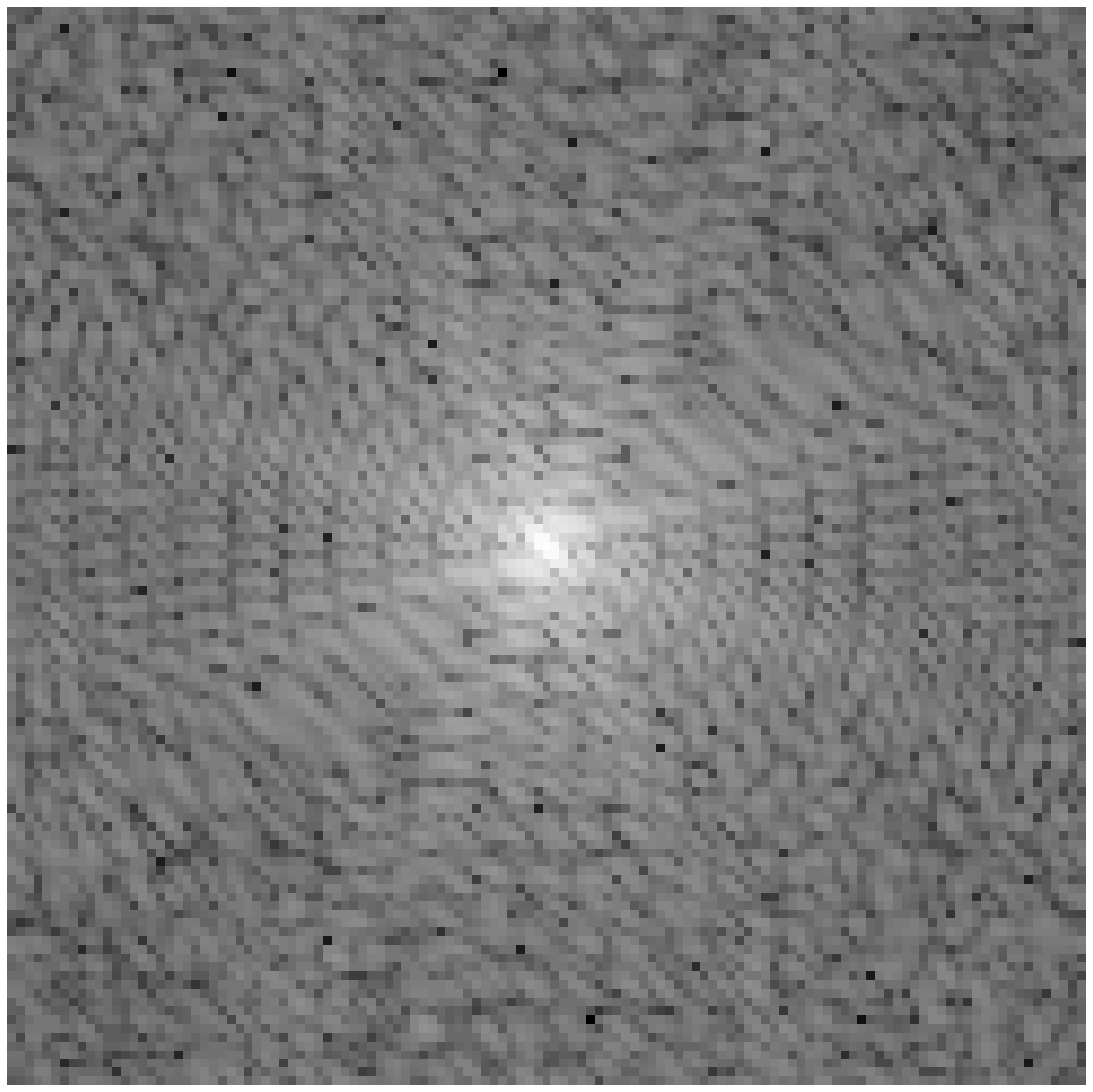}
&\includegraphics[width=0.19\linewidth]{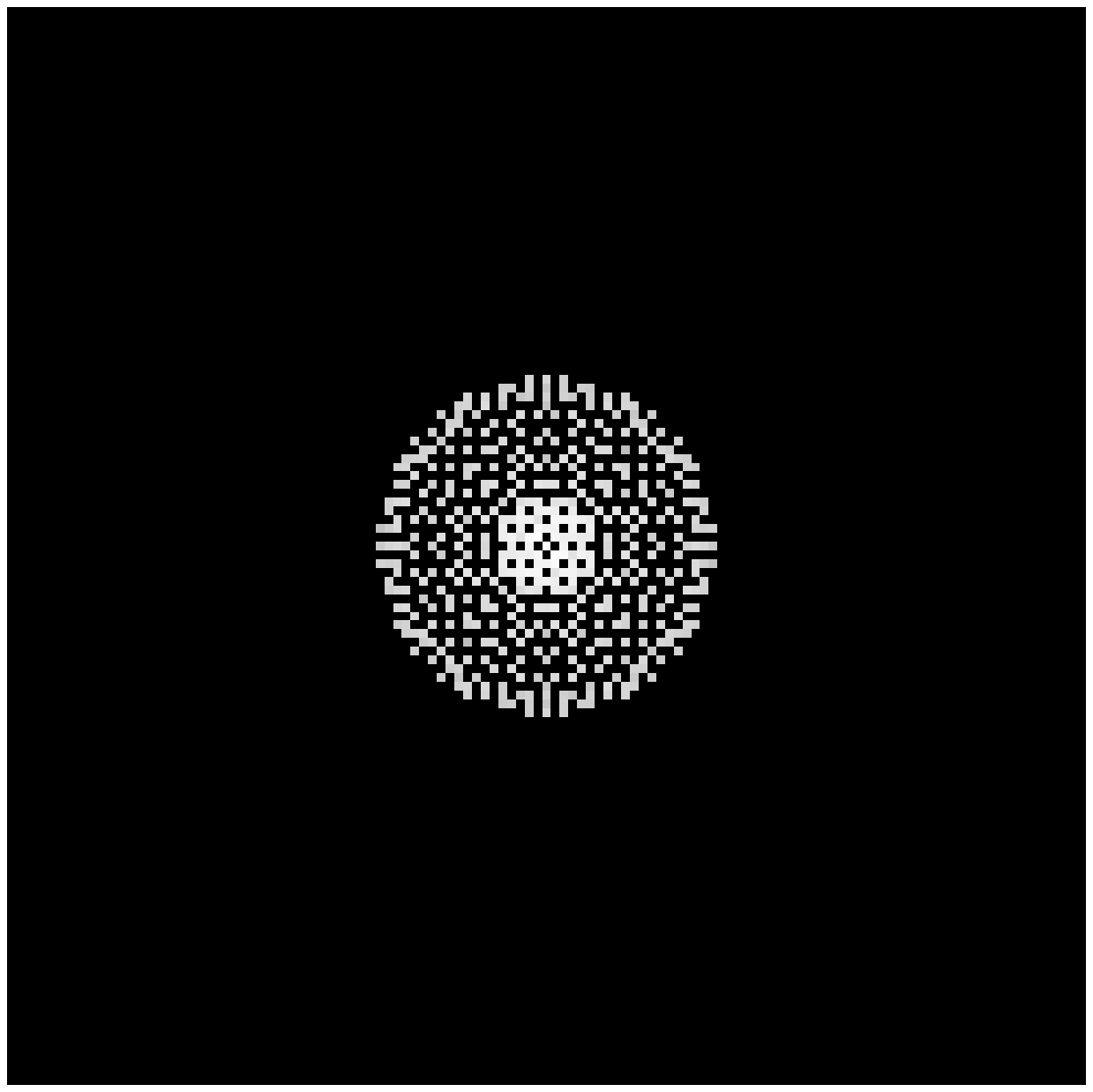}
&\includegraphics[width=0.19\linewidth]{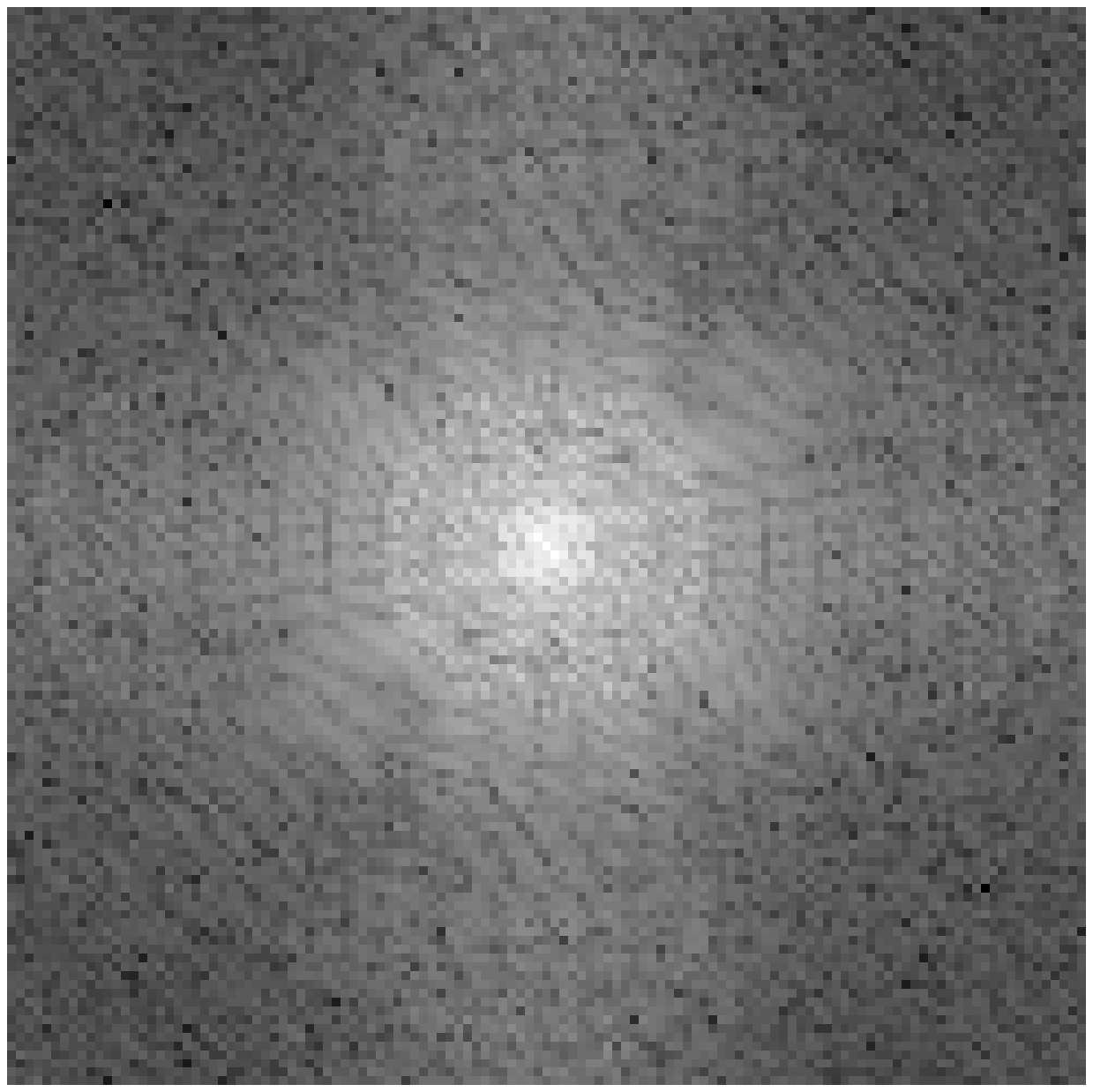}
&\includegraphics[width=0.19\linewidth]{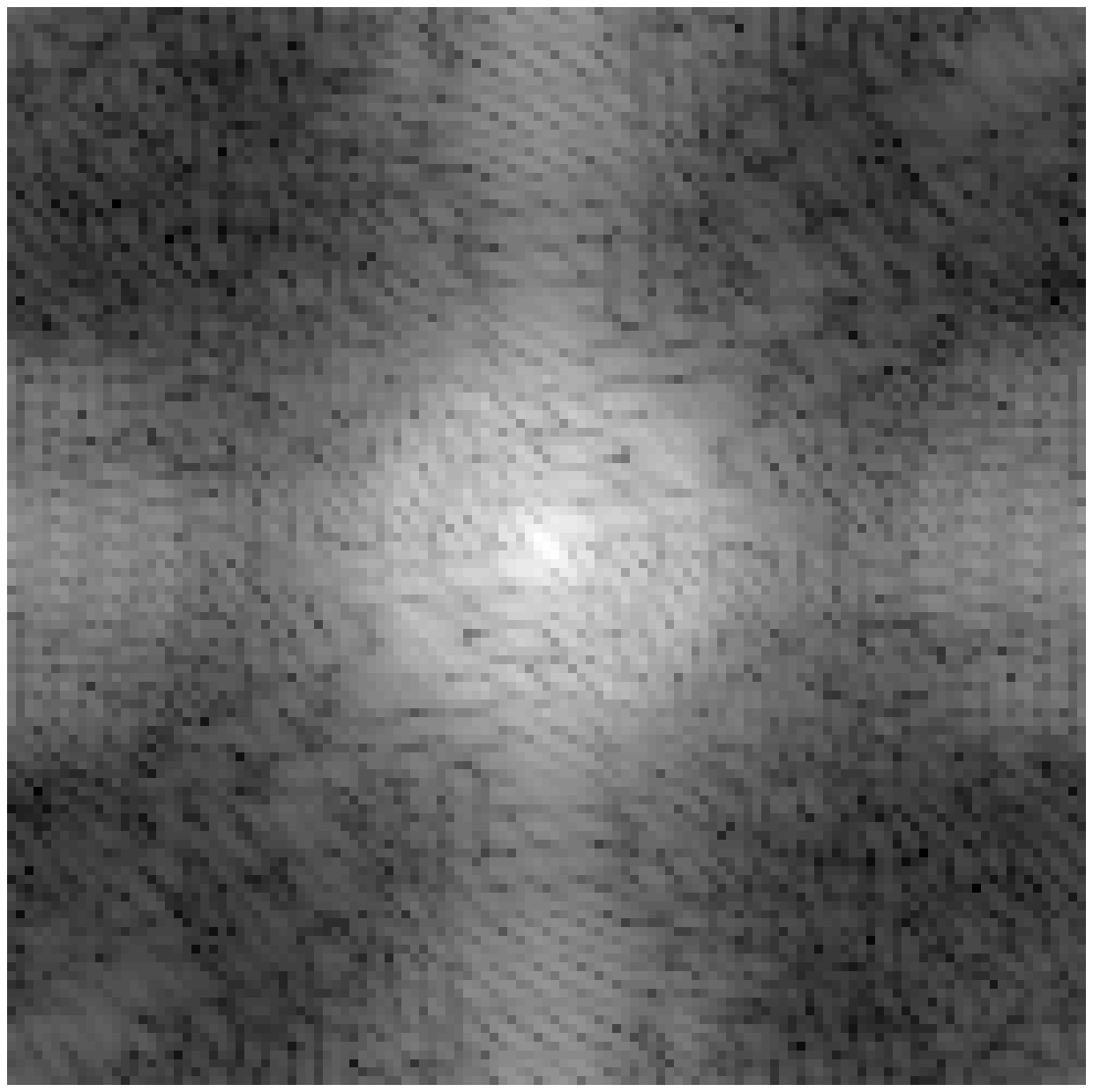}
&\includegraphics[width=0.19\linewidth]{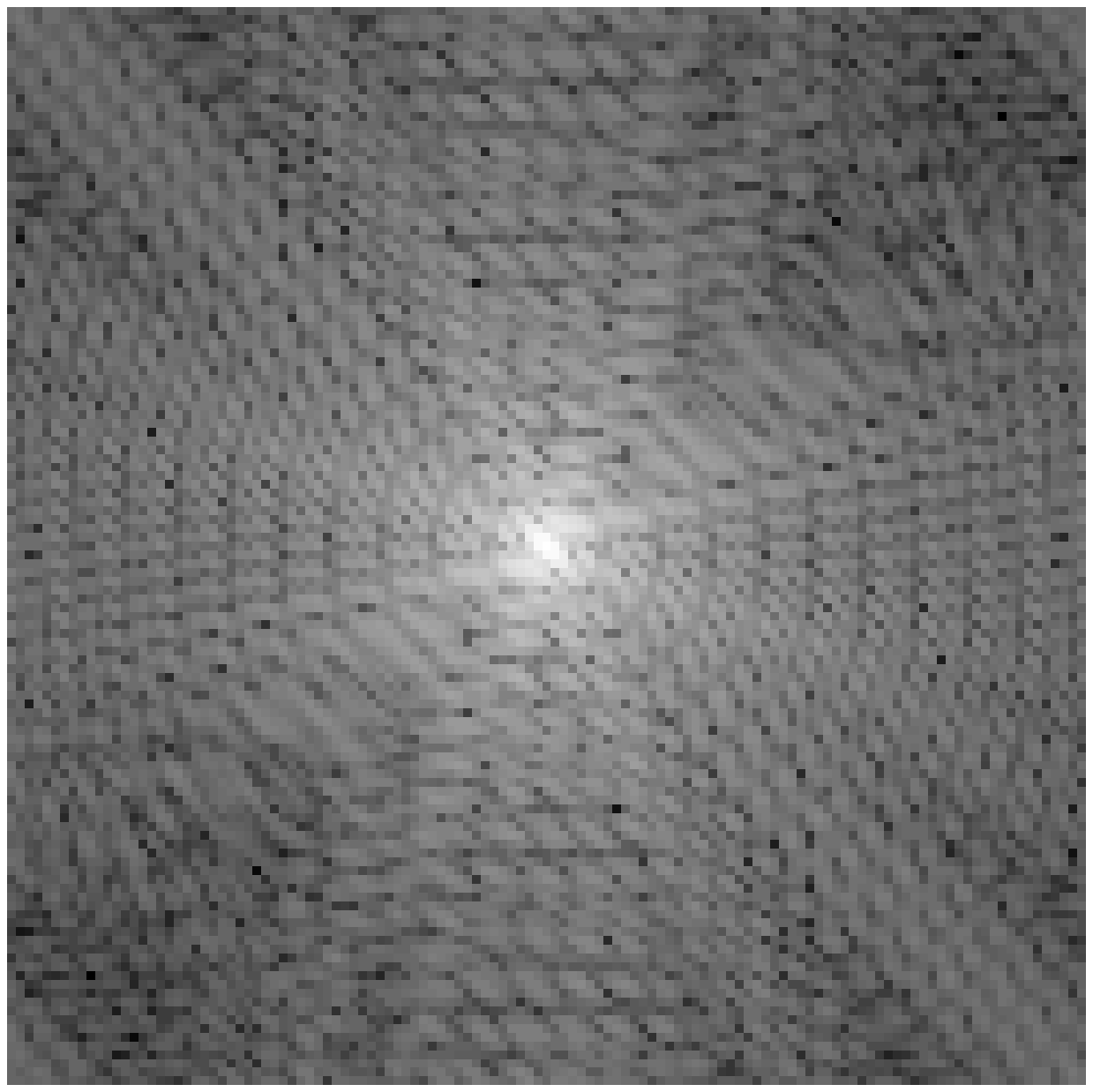}\\
%Spectrum of $u+v_{\delta^3}$ & Spectrum of $u+v_{\delta^1}$ & Spectrum of TV restored\\
\end{tabular}
\vspace{-0.4cm}
\caption%[]
{Corrupted and restored scatterers and their respective spectra.}\label{figScat}
% \end{minipage}
\end{figure}

\begin{figure}[H]
\centering
% \begin{minipage}[c]{\linewidth}
\begin{tabular}{ccccc}
Original $g_0$ & Corrupted $g$ & SSD - $\delta^3$ & NL-Atom - $\delta^1$ & TV restored\\
\includegraphics[width=0.19\linewidth]{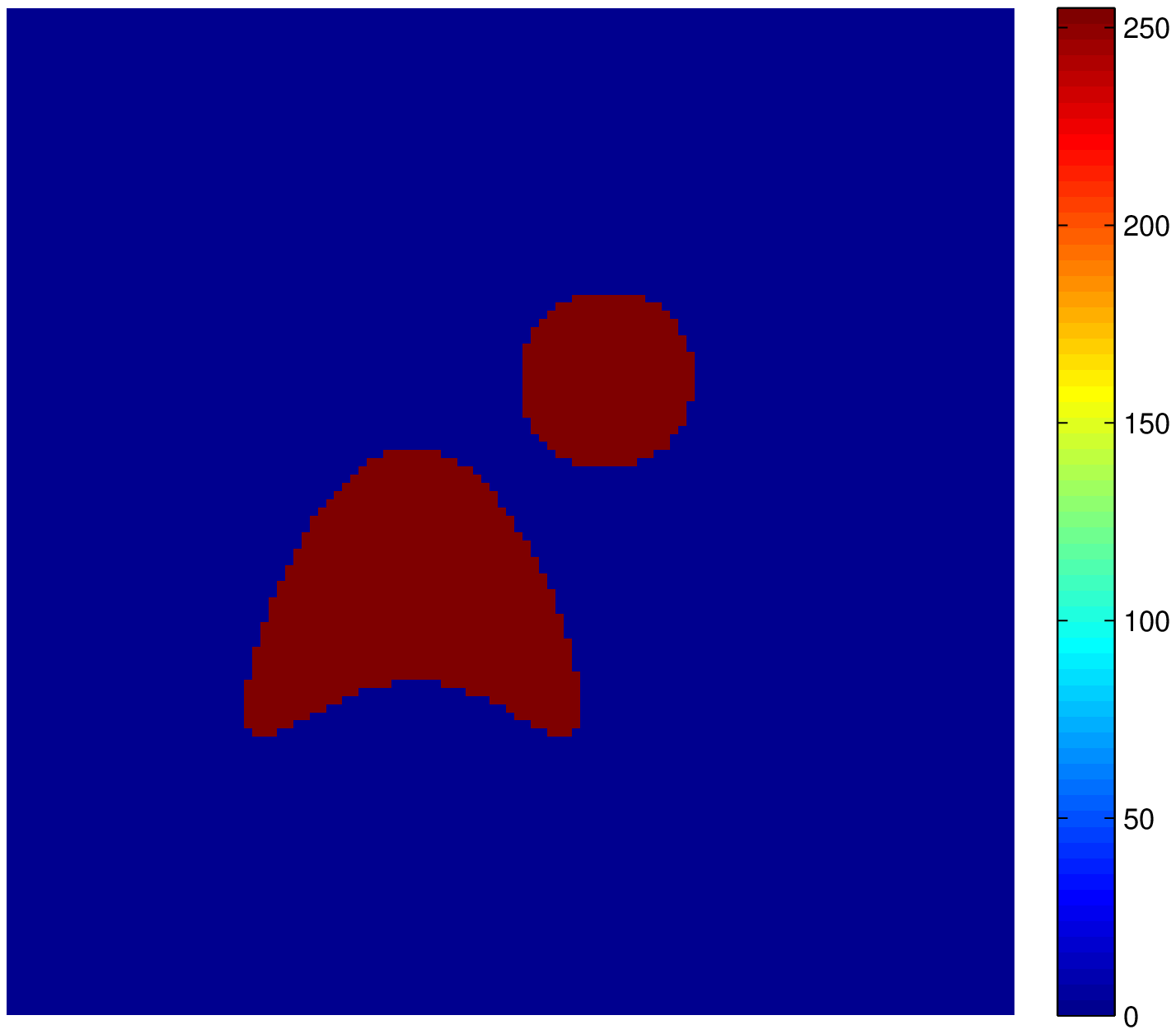}
&\includegraphics[width=0.19\linewidth]{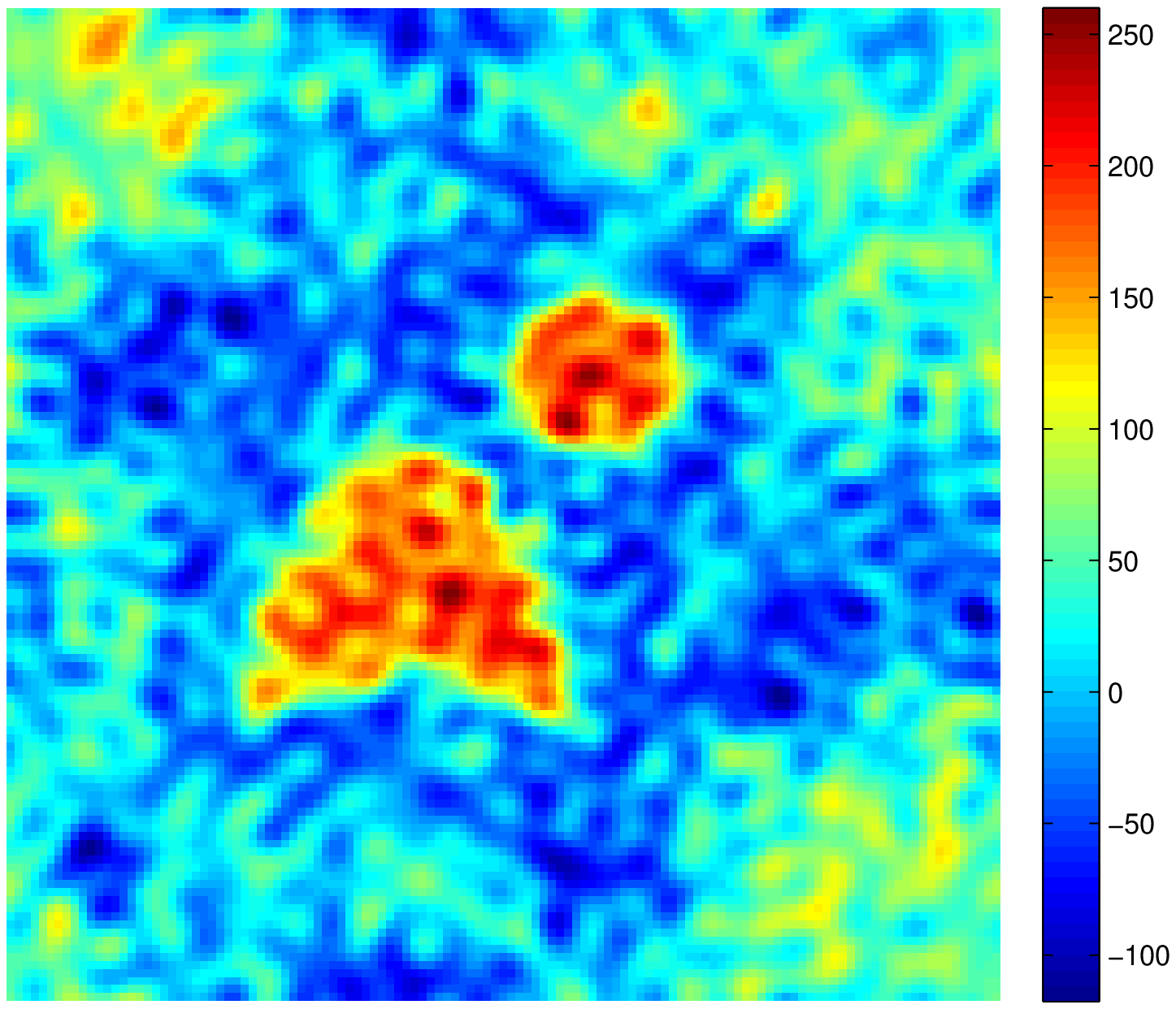}
&\includegraphics[width=0.19\linewidth]{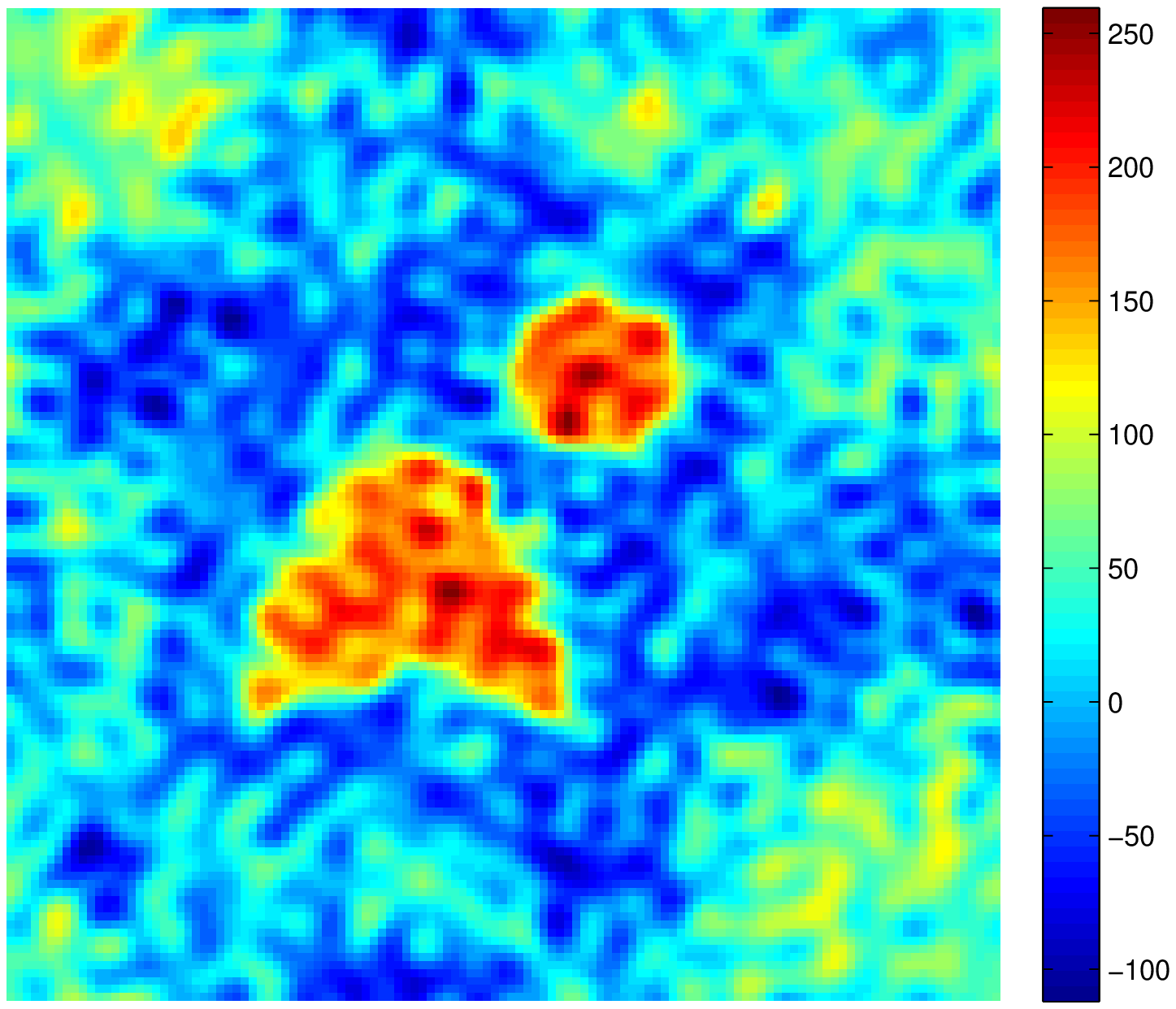}
&\includegraphics[width=0.19\linewidth]{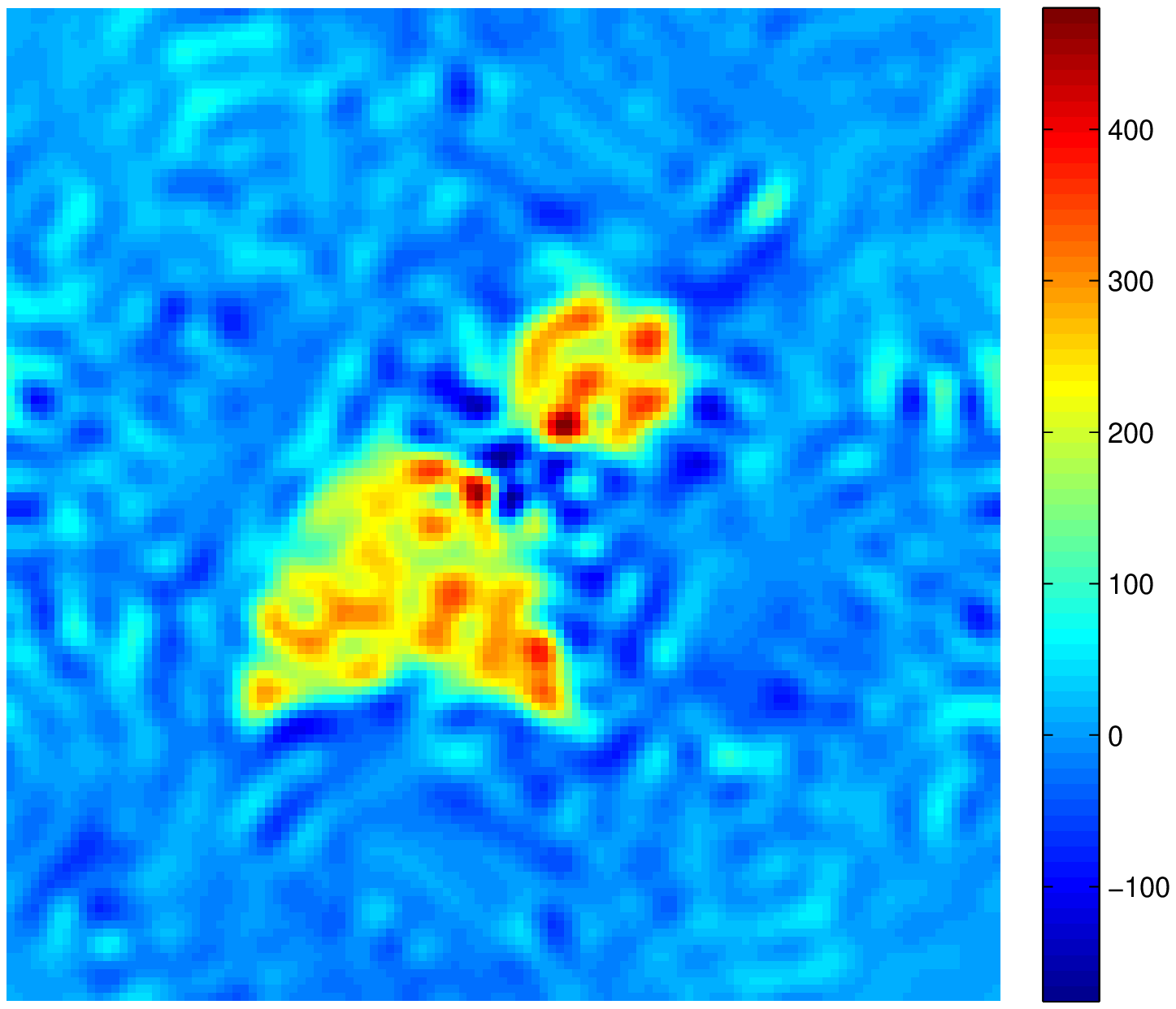}
&\includegraphics[width=0.19\linewidth]{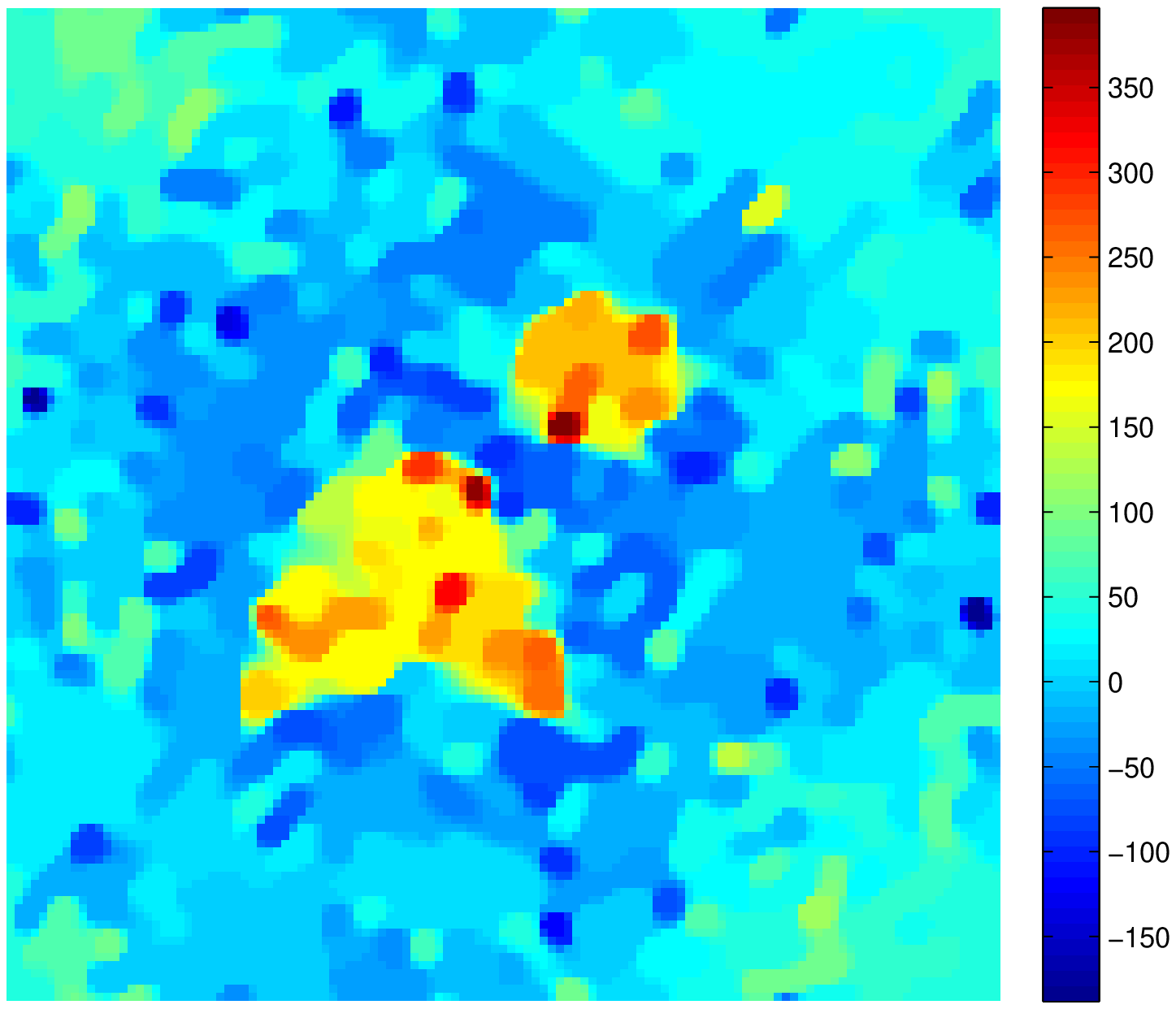}\\
%Corrupted image $u$ & $\delta^3$-NL restored & Atom-based restoration & TV restored\\
 &PSNR=13.4dB &PSNR=13.6dB & PSNR=17.1dB & PSNR=14.4dB\\
\includegraphics[width=0.19\linewidth]{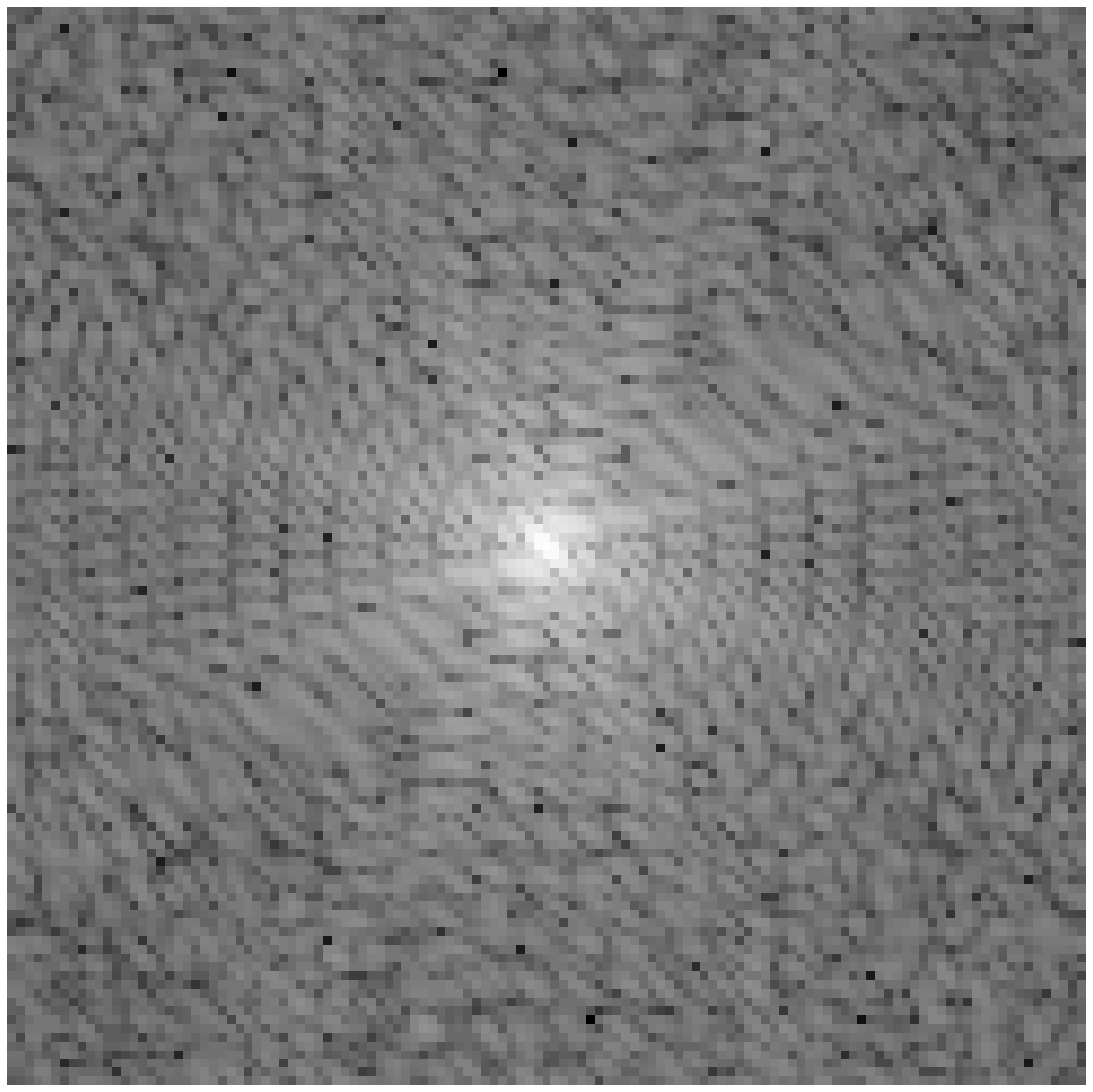}
&\includegraphics[width=0.19\linewidth]{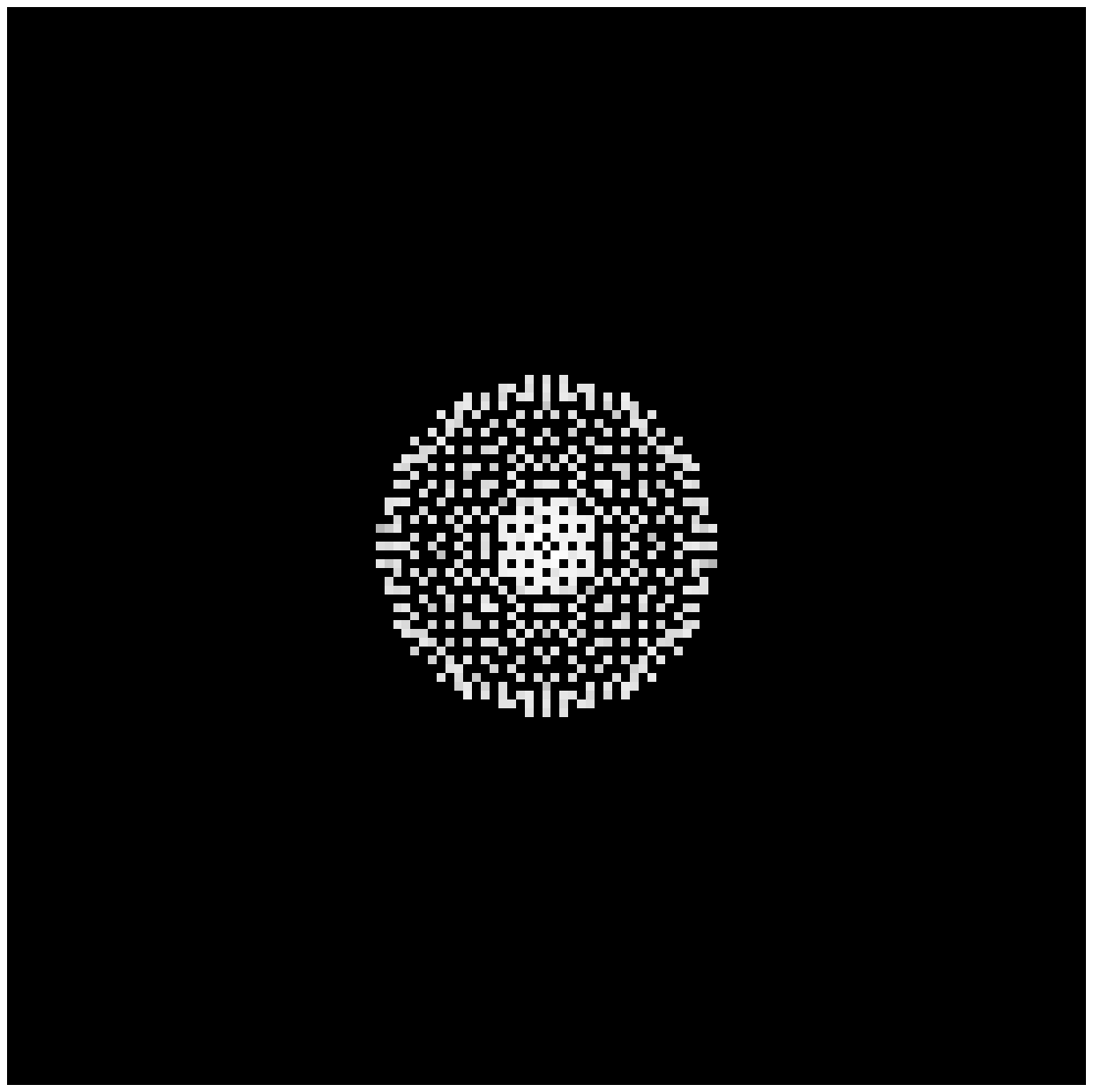}
&\includegraphics[width=0.19\linewidth]{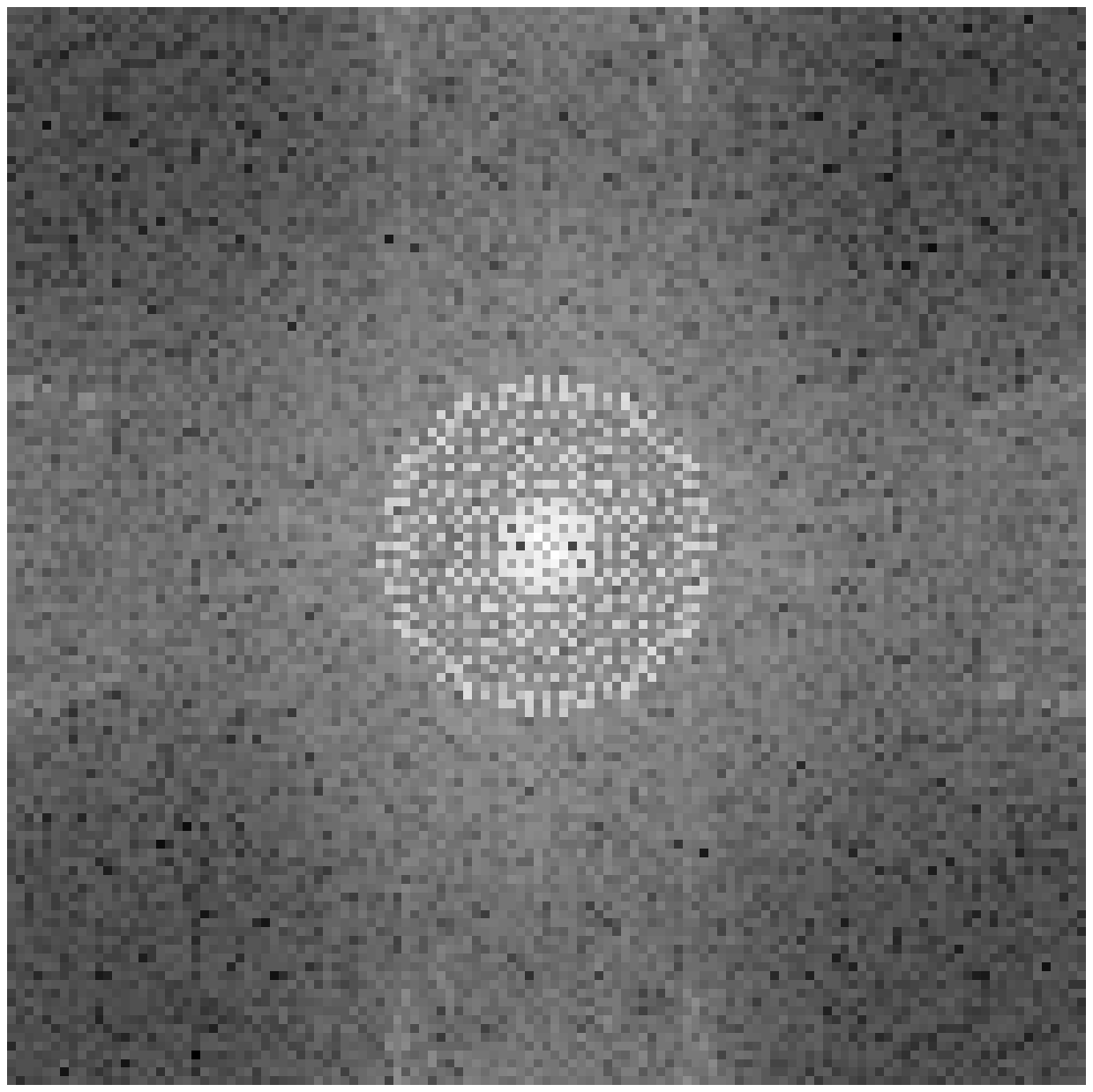}
&\includegraphics[width=0.19\linewidth]{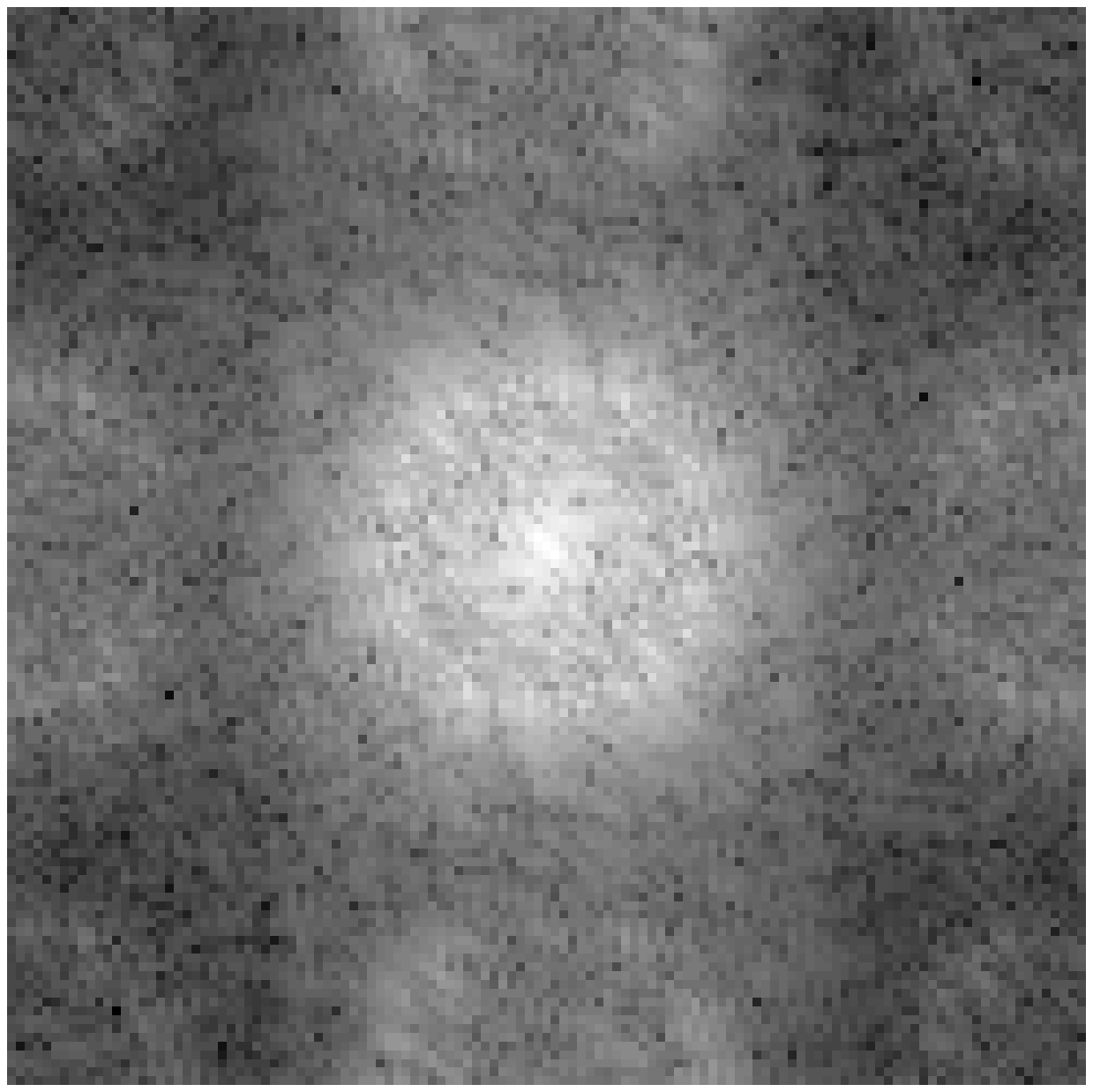}
&\includegraphics[width=0.19\linewidth]{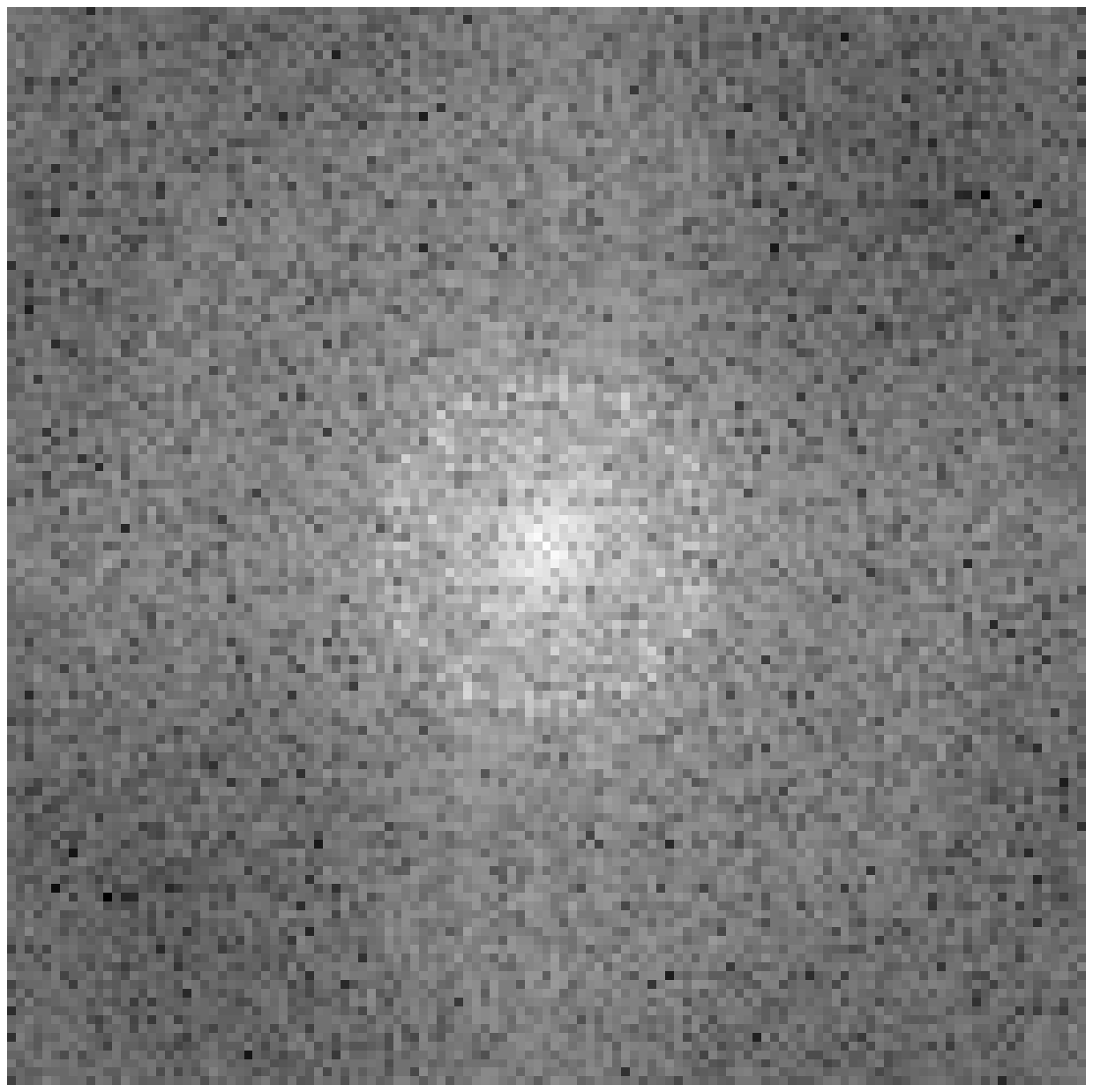}\\
%Spectrum of $u+v_{\delta^3}$ & Spectrum of $u+v_{\delta^1}$ & Spectrum of TV restored\\
\end{tabular}
\vspace{-0.4cm}
\caption%[]
{Corrupted and restored scatterers and their respective spectra.}\label{figScatNoisy}
% \end{minipage}
\end{figure}

\begin{figure}[H]
\centering
% \begin{minipage}[c]{\linewidth}
\begin{tabular}{ccccc}
Original $g_0$ & Corrupted $g$ & SSD - $\delta^3$ & NL-Atom - $\delta^1$ & TV restored\\
\includegraphics[width=0.19\linewidth]{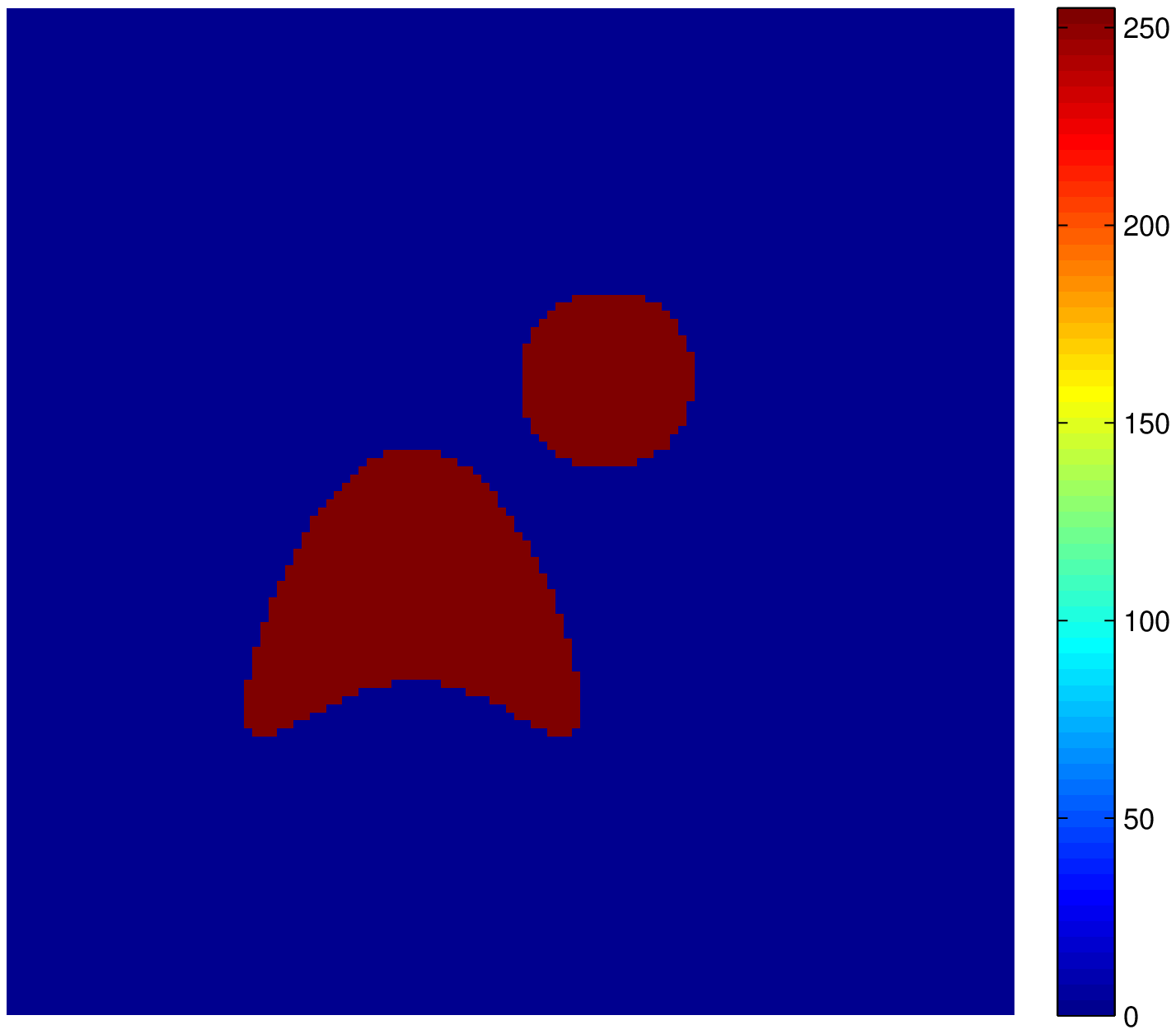}
&\includegraphics[width=0.19\linewidth]{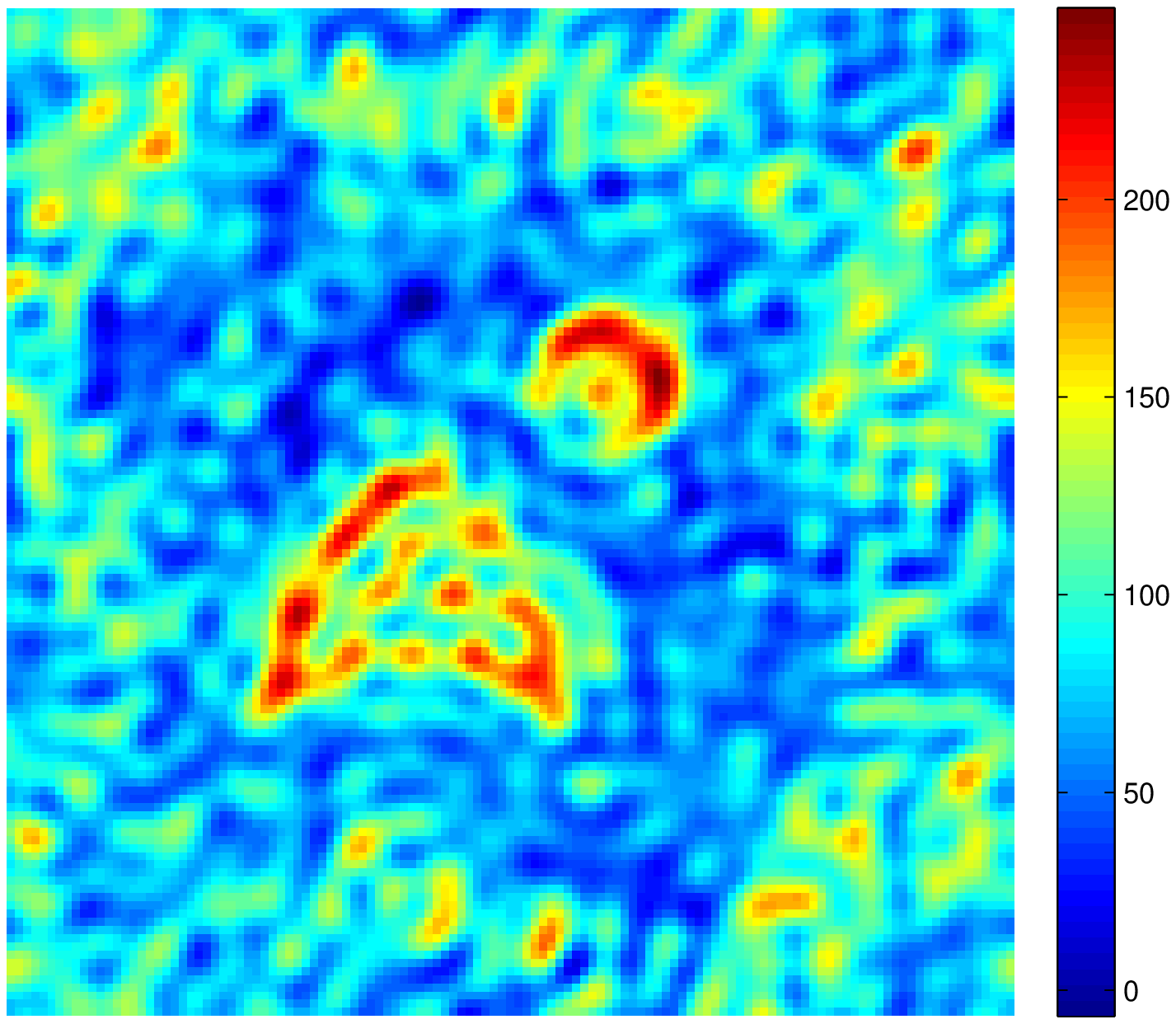}
&\includegraphics[width=0.19\linewidth]{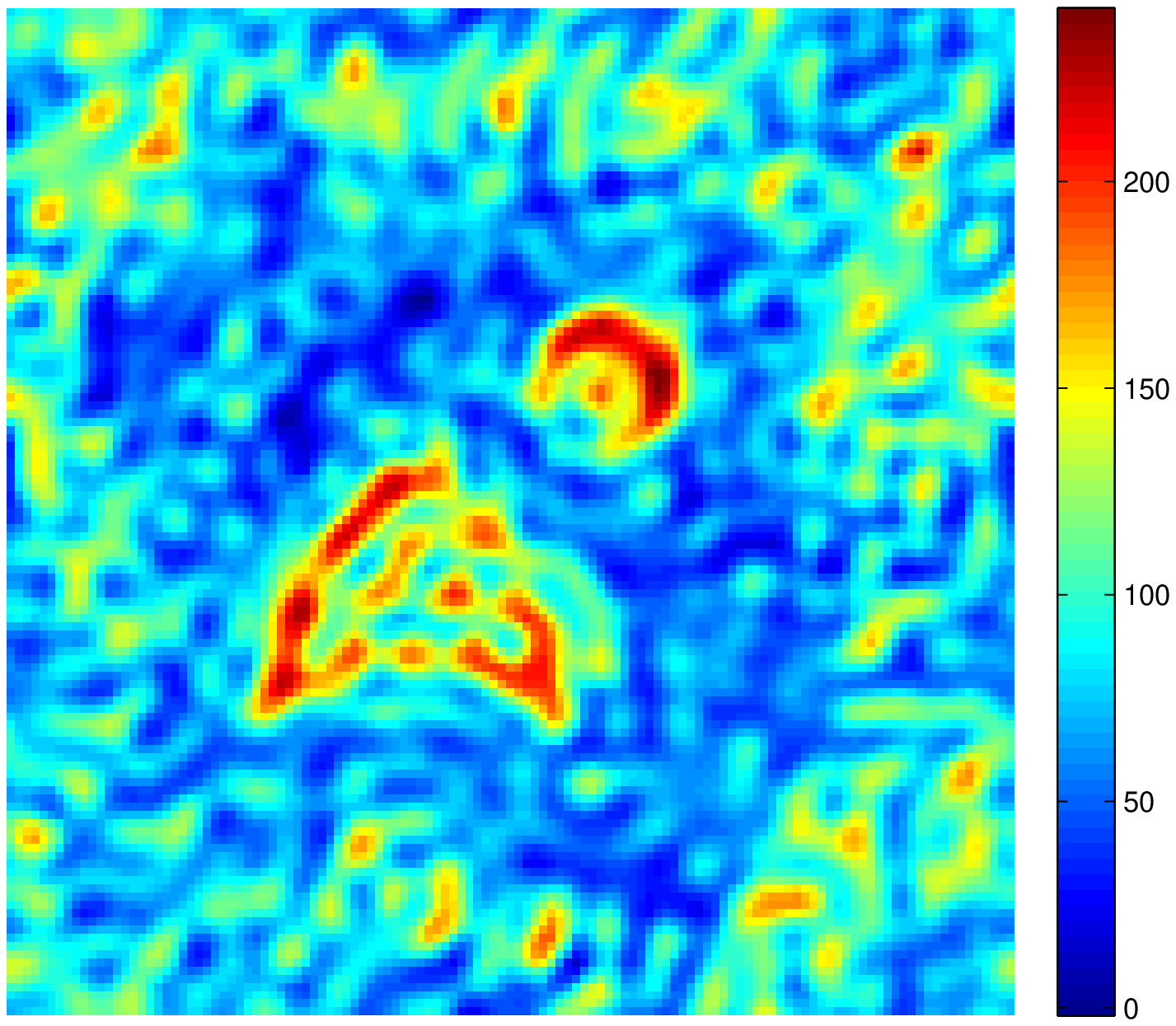}
&\includegraphics[width=0.19\linewidth]{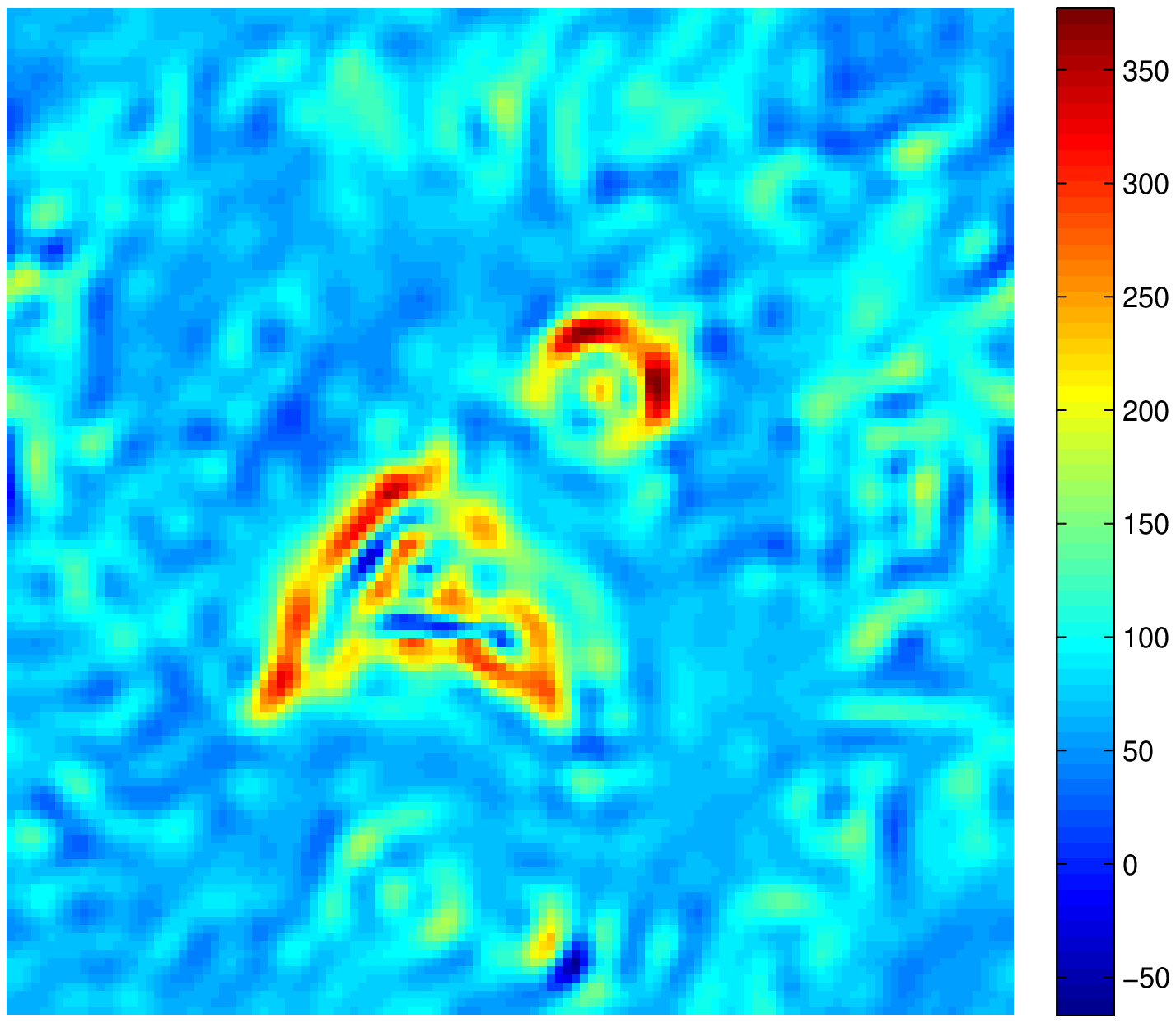}
&\includegraphics[width=0.19\linewidth]{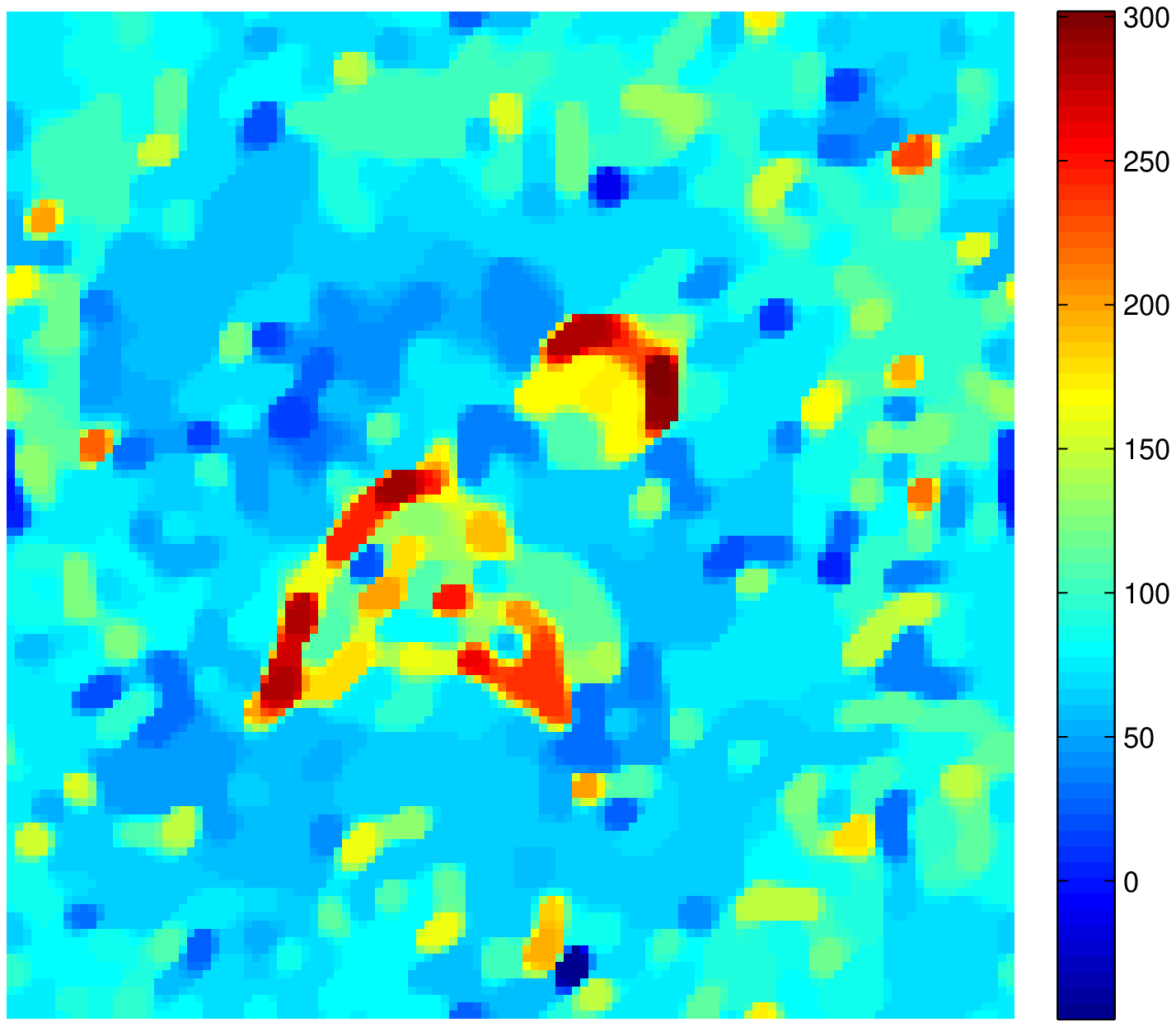}\\
%Corrupted image $u$ & $\delta^3$-NL restored & Atom-based restoration & TV restored\\
& PSNR=9.24dB &PSNR=9.25dB & PSNR=9.81dB & PSNR=9.42dB\\
\includegraphics[width=0.19\linewidth]{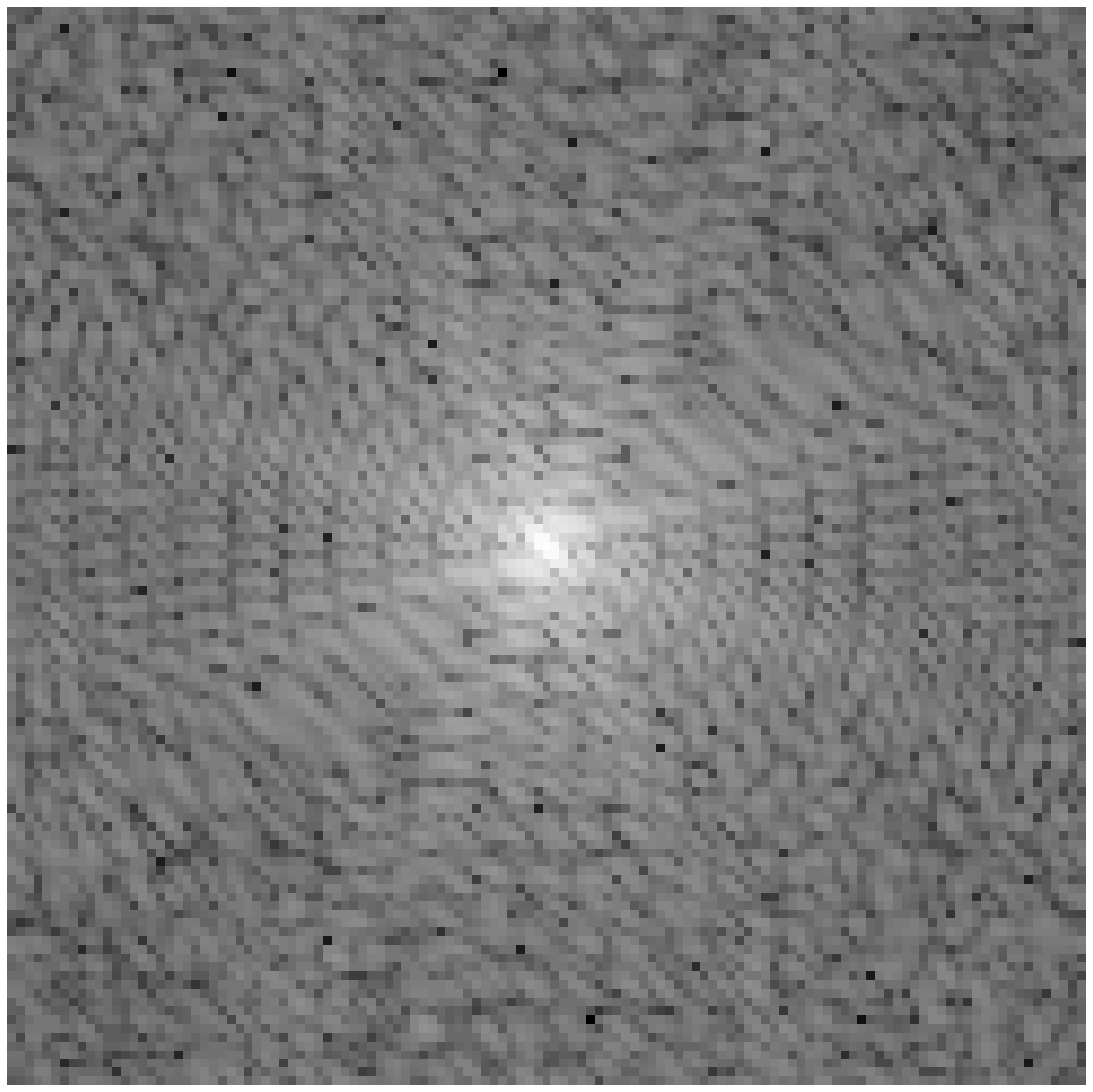}
&\includegraphics[width=0.19\linewidth]{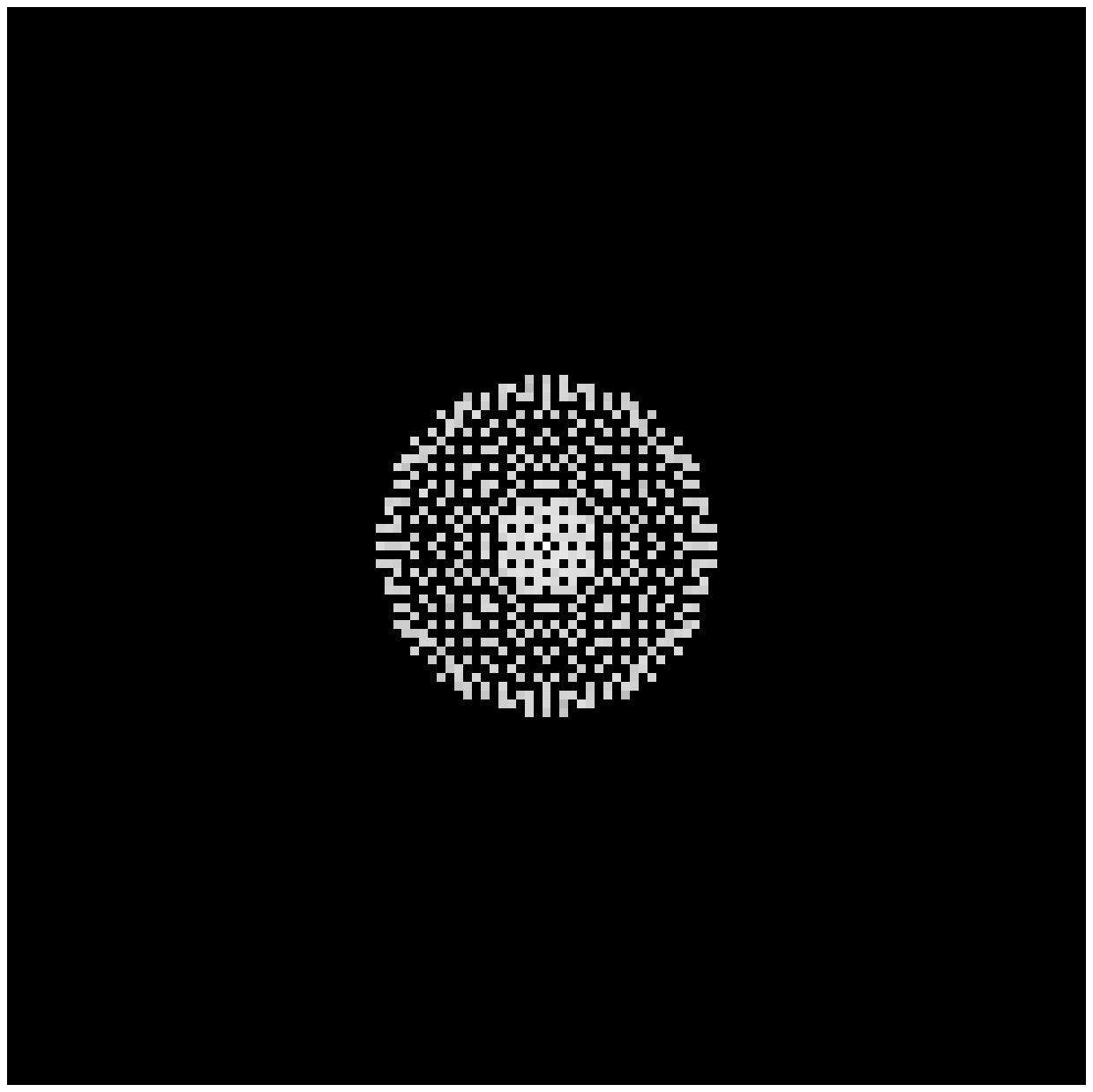}
&\includegraphics[width=0.19\linewidth]{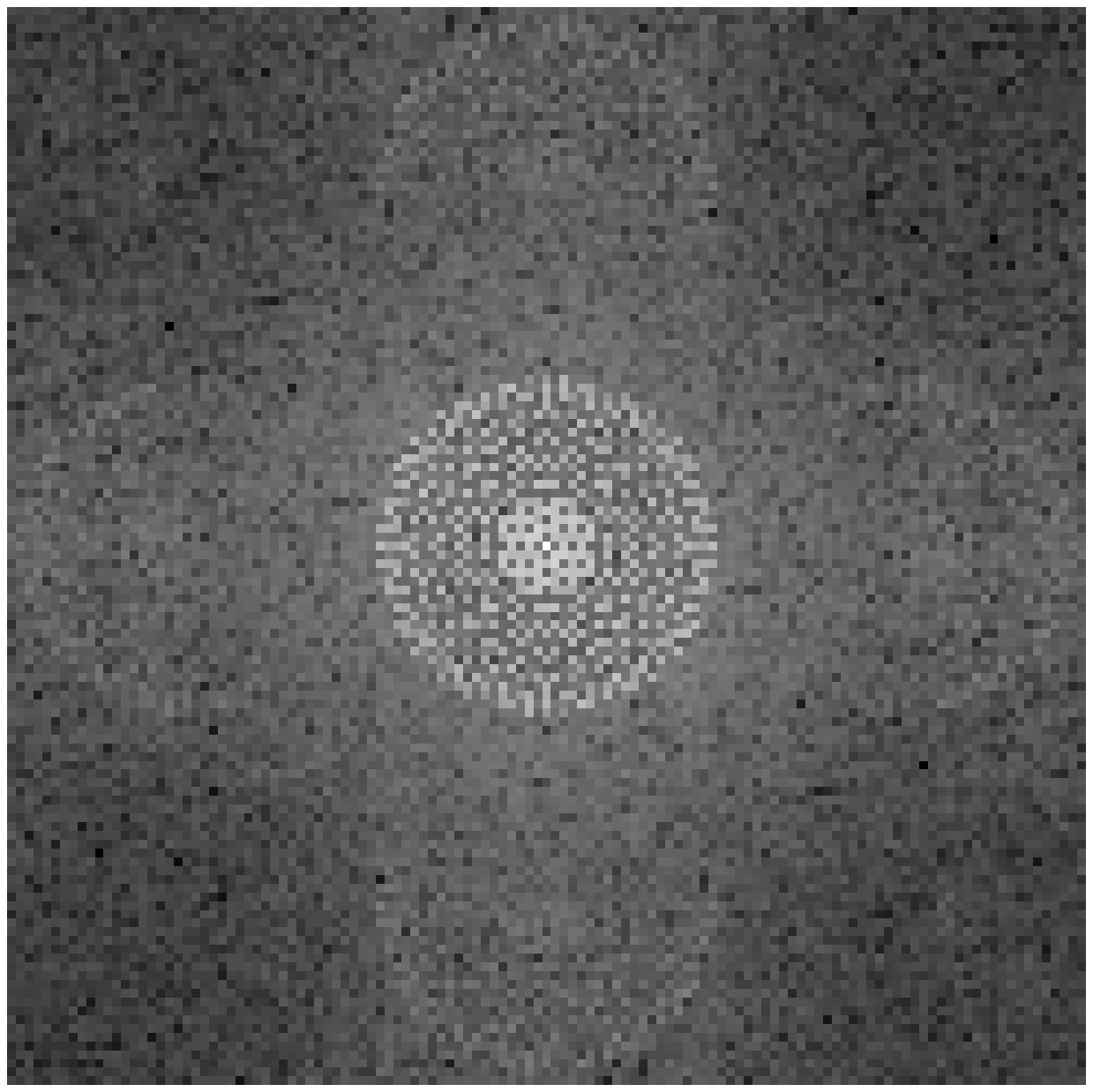}
&\includegraphics[width=0.19\linewidth]{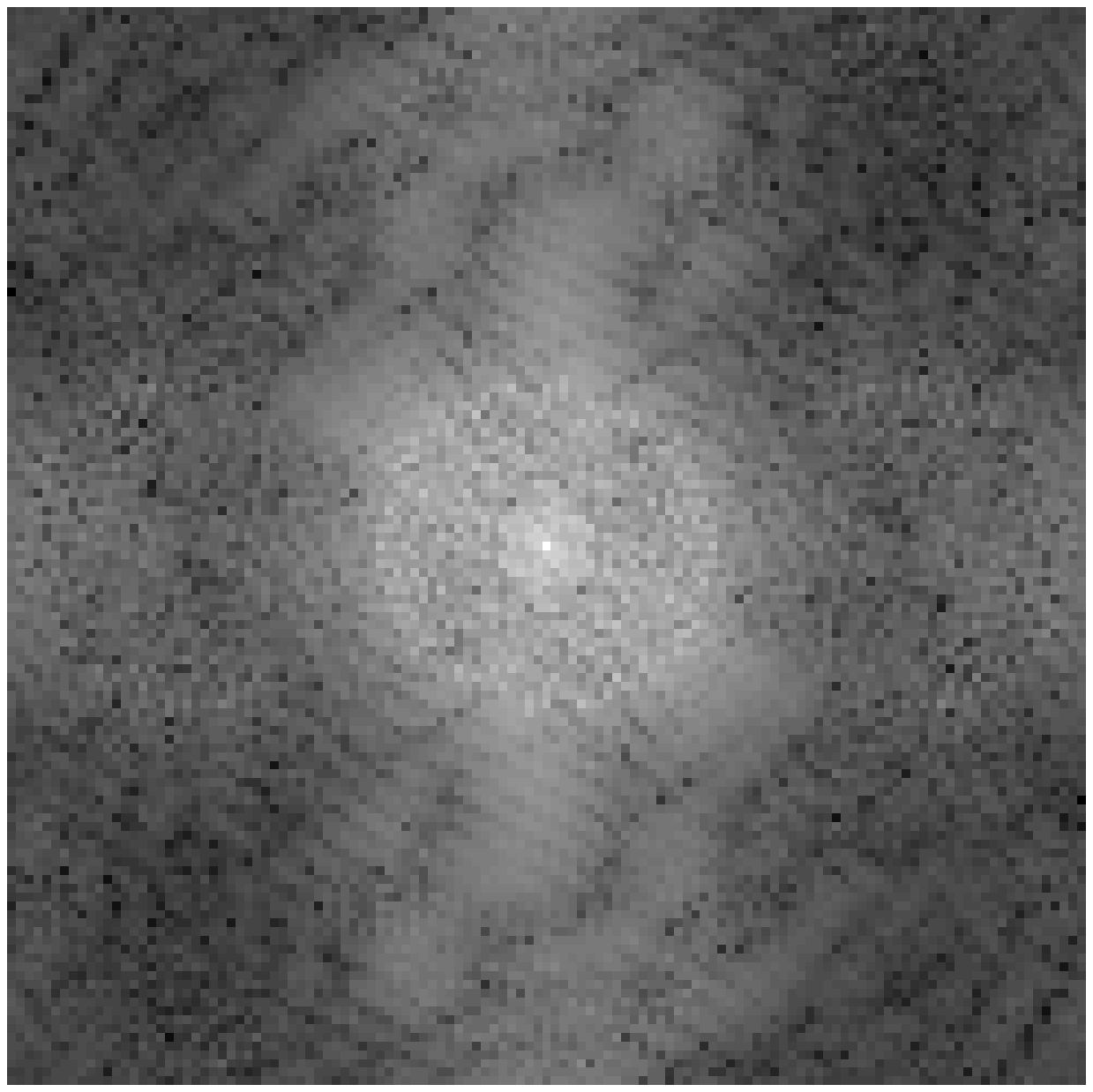}
&\includegraphics[width=0.19\linewidth]{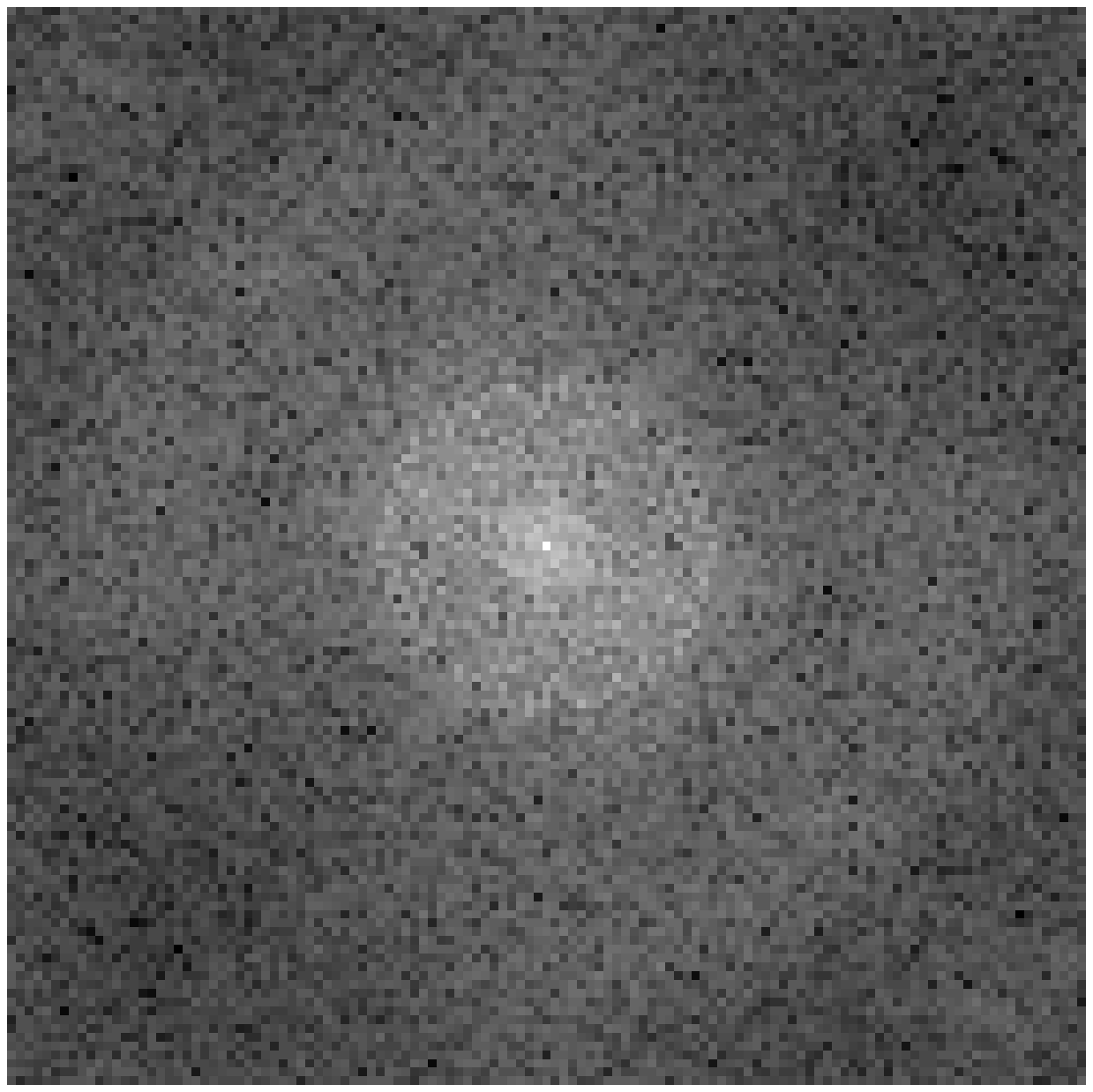}\\
%Spectrum of $u+v_{\delta^3}$ & Spectrum of $u+v_{\delta^1}$ & Spectrum of TV restored\\
\end{tabular}
%\vspace{-0.4cm}
\caption%[]
{Corrupted and restored scatterers and their respective spectra.}\label{figScatBorn}
% \end{minipage}
\end{figure}

\begin{figure}
% \hfill
\begin{minipage}[c]{\linewidth}
\centering
\begin{tabular}{ccc}
\includegraphics[width=0.18\linewidth]{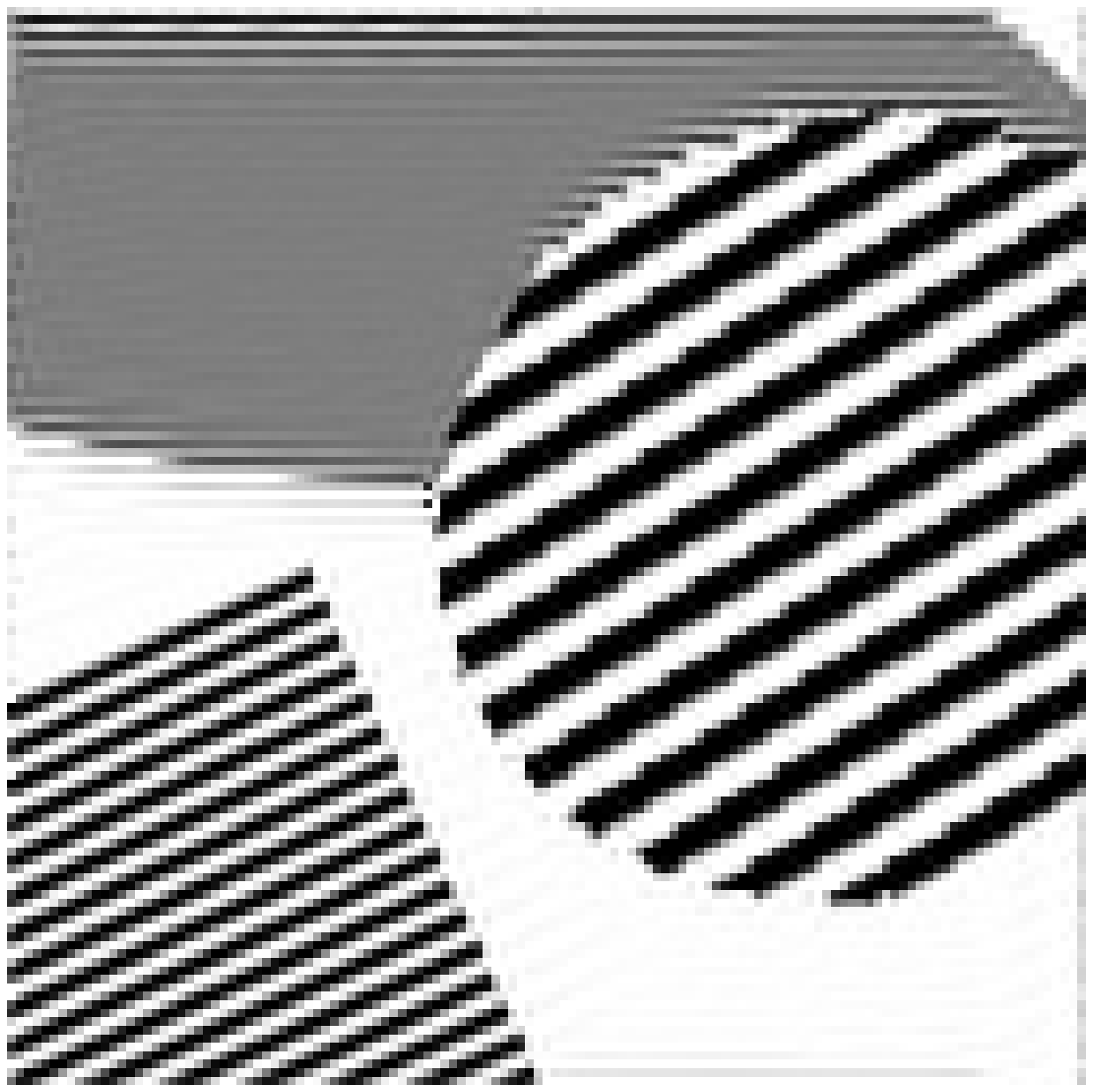}
&\includegraphics[width=0.18\linewidth]{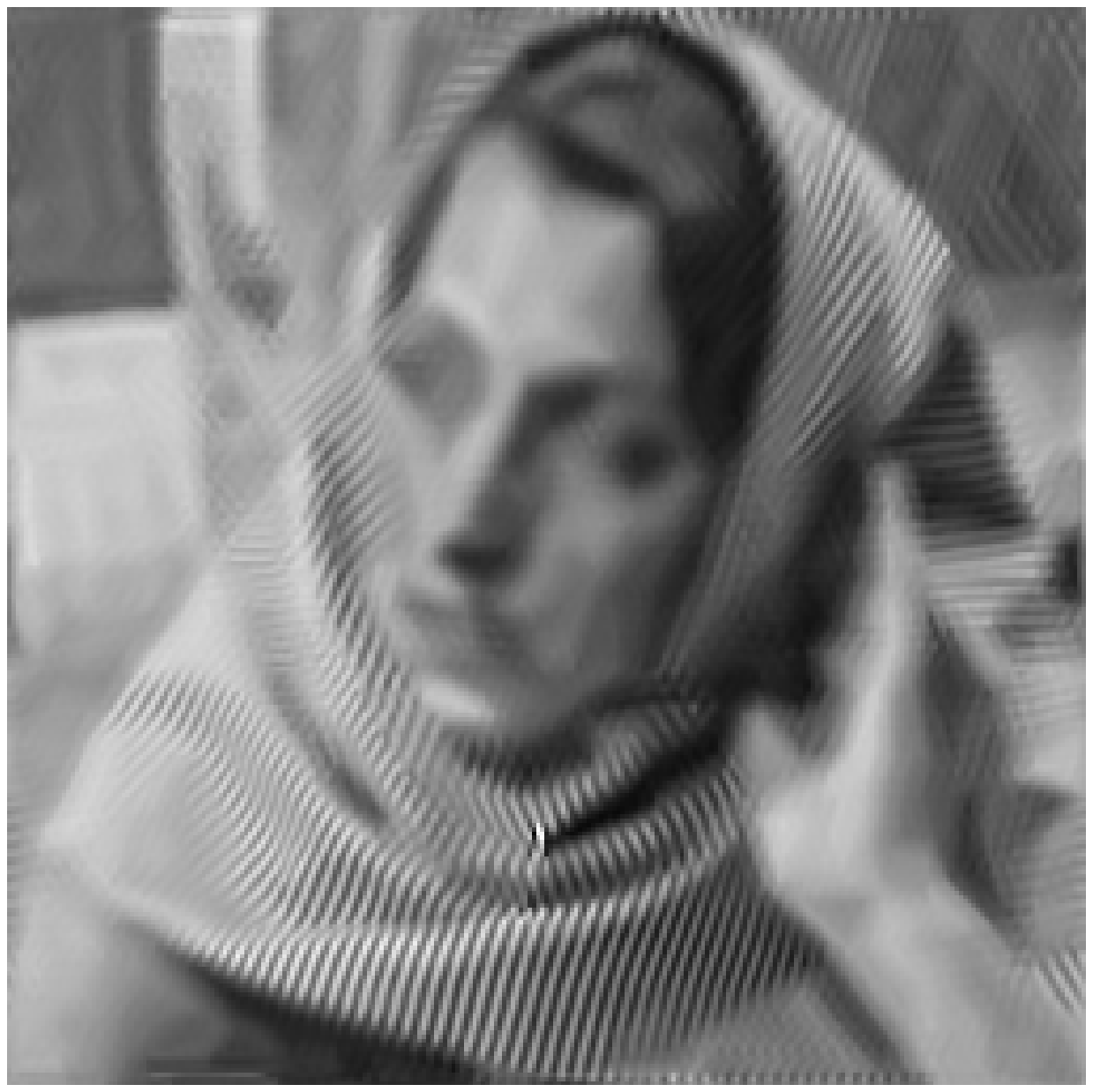}
&\includegraphics[width=0.205\linewidth]{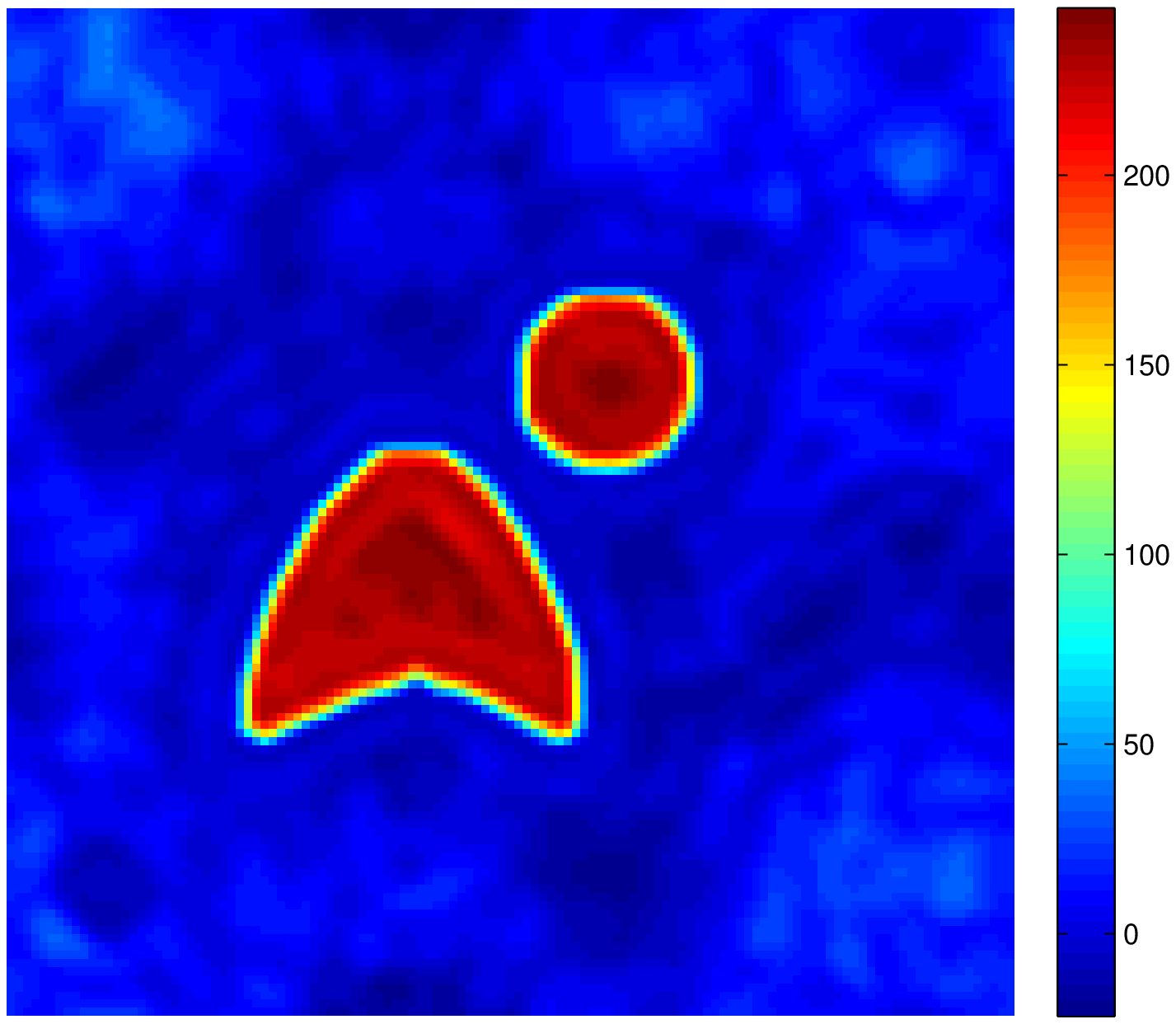}\\
%PSNR=21.6dB &
PSNR=12.1dB & PSNR=24.9dB&PSNR=22.3 dB \\
%$(\alpha=1)$& & &\\
%\vspace{0.2cm}
%\includegraphics[width=0.18\linewidth]{NL_7sc_resc_nc_1_recompute}
\includegraphics[width=0.205\linewidth]{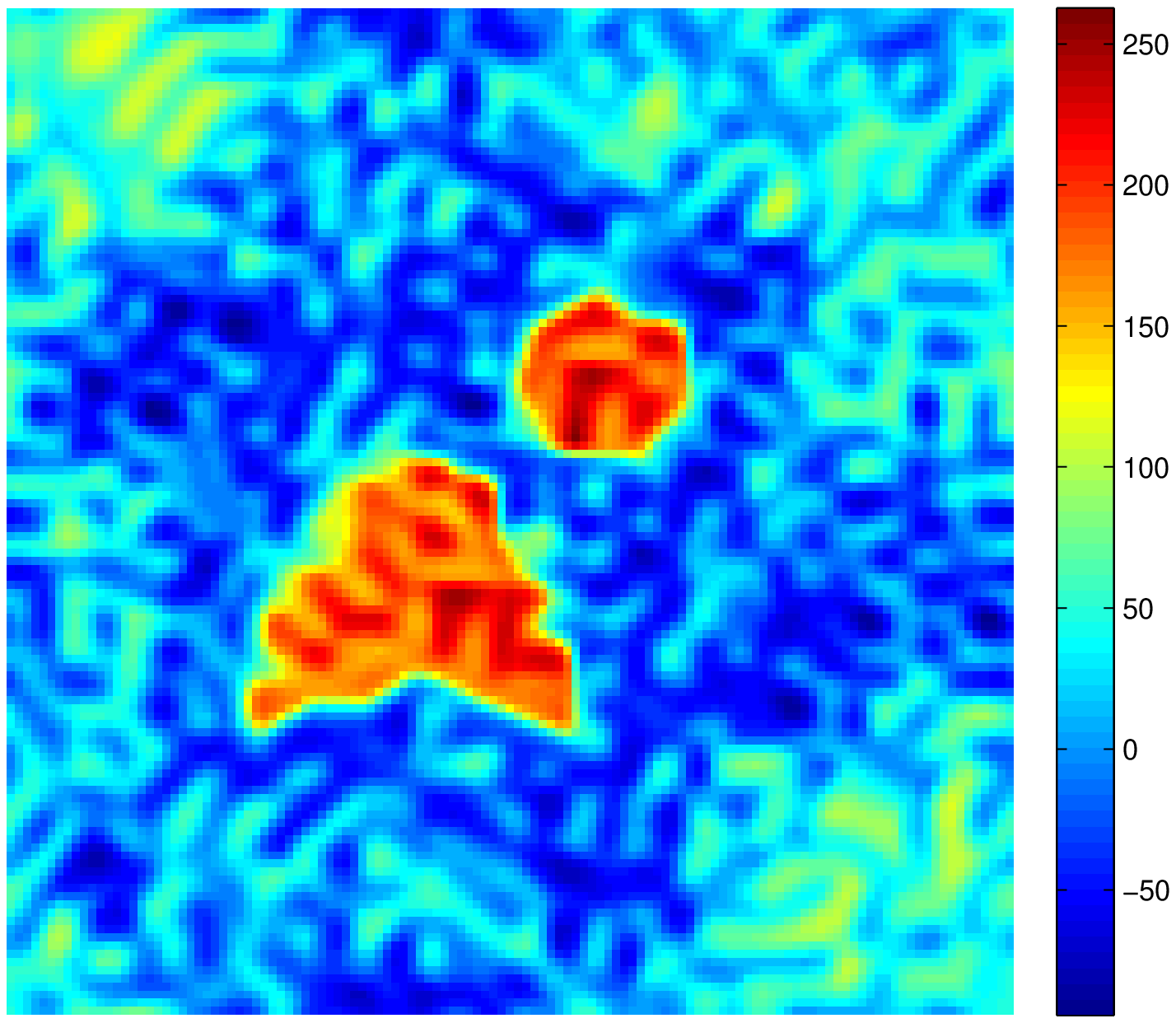}
&\includegraphics[width=0.205\linewidth]{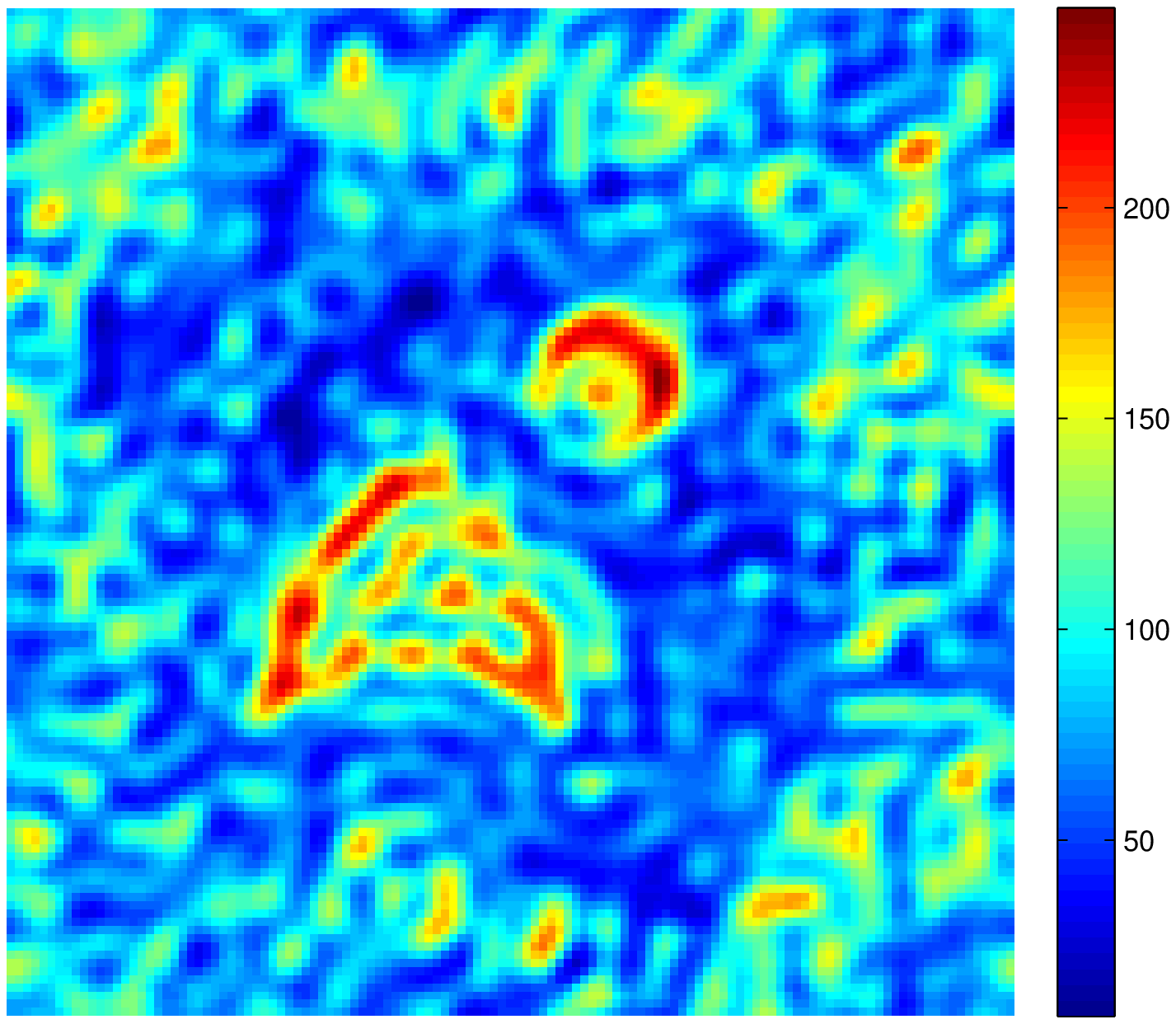}
&\includegraphics[width=0.205\linewidth]{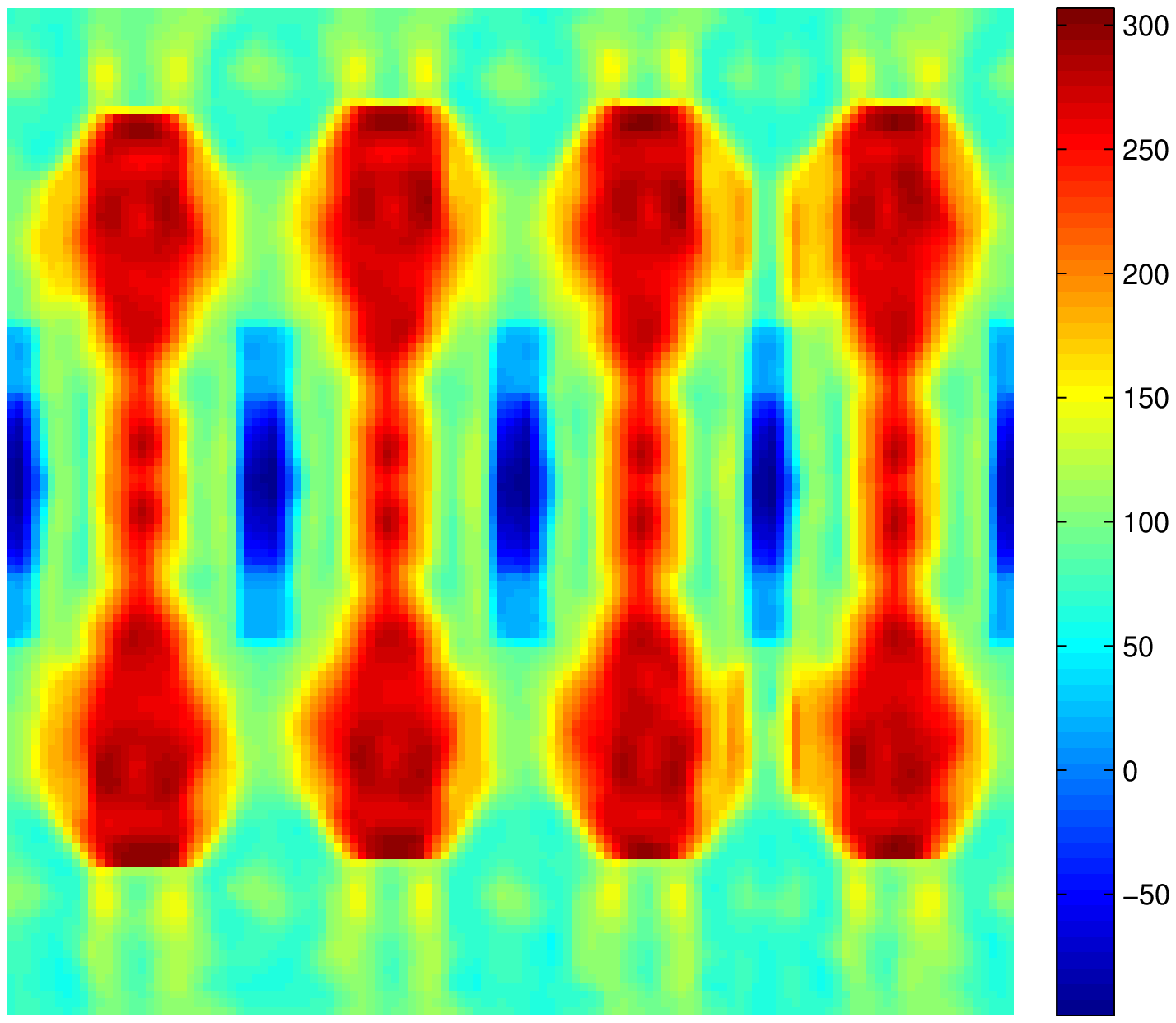}\\
%PSNR=21.9dB &
PSNR=14.6dB &PSNR=9.27dB &PSNR=8.9dB\\
%$(\alpha=2)$& & &
% \includegraphics[width=0.18\linewidth]{rayures_poids_NL_means_recomputed_NL_13}
% &\includegraphics[width=0.18\linewidth]{barbara_NLmeans_recomputedNL_13}
% &\includegraphics[width=0.18\linewidth]{scatt2_NL_means_recomputed_NL_13}
% &\includegraphics[width=0.18\linewidth]{scatt3_NLmeansNL_13}
% &\includegraphics[width=0.18\linewidth]{scatt_NL_means_recomputed_NL_13}
\end{tabular}
\end{minipage}
\vspace{-0.4cm}
\caption
{Restored images after many SSD weight recomputations.}
\vspace{-0.4cm}
\label{weight_recompute}
\end{figure}

\noindent\textbf{Numerical results:}
In the following tests we consider three different energies:% to be minimized to restore the 
\begin{itemize}
 \item The constrained total variation minimization problem
\begin{align*}
\min_{v\in\M^\perp} TV(g+v),
\end{align*}
%\newbf{This way, there is no Lagrange multiplier to tune.}
which can be solved for instance using \newbf{\cite{Beck,Nesterov07} or even \cite{ChamPD}}.

 \item The constrained minimization problem (\ref{E_NL}) where we consider the SSD measure of similarity $\delta^3$. In the simulations of Fig. \ref{figToy} and Fig. \ref{figBarb}, we take the window size $\eta=20$ pixels, the patch size $\rho=7$ pixels and $\e=5$ pixels. For a fixed patch, we only keep the $m_0=10$ best matches and $h=100$. For the tests of Fig.~\ref{figScat}-\ref{weight_recompute}, %% 13-17??
we set $\eta=100$, $\rho=5$, $\e=3$, $m_0=6$, $h=100$. For the smaller example of Fig. \ref{figLen} we chose $\eta=20$, $\rho=5$, $\e=1$, $m_0=8$, $h=100$.
%, $p=4$.
In Fig.~\ref{nl+atom} %% 18?
we chose $\eta=60$, $\rho=9$, $\e=3$, $m_0=10$, $h=100$. %, $p=4$.
 \item The constrained minimization problem (\ref{E_NL}) where the weight is computed with the atom-based distance $\delta^1$. In the tests of Fig. \ref{figToy} and Fig. \ref{figBarb}, we consider the $n_0=25$ first atoms with a moment $p=4$. As above we take $\e=5$ pixels, $h=100$ and for a fixed patch, we only consider the $m_0=10$ best matches. For Fig.~\ref{figScat}-\ref{figScatParallel}, %% 13-16
we set $\eta=100$, $\e=3$, $m_0=6$, $n_0=18$, $h=100$, $p=4$. To produce Fig.~\ref{nl+atom} %% 18 ?
 we took $\eta=60$, $\e=3$, $m_0=10$, $n_0=18$, $h=100$, $p=4$. For the small image of Fig. \ref{figLen} we picked $\eta=20$, $\e=1$, $m_0=8$, $n_0=18$, $h=100$, $p=20$.
\end{itemize}

\newbfb{In our first example (Fig. \ref{figLen}), we consider a \emph{small} $64\times64$ image and both exponents $\alpha=1$ and $\alpha=2$ in equation (\ref{E_NL}):}

\newbfb{For this first example, both choices for the parameter $\alpha$ provide improved results. The minimization is possible either by approximating the $\ell^1$ norm or by using modern splitting algorithms as introduced in \cite{Raguet,Combettes2011}}. 

In this experiment, it is interesting to observe that the robust choice $\alpha=1$ is particularly useful with the SSD distance which is likely to produce outliers (spurious matches), and almost useless with the oracle distance (see Fig. \ref{simul_ideal}) which produces perfect matches. Our distance seems to lie in-between, in a region where the choice is apparently less decisive.

Henceforth, we shall only consider the exponent $\alpha=2$ which is numerically tractable for larger images (the case $\alpha=1$ being in 
practice, as expected, quite slower than the quadratic case). The experiment that follows (Fig. \ref{figToy}) is a simple but larger $128\times128$ toy example.

\vspace{-0.7cm}
Fig.~\ref{figToy} and~\ref{figBarb} show an example with the
Mask of Fig.~\ref{Atom_mask}. In this case, the best results
are obtained with the atom-based distance. In particular in Fig.~\ref{figToy},
the oscillations with frequencies in the mask are almost perfectly
recovered and the spurious contamination is very low. %% commentaire
%% sur alpha=1??
\newline
\indent Let us see how these methods perform for the acoustic scattering problem of Section \ref{sec_scatt}. To do so, we consider the mask introduced in Fig.~\ref{spectrum_scatt}.
\vspace{-0.7cm}
\vspace{-0.7cm}
In this case the best results are obtained with the standard Total
Variation: this is consistent with the fact that the object to
recover is piecewise constant and the loss is mostly in the
high frequencies (for further details on these qualitative properties see \cite{JalalzaiJump,JalalzaiStairc,JalalzaiPhD}). 
%Since the Total Variation is well-behaved with functions that are piecewise constant, we could have expected that the TV restoration would perform this well.

In the next simulation, we are going to assume that the Fourier coefficients that we kept are contaminated by a Gaussian noise of magnitude $0.03{\|g_0\|}_2$.
\vspace{-0.7cm}
% \vspace{-0.7cm}
In this case, the Total Variation based optimization find it much harder
to recover the data and the Atom-based distance yields satisfactory results.
Now considering the Born approximation (\ref{Born_approx}), let us use the data that comes out of the direct problem. In a sense, this amounts to adding to the Fourier coefficients a noise whose distribution is unknown. Our method is able to get rid of all the spurious objects.
%\vspace{-0.6cm}
%\vspace{-0.6cm}
In this kind of problems it is usually important to distinguish objects that are very close. Let us consider such a situation:
%\vspace{-0.6cm}
\begin{figure}[H]
\centering
% \begin{minipage}[c]{\linewidth}
\begin{tabular}{ccccc}
Original $g_0$ & Corrupted $g$ & SSD - $\delta^3$ & NL-Atom - $\delta^1$ & TV restored\\
\includegraphics[width=0.19\linewidth]{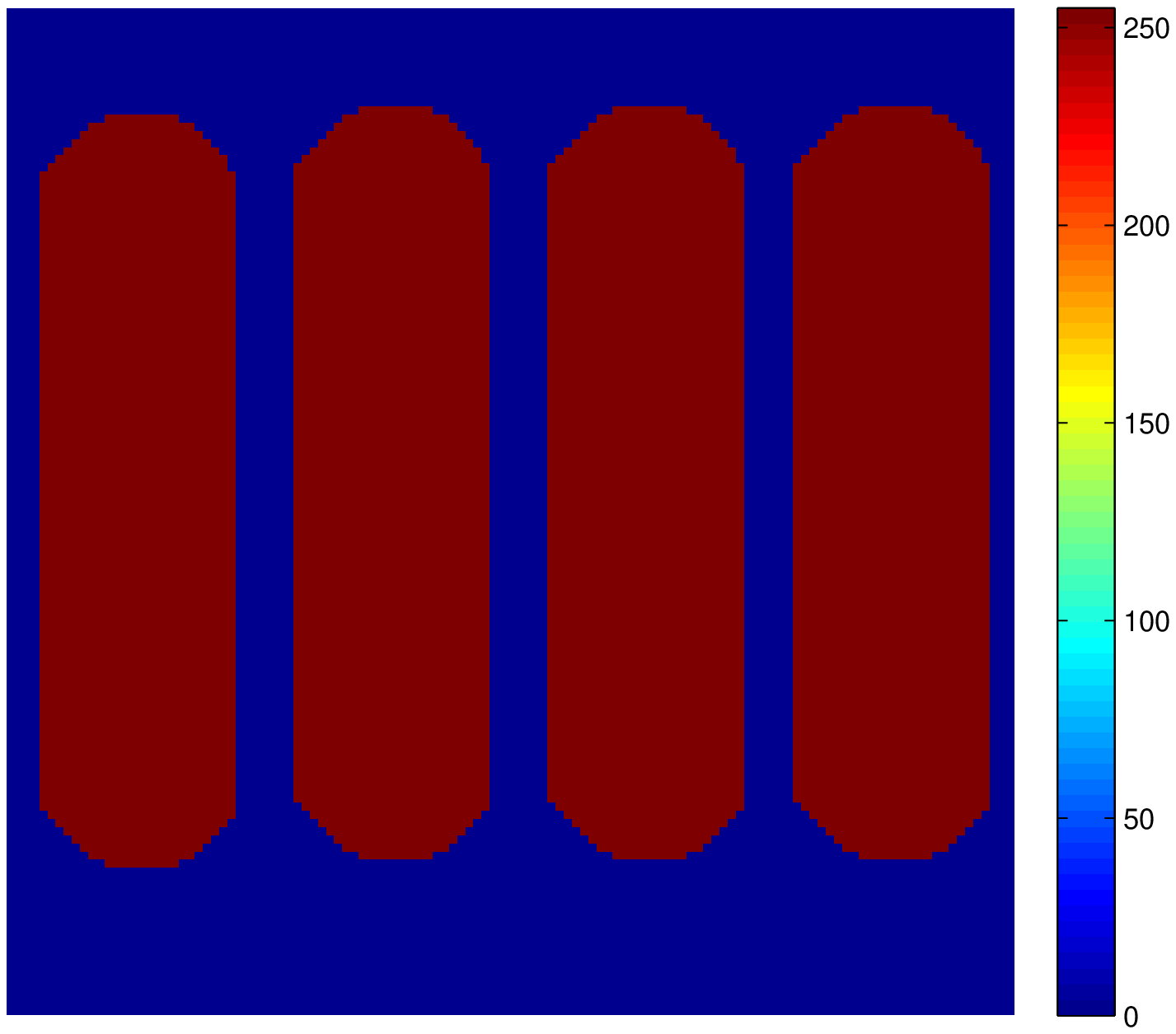}
&\includegraphics[width=0.19\linewidth]{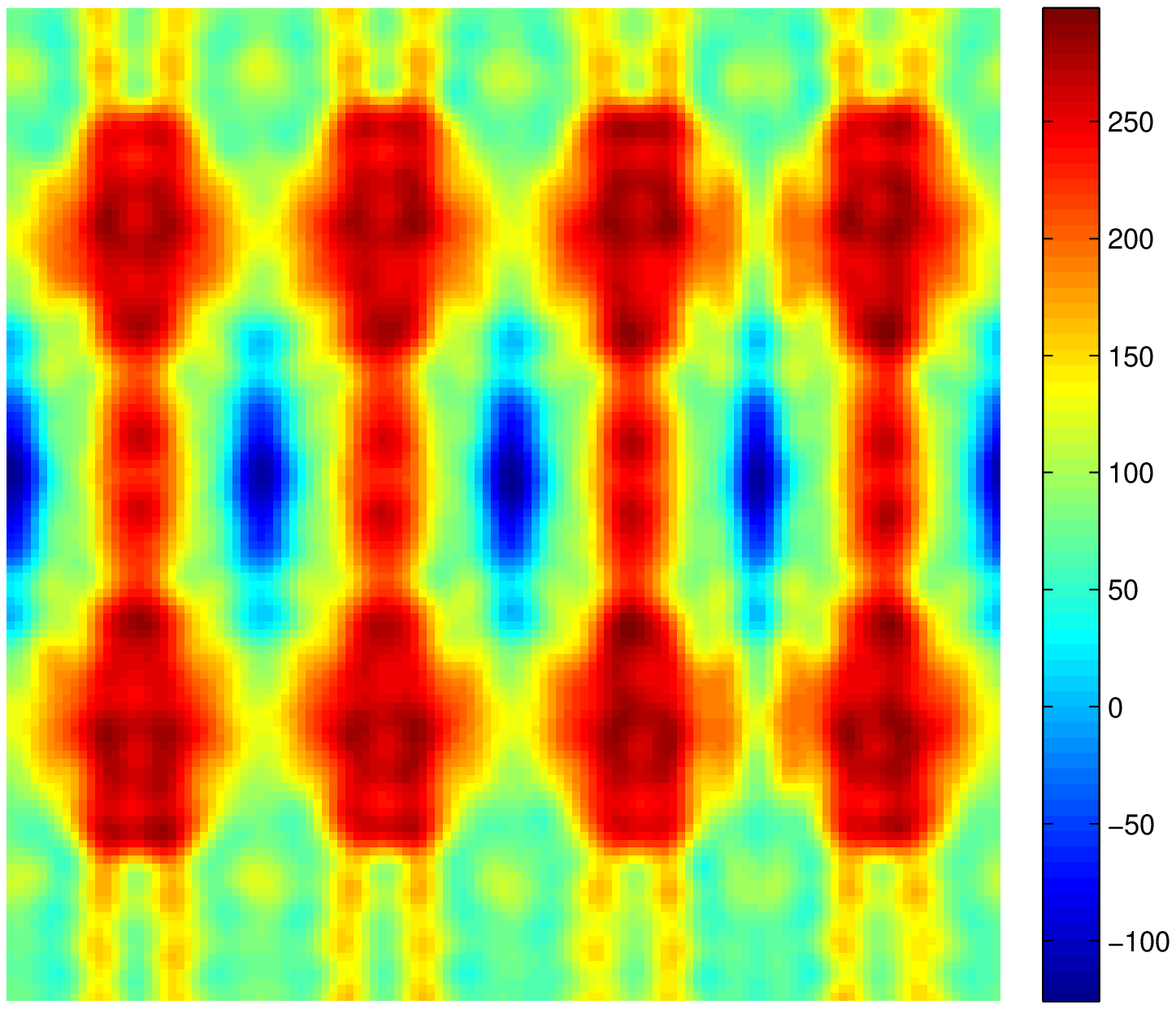}
&\includegraphics[width=0.19\linewidth]{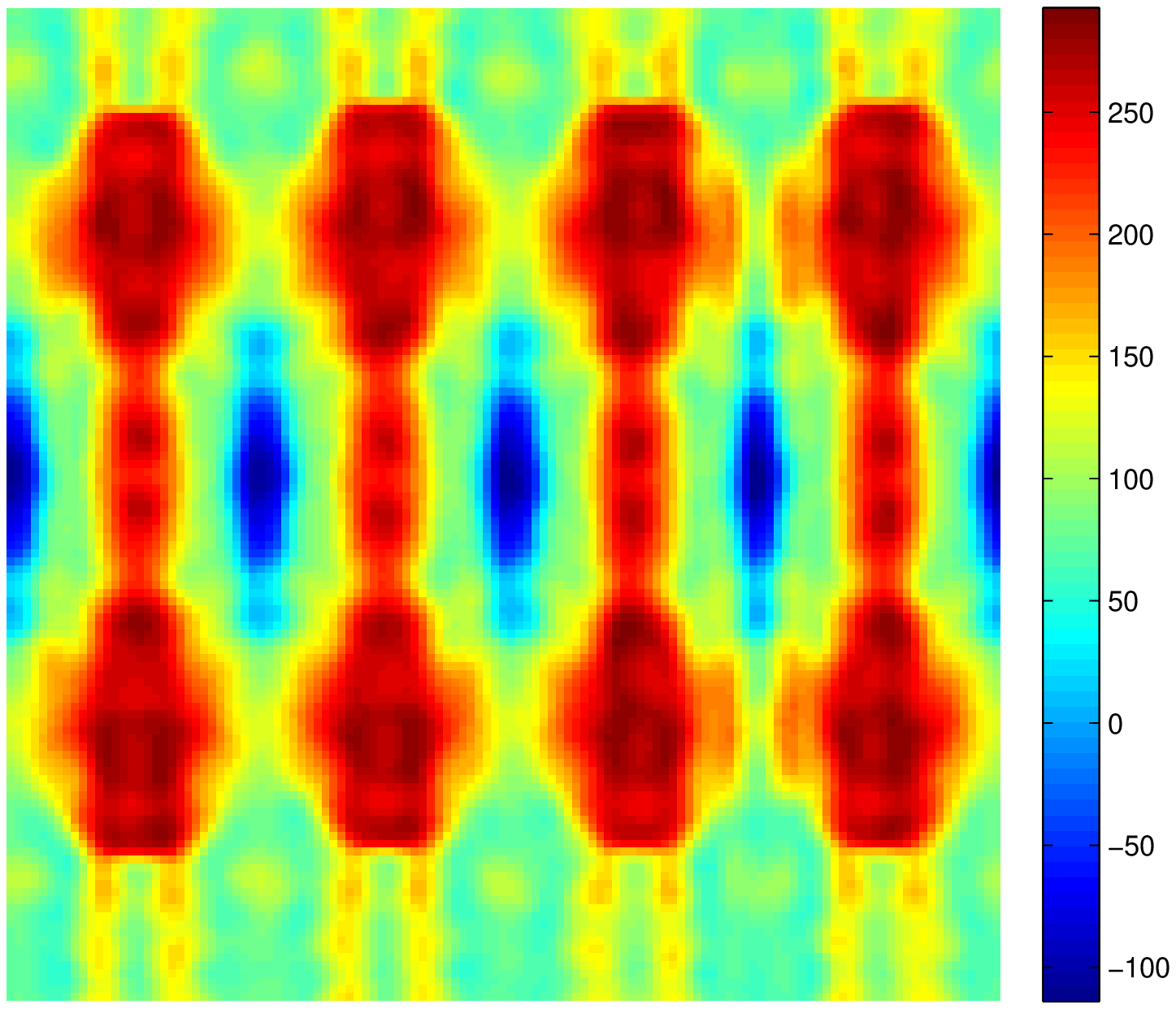}
&\includegraphics[width=0.19\linewidth]{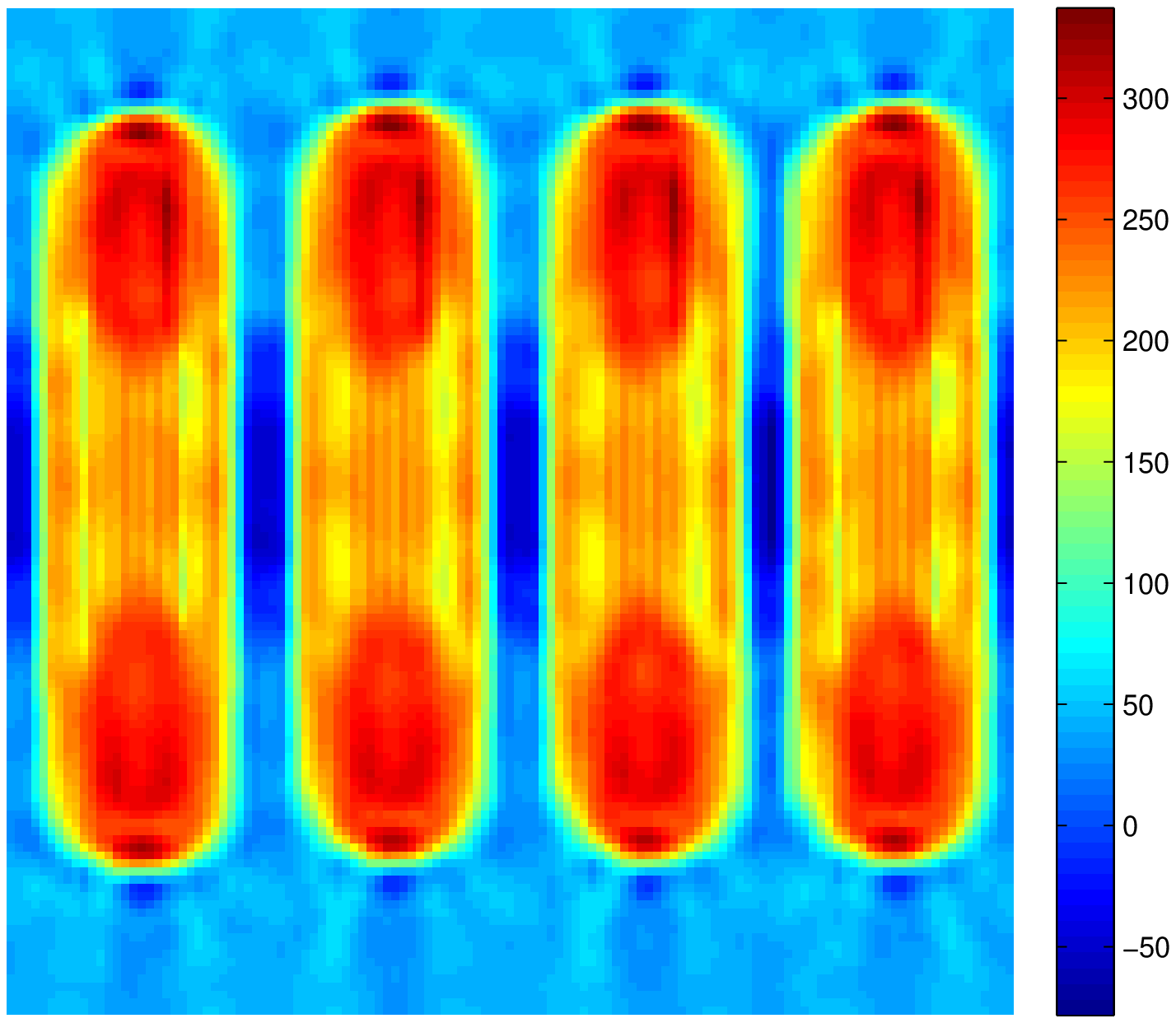}
&\includegraphics[width=0.19\linewidth]{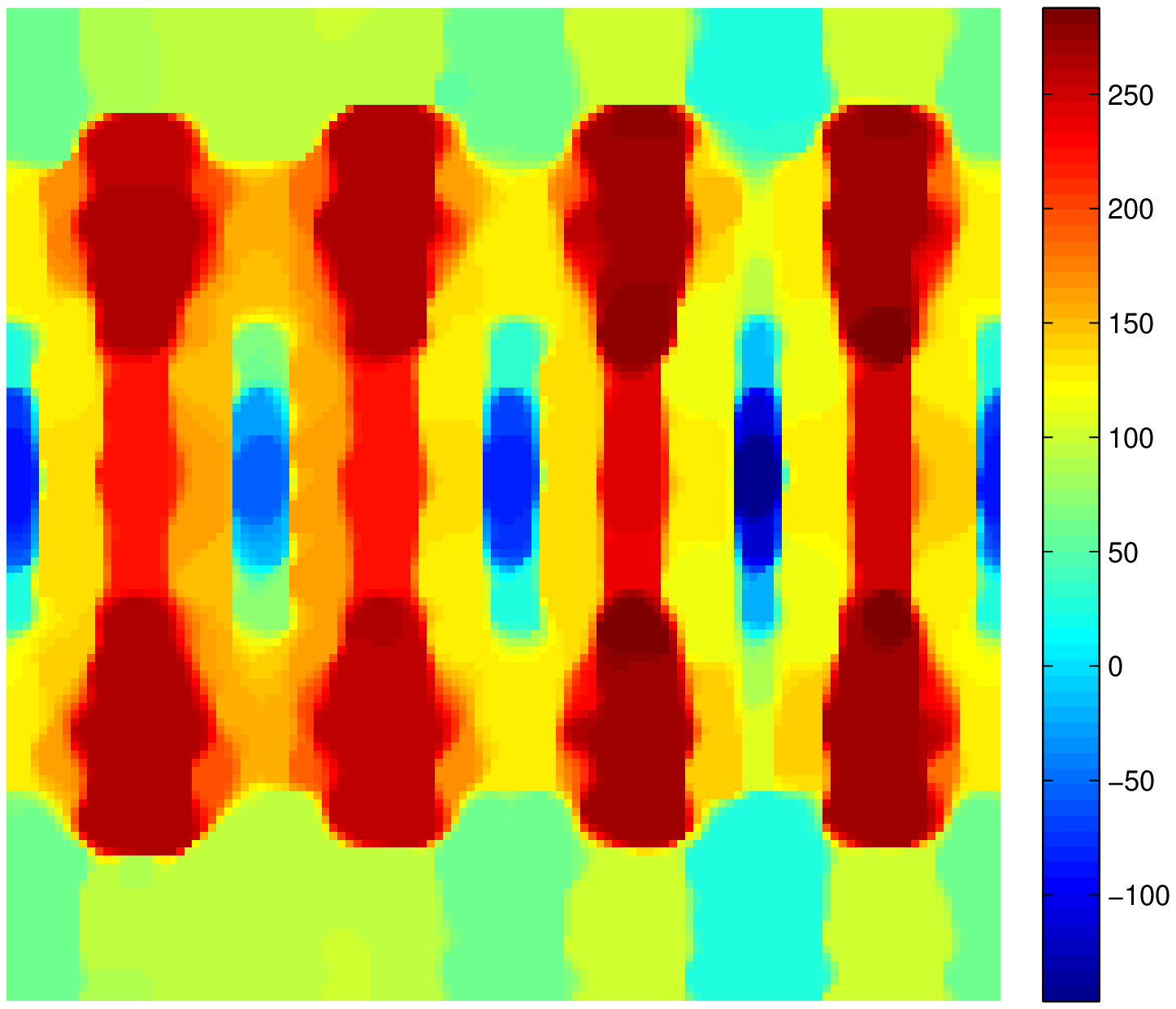}\\
%Corrupted image $u$ & $\delta^3$-NL restored & Atom-based restoration & TV restored\\
& PSNR=8.2dB &PSNR=8.4dB & PSNR=14.5dB & PSNR=9.6dB\\
\includegraphics[width=0.19\linewidth]{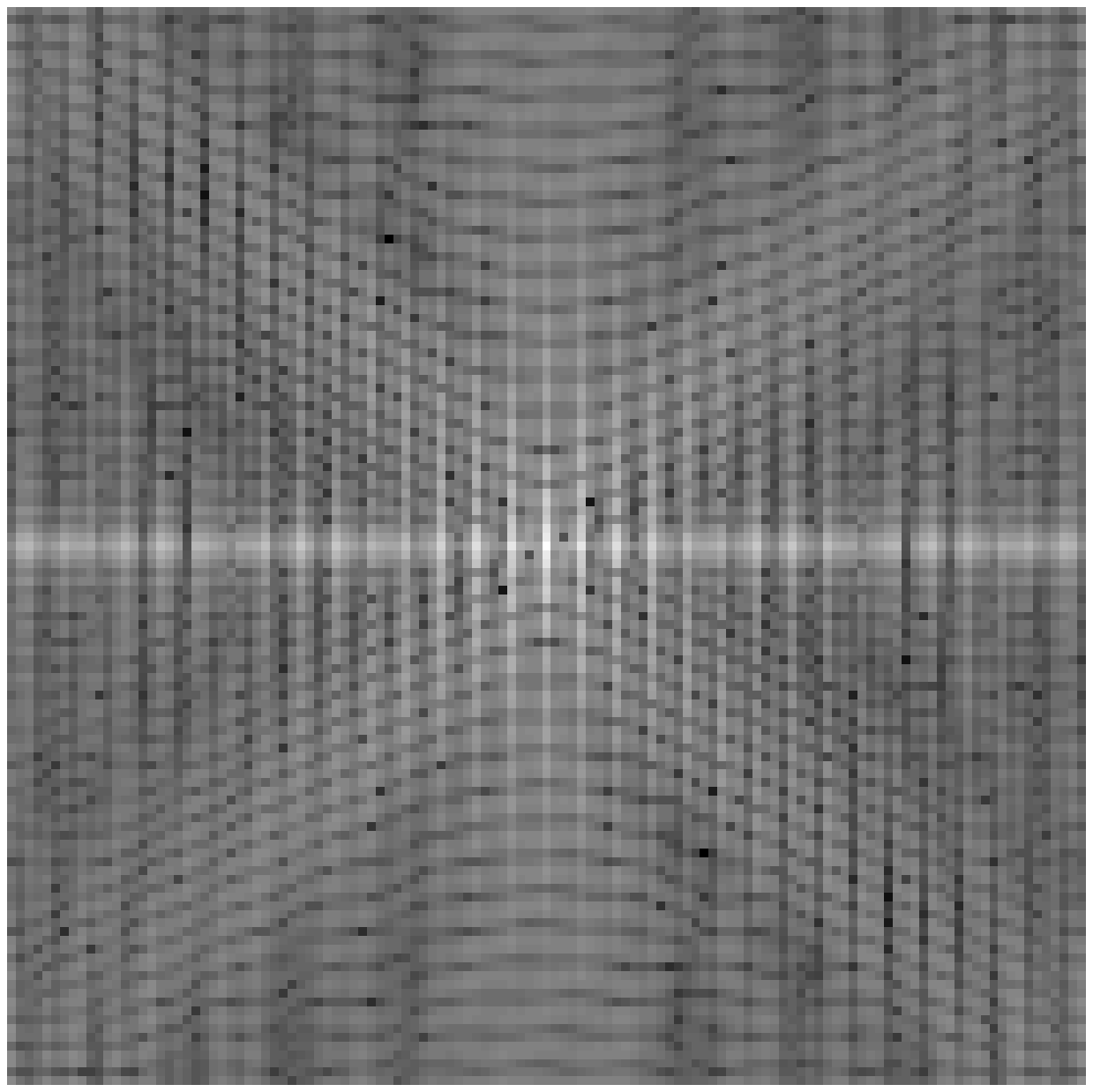}
&\includegraphics[width=0.19\linewidth]{scatt_NL_14}
&\includegraphics[width=0.19\linewidth]{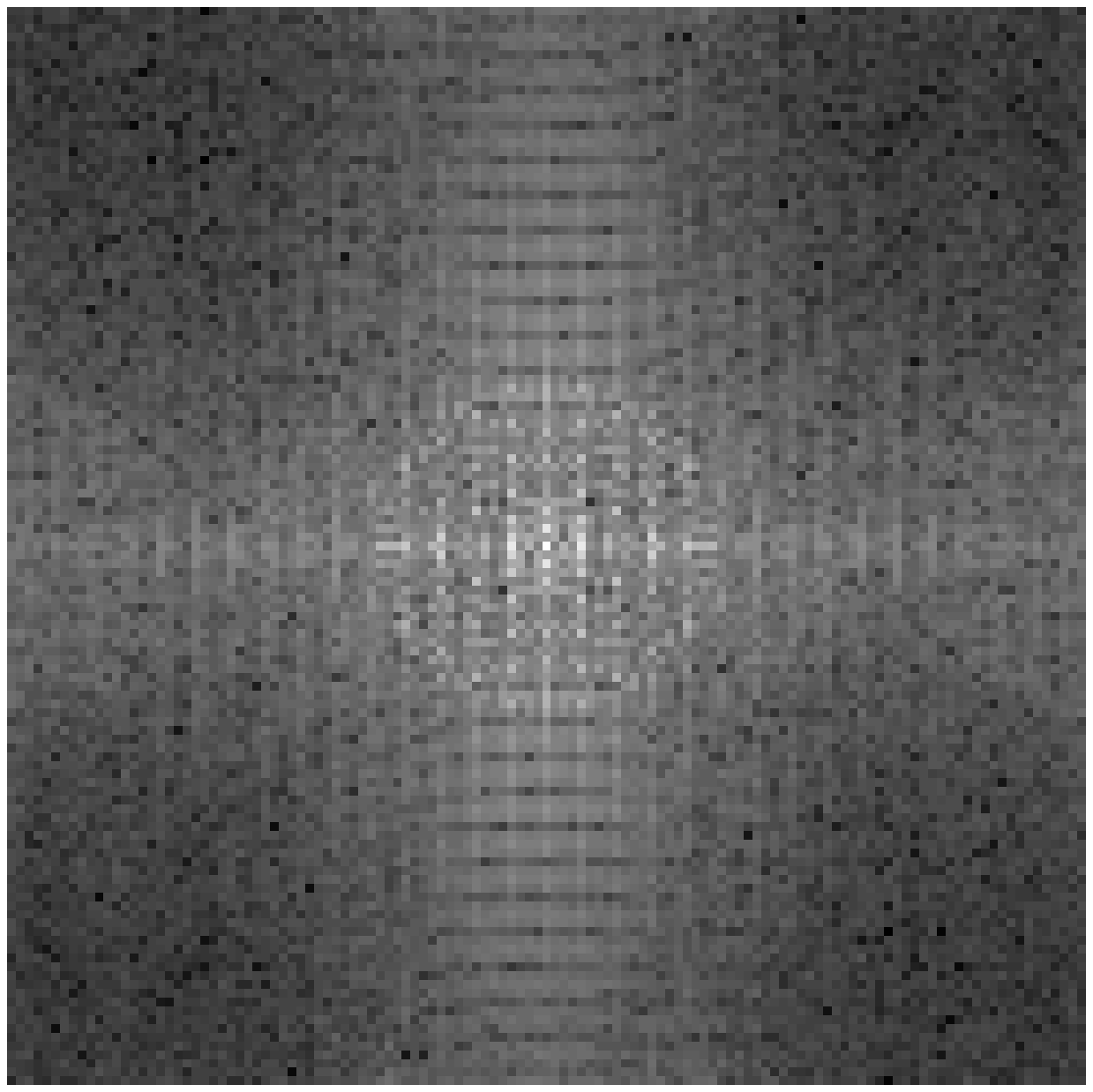}
&\includegraphics[width=0.19\linewidth]{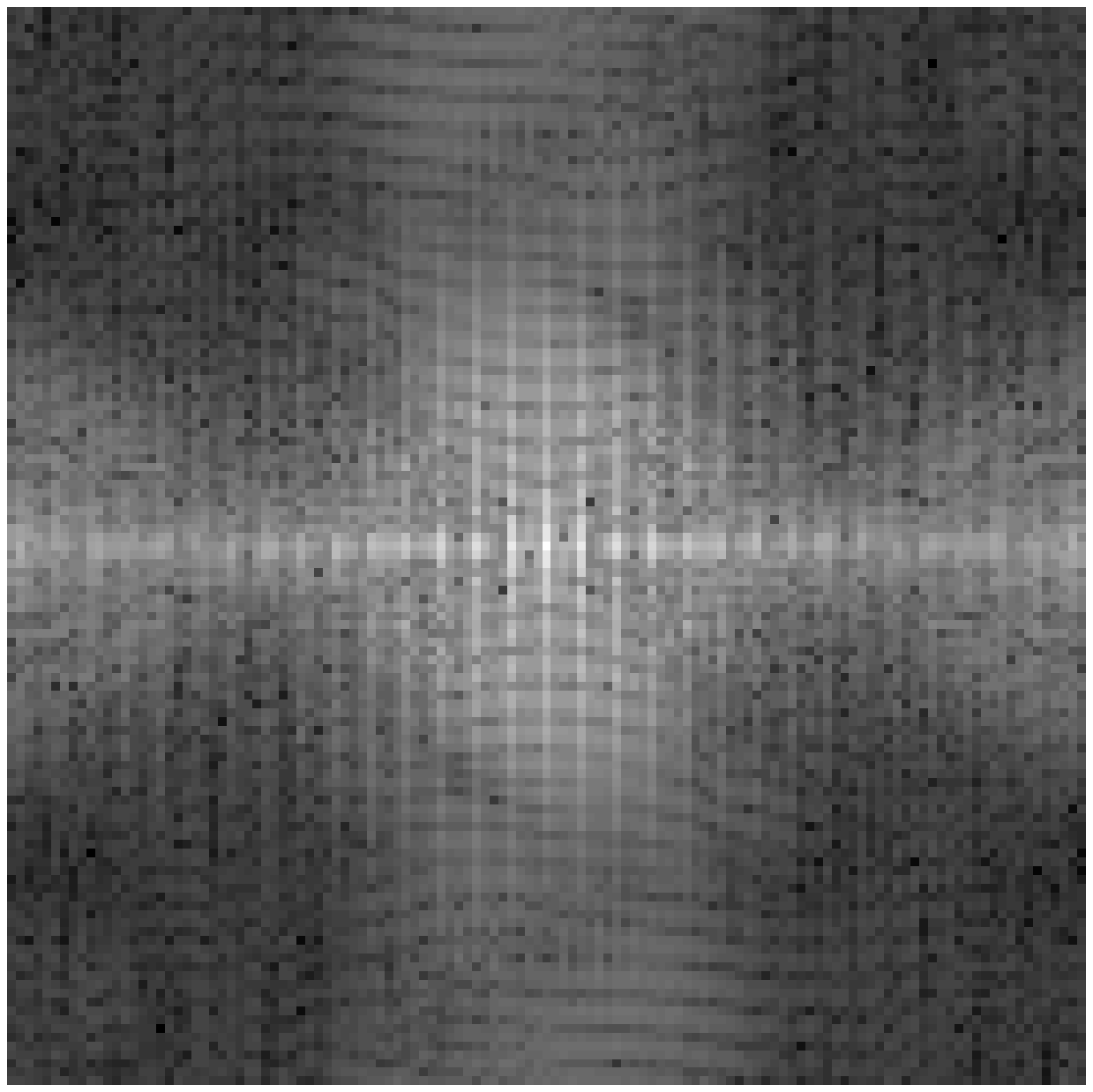}
&\includegraphics[width=0.19\linewidth]{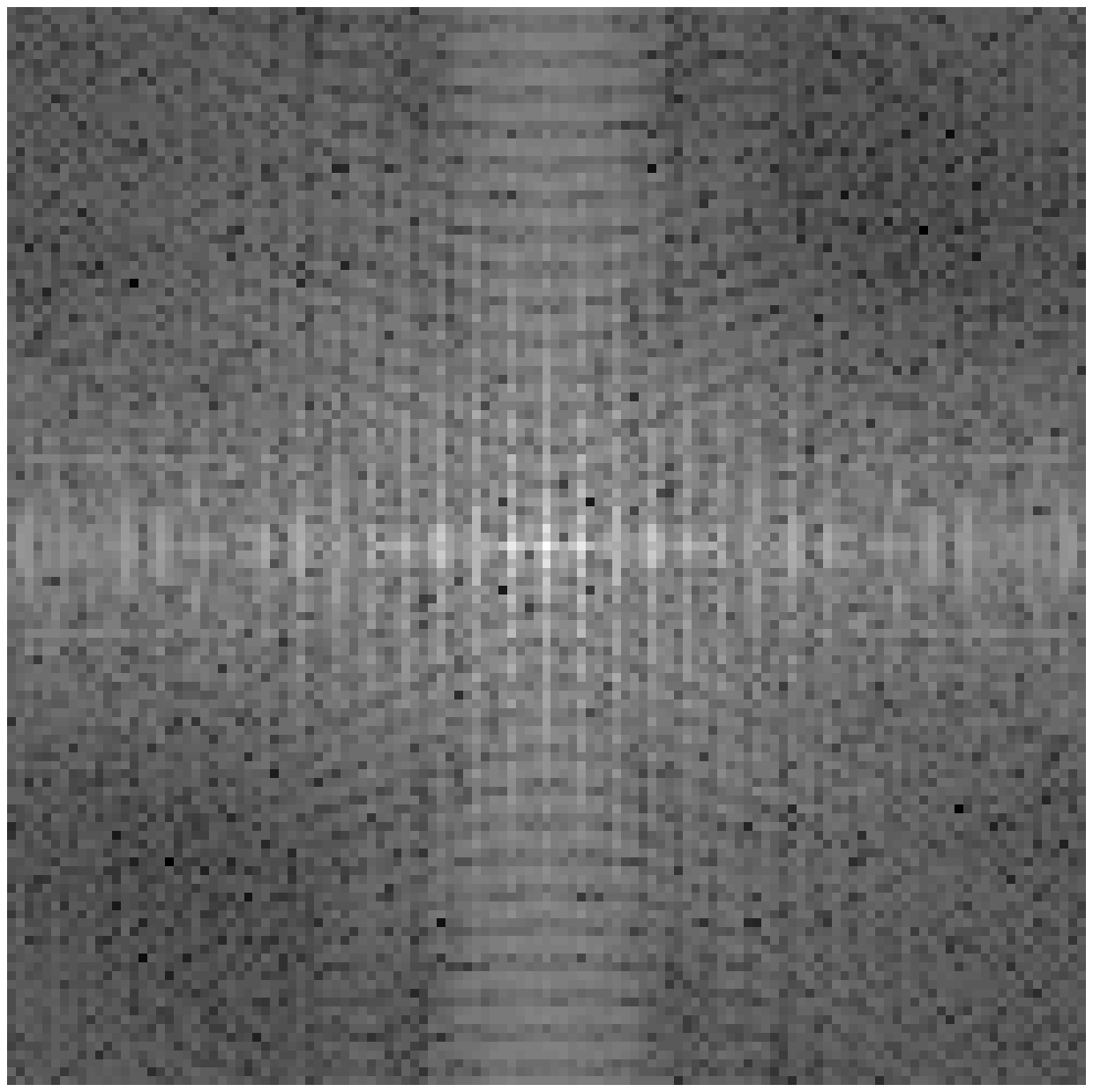}\\
%Spectrum of $u+v_{\delta^3}$ & Spectrum of $u+v_{\delta^1}$ & Spectrum of TV restored\\
\end{tabular}
%\vspace{-0.4cm}
\caption%[]
{The original image $128\times128$ is made of 4 scatterers separated by 6 pixels. We then depicted respectively the corrupted and the restored scatterers and their spectra.}\label{figScatParallel}
% \end{minipage}
\label{close_objects}
\end{figure}
%\vspace{-0.7cm}
Actually, in the latter experiment, we were able to distinguish objects that are separated by $0.56\la$ which is quite close to the theoretical limit $0.5\la$. \newbf{We were able to achieve this result since, we recall, two identical patches in the original image $g_0$ do not necessarily match according to the SSD but do match according to the atom-based distance. This explains why in Fig.~\ref{close_objects}, the atom-based distance improves the matching step (as in Fig.~\ref{best_match}) and provides logically better results as observed in the different tests. In this case, there is to our knowledge no method but the one we propose to improve the results over the raw data.}\\

% \vspace{-0.7cm}
\noindent\textbf{Weight recomputation:}
In \cite{Gilboa06nl,Gilboa,Peyre08,Peyre09,Faccolio}, it was pointed out that one can get improved results (especially for inpainting problems) by allowing recomputation of the SSD-based weight on the being restored image. We depict in Fig.~\ref{weight_recompute} the restorations one can expect after several recomputations.

% \vspace{-0.5cm}
% \begin{figure}
% \begin{minipage}[c]{\linewidth}
% \centering
% \begin{tabular}{cccc}
% \includegraphics[width=0.19\linewidth]{rayures_poids_NL_means_recomputed_NL_7}
% &\includegraphics[width=0.19\linewidth]{rayures_poids_NL_means_recomputed_NL_13}
% &\includegraphics[width=0.19\linewidth]{barbara_NLmeans_recomputedNL_7}
% &\includegraphics[width=0.19\linewidth]{barbara_NLmeans_recomputedNL_13}\\
% PSNR=12.0dB & & PSNR=24.9dB &
% 
% \end{tabular}
% \end{minipage}
% % \vspace{-0.5cm}
% \caption
% {Restored images and their respective spectra after weight recomputations.}
% \end{figure}

% \vspace{-0.5cm}

%To be able to recompute the atom-based distance one should understand how to extend the confidence regions on the Fourier side.

This procedure is really cumbersome and does not always improve results over the atom-based method. The weight recomputation is not possible for the atom-based distance we introduced since the distance computed on the restored image is exactly the same as the one computed on the corrupted image. However, our method can be used as an initialization for the classical weight recomputation to improve results further. This is the strategy we adopt in the following tomography problem where the Fourier coefficients got corrupted by a Gaussian noise of magnitude $0.3{\|g_0\|}_2$:
\vspace{-0.6cm}
\begin{figure}[H]
\begin{minipage}[c]{\linewidth}
\centering
\begin{tabular}{cccc}
Original $g_0$ & Spectrum of $g_0$ &Corrupted $g$ & Spectrum of $g$\\
\includegraphics[width=0.19\linewidth]{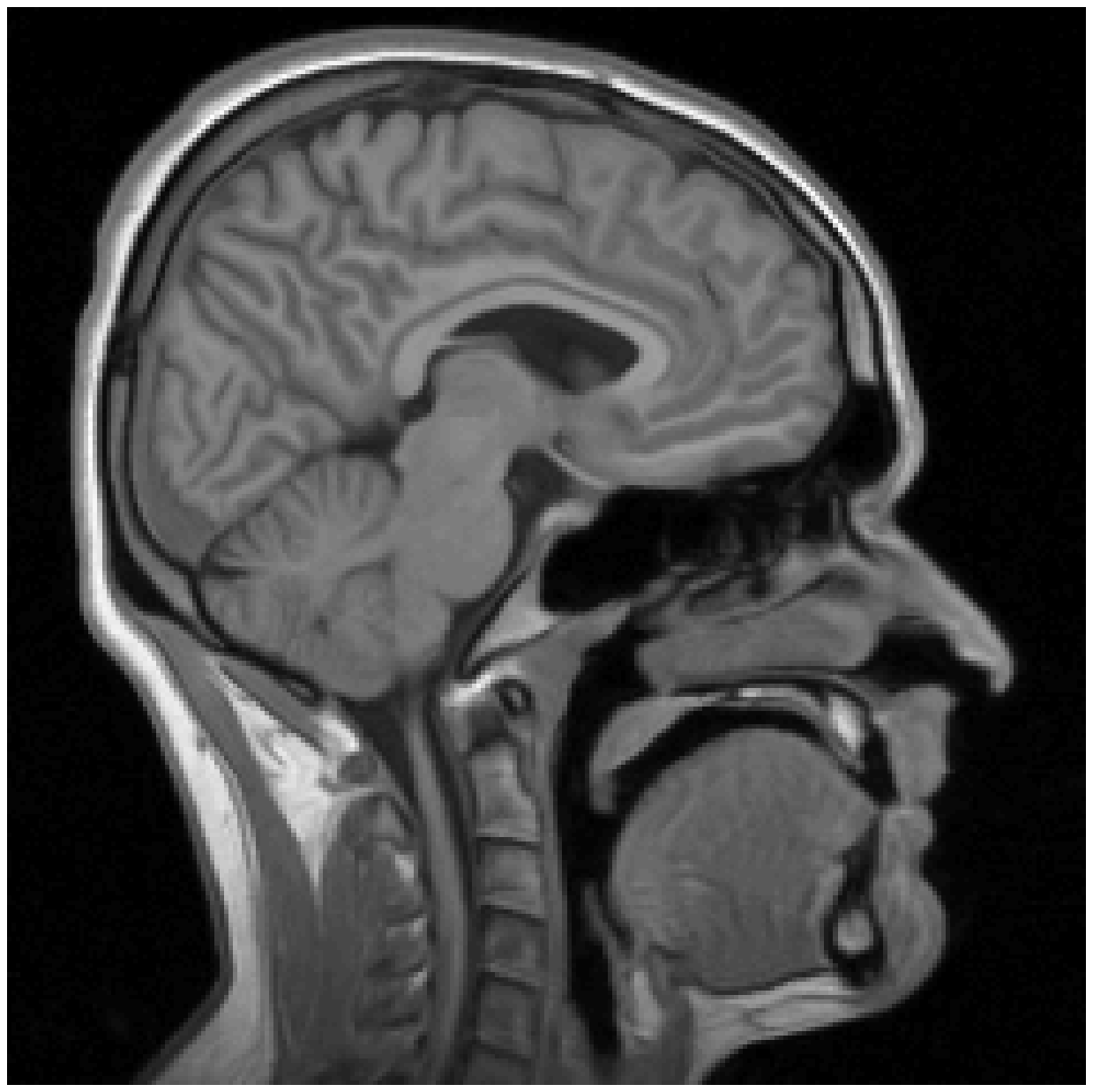}
&\includegraphics[width=0.19\linewidth]{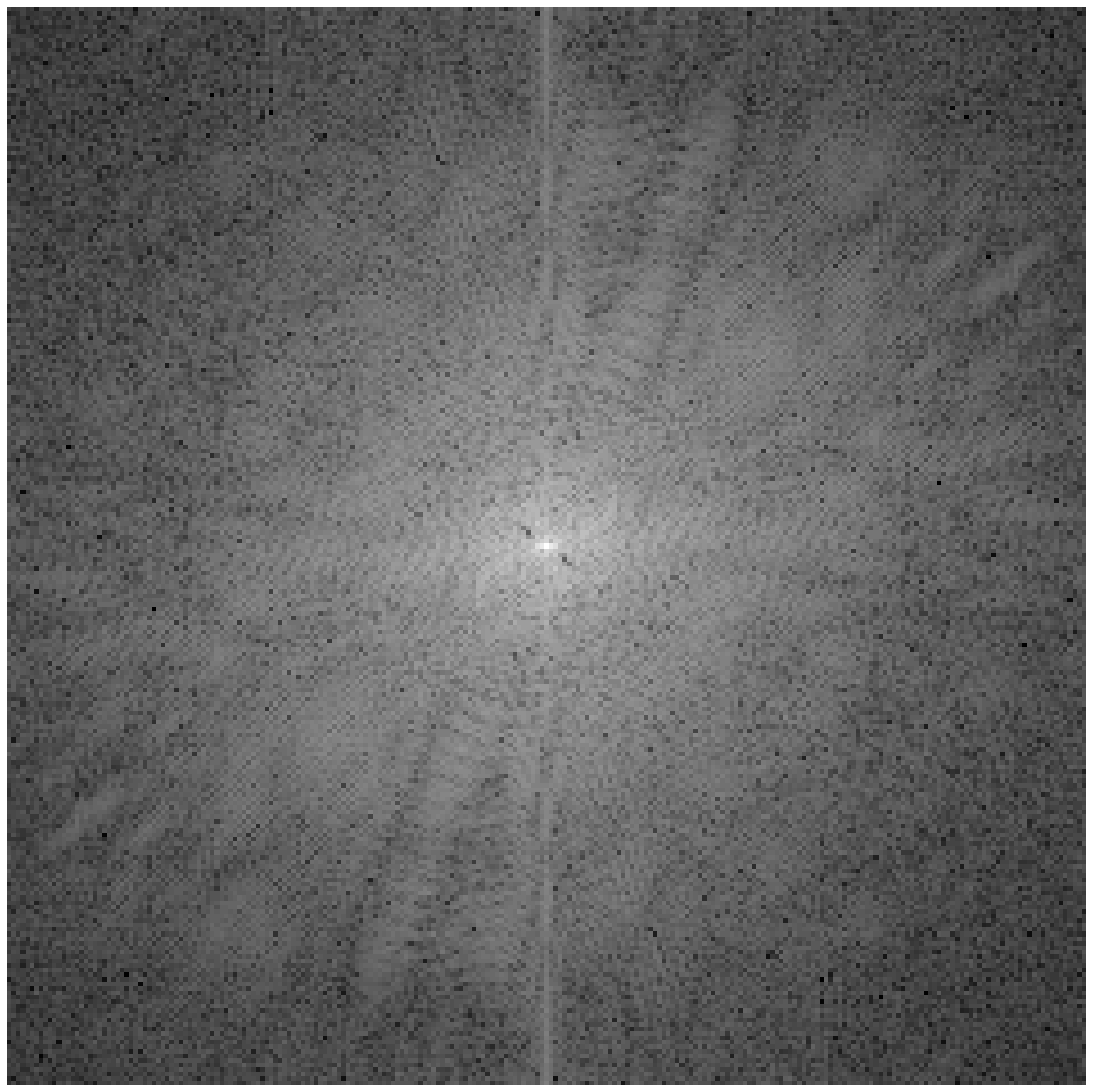}
&\includegraphics[width=0.19\linewidth]{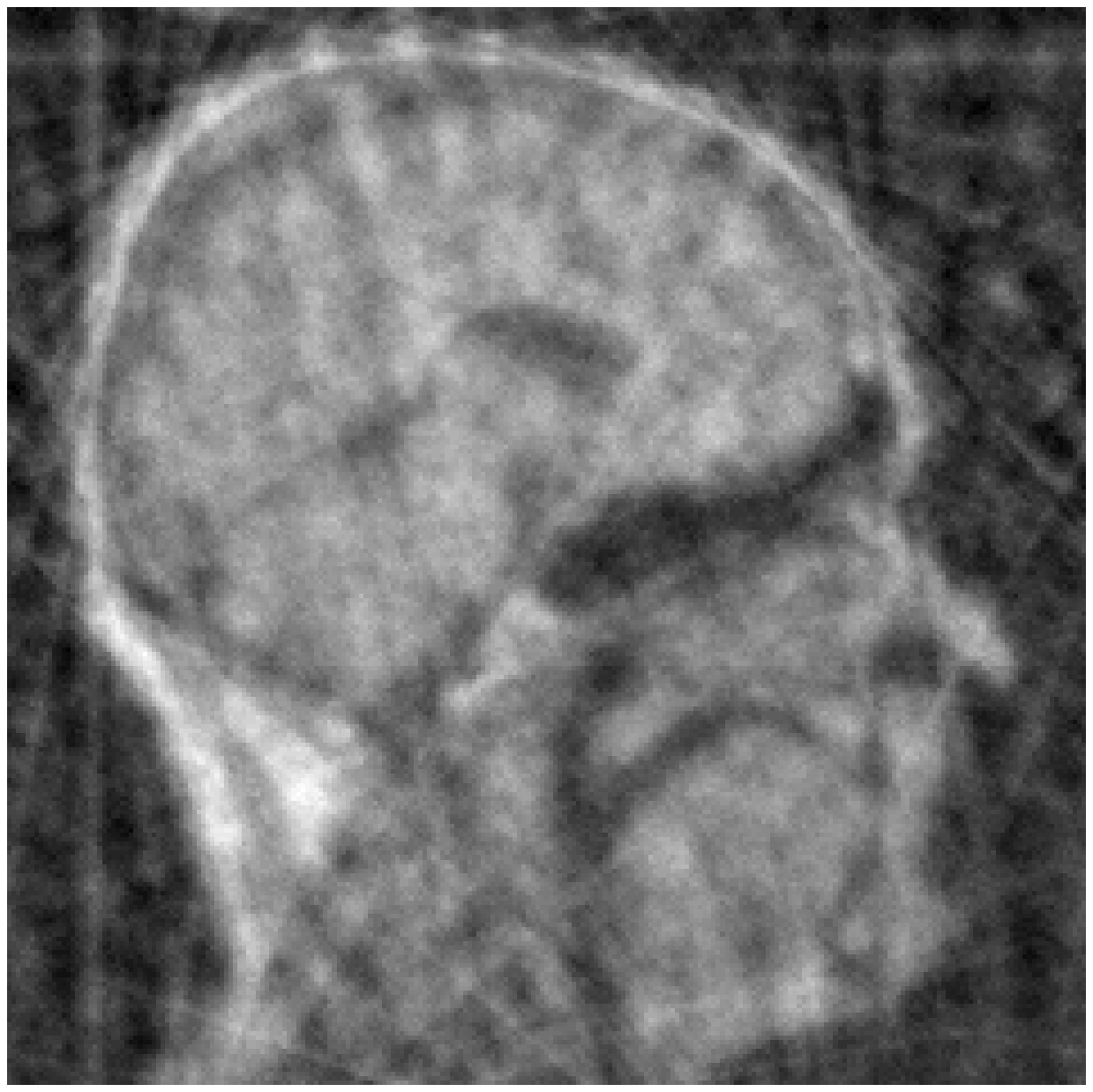}
&\includegraphics[width=0.19\linewidth]{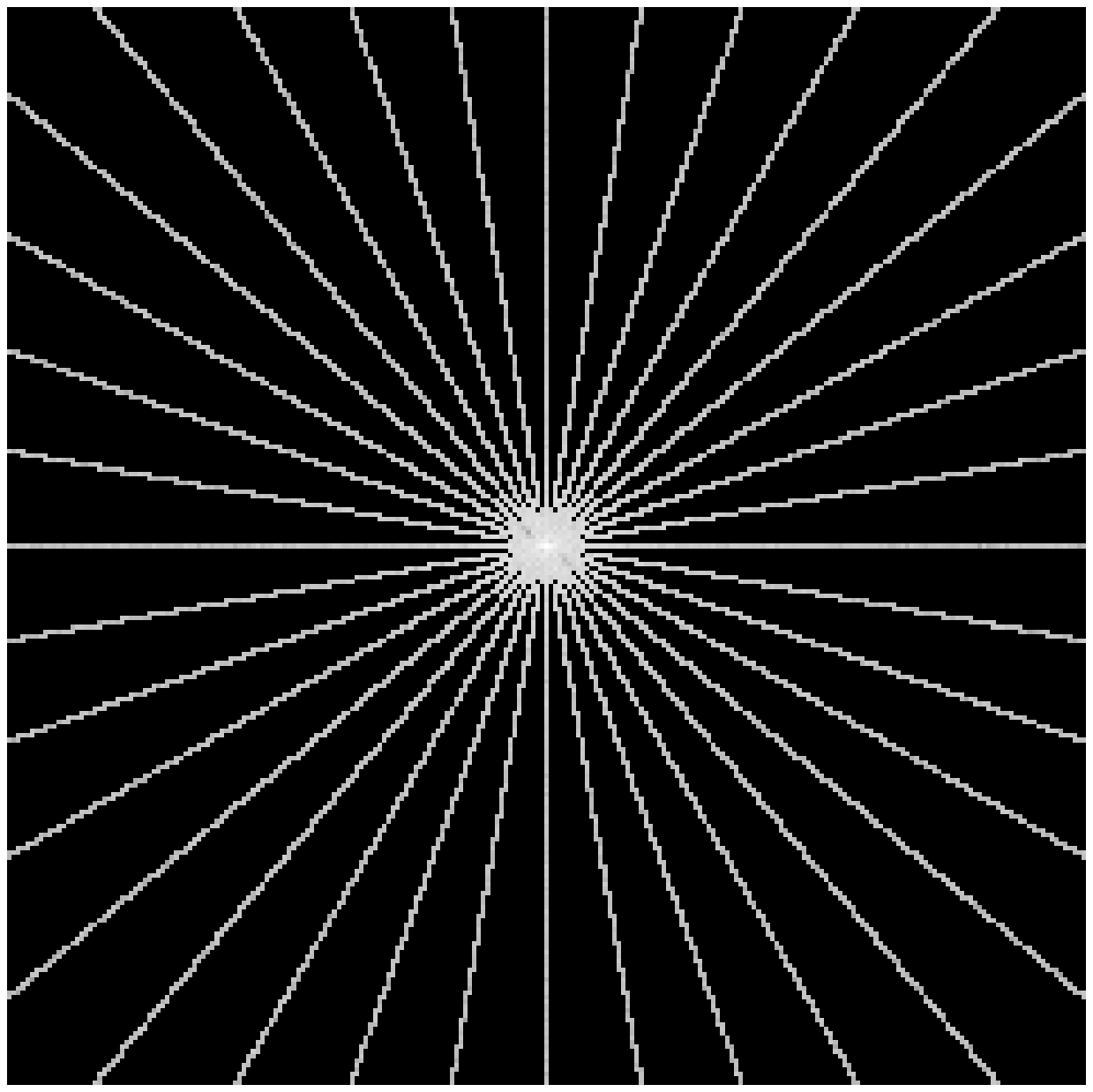}\\
&&PSNR=22.4dB&\\
\end{tabular}
\end{minipage}
\begin{minipage}[c]{\linewidth}
\centering
\begin{tabular}{ccccc}
% Original $g_0$ & Corrupted $g$ & 
&  & $\delta^1$ then one & $\delta^3$ recomputed& \\
SSD - $\delta^3$ &NL-Atom - $\delta^1$&computation of $\delta^3$  & 20 times &TV restored\\
% \includegraphics[width=0.19\linewidth]{tomo_NL_4sc_nc}
% &\includegraphics[width=0.19\linewidth]{tomo_NL_6sc_nc}
\includegraphics[width=0.19\linewidth]{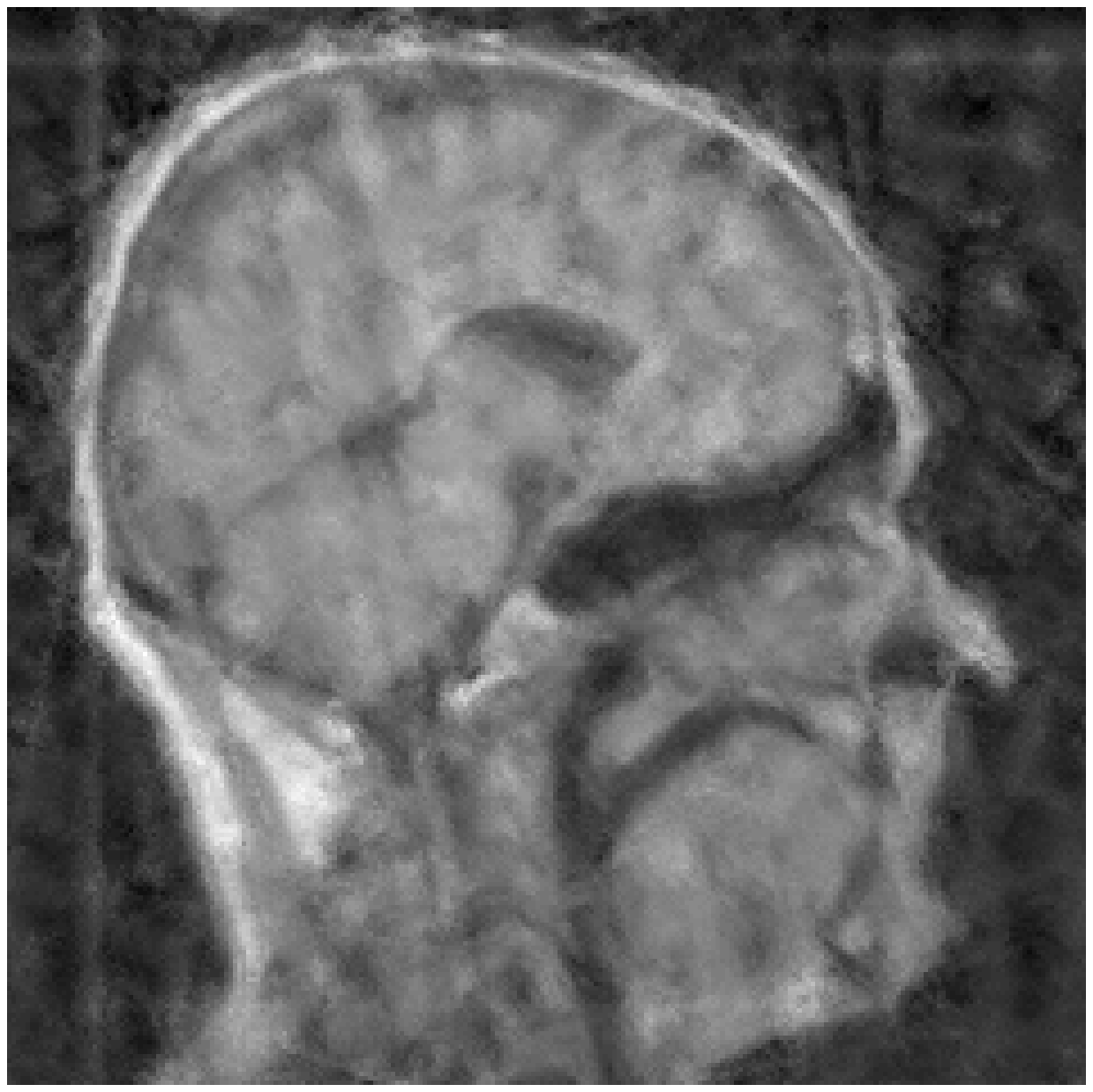}
&\includegraphics[width=0.19\linewidth]{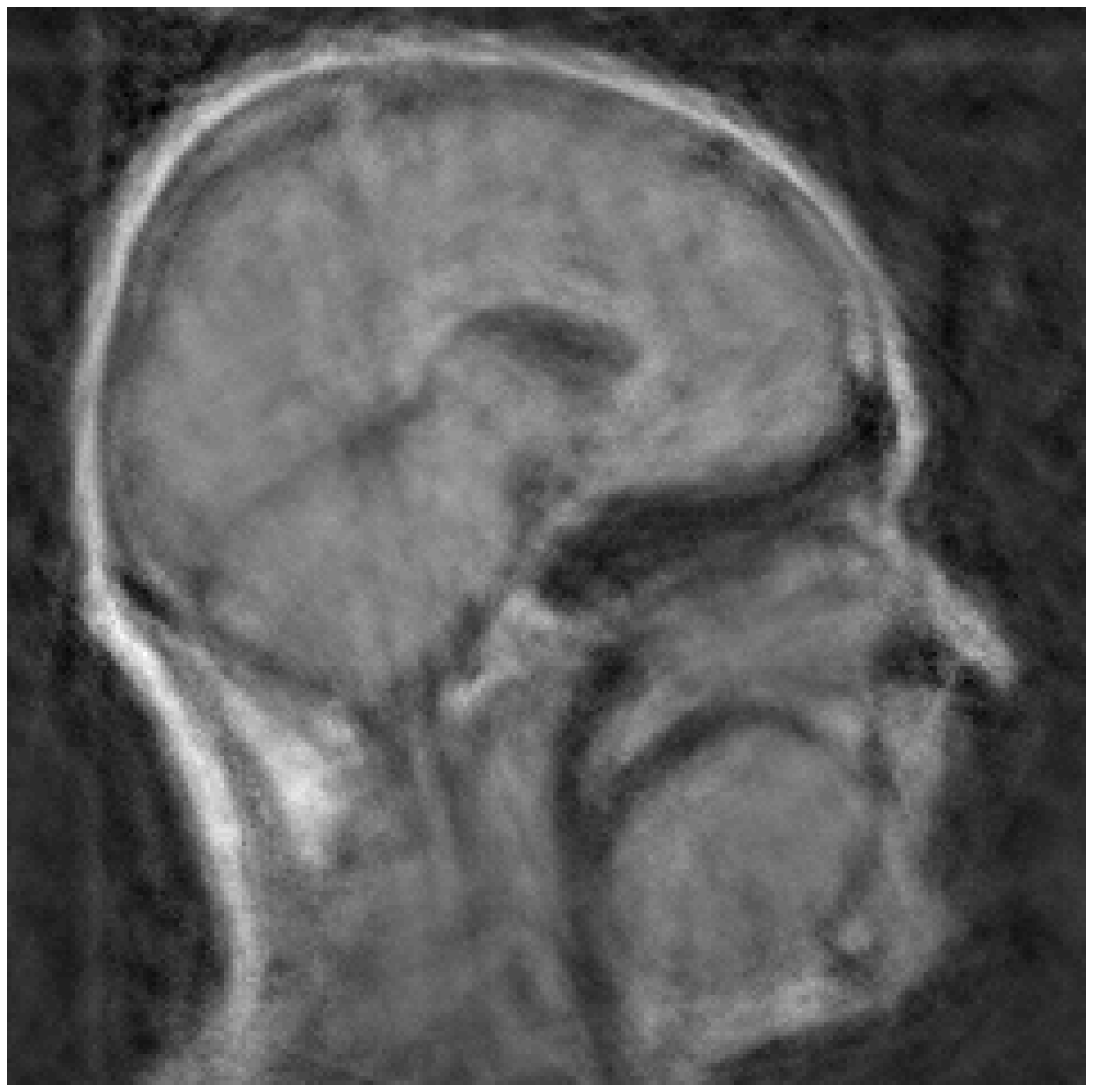}
&\includegraphics[width=0.19\linewidth]{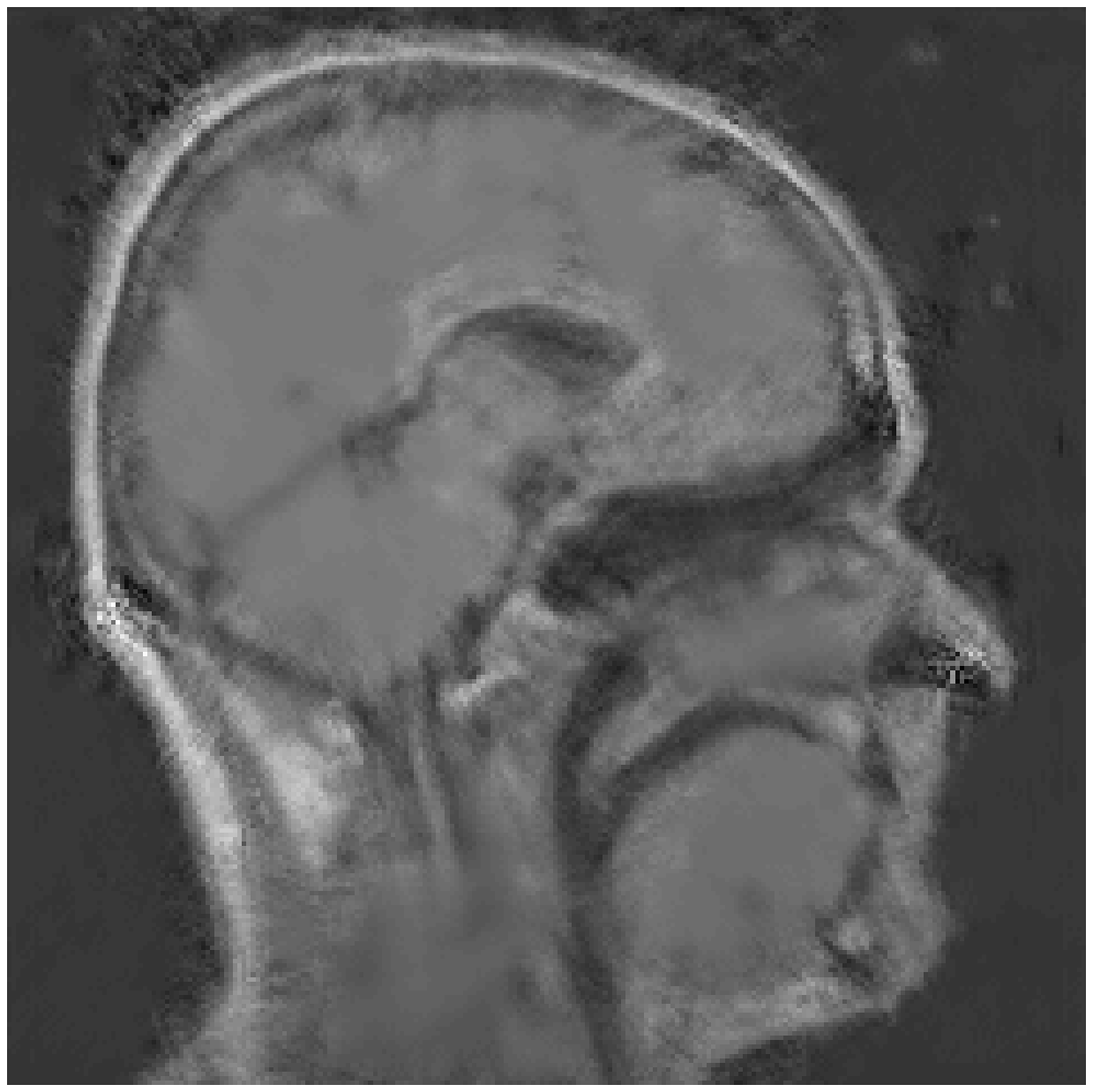}
&\includegraphics[width=0.19\linewidth]{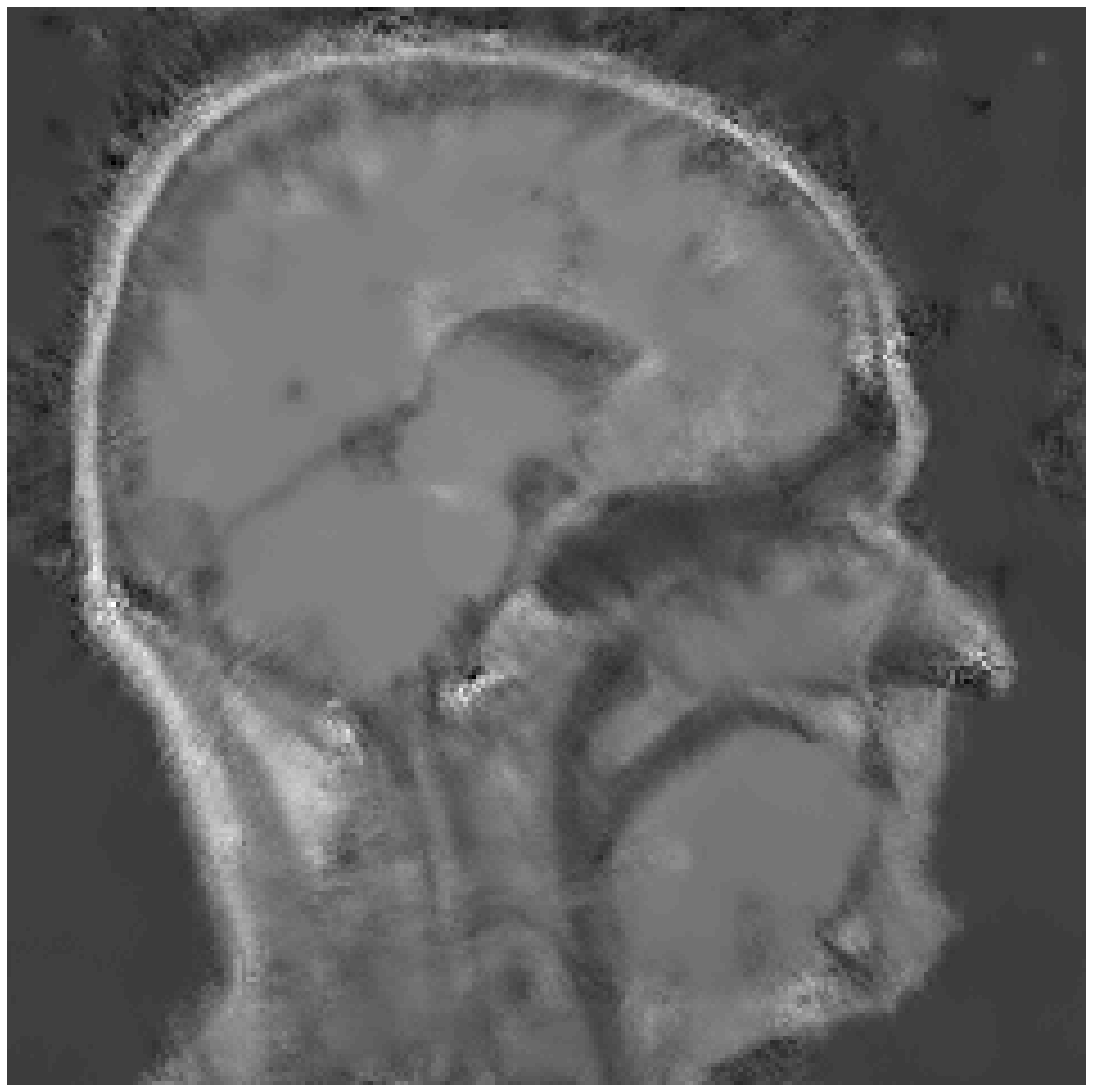}
&\includegraphics[width=0.19\linewidth]{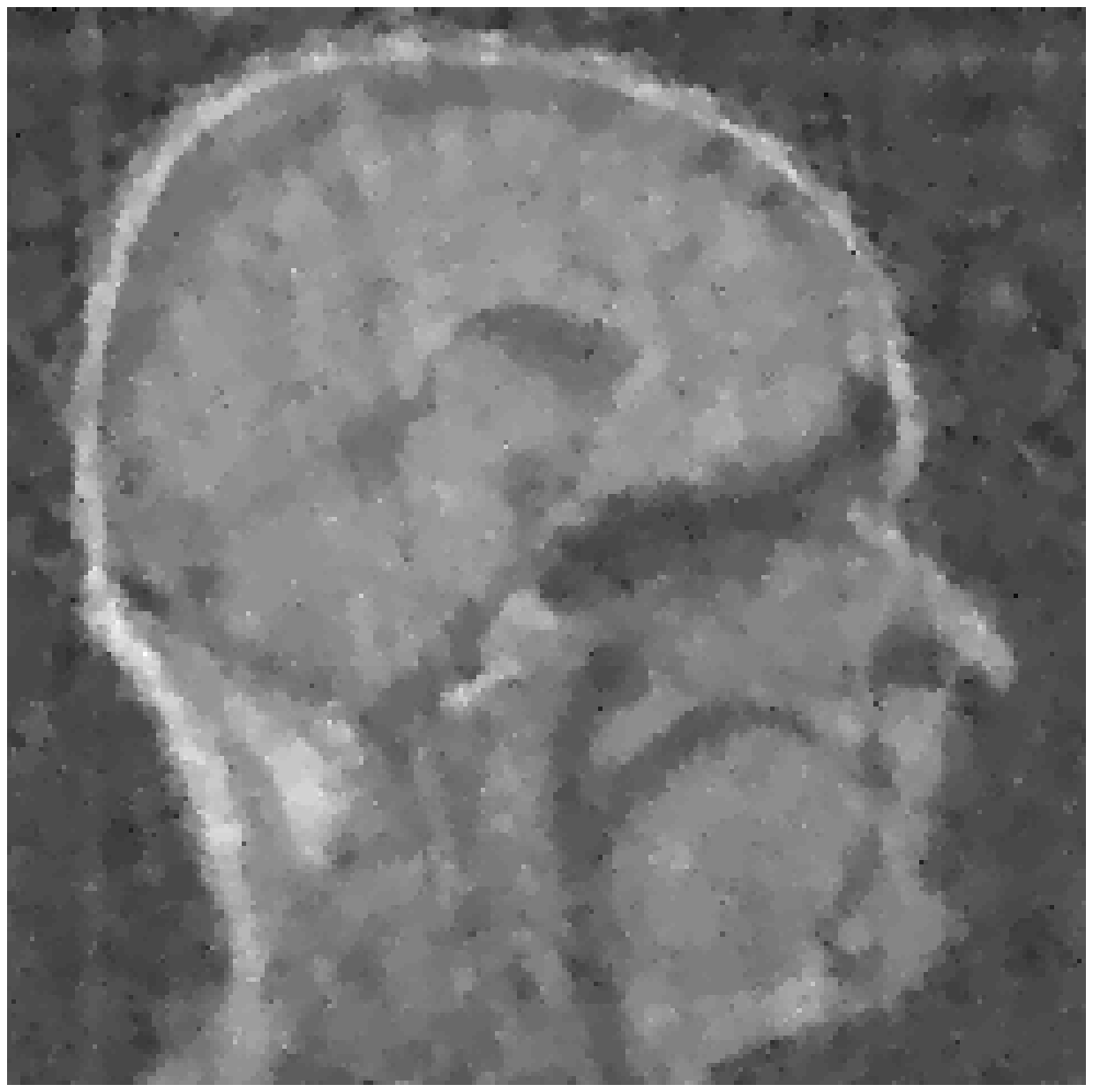}\\
%Corrupted image $u$ & $\delta^3$-NL restored & Atom-based restoration & TV restored\\
% & PSNR=22.4dB &
PSNR= 23.8dB & PSNR=24.9dB & PSNR=25.8dB & PSNR=24.8dB &PSNR=23.6dB\\
% % \includegraphics[width=0.18\linewidth]{tomo_NL_11}
% % &\includegraphics[width=0.18\linewidth]{tomo_NL_14}
% \includegraphics[width=0.18\linewidth]{tomo_NLmeansNL_13}
% &\includegraphics[width=0.18\linewidth]{tomo_NL_13}
% &\includegraphics[width=0.18\linewidth]{tomo_NL_15}\\
% %Spectrum of $u+v_{\delta^3}$ & Spectrum of $u+v_{\delta^1}$ & Spectrum of TV restored\\
\end{tabular}
\end{minipage}
\vspace{-0.4cm}
\caption
{Restoration for a $240\times240$ tomography image.}
\vspace{-0.4cm}
\label{nl+atom}
\end{figure}

\noindent\textbf{Computation Times:} {As we claimed previously in our discussion on complexity, the computation time with our approach (which in a first step reduces the dimensionality to retain only the information which is useful for the matching) is quite reduced. For instance, in Fig.~\ref{nl+atom} we observed the following computation times
\begin{center}
\begin{tabular}{|c|c|}
\hline
Approach & Computation time\\
\hline
SSD - $\delta^3$ & 83s (dist) + 15s (iter) = 98s \\
\textbf{NL-Atom - $\delta^1$} &\textbf{8s (dist) + 25s (iter) = 33s}\\
$\delta^1$ then one $\delta^3$ &143s\\
$\delta^3$ recomputed 20 times & $>$30min\\
$TV$ restored &82s\\
\hline
\end{tabular}
\end{center}
We performed our experiments in Matlab, so the absolute CPU times are not really relevant.}

\section{Conclusion}

\newbfa{In this paper, we have considered the problem of reconstructing
an image with a known pattern of missing Fourier coefficients, by
means of a non-local method which assumes and exploits some
spatial redundancy of the original image. In order to detect
similar patches on the original image, we have introduced
an original similarity criterion which is different from a standard
quadratic distance, and is insensitive to the image degradation.
This distance is based on precomputed atoms which are then used to
filter out the degradation of the image, without destroying the local features.}

\newbfa{Then, by minimizing
a simple variational model (based on a quadratic energy, which
we claim would produce the exact solution if the similar patches were
exactly known---and in large enough quantity), we have experimentally
shown the efficiency of this non-local framework, for some Fourier
pattern. If the loss of coefficients is only in the high frequencies,
then our approach does not produce any improvement over standard zooming
techniques (local or nonlocal). On the other hand, for more complex
masks which miss low frequency Fourier coefficients, our results are
far superior to results obtained with more usual measures of
redundancy.}

%\newpage
\bibliographystyle{siam}
\bibliography{biblio_atom}

\end{document}